%% file: master.tex
\documentclass[twoside,reqno,a4paper]{myamsbook}
\usepackage{eucal,xspace,amssymb}
\usepackage[german,english]{babel}

\numberwithin{equation}{section}
\input{nc}

\begin{document}
\frontmatter
\include{title}
\mainmatter
\include{intro}

\include{w}
\include{cu}
\include{ca}
\include{cc}
\include{sectors}
\include{not}
\backmatter
\bibliographystyle{amsalpha}
\bibliography{abk,mrabbrev,bib}
\end{document}

%% file: nc.tex
\newcommand{\nc}{\newcommand}
\nc{\rnc}{\renewcommand}
\nc{\DMO}{\DeclareMathOperator}
\nc{\nt}{\newtheorem}

\rnc{\thefootnote}{\alph{footnote}}

\nt{Th}{Theorem}[section]
\nt{Lem}[Th]{Lemma}
\nt{Prop}[Th]{Proposition}
\nt{Cor}[Th]{Corollary}

\theoremstyle{definition}
\nt{Def}[Th]{Definition}
\nt{Ex}{Example}

\theoremstyle{remark}
\nt*{Rem}{Remark}
\nt*{Dank}{\Large\sc Danksagung}

\nc{\ADD}[1]{\begin{quote}\tt #1\end{quote}}
\nc{\ART}{and references therein}
\nc{\BH}{bilinear Hamiltonian\xspace}
\nc{\BHs}{{\BH}s\xspace}
\nc{\BP}{basis projection\xspace}
\nc{\BPs}{{\BP}s\xspace}
\nc{\CA}{Cuntz algebra\xspace}
\nc{\CAs}{{\CA}s\xspace}
\nc{\CAL}{$C^*$--algebra\xspace}
\nc{\CALs}{{\CAL}s\xspace}
\nc{\CAR}{CAR algebra\xspace}
\nc{\CARs}{{\CAR}s\xspace}
\nc{\CCR}{CCR algebra\xspace}
\nc{\CCRs}{{\CCR}s\xspace}
\nc{\CE}{conditional expectation\xspace}
\nc{\CEs}{{\CE}s\xspace}
\nc{\COC}{continuous curve\xspace}
\nc{\CR}{commutation relation\xspace}
\nc{\CRs}{{\CR}s\xspace}
\nc{\DHR}{Doplicher, Haag and Roberts\xspace}
\nc{\DR}{Doplicher and Roberts\xspace}
\nc{\FA}{field algebra\xspace}
\nc{\FAs}{field algebras\xspace}
\nc{\FS}{Fock state\xspace}
\nc{\FSs}{{\FS}s\xspace}
\nc{\GG}{gauge group\xspace}
\nc{\HK}{Haag and Kastler\xspace}
\nc{\HSP}{Hilbert space\xspace}
\nc{\HSPs}{{\HSP}s\xspace}
\nc{\IFF}{if and only if }
\nc{\ONB}{orthonormal basis\xspace}
\nc{\PO}{Poincar\'e\xspace}
\nc{\QF}{quadratic form\xspace}
\nc{\QFs}{{\QF}s\xspace}
\nc{\QFT}{quantum field theory\xspace}
\nc{\SSS}{superselection sector\xspace}
\nc{\SSSs}{{\SSS}s\xspace}
\nc{\VN}{von Neumann\xspace}
\nc{\WRT}{with respect to }

\nc{\+}{^\dagger}
\nc{\0}{\mspace{-1.5mu}\circ\mspace{-0.5mu}}
\nc{\1}{\ensuremath{{\boldsymbol1}}}
\nc{\2}{{\frac{1}{2}}}
\nc{\3}{{\tfrac{1}{2}}}
\nc{\4}[1]{{\overline{#1}}}
\nc{\ABS}[1]{\lvert{#1}\rvert}
\nc{\bigABS}[1]{\bigl\lvert{#1}\bigr\rvert}
\nc{\BigABS}[1]{\Bigl\lvert{#1}\Bigr\rvert}
\nc{\biggABS}[1]{\biggl\lvert{#1}\biggr\rvert}
\rnc{\AA}{\ensuremath{\mathfrak A}\xspace}
\nc{\BB}{{\mathfrak B}}
\nc{\CC}{\ensuremath{\mathbb C}\xspace}
\nc{\CK}{\ensuremath{\mathfrak{C}}(\KK)\xspace}
\nc{\CKG}{\ensuremath{\mathfrak{C}(\KK,\kappa)}\xspace}
\nc{\DD}{\ensuremath{\mathfrak{D}}\xspace}
\nc{\DEF}{\equiv}
\nc{\EE}{{\mathfrak E}}
\nc{\EH}{\ensuremath{\WO{\exp\bigl(\3b(H)\bigr)}}}
\nc{\EHV}{\ensuremath{\WO{\exp\bigl(\3b(H_V)\bigr)}}}
\nc{\EHW}{\ensuremath{\WO{\exp\bigl(\3b(H_{W_V})\bigr)}}}
\nc{\HH}{\ensuremath{\mathcal H}\xspace}
\nc{\HSNORM}[1]{{\lVert{#1}\rVert}_{\text{\rm HS}}}
\nc{\IDIAG}{{\II_{\text{\rm diag}}}}
\nc{\IFIN}{{\II^{\text{\rm fin}}}}
\nc{\IHS}{{\II_{\text{\rm HS}}}}
\nc{\II}{\ensuremath{\mathcal{I}}\xspace}
\DMO{\IND}{ind}          
\nc{\KK}{\ensuremath{\mathcal K}\xspace}                    
\nc{\MAT}[4]{\bigl(\begin{smallmatrix}#1&#2\\#3&#4\end{smallmatrix}\bigr)}
\nc{\NN}{\ensuremath{\mathbb N}\xspace}                  
\nc{\NORM}[1]{\lVert{#1}\rVert}
\nc{\bigNORM}[1]{\bigl\lVert{#1}\bigr\rVert}
\nc{\BigNORM}[1]{\Bigl\lVert{#1}\Bigr\rVert}
\nc{\biggNORM}[1]{\biggl\lVert{#1}\biggr\rVert}
\nc{\PLV}{{\mathcal{P}_{L_V}}}
\nc{\PLVT}{{\tilde\mathcal{P}_{L_V}}}
\nc{\PME}{{\Psi(-\1)}}
\nc{\PP}{{\mathfrak P}}
\nc{\PU}{\ensuremath{\Psi(U)}\xspace}
\nc{\PUV}{\ensuremath{\Psi(U_V)}\xspace}
\nc{\PV}[1]{\ensuremath{{\Psi_{#1}(V)}}}
\nc{\PW}[1]{\ensuremath{{\Psi_{#1}(W_V)}}}
\nc{\QK}{\mathcal{Q}(\KK)}
\rnc{\rho}{\ensuremath{\varrho}\xspace}
\DMO{\RAN}{ran}
\DMO{\RE}{Re}
\nc{\RR}{\ensuremath{\mathbb R}\xspace}
\nc{\SC}{^\sharp}
\nc{\SET}[2]{\{#1\ |\ #2\}}
\DMO{\SGN}{sign}
\nc{\SLIM}[1]{\underset{#1\to\infty}{\text{\rm s--lim}}\,}
\nc{\SP}[2]{{{\mathcal S}_{{#1}}^{#2}({\KK},\kappa)}}
\DMO{\SPAN}{span}          
\rnc{\SS}{{\mathfrak S}}
\DMO{\SUPP}{supp}          
\DMO{\TR}{tr}          
\nc{\TT}{\ensuremath{{\mathbb T}}\xspace}
\nc{\WK}{{\WW(\KK,\kappa)}}      
\nc{\WO}[1]{{\;:\!{#1}\!:\;}} 
\nc{\WT}[1]{{\widehat{#1}}}
\nc{\WW}{\mathfrak{W}}
\nc{\ZZ}{\ensuremath{\mathbb Z}\xspace}

%% file: title.tex
\title[Endomorphisms of CAR and CCR Algebras]{Charged Quantum Fields
  Associated with Endomorphisms of CAR and CCR Algebras\vspace*{3cm}
  \begin{otherlanguage}{german}
    {\rm Dissertation zur Erlangung der Doktorw"urde, beim Fachbereich
      Physik der Freien Universit"at Berlin eingereicht im Mai 1998,
      von}
  \end{otherlanguage}\vspace*{2cm}}
\author[C.~Binnenhei]{{\huge Carsten Binnenhei}}
\thanks{\begin{otherlanguage}{german}
    \begin{tabular}{ll}
      Author address: & Institut f\"ur Theoretische Physik der Freien
      Universit"at Berlin,\\
      & Arnim\-allee 14, D--14195 Berlin\\
      E--mail address: & {\tt Carsten.Binnenhei@physik.fu-berlin.de}
    \end{tabular}
  \end{otherlanguage}}
\maketitle
\newpage  
\thispagestyle{plain}
\hspace{14pt}{\sc\begin{tabular}{ll}
    Erstgutachter: & Prof.~Bert Schroer\\
    Zweitgutachter: & Prof.~Robert Schrader\\[1ex]
    Tag der Disputation: & 11.\ Juni 1998
  \end{tabular}}
\vspace*{\fill}
\begin{abstract}
  The appearance of the \CAs ${\mathcal O}_d$ is a generic feature of
  local \QFT. This fact has been discovered by S.~Doplicher and
  J.~E.~Roberts within the algebraic theory of \SSSs. Generators of
  the \CAs arise as charged field operators which implement localized
  endomorphisms of the observable algebra. Whereas the existence of
  such operators can be derived from first principles, little is known
  about the actual construction of these fields in concrete models.
  
  In view of this apparent discrepancy, we develop a comprehensive
  theory of quasi-free endomorphisms of the CAR and \CCRs which give
  rise to representations of the \CAs $\mathcal{O}_d$ on Fock space.
  The number $d$ is the statistics dimension of the endomorphism. It
  can be any power of $2$ (including $1$ and $\infty$) in the CAR
  case, but takes only the values $1$ and $\infty$ in the CCR case.
  
  We obtain necessary and sufficient conditions for implementability
  of quasi-free endomorphisms. By studying extensions of partial \FSs,
  we find that the semigroup of implementable quasi-free endomorphisms
  has a simple structure: It can be written as a product of a group of
  automorphisms which are close to the identity, and the semigroup of
  endomorphisms which leave the given \FS invariant.
  
  We describe the construction of the implementers of a quasi-free
  endomorphism in terms of annihilation and creation operators in
  detail. These implementers span a $d$--dimensional \emph{Fock space
    of isometries.} The Fock structure of the space of implementers is
  the key to the determination of the charge quantum numbers of the
  endomorphism. It entails that implementable endomorphisms with
  statistics dimension $d\neq1$ are always reducible.
  
  We compare the structure of the semigroup of (gauge invariant)
  quasi-free endomorphisms with the generic superselection structure
  of \QFT.\\[.6cm]
  \begin{otherlanguage}{german}
    {\sc Zusammenfassung.}
    Wie S.~Doplicher und J.~E.~Roberts gezeigt haben, treten die
    Cuntz--Algebren $\mathcal{O}_d$ ganz allgemein in der lokalen
    Quantenfeldtheorie auf. Darstellungen der Cuntz--Algebren werden
    erzeugt von ladungstragenden Quantenfeldern, welche lokalisierte
    Endomorphismen der Observablenalgebra implementieren. Es sind aber
    bisher keine Beispiele f"ur die Konstruktion solcher geladenen
    Felder in konkreten Modellen bekannt.
    
    Angesichts dieser Situation entwickeln wir eine vollst"andige
    Theorie derjenigen quasifreien Endomorphismen der CAR-- und
    CCR--Algebren, die zu Darstellungen der Cuntz--Algebren
    $\mathcal{O}_d$ auf dem Fockraum f"uhren. Dabei ist $d$ die
    statistische Dimension des Endomorphismus. Wie sich zeigt, kann
    $d$ im CAR--Fall eine beliebige Potenz von zwei sein, im CCR--Fall
    jedoch nur eins oder unendlich.

    Wir beweisen ein notwendiges und hinreichendes Kriterium f"ur die
    Implementierbarkeit quasifreier Endomorphismen. Die Halbgruppe der
    implementierbaren Endomorphismen hat eine einfache Struktur: Jeder
    implementierbare Endomorphismus l"a"st sich zerlegen in ein
    Produkt eines Automorphismus, der "`nahe"' bei der Identit"at
    liegt, und eines Endomorphismus, der den gegebenen Fockzustand
    invariant l"a"st.

    Wir finden explizite Formeln f"ur die geladenen Felder, die einen
    gegebenen Endomorphismus implementieren. Aus diesen Formeln l"a"st
    sich ableiten, da"s der von den Feldern aufgespannte
    $d$--dimensionale Hilbertraum selbst eine Fockraumstruktur
    tr"agt. Weiter lassen sich die Ladungsquantenzahlen der
    Endomorphismen bestimmen. Endomorphismen mit statistischer
    Dimension ungleich eins sind stets reduzibel.

    Wir diskutieren "Ahnlichkeiten und Unterschiede zwischen der
    Halbgruppe der (eichinvarianten) quasifreien Endomorphismen und
    der Halbgruppe der lokalisierten Endomorphismen in der Theorie der
    Superauswahlsektoren. 
  \end{otherlanguage}
\end{abstract}

\tableofcontents

\thispagestyle{plain}
\begin{otherlanguage}{german}
  \begin{quotation}
    \large{\it Da"s mit den steigenden Anspr"uchen an die Genauigkeit
      der Messungen auch die Kompliziertheit der Instrumente immer
      gr"o"ser wird, findet ohne weiteres Verst"andnis und Billigung.
      Aber da"s bei der fortgesetzten Verfeinerung der gesetzlichen
      Zusammenh"ange zu ihrer Formulierung Definitionen und Begriffe
      benutzt werden, die sich immer weiter von altgewohnten Formen
      und anschaulichen Vorstellungen entfernen, macht man
      stellenweise der theoretischen Forschung zum Vorwurf, ja man
      will darin Anzeichen daf"ur erblicken, da"s sie sich auf einem
      Irrweg befindet.

      \noindent Nichts kann kurzsichtiger sein als eine derartige 
      Vermutung.}  
    \\[1ex] Max Planck, {\em Sinn und Grenzen der exakten 
      Wissenschaft}, 1941
  \end{quotation}
\end{otherlanguage}


%% file: intro.tex
\chapter*{Introduction} 
\noindent
{\sc This} thesis is essentially concerned with three different \CALs:
the algebra of the canonical anti\CRs (CAR); the Weyl algebra, the
exponentiated version of the canonical \CRs (CCR); and the Cuntz
algebra\footnote{The reader who is unfamiliar with operator algebras
  may think of a \CAL as an algebra of bounded linear operators on
  some \HSP which is closed under taking adjoints and uniform limits.
  Textbooks on operator algebras, with applications to physics, are
  \cite{BR,EK}.}. Common to all three is the fact that each is
associated, in a specific way, with an underlying vector space: The
\CAR is the $C^*$--Clifford algebra over a real \HSP, the Weyl algebra
is the \CAL generated by a projective unitary representation of a real
symplectic space, and the \CA is the universal \CAL generated by a
complex \HSP. (The vector spaces belonging to the CAR and \CCRs will
always be assumed to be infinite dimensional, as we will be
exclusively dealing with systems possessing infinitely many degrees of
freedom.)

The CAR and \CCRs are the most prominent algebras in quantum physics,
due to their distinguished r\^{o}le in describing systems of Fermions
and Bosons.
The {\em canonical anti\CRs} have been introduced by Jordan and Wigner
in 1928 in their analysis of the implications of Pauli's exclusion
principle for the Fermi gas \cite{JW}.
As an abstract \CAL, the \CAR has a very simple structure: It is an
approximately finite dimensional algebra, in a sense a non--commutative
analogue of a zero--dimensional topological space, and is isomorphic to
an infinite tensor product of copies of the two by two matrices.

Heisenberg's {\em canonical \CRs} were found by Born in 1925.  They
first appeared in Born's joint work with Jordan on the matrix
formulation of quantum mechanics \cite{BJ}, but were independently
obtained two months later by Dirac \cite{Dir26}.
Inspired by group theoretic considerations, H. Weyl discovered the
usefulness\footnote{He regarded his relations as the answer to the 
\begin{otherlanguage}{german}
  ``Frage nach dem Wesen und der richtigen Definition der kanonischen
  Variablen''
\end{otherlanguage} \cite{Wey28}.}
of replacing Heisenberg's \CRs (which do not have bounded \HSP
realizations) with their exponential form \cite{Wey28}.  In contrast
to the CAR case, the Weyl algebra is a very ``large'' object (it is
not separable), and not much seems to be known about its abstract
properties.

The {\em\CAs} are the basic examples of infinite \CALs (those
containing non--unitary isometries) and are of great importance in the
general structure theory of \CALs. They have been introduced by Cuntz
in 1977 \cite{C77}, but their generators (``\HSPs of isometries'') had
been studied before by \DR in the context of general \QFT
\cite{DR72,R76}.  The generic appearance of the \CAs in \QFT has been
established rather recently by \DR \cite{DR90}.

Each of these algebras possesses a natural class of structure
preserving transformations (``*--endomorphisms''), namely those which
arise from linear symmetries of the underlying vector spaces. These
transformations will in all three cases be termed \emph{quasi--free
  endomorphisms,} although in the CAR and CCR cases several other
names are also used in the literature (linear canonical
transformations, Bogoliubov transformations, one--particle
transformations, \dots).  Thus quasi--free endomorphisms correspond to
real orthogonal transformations in the CAR case, to real symplectic
transformations in the CCR case, and to complex isometries in the \CA
case. Note that all these transformations need not be invertible
(their ranges need not be the whole space) if the underlying vector
spaces are infinite dimensional. If they are invertible, then the
transformations of the algebras will be called quasi--free {\em
  automorphisms}.

It seems that quasi--free automorphisms (of the \CCR) first appeared
in Bogoliubov's treatment of the Bose gas, although general (not
necessarily linear) canonical transformations already occur in the
famous ``Dreim\"annerarbeit'' of Born, Heisenberg and Jordan
\cite{BHJ}. In the meantime there has been a tremendous work on
quasi--free automorphisms of the CCR and \CARs, so that it would be
hard to say anything new about them.
Our interest is mainly in genuine (non--invertible) endomorphisms,
which have not been treated systematically in the literature so far.
\emph{We develop a complete theory of those quasi--free endomorphisms
  of the CAR and \CCRs which are, in the widest sense, related to
  second quantization}\footnote{It should be noted that parts of our
  results have already been published \cite{CB1,CB2,CB3}.}. The second
quantization of a genuine endomorphism is however no longer a single
unitary operator, but, as follows from the work of \DR, a whole \HSP
of isometries on Fock space.  This means that there is a
representation of a \CA associated with each such endomorphism.

But why should a physicist care about endomorphisms of \CALs?  Let us
give an answer to this question by sketching the history of the theory
of \SSSs.
\specialsection*{The Algebraic Theory of Superselection Sectors}
\label{sec:SSS}
The theory of \SSSs is an important and particularly successful branch
of local quantum field theory\footnote{We refer 
  to Haag's beautiful book \cite{H}
  for a comprehensive introduction into the subject. We further
  recommend the lecture notes of Fredenhagen \cite{F}, Roberts
  \cite{R90}, and Schroer \cite{BSch}.  The early history of
  superselection rules has been nicely reviewed by Wightman
  \cite{W95}.}. It was initiated by the observation of Wick, Wightman
and Wigner in 1952 that the validity of the superposition principle in
quantum physics is limited by what they called \emph{superselection
  rules} \cite{WWW}. For instance, there is no interference between a
single--electron state $\psi_-$ and a single--positron state $\psi_+$.
In the state $\psi=\alpha_+\psi_++\alpha_-\psi_-$, the relative phase
between the $\psi_\pm$--components cannot be measured, but can be
arbitrarily changed by applying
global gauge transformations $\psi\mapsto\psi'=
\alpha_+e^{i\lambda}\psi_++\alpha_-e^{-i\lambda}\psi_-$.  Such $\psi$
is not a coherent superposition, but a mixture of the pure states
$\psi_\pm$, with weights $\ABS{\alpha_\pm}^2$.  Accordingly, matrix
elements of physical observables between $\psi_+$ and $\psi_-$ must
vanish, observables are gauge invariant, and the physical \HSP splits
up into invariant ``coherent'' subspaces, each carrying a definite
value of the electric charge. The unobservability of relative phases
in such situations led Wick, Wightman and Wigner to the conclusion
that the parities of elementary particles with different charges
cannot be compared.

Significant progress towards a deeper understanding of the general
structure of \QFT, and in particular of the concept of superselection
rules,
was achieved by \HK in 1964 \cite{HK}. Building on earlier ideas of
Haag \cite{H57}, they proposed a {\CAL}ic treatment of \QFT. Whereas
``global'' \CAL approaches to quantum field theory had previously been
advocated by several authors, notably by Araki, Haag, Schroer, and
Segal, \HK emphasized the importance of the \emph{principle of
  locality} in field theory, i.e.\ the absence of actions at a
distance.  Locality allows to assign to any bounded region $O$ in
space--time the \CAL $\AA(O)$ generated by all observables which can
be measured within that region, such that the natural partial ordering
between regions is preserved. The local algebras provide in some sense
a ``coordinate--free'' description of \QFT as opposed e.g.\ to
Wightman's approach \cite{PCT} where one has to make a specific choice
among all fields belonging to the same Borchers class.  \emph{Einstein
  causality} requires that algebras belonging to spacelike separated
regions have to commute with each other. \HK postulated that such a
correspondence between space--time regions and algebras of local
observables should fix the content of the theory completely.  This
point of view is plausible because ``ultimately all physical processes
are analyzed in terms of geometric relations''
\cite{HK}, and was supported by the general theory of collision
processes that had been developed earlier by Haag and Ruelle
\cite{H58,Ru}.

The formulation of \QFT in terms of local algebras permitted a new
look at superselection rules. \HK introduced the \emph{quasilocal
  algebra} \AA as the \CAL generated by all local observables. The
algebra \AA does not contain global quantities such as total charge or
total energy; these can only be obtained as strong limits in specific
representations. Now \AA is expected to possess an abundance of
inequivalent irreducible representations, e.g.\ representations
associated with states having different behaviour at spacelike
infinity. (\HK believed that this was the only mechanism to produce
inequivalent representations in \QFT, but it was soon recognized that
there exist inequivalent representations even among the states with
the same asymptotic behaviour; cf.\ \eqref{DHR} below.) That generic
\CALs have lots of inequivalent representations had already been
discovered by \VN in 1939 \cite{vN}\footnote{He showed this on the
  example of an infinite tensor product of $2\times2$ matrix algebras
  (isomorphic to the \CAR) by exhibiting representations of type ${\rm
    I}_\infty$, ${\rm II}_1$ and ${\rm II}_\infty$. The field
  theoretic examples of inequivalent representations found in the
  1950s also came from CAR and \CCRs, e.g.\ from non--implementable
  quasi--free automorphisms, and thus are closely related to the
  subject of this thesis.}.  Actually, any simple infinite dimensional
\CAL (besides the algebra of compact operators on a separable \HSP,
which corresponds to the \CCR for finitely many degrees of freedom,
i.e.\ to ordinary quantum mechanics) possesses uncountably many
inequivalent irreducible representations. It was Haag who realized in
the mid--fifties the \emph{need} for inequivalent representations in
order to obtain interacting quantum fields (``Haag's Theorem'').

The different coherent subspaces are now interpreted as inequivalent
irreducible representation spaces of the single algebra \AA. \HK
called (the unitary equivalence class of) an irreducible
representation of \AA a {\em\SSS}. However, only a small subclass of
all representations of \AA can be expected to have a physical
interpretation. In \QFT one is mainly interested in states which
describe local finite--energy excitations of the vacuum.  The
corresponding sectors are called \emph{charge} \SSSs (and from now on,
a sector will always mean a charge \SSS).  Here the term ``charge'' is
used in a very broad sense: It applies to any quantity which can be
used to label the various sectors.

\HK argued that already a single sector should comprise all relevant
physical information. If one starts e.g.\ with a state in the vacuum
sector, one can create a particle of unit charge together with its
antiparticle and then send the antiparticle ``behind the moon''. The
resulting state will deviate, \WRT local measurements, arbitrarily
little from a state in the charge--one sector. Thus any given state
belonging to some sector can be approximated by states in each other
sector. In the terminology of \HK, all sectors are ``physically
equivalent'', which is, by a theorem of Fell, tantamount to all
sectors being ``equally faithful'' (they all have the same kernel).
All \SSSs thus determine the same ``abstract'' \CAL \AA, and the
choice of a particular representation of \AA appears essentially as a
matter of convenience. This is the solution of the (at that time much
discussed) problem of inequivalent representations in \QFT offered by
\HK.

If one takes this philosophy seriously, one faces the basic problem:
Given the quasilocal algebra \AA (together with its local structure)
in a, say, vacuum sector, how can one extract the interesting physical
information? In particular, one would like to determine all charge
\SSSs together with a set of unobservable fields generating these
sectors from the vacuum. The charge quantum numbers ascribed to the
various sectors and fields are expected to be related to some sort of
inner symmetries which act covariantly on the fields and trivially on
the observables, and one would like to understand the laws of
composition and exchange of charges (``statistics'').

In the first step of such an investigation one has to specify which
representations of \AA are to be regarded as ``local excitations of
the vacuum'', i.e.\ as charge \SSSs. In his pioneering work
\cite{B65,B67} Borchers proposed to consider all irreducible positive
energy representations $\pi$ which are ``strongly locally equivalent''
to a given vacuum representation $\pi_0$ and fulfill a certain ``weak
duality'' condition\footnote{\label{foo:PER}$\pi$ is a \emph{positive
    energy representation} if the space--time translations are
  unitarily implemented in $\pi$ such that the relativistic spectrum
  condition holds (see \cite[\ART]{BFK} for a local version of the
  spectrum condition).  If one assumes that the local algebras are
  weakly closed (this can be done in presence of a distinguished
  vacuum representation), then the restrictions of positive energy
  representations to the local algebras are known to be unitarily
  equivalent. It is commonly believed that the \VN algebras associated
  with double cones (see below) are all isomorphic to the unique
  hyperfinite type ${\rm III}_1$ factor \cite{A64,F85,BAF}. --- $\pi$
  is \emph{strongly locally equivalent} to $\pi_0$ if the restrictions
  of $\pi$ and $\pi_0$ to the relative commutant of any local algebra
  are equivalent, and \emph{weak duality} means that, in the
  representation $\pi$, the relative commutant is weakly dense in the
  commutant.  Unfortunately, the DHR criterion \eqref{DHR} is also
  sometimes referred to as ``strong local equivalence''.}. Under these
assumptions the unitary operators implementing the strong local
equivalence could be interpreted as charged local fields.

However, it soon became clear that Borchers' assumptions were violated
in typical examples. Positivity of the energy should of course hold in
any reasonable theory, but is in general, due to infrared problems,
too weak a condition to allow a complete classification of
representations (such a classification is however possible under
certain circumstances, as e.g.\ in conformal field theory \cite{BMT}).
In order to clarify the r\^{o}le of Borchers' other assumptions, \DHR
started a careful analysis of the superselection structure of
elementary particle physics in 1969 which led to a series of seminal
papers \cite{DHR1,DHR2,DHR3,DHR4} and was to a certain extent
completed some 20 years later by Doplicher and Roberts \cite{DR90}.

In \cite{DHR1} \DHR considered theories where a compact global gauge
group acts on a \emph{given} \FA such that the observables are
precisely the gauge invariant fields. They studied the representations
of \AA that are contained in a fixed vacuum representation of the
larger \FA. Among other things, they found that the sectors of \AA
occurring in this manner are in a natural way labelled by the
inequivalent irreducible representations of the \GG, and that weak
duality implies the failure of strong local equivalence (see also
\cite{R70}). Therefore Borchers' fields do not exist in this case.
Instead \DHR obtained two properties of these sectors which are
closely related to, but much more significant than, strong local
equivalence and weak duality. The first is that the various
representations of the algebras belonging to the spacelike complement
of any bounded region are unitarily equivalent, in symbols
\begin{equation}
  \label{DHR}
  \pi|_{\AA(O')}\simeq\pi_0|_{\AA(O')}.
\end{equation}
Here $\pi_0$ is the vacuum representation of \AA (corresponding to
the trivial representation of the \GG), $\pi$ is some superselection
sector, and $\AA(O')$ is the $C^*$--subalgebra of \AA generated by
all local observables which can be measured within the spacelike
complement $O'$ of the bounded space--time region $O$. This condition
is weaker than Borchers' strong local equivalence.  The second
property is a strengthening of Einstein causality
\begin{equation}
  \label{HD}
  \pi(\AA(O))=\pi(\AA(O'))'.
\end{equation}
The prime on the right denotes the commutant (the set of all bounded
operators on the representation space of $\pi$ which commute with all
elements of $\pi(\AA(O'))$). Here $\pi(\AA(O))$ is assumed to be
weakly closed, cf.\ footnote (\ref{foo:PER}).  Eq.~\eqref{HD} can only
be expected to hold for particularly simple regions $O$, e.g.\ for
double cones\footnote{Double cones are non--void intersections of
  suitably situated open backward with forward light cones. They
  constitute the simplest (yet sufficiently rich) class of causally
  complete bounded regions in Minkowski space which is closed under
  \PO transformations.}, and henceforth, we will generically take
double cones as localization regions. \DHR showed that Eq.~\eqref{HD}
holds in this form in all \emph{simple} sectors, i.e.\ in all sectors
corresponding to one--dimensional representations of the \GG, but not
in non--simple sectors. The condition \eqref{HD} is stronger than
Borchers' weak duality, but the point is that it is not assumed to
hold in \emph{all} sectors.

Eq.~\eqref{DHR} is commonly called the \emph{DHR selection criterion.}
It is supposed to single out most sectors of interest in theories with
short range forces.
It implies that the states in the representation
$\pi$ look asymptotically like the vacuum, and that the charges
distinguishing between $\pi$ and $\pi_0$ can be localized in any
bounded region. The DHR criterion is known to hold for all irreducible
positive energy representations in conformal field theory \cite{BMT},
but it excludes ``topological charges'' which appear even in purely
massive theories \cite{BF,FM}, and charges which can be measured at
arbitrary distances, e.g.\ by virtue of Gau\ss' law. Sectors
satisfying the DHR criterion are automatically \PO covariant with
positive energy under rather general assumptions \cite{GL}, so that
covariance does not have to be assumed from the outset.

Following a proposal of Schroer in \cite{FRS1}, \eqref{HD} is called
\emph{Haag duality} in order to distinguish it from several other
concepts of ``duality'' occurring in physics.  It means that the local
algebras cannot be enlarged in the representation $\pi$ without
violating Einstein causality. It was originally invented by Haag and
Schroer as an expression of underlying relativistic dynamics \cite{HS,
  H63}. The failure of Haag duality for double cones in the vacuum
sector indicates that the theory is in some sense incomplete. It is
typically caused by spontaneous breakdown of inner symmetries
\cite{R74,BDLR}, but is also generic in two--dimensional \QFT (cf.\ 
\cite{MM1}). This phenomenon can be traced back to the existence of
operators which are only invariant under the unbroken symmetries, but
not under the broken ones, in the first case; and to the existence of
``kink operators'' in the second case.
\medskip{}

In \cite{DHR2} the converse problem of reconstructing the field
algebra and the \GG from \AA and $\pi_0$ was solved for the set of all
\PO covariant sectors satisfying the DHR criterion and Haag duality
(simple sectors). This set of sectors has the structure of a discrete
Abelian group, and its Abelian dual can be viewed as the \GG.  There
is a unique \FA consisting of Bose and Fermi fields which has \AA as
its gauge invariant part and acts irreducibly on the ``physical'' \HSP
which contains each simple sector with multiplicity one.

\DHR extended this analysis to the class of \emph{all} charge \SSSs
conforming with the DHR criterion \eqref{DHR} in \cite{DHR3}. It is
remarkable how little input is needed for their methods to apply, the
deeper reason for that being an underlying general duality theory for
compact groups. (The striking analogy of the superselection theory with
the representation theory of compact groups had been fully recognized
in \cite{DR72}, but it took almost 20 years to finish the proof that
both structures are really identical \cite{DR90}).
It suffices to start with the quasilocal algebra \AA, together with a
faithful irreducible representation $\pi_0$ on a separable \HSP
$\HH_0$, such that Haag duality \eqref{HD} for double cones and
Borchers' ``property B''\footnote{\label{foo:B}Property B means that
  the local algebras are ``almost type III'' (any non--trivial local
  projection is, at least within a slightly larger algebra, equivalent
  to \1).  This property was derived from standard assumptions by
  Borchers \cite{B67a}. It implies that any local algebra contains a
  subalgebra isomorphic to the \CA $\mathcal{O}_2$.
} hold in $\pi_0$.  Let us sketch some of their results.

First of all, it is an immediate consequence of \eqref{DHR} and
\eqref{HD} that any representation $\pi$ satisfying the DHR criterion
is unitarily equivalent to a representation on $\HH_0$ of the form
$\pi_0\0\rho$ where \rho is a \emph{localized endomorphism} of \AA.
Such \rho has the following properties: \label{page:LOCEND}
\begin{itemize}
\item \rho is a unital map from \AA into itself which preserves the
  algebraic structure, the star ``$*$'' and the norm.
\item \rho is \emph{localized} in some double cone $O$ (and then also
  in every larger double cone) in the sense that
  \begin{equation}
    \label{LOCEND}
    \rho(a)=a, \quad a\in\AA(O'),
  \end{equation}
  and it maps the algebras belonging
  to larger regions than $O$ into themselves.
\item \rho is \emph{transportable}: If $\hat O$ is another double
  cone, then there is an endomorphism $\hat\rho$ localized in $\hat O$
  such that the representations $\pi_0\0\rho$ and $\pi_0\0\hat\rho$
  are equivalent. Such $\hat\rho$ has the form
  \begin{equation}
    \label{TRANS}
    \hat\rho(a)=u\rho(a)u^*, \quad a\in\AA,
  \end{equation}
  where $u$ is a unitary element contained in any $\AA(\check{O})$
  with $\check{O}\supset O\cup\hat O$ by Haag duality.
\end{itemize}
Conversely, any localized endomorphism \rho gives rise to a
representation $\pi_0\0\rho$ fulfilling \eqref{DHR}. Thus the sectors
fulfilling the DHR criterion are in one--to--one correspondence with
the equivalence classes $[\rho]$ of irreducible localized
endomorphisms \rho. Here two endomorphisms \rho, $\hat\rho$ are called
\emph{equivalent} if they are related to each other as in
\eqref{TRANS}, and \rho is \emph{irreducible} if the corresponding
representation is.  Irreducible endomorphisms are not necessarily
invertible, it can happen that $\pi_0(\rho(\AA))''=\pi_0(\AA)''$ but
$\rho(\AA)\subsetneq\AA$. Invertible endomorphisms (= automorphisms)
correspond to simple sectors and are characterized by the property
that $\rho^2$ is irreducible.

Several operations can be performed within the set of localized
endomorphisms. The \emph{direct sum} of localized endomorphisms
$\rho_1,\dots,\rho_n$ is defined as follows.  Take local observables
$v_1,\dots,v_n$ with the properties
\begin{subequations}
  \label{CUNTZ}
  \begin{align}
    v_j^*v_l &=\delta_{jl}\1,\label{CUNTZa}\\
    \sum_{j=1}^nv_jv_j^* &=\1\label{CUNTZb}
  \end{align}
\end{subequations}
(by the way, these are the defining relations of the \CA
$\mathcal{O}_n$), and set
\begin{equation}
  \label{DIRSUM}
  \bigoplus_{j=1}^n\rho_j(a)\DEF\sum_{j=1}^nv_j\rho_j(a)v_j^*, \quad
  a\in\AA.
\end{equation}
Such $v_j$ exist due to property B. The direct sum $\oplus_j\rho_j$ is
again a localized endomorphism, and its equivalence class does not
depend on the choice of the $v_j$ (if $\hat v_1,\dots,\hat v_n$ is
another collection of local observables fulfilling \eqref{CUNTZ}, then
the equivalence between the two direct sums is established, in the
sense of \eqref{TRANS}, by the unitary $u=\sum \hat v_jv_j^*$).
\label{page:SUB} 
Similarly, if \rho is reducible, and $p\in\pi_0(\rho(\AA))'$ a
non--trivial projection, then there exists (by Haag duality and
property B) a local observable $v$ such that 
\begin{equation}
  \label{BORCHERS}
  vv^*=p,\qquad v^*v=\1,
\end{equation}
and the \emph{``subobject''} of \rho corresponding to $p$ can be
defined by $\rho_p(a)=v^*\rho(a)v$.  The equivalence class of $\rho_p$
is again independent of the choice of $v$.  Finally, the {\em
  composition} $\rho_1\0\rho_2$ of two localized endomorphisms is a
localized endomorphism whose equivalence class depends only on the
classes of $\rho_1$ and $\rho_2$, and in particular not on the order
of the factors, because endomorphisms which are localized in mutually
spacelike double cones commute with each other.

Thus direct sums and subrepresentations of representations fulfilling
the DHR criterion also fulfill this criterion, and one has a
well--defined commutative product (``fusion'') of equivalence classes
of such representations, corresponding to the composition of charges
and given by the composition of the associated endomorphisms.  The
availability of this product of sectors is the main advantage of
working with endomorphisms.

These observations provide the basis for an intrinsic understanding of
statistics, which is independent of a possible particle interpretation
of the theory, and for the reconstruction of gauge symmetries and
charged fields from observable data only.

The \emph{statistics} of a sector $[\rho]$ describes the effect of
exchanging identical charges (remember that every sector carries a
specific charge). It is determined by the \emph{statistics operator}
$\varepsilon_{\rho}$, which can be defined as follows. Pick a unitary
``charge transporter'' $u$ as in \eqref{TRANS} such that \rho and the
corresponding $\hat\rho$ are localized in spacelike separated double
cones. Then
\begin{equation}
  \label{STOP}
  \varepsilon_{\rho}\DEF u^*\rho(u)
\end{equation}
is a unitary operator which commutes with all elements of
$\rho^2(\AA)$, and its definition is independent of the particular
choice of $u$.  (This is not true in two--dimensional Minkowski space,
where the spacelike complement of a double cone has two connected
components.  There one can have two different choices of
$\varepsilon_{\rho}$, one the other's inverse, depending on whether
$\hat\rho$ is localized to the left or right of \rho. This possibility
had already been observed by Streater and Wilde in 1970 \cite{StrW}.)
The statistics operator $\varepsilon_{\rho}$ fulfills the algebraic
relations
\begin{align}
  \varepsilon_{\rho}\rho\bigl(\varepsilon_{\rho}\bigr)
    \varepsilon_{\rho} &=\rho\bigl(\varepsilon_{\rho}\bigr)
    \varepsilon_{\rho}
    \rho\bigl(\varepsilon_{\rho}\bigr),\label{YB}\\
  \varepsilon_{\rho}^2 &=\1,\label{E1}
\end{align}
so that the operators $\rho^n(\varepsilon_{\rho}),\ n\geq0$, fulfill
the characteristic relations of elementary permutations
(transpositions) which exchange $n$ and $n+1$. Thus canonically
associated with any sector is a unitary representation of the infinite
permutation group. (Relation \eqref{E1} gets in general lost in two
dimensions, so that one obtains representations of the infinite
\emph{braid group} instead. The first examples of sectors with Abelian
braid group statistics were again given in \cite{StrW}.) This
permutation group representation is
analogous to the action of the permutation group on wave functions in
quantum mechanics. It describes permutations of factors in products of
localized state vectors (or permutations of identical particles, if
the theory has a particle content; see \cite{DHR4}).

The analysis of the statistics of the sector $[\rho]$ now proceeds
with the help of left inverses. A \emph{left inverse} $\phi_\rho$ of a
non--invertible endomorphism \rho is a substitute for the inverse of
an automorphism. It is a unital positive linear map from \AA into
itself which is not multiplicative on the whole algebra \AA, but
satisfies
$$\phi_\rho(a\rho(b))=\phi_\rho(a)b,\qquad a,b\in\AA.$$
It follows
from this that $\phi_\rho\0\rho=\text{id}$, that $\rho\0\phi_\rho$ is
a \CE from \AA onto $\rho(\AA)$, and that $\phi_\rho$ enjoys the same
localization properties \eqref{LOCEND} as \rho. Such $\phi_\rho$
always exists, and corresponds to the physical operation of
transferring the charge of \rho to spacelike infinity. 
By applying a left inverse $\phi_\rho$ to the statistics operator
$\varepsilon_{\rho}$ one obtains the \emph{statistics parameter}
$\lambda_{[\rho]}$,
$$\lambda_{[\rho]}\DEF\phi_\rho(\varepsilon_{\rho}).$$
The
statistics parameter is a scalar because $\phi_\rho$ maps
$\rho^2(\AA)'$ into $\rho(\AA)'$, and the latter contains only scalars
by Schur's Lemma. $\lambda_{[\rho]}$ is a numerical
invariant\footnote{See \cite{FRS1,FRS2,KMR}
  for a discussion of uniqueness of $\lambda$ 
  in two dimensions.}
of the sector $[\rho]$. It can be used to classify the statistics of
$[\rho]$. In the case of ``finite statistics''\footnote{The case
  $\lambda_{[\rho]}=0$ does not occur in theories with
  particle--antiparticle symmetry and is usually disregarded.
  Reasonable examples of sectors with infinite statistics dimension
  have been given by Fredenhagen \cite{F94}. See also \cite{BCL}.}
(i.e.\ $\lambda_{[\rho]}\neq0$), the left inverse $\phi_\rho$ is
unique,
and one obtains the \emph{statistics phase} $\eta_{[\rho]}$, a complex
number of modulus one, and the \emph{statistics dimension}
\label{page:D} $d_{[\rho]}\geq1$ by polar decomposition:
$$\lambda_{[\rho]}=\frac{\eta_{[\rho]}}{d_{[\rho]}}.$$
Simple sectors
are precisely the sectors with $d_{[\rho]}=1$. The statistics
dimension can be viewed as a measure for the deviation from Haag
duality in the sector $[\rho]$. It can also be defined for reducible
endomorphisms and coincides with the square root of the minimal index
\cite{J,Ko} of the inclusion $\rho(\AA(O))\subset\AA(O)$ \cite{L}. It
is additive on direct sums and multiplicative on products of
endomorphisms. Any localized endomorphism with finite statistics is a
finite direct sum of irreducible endomorphisms with finite statistics.

\DHR classified the possible statistics in Minkowski space of
dimension greater than 2 \cite{DHR3}. There the statistics phase is
just a sign $\eta_{[\rho]}=\pm1$, and the statistics dimension
$d_{[\rho]}$ is a natural number (in the infinite statistics case one
sets $d_{[\rho]}=\infty$). These numbers characterize the unitary
representation of the permutation group: 
Depending on the sign $\eta_{[\rho]}$, a sector obeys either {\em
  para--Bose} or \emph{para--Fermi} statistics of order $d_{[\rho]}$.
This means that all representations of the permutation group occur
whose Young tableaux have columns resp.\ rows up to length
$d_{[\rho]}$. In a sector with infinite statistics, all irreducible
representations of the permutation group occur. Moreover, for every
sector $[\rho]$ with finite statistics, there exists a unique {\em
  conjugate sector} $[\bar\rho]$ which is determined by the property
that the product $[\rho\bar\rho]$ contains the vacuum sector as a
subrepresentation. \label{page:CONJ} A sector and its conjugate have
the same statistics: $\lambda_{[\rho]}=\lambda_{[\bar\rho]}$
(``particle--antiparticle symmetry''). The conjugate sector can be
viewed as arising from the state induced by applying the left inverse
to the vacuum state.
\medskip{}

For the program of reconstructing the \FA and the \GG from the
observables, it proved instructive to reinvestigate the situation with
given \FA and \GG \cite{DR72}. As mentioned above, the sectors
occurring under these circumstances correspond to the various
irreducible representations of the \GG and satisfy the DHR
criterion~\eqref{DHR}, hence can be described by localized
endomorphisms (one also needs a duality assumption on the fields for
this). They all have finite statistics. In \cite{DR72} \DR assigned to
a localized endomorphism \rho the closed linear subspace $H(\rho)$ of
local field operators $\Psi$ which induce \rho:
\begin{equation}
  \label{HRHO}
  H(\rho)\DEF\{\Psi\ |\ \Psi a=\rho(a)\Psi\quad\text{for all }a\in\AA\}.
\end{equation}
Because any field which commutes with all quasilocal observables is a
multiple of the identity, each $\Psi\in H(\rho)$ is a multiple of an
isometry
$$\Psi^*\Psi=\NORM{\Psi}^2\1.$$
Likewise, $\Psi^*\Psi'$ is proportional to \1 for any two
$\Psi,\Psi'\in H(\rho)$, defining an inner product in $H(\rho)$:
$$\langle\Psi,\Psi'\rangle\1\DEF\Psi^*\Psi'.$$
This inner product induces the usual operator norm. Thus $H(\rho)$ is
a \emph{\HSP of isometries} \cite{R76} inside the \FA. It has
the property that the joint kernel of all $\Psi^*$ vanishes:
$\cap\ker\Psi^*=\{0\}$. Moreover, the dimension of $H(\rho)$ is equal
to the statistics dimension of \rho: 
\begin{equation}
  \label{DIMHR}
  \dim H(\rho)=d_{[\rho]},
\end{equation}
and the gauge action restricts to a
continuous unitary representation of the \GG on $H(\rho)$. In this way
one obtains a concrete equivalence between the ``category'' of
localized endomorphisms with finite statistics (whose morphisms are
the intertwiners between endomorphisms), and the category of finite
dimensional continuous unitary representations of the compact \GG
(with morphisms the intertwiners between representations). This
equivalence preserves irreducibility and direct sums. Since
$H(\rho_1)\otimes H(\rho_2)$ is canonically isomorphic to
$H(\rho_1\rho_2)=H(\rho_1)\cdot H(\rho_2)$, the composition of
endomorphisms corresponds to taking tensor products of
representations. The permutation symmetry is related to changing the
order of factors in tensor powers, and charge conjugation corresponds
to passing to the complex conjugate representation.

Since our own work will be concerned with the description of the \HSPs
$H(\rho)$ associated with quasi--free endomorphisms \rho, let us add
the remark that any \ONB $\Psi_1,\dots,\Psi_{d_{[\rho]}}$ in $H(\rho)$
fulfills Cuntz' relations \eqref{CUNTZ} and implements \rho in the
sense that
\begin{equation}
  \label{IMP}
  \rho(a)=\sum_{j=1}^{d_{[\rho]}}\Psi_ja\Psi_j^*,\qquad a\in\AA.
\end{equation}
This concept of \emph{implementation of endomorphisms by \HSPs of
  isometries} reduces, in the case $d_{[\rho]}=1$, to the familiar
unitary implementation of automorphisms. Since $H(\rho)$ is a
representation space of the \GG, the $\Psi_j$ transform like a tensor
under gauge transformations. Indeed, \DR have shown that the elements
of $H(\rho)$ are the ``typical elements'' of the \FA in the sense that
any irreducible tensor $\Phi_1,\dots,\Phi_d$ of local fields is of the
form $\Phi_j=a\Psi_j$ with $a\in\AA$ and $\Psi_j\in H(\rho)$, for some
irreducible localized endomorphism \rho. It follows that the linear
span of all $\Psi\in H(\rho)$, where \rho runs through all
endomorphisms localized in a double cone $O$, is weakly dense in the
\VN algebra of fields localized in $O$. We would also like to note
that, associated with $H(\rho)$, there is a ``Bosonized'' version
$\hat\varepsilon_{\rho}$ of the statistics operator, obtained by
setting
\begin{equation}
  \label{BSO0}
  \hat\varepsilon_{\rho}\DEF\sum_{j,l=1}^{d_{[\rho]}}\Psi_j\Psi_l
  \Psi_j^*\Psi_l^*.
\end{equation}
This is the operator which effects the exchange of factors in a tensor
product, since it fulfills
$$\hat\varepsilon_{\rho}\Psi\Psi'=\Psi'\Psi,\qquad
\Psi,\Psi'\in H(\rho).$$
There is also a simple formula for the left inverse of \rho in terms
of the $\Psi_j$:
\begin{equation}
  \label{LI0}
  \phi_\rho(a)=\frac{1}{d_{[\rho]}}\sum_{j=1}^{d_{[\rho]}}
  \Psi_j^*a\Psi_j,\qquad a\in\AA.  
\end{equation}
The corresponding ``Bosonized'' statistics parameter is of course
given by
\begin{equation}
  \label{BSP}
  \phi_\rho\bigl(\hat\varepsilon_{\rho}\bigr)
  =\frac{1}{d_{[\rho]}};
\end{equation}
that is, the information about the statistics
phase $\eta_{[\rho]}$ is lost. The full field theoretic statistics
operator $\varepsilon_{\rho}$ has instead the form
$$\varepsilon_{\rho}=\sum_{j,k,l=1}^{d_{[\rho]}}\Psi_j\hat\Psi_j^*
\Psi_k\hat\Psi_l\Psi_l^*\Psi_k^*,$$
where the $\hat\Psi_j\DEF u\Psi_j$
are an \ONB in the \HSP $H(\hat\rho)$ implementing the endomorphism
$\hat\rho$ (cf.\ \eqref{STOP} and \eqref{TRANS}). This expression is
obtained by substituting \rho in \eqref{STOP} by the formula
\eqref{IMP}, and by writing the unitary charge transporter $u$ in
terms of the $\Psi_j,\hat\Psi_j$ as
$u=\sum_j\hat\Psi_j\Psi_j^*$. Using the asymptotic \CRs of the fields,
which are in the present case of Bose or Fermi type
($\Psi_k\hat\Psi_l=\pm\hat\Psi_l\Psi_k$), one gets that
$\varepsilon_{\rho}=\pm\hat\varepsilon_{\rho}$. 
\medskip{}
  
These observations led \DR to the conjecture that the (finite
statistics) superselection structure described in \cite{DHR3}, and
valid in at least 3 space--time dimensions, should always be
equivalent to the representation theory of a unique compact group.
The proof of this conjecture was completed in the late 1980s
\cite{DR90} via an extension \cite{DR88,DR89a,DR89} of the
Tannaka--Krein duality theory of compact groups. The
\emph{Tannaka--Krein theory} allows to recover a compact group from
its ``concrete dual'', i.e.\ from the collection of
finite dimensional unitary representation spaces together with the
intertwiners between representations. The group elements can then be
identified with certain functions assigning to each representation
space a unitary operator on that space. \DR instead characterized the
\emph{abstract} duals of compact groups. They found that any category
which has essentially all the properties shared by the category of
localized endomorphisms, namely a composition law with permutation
symmetry, and the existence of subobjects, direct sums and conjugates,
is equivalent to a category of finite dimensional continuous unitary
representations of a unique compact group.

The construction of \FA and \GG from the observables and
localized endomorphisms now amounts to the construction of a concrete
group dual from an abstract one. The \FA can be obtained as
the ``cross product'' of \AA by the semigroup of localized
endomorphisms. This is a \CAL which contains \AA and, for each
endomorphism \rho, a finite dimensional \HSP $H(\rho)$ inducing
\rho, with certain relations between the elements of \AA and the
elements of the algebras generated by the $H(\rho)$. The \GG can be
identified with the compact group of all automorphisms of the field
algebra which leave \AA pointwise fixed. 

Summarizing, the main result of \DR states that, in Minkowski space of
dimension greater than two, the observable algebra can always be
embedded into a larger \FA on which a compact \GG acts in such a way
that the observables are precisely the gauge invariant fields. The
fields are local relative to the observables and act irreducibly on a
\HSP which contains each \SSS \rho with multiplicity $d_{[\rho]}$. The
charge quantum numbers are in one--to--one correspondence with the
equivalence classes of irreducible representations of the \GG. This
construction is unique up to unitary equivalence, if one requires that
fields commute or anticommute at spacelike separations (``normal
\CRs'').

Thus the program of reconstructing charged fields and gauge symmetries
from local observables has been carried through successfully, at least
in the case of strictly localizable charges (cf.\ \eqref{DHR}) in
space--time of dimension greater than two. This is of course a strong
confirmation of the central idea behind local \QFT, that all physical
information should be encoded in the relative position of the algebras
of local observables.
\medskip{}

The picture is however less complete in low dimensional \QFT. An
analysis of the low dimensional superselection structure based on the
DHR selection criterion \eqref{DHR} and on Haag duality \eqref{HD} in
the vacuum sector has been worked out by Fredenhagen, Rehren and
Schroer \cite{FRS1,FRS2}.  The methods of \DHR can be adapted to this
situation, and, as already indicated above, one finds a somewhat
richer structure in this case.  The statistics phase can be an
arbitrary element of the circle group \TT, and the statistics
dimension need not be an integer, but can take almost any value that
is allowed by Jones' list \cite{J} of indices of subfactors.
Statistics is governed by the braid group instead of the permutation
group, and the gauge symmetries (the ``dual object'' of the
superselection structure) do not form a group in general.
However, there are only partial results concerning the quantum
symmetry problem (e.g.\ \cite{MS92,FK,Re96,BNS,NSW}), the
classification of the occurring braid group representations (see e.g.\ 
\cite{FRS1,L2}), and the reconstruction of charged fields (e.g.\ the
``reduced field bundle'' \cite{FRS2}, a bounded version of the
conformal ``exchange algebras'' of Rehren and Schroer \cite{RSch}, or
\cite{VSch}). It is interesting to note that many structural
peculiarities such as braid group statistics, non--integer dimensions,
Verlinde's modular algebra etc., which had been found previously in
conformal field theory, and which were thought to be consequences of
conformal invariance, could be shown to be generic features of low
dimensional \QFT, independent of conformal invariance; see e.g.\ 
\cite{Re90,FRS2}.
\medskip{}

Whereas Haag duality automatically holds in conformally invariant
theories on the circle (due to space--time compactification, see
\cite{BSM,BGL,FG}), it is not such a reasonable assumption in
two--dimensional Minkowski space.  The basic mechanism for the
violation of Haag duality in the vacuum sector is the following: If
$O$ is a given double cone, then there exist operators in $\AA(O)'$
which create a charge in the left spacelike complement of $O$ and
annihilate a charge of the same type in the right spacelike
complement. Such operators cannot be approximated by observables in
$\AA(O')$, so that \eqref{HD} fails.

A preliminary analysis of the two--dimensional situation without
assuming Haag duality from the outset has recently been attempted by
M\"uger \cite{MM1}. Similar as \DHR in \cite{DHR1}, he starts from a
\FA with normal \CRs from which the observables are selected by a
gauge principle. The fields are assumed to satisfy a certain duality
property, which would imply Haag duality for the observables in higher
dimensions, but entails only a weaker form of duality (``essential
duality'') in two dimensions; and a specific form of causal
independence, which is believed to hold in massive theories. One can
then enlarge the local \FAs by adding certain non--local ``kink'' or
``disorder'' operators to the fields which act like the identity on
one half of the spacelike complement of some double cone, and like a
global gauge transformation on the other half. One can also introduce
a system of enlarged local observable algebras on the vacuum \HSP of
the original observables, which fulfills Haag duality and is therefore
called the ``dual net''\footnote{``Essential duality'' means that the
  so--defined dual net satisfies Einstein causality \cite{R74}:
  $\AA^d(O')\subset\AA^d(O)'$. Essential duality is known to hold if
  the local algebras are generated by Wightman fields \cite{BW}. The
  passage from a non--Haag dual theory to the dual net is the
  customary way of restoring Haag duality. Its value lies in the fact
  that, in higher dimensions and under the assumption of essential
  duality, the dual net possesses precisely the same \SSSs as the
  original theory \cite{R80}. This is however not true in two
  dimensions, so that the significance of the dual net is somewhat
  limited in this case. Under M{\"u}ger's assumptions, the dual net
  has in fact \emph{no} \SSSs fulfilling the DHR criterion. It seems
  that \SSSs of the original theory extend at best to soliton sectors
  of the dual net \cite{MM2}. --- Similar questions have been
  investigated for conformally invariant theories on the real line in
  \cite{GLW}.}, by setting
$$\AA^d(O)\DEF\AA(O')'.$$
It turns out that this dual net is just the
fixed point net of the enlarged \FA under the \GG $G$; and that, if
$G$ is finite, there is a natural action of a certain Hopf algebra
containing $G$ (the ``quantum double'' of $G$) on the enlarged \FA
such that the original observables are the fixed points under the
action of the whole quantum double. In this sense the violation of
Haag duality in two--dimensional Minkowski space is also related to
symmetry breaking.
\medskip{}

Let us finally comment on the status of the DHR selection criterion
\eqref{DHR}. As we already pointed out, all \SSSs which can be reached
by applying local fields to the vacuum, and all positive energy
representations in conformal field theory, fulfill this criterion. But
even in purely massive theories it is not true that all positive
energy representations can be localized in bounded space--time
regions. Buchholz and Fredenhagen proved that for a primary positive
energy representation $\pi$ whose energy--momentum spectrum starts
with an isolated mass shell, there is a unique vacuum representation
relative to which $\pi$ can be localized around ``semi--infinite
strings'' extending from one point to spacelike infinity \cite{BF}.
(In two dimensions there are possibly two inequivalent vacua
associated with $\pi$, and one can have \emph{soliton sectors}
interpolating between these vacua; see
\cite{Froe,Fre90,Fre93,Schl,Re97}.) One can perform an analysis of
\SSSs having this weaker localization property relative to a fixed
vacuum. Such an analysis is technically more involved than in the case
of the DHR criterion, but the resulting structure is very similar. In
particular, \SSSs can still be described with the help of
endomorphisms. Braid group statistics arises already in three
dimensions; see \cite{BF,DR90,FGM,MSch,FGR} for details.

Even weaker localization properties must be expected for charged
sectors in the presence of long range forces, e.g.\ in QED, where the
asymptotic direction of the electric flux at spacelike infinity has to
be taken into account. There is some hope that the localization of
charged states can be improved by comparing them with an
``infravacuum'' state (a certain radiation background) instead of a
vacuum state, so that the criterion of Buchholz and Fredenhagen would
apply \cite{Buch82,Walter}. Another mechanism which could make the
methods of superselection theory applicable to charges obeying Gau\ss'
law has recently been proposed in \cite{BDMRS}. These authors showed
on the example of a free massless scalar field that such charges can
be described by automorphisms \rho which violate the
Buchholz--Fredenhagen condition on \AA, but fulfill the stronger DHR
criterion \eqref{DHR} relative to a smaller subalgebra. Although one
can no longer conclude from Haag duality that the charge transporters
$u$ entering the definition \eqref{STOP} of the statistics operators
$\varepsilon_{\rho}$ are contained in \AA (cf.\ \eqref{TRANS}), it is
nevertheless possible to define $\rho(u)$ unambiguously. Then the
statistics operators are well--defined, and one can discuss the
statistics of the model in the usual way.  
\medskip{}

There are many other important topics in the theory of \SSSs that we
cannot touch upon in this introductory survey. Let us mention the deep
connection with the theory of \VN algebras, whose powerful tools as
e.g.\ the modular theory of Tomita and Takesaki (leading for instance
to a new approach to the just mentioned localization problem
\cite{BSch97a,BSch97b}) and the techniques of the theory of subfactors
have found a lot of applications in local quantum physics. It is even
true that some of these developments had been anticipated in physics
before they were established in greater generality in mathematics. As
these techniques will not be applied in our work, we refer to the
original literature \cite{L,L2,FRS1} and to the reviews
\cite{B95,HWW96,BSch96,BSch98}, and references quoted there, for these
matters.


%% file: w.tex
\specialsection*{What We Have Done}
\label{sec:overview}
\subsection*{Review of the perspective}
It should have become clear from the preceding excursion to
non--perturbative \QFT that endomorphisms of \CALs do play an
important r\^ole in quantum physics. Starting from a theory of local
observables, for which clear physical principles can be formulated,
localized endomorphisms are the basic tools for an intrinsic
construction of unobservable charged fields, for the derivation of
their localization properties, \CRs and symmetries. Since one has less
intuition concerning the physical properties of unobservable
quantities, it is gratifying that one does not have to postulate e.g.\ 
the \CRs of fields, but can deduce them from the principles of
locality and causality. It is equally remarkable that the gauge
symmetries can be derived from the interrelations between local
observables, which are by their very nature gauge invariant.

But, as it sometimes happens if one tries to put a physical theory on
a sound mathematical basis, it is difficult to identify the general
structures, whose existence is predicted by the abstract mathematical
analysis, in concrete models. This applies in particular to \QFT where
mathematically rigorous models beyond free fields are still lacking in
four--dimensional space--time. Let us have a look at some field
theoretical models whose superselection structure has been puzzled
out.

Localized automorphisms of the \emph{free charged Klein--Gordon
  field,} with \GG \TT, can be constructed as follows \cite{F73,BLOT}.
One smears the field with a smooth real test function which has
support in a compact region $O$ and whose Fourier transform does not
vanish identically on the positive mass shell. Then one obtains a
unitary operator by polar decomposition of the smeared field which
implements an automorphism localized in $O$ and which carries one unit
charge. (One cannot use quasi--free automorphisms of the \CCR for this
purpose because they are all neutral.)

Localized automorphisms of the \emph{free Majorana field,} with \GG
$\ZZ_2$, are even simpler to get. Here one can take the Majorana field
itself, smeared out with a suitable localized real test function, as
unitary implementer. This amounts to an especially simple choice of a
quasi--free automorphism (a reflection).

The situation is a bit different in the case of the \emph{free Dirac
  field}, with \GG \TT. Since the field operators are no longer
injective, one cannot build automorphisms and unitary implementers
directly out of the field operators. It is then natural to look for
localized automorphisms among the class of quasi--free automorphisms
of the \CAR. In \cite{CB0} we exhibited a family of charge--carrying
localized quasi--free automorphisms which are induced by certain
unitary multipliers on the single particle space. This construction
works however only in two dimensions, and it is unlikely that unitary
multipliers yield charge--carrying implementable automorphisms in
higher dimensions. Our construction has been generalized by admitting
kink--like multipliers by Adler \cite{Ad}.  The corresponding
automorphisms then show Abelian braid group statistics and extend to
solitons of the dual net.

The \emph{current algebra} derived from the massless free scalar field
has \emph{no} \SSSs in dimension greater than two, but exhibits
spontaneous symmetry breaking \cite{Str74,BDLR}. On the contrary, it
possesses a \emph{continuum} ($\cong\RR^2$) of \SSSs in two
dimensions, due to its peculiar infrared properties. The associated
localized automorphisms correspond to displacements of the fields
\cite{StrW}.  The superselection structure of the conformal current
algebra on the circle and of its local extensions (\GG $\ZZ_{2N}$)
has been studied in \cite{BMT}.

Of greater interest is the case of genuine endomorphisms,
corresponding to \emph{non--simple} sectors. The only explicit
examples of genuine localized endomorphisms constructed so far are, to
the best of our knowledge, the ones leading to the non--simple sectors
of the chiral conformal $so(N)$ WZW models at level one.
\label{page:WZW} They all belong to the class of quasi--free
endomorphisms of the \CAR. The first example appeared in the treatment
of the conformal Ising model by Mack and Schomerus \cite{MS}. This
model has one non--simple sector, with highest weight $\tfrac{1}{16}$
and statistics dimension $\sqrt{2}$, and Mack and Schomerus offered a
candidate of a localized endomorphism which was conjectured to
describe this sector.  They actually did not use this localized
endomorphism in their computations, but worked with a global
endomorphism throughout. This makes the analysis technically simpler,
but also somewhat questionable, because e.g.\ the concept of
statistics depends crucially on locality. The ideas of Mack and
Schomerus were soon generalized by Fuchs, Ganchev and Vecserny{\'e}s
to the level one $so(N)$ WZW models, which also have a Fermionic
realization \cite{FGV}.  These models possess one non--simple sector,
with highest weight $\frac{N}{16}$ and statistics dimension
$\sqrt{2}$, if $N$ is odd; the case $N=1$ reproduces the Ising model.
But these authors also used global endomorphisms.  This state of
affairs was subsequently improved by B{\"o}ckenhauer who constructed
localized endomorphisms, among them the candidate of Mack and
Schomerus, which are equivalent to the global ones of \cite{MS,FGV},
and which imply the same fusion rules \cite{JMB1,JMB2}.  Note that,
due to the non--integer statistics dimension, none of these
endomorphisms can have the properties predicted by \DR (cf.\ 
\eqref{HRHO}--\eqref{IMP}).  These endomorphisms are not related to
group symmetries, but to genuine quantum symmetries.

In this connection, we should also mention the attempts of
A.~Wassermann \cite{Was} and Recknagel \cite{Rec93,Rec96} to
substitute localized endomorphisms by certain other structures. In
\cite{Was} the fusion of positive energy representations of the loop
groups $LSU(N)$ is described. These representations can be constructed
using implementers of certain quasi--free automorphisms of the CAR and
\CCRs over $L^2(\TT)$ \cite{PrSe}, and are closely related to the
$su(N)$ WZW models. In \cite{Was}, their fusion is not performed with
the help of endomorphisms, but uses an equivalent technique, the
tensor product of bimodules over \VN algebras \cite{Con} (``Connes
fusion''). In \cite{Rec93} it is proposed to replace endomorphisms of
algebras by endomorphisms of some associated $K_0$--groups which are
in principle much easier to handle.  Though this heuristic approach is
plagued with some serious shortcomings, it was possible on its basis
to reproduce the fusion rules of the $su(2)$ WZW model. In
\cite{Rec96} it is tried to reach the sectors of some minimal models
by ``amplimorphisms'' of certain associated path algebras. A
characteristic feature of \cite{Rec93,Rec96} is the complete lack of
locality. (There is also the reformulation of the DHR theory given by
Fredenhagen, where representations and endomorphisms are replaced by
states and completely positive maps \cite{F92}. The usefulness of this
approach has been demonstrated on some subtheories of the conformal
current algebra \cite{F94}.)

Summing up, one is confronted with a scarcity of field theoretic
models whose localized endomorphisms are explicitly known, and there
are in fact \emph{no} known examples of endomorphisms which fit
completely into the scheme of \DR.
\subsection*{Quasi--free endomorphisms of CAR and \CCRs}
Ever since the invention of \QFT in the late 1920s \cite{Dir27,JW},
the CAR and \CCRs have been the dominating algebras in this field.  In
view of the above remarks it is natural to ask whether one can find
endomorphisms of these algebras which share all the properties
predicted by the theory of \DR (see the discussion after
Eq.~\eqref{HRHO}).

Specifically, the questions we are facing are the following: Do there
exist quasi--free endomorphisms of the CAR and \CCRs which can be
implemented on Fock space by \HSPs of isometries in the sense of
Eq.~\eqref{IMP}?  If yes, how can such endomorphisms be characterized?
What are their possible charge quantum numbers? And can one construct
the corresponding \HSPs of isometries (the ``charged fields'')
explicitly?

Very briefly, the answers implied by our work can be summarized as
follows. Each algebra possesses a rich semigroup of quasi--free
endomorphisms having the desired properties. These semigroups are the
natural generalizations of the well--known restricted orthogonal and
symplectic groups. A quasi--free endomorphism belongs to one of them
\IFF its associated one--particle operator fulfills a certain
Hilbert--Schmidt condition. There are detailed formulas for the
corresponding charged fields on Fock space which have a well--defined
meaning as infinite sums converging strongly on a dense domain. The
\HSPs of isometries spanned by these fields can in a natural way be
regarded as Fock spaces over some auxiliary space. (These auxiliary
spaces can have finite or infinite dimension. Be aware that these
``Fock spaces of isometries'' are not contained in the original Fock
space, but consist of operators acting on the latter.) This Fock space
structure is compatible with the action of the gauge symmetries, and
provides the key to the determination of the charge quantum numbers.
Genuine endomorphisms are always reducible; the possible values of
their statistics dimensions are the powers of 2 (CAR) resp.\ $\infty$
(CCR). They induce representations of the \CAs $\mathcal{O}_{2^n}$ and
$\mathcal{O}_\infty$. 
\medskip{}

For the convenience of the reader who finds the preceding remarks too
condensed we would like to give now a detailed exposition of the
material contained in the central chapter of this thesis.
\bigskip{}

\emph{Section~\ref{sec:CUN}.} Here we review the \CAs $\mathcal{O}(H)$
and their basic properties. This section is intended as a supplement
to the main text, and its content is not essential for an
understanding of the remainder.

After stating the definition of $\mathcal{O}(H)$, we quote Evans' Fock
space construction of $\mathcal{O}(H)$ as a quotient of the
Cuntz--Toeplitz algebra. Some general properties of $\mathcal{O}(H)$,
mostly due to Cuntz as e.g.\ its $K$--theory, are mentioned, and
endomorphisms of $\mathcal{O}(H)$ are discussed.  Quasi--free
endomorphisms and quasi--free states of $\mathcal{O}(H)$ are closely
related to the structure of the gauge invariant subalgebra of
$\mathcal{O}(H)$, and have been studied by Evans et al. Quasi--free
group actions on $\mathcal{O}(H)$ are an important element in the
theory of \DR.

There has recently been some interest, e.g.\ in connection with
Powers' $E_0$--semigroups, in the relation between representations of
$\mathcal{O}(H)$ on a \HSP \HH and endomorphisms of $\BB(\HH)$. Note
that the \CAs enter our analysis in exactly the same way: We study
endomorphisms of the CAR and \CCRs which give rise to representations
of the \CAs on Fock space via Eq.~\eqref{HRHO} and \eqref{IMP}. Thus
our results also provide interesting examples for the representation
theory of $\mathcal{O}(H)$.

Finally some remarks concerning the r\^ole of the \CAs in the general
theory of \CALs are made, but there is no room to discuss these
important topics in greater detail.
\bigskip{}

\emph{Section~\ref{sec:CAR}.} This section contains a thorough
analysis of the semigroup of all quasi--free endomorphisms of the \CAR
which can be implemented by \HSPs of isometries in a fixed Fock
representation. Essentially all results which go beyond the case of
automorphisms are new. Most of them have already been published in
\cite{CB1,CB2}, but the presentation given here is in several respects
superior to the one in \cite{CB1}.  This analysis is completely
general in that we do not make specific assumptions on the structure
of the real \HSP underlying the \CAR. The price that one has to pay
for this generality is that the concept of locality is not
incorporated at this level, but has to be discussed separately in a
more restrictive setting. As a consequence, we cannot apply the
methods of the theory of \SSSs which hinge upon the principle of
locality. Our methods are instead taken from the representation theory
of the \CAR and from general functional analysis (e.g.\ Fredholm
theory), and are largely independent of the Doplicher--Roberts theory.
(Analogous remarks are valid for the treatment of the CCR case in
Section~\ref{sec:CCR}.)  \medskip

\emph{Section~\ref{sec:QFCAR}.} The basic objects and facts that will
be needed later on are introduced. Araki's ``selfdual'' \CAR formalism
is used throughout which amounts to complexification of the underlying
real \HSP. Quasi--free endomorphisms are in one--to--one
correspondence with their restrictions (``Bogoliubov operators'') to
this space. Bogoliubov operators are isometric. The fundamental
invariant of a Bogoliubov operator $V$ is its Fredholm index, a
non--positive integer (or $\infty$), and we find that this index is
related to the statistics dimension $d_V$ of the corresponding
endomorphism $\rho_V$ by the formula
\begin{equation}
  \label{INDFOR}
  d_V=2^{-\3\IND V}.
\end{equation}
Here the statistics dimension $d_V$ is defined to be the square root
of the Watatani index of the range of $\rho_V$.  This is a purely
{\CAL}ic notion which does not depend on representations.
Automorphisms are characterized by $d_V=1$. By a somewhat surprising
result in \cite{CB2}, any $\rho_V$ with statistics dimension
$d_V=\sqrt{2}$ induces in a canonical way an isomorphism from the \CAR
onto its even subalgebra. These isomorphisms can be used to study the
even subalgebra, which models the algebra of observables in various
physical systems.

The class of quasi--free states is defined, and the main technical
tool to be used later, the quasi--equivalence criterion for
quasi--free states of Powers and St{\o}rmer and Araki, is stated. Fock
states are the pure quasi--free states. The Powers--St{\o}rmer--Araki
criterion can be simplified if one of the states involved is a Fock
state. This has been observed by Powers, but we arrived independently
at the same conclusion, by an argument which can be found in the
preprint version of \cite{CB1}.

Associated with every Fock representation is a second (equivalent)
representation which we call the ``twisted Fock representation''. The
twisted Fock representation provides a convenient way to describe
``twisted duality'', the analogue of Haag duality in the presence of
Fermi fields, and is always useful if one has to deal with commutants
of ``local'' subalgebras. The consequent use of the twisted Fock
representation will lead to some simplifications in the construction
of the charged fields in Section~\ref{sec:FORM}.
\medskip

\emph{Section~\ref{sec:REP}.} This section is concerned with the
representations $\pi\0\rho$ that are obtained by composing a Fock
representation with a quasi--free endomorphism. (Recall from
p.~\pageref{page:LOCEND} that the representations occurring in the
theory of \SSSs have such a form.) We start with a discussion of the
implementation problem for endomorphisms of arbitrary \CALs. The
conclusion to be drawn from this general discussion is that an
endomorphism \rho is implementable in an irreducible representation
$\pi$ \IFF the representations $\pi$ and $\pi\0\rho$ are
quasi--equivalent.

If $\pi$ is the Fock representation of the \CAR induced by a \FS
$\omega$, and if \rho is a quasi--free endomorphism, we show that
$\pi\0\rho$ is a multiple of the GNS representation of the quasi--free
state $\omega\0\rho$. (We do actually a little bit more because we
give an explicit decomposition of Fock space into invariant subspaces.
This detailed description will be used in Section~\ref{sec:FORM} to
prove the completeness relation \eqref{CUNTZb} for the charged
fields.) The multiplicity is some power of two or infinite.  The
question of implementability is thereby reduced to the question of
quasi--equivalence of quasi--free states, and we can derive our basic
implementability condition from the Powers--St{\o}rmer--Araki
criterion. This condition generalizes the well--known
Shale--Stinespring condition which is restricted to the case of
automorphisms. Both conditions are formally the same. The Fredholm
index of the Bogoliubov operator corresponding to an implementable
endomorphism is always even, so that the statistics dimensions of
implementable endomorphisms are integers (or $\infty$) as they should
(cf.\ \eqref{DIMHR}).

The second half of this section deals with representations $\pi\0\rho$
where $\pi$ is a given Fock representation, but \rho an arbitrary
(non--implementable) quasi--free endomorphism with finite index. (The
analysis would be trivial in the implementable case, because
$\pi\0\rho$ is then unitarily equivalent to $d_V\cdot\pi$.) This
analysis is e.g.\ relevant for the endomorphisms describing the
non--simple sectors of the WZW models (cf.\ the remarks made on
p.~\pageref{page:WZW}). We derive criteria for unitary equivalence of
two such representations, describe the quasi--free states of the form
$\omega\0\rho$ where $\omega$ is the Fock state corresponding to
$\pi$, and characterize the endomorphisms \rho for which the states
$\omega\0\rho$ are pure or ``almost pure'' (mixtures of two
inequivalent pure states. This is in some sense the best one can get
for endomorphisms $\rho_V$ with $\IND V$ odd.) These preparatory
results are then applied to give an alternative proof a theorem of
B\"ockenhauer, namely that $\pi\0\rho_V$ is a multiple of another Fock
representation, with multiplicity $d_V$, if $\IND V$ is even; and that
it is a multiple of two ``pseudo Fock representations'', with
multiplicity $d_V/\sqrt{2}$, if $\IND V$ is odd. Invoking our
isomorphism onto the even subalgebra from Section~\ref{sec:QFCAR}, we
obtain analogous results for the restrictions of these representations
to the even subalgebra.  It was observed by Szlach\'anyi and
B\"ockenhauer, but regarded as a curiosity, that the restriction of
$\pi\0\rho_V$ to the even subalgebra behaves like a representation
$\pi\0\rho_{V'}$ of the whole \CAR, where $V'$ is a Bogoliubov
operator with $\IND V'=\IND V-1$. Our approach gives a natural
explanation for this phenomenon.  \medskip

\emph{Section~\ref{sec:IP}.} Having established the necessary and
sufficient condition for implementability in Section~\ref{sec:REP}, we
study here the structure of the topological semigroup of all
quasi--free endomorphisms which can be implemented in a fixed Fock
representation. This semigroup is an extension of the restricted
orthogonal group. One of our achievements is the proof that this
semigroup can be written as a product of a small subgroup consisting
of automorphisms which are close to the identity, and the
sub--semigroup of endomorphisms which leave the given \FS invariant.
What is more, we are able to make a definite choice of the two factors
in which an implementable endomorphism decomposes.

There are some results involved in the proof of this product
decomposition which are of independent interest. The first is a useful
parameterization of the class of all Fock states which are unitarily
equivalent to the given one. This parameterization is done in terms of
certain pairs $(T,\mathfrak{h})$ consisting of an antisymmetric
Hilbert--Schmidt operator $T$ and a finite dimensional subspace
$\mathfrak{h}$ of the kernel of $T$, and is adapted to the structure
of the cyclic vectors in Fock space which induce the states. The next
result is a canonical (up to a finite dimensional part related to
$\mathfrak{h}$) choice of a quasi--free automorphism belonging to the
small subgroup mentioned above which transforms a given equivalent \FS
into the original one.

Associated with any quasi--free endomorphism $\rho_V$ there is a
``partial'' \FS, viz.\ a \FS of the \CAR over the range of $V$. Using
the above parameterization, we can extend this partial \FS (in the
implementable case) to a proper \FS, say $\omega_V$, of the whole
\CAR.  (This procedure is reminiscent of the construction of the
conjugate sector, with the help of the left inverse, in \QFT; cf.\ 
p.~\pageref{page:CONJ}.) The choice of $\omega_V$ is made definite by
minimizing both ``parameters'' $T$ and $\mathfrak{h}$ in an
appropriate way. We can thus assign, in an unambiguous way, to each
implementable quasi--free endomorphism $\rho_V$ a \FS $\omega_V$ which
is equivalent to the original \FS $\omega$, and such that
$\omega_V\0\rho_V=\omega$. By the above, we get in addition a
quasi--free automorphism $\rho_U$ in the small subgroup which also has
the property that $\omega_V\0\rho_U=\omega$. The announced product
decomposition of $\rho_V$ is then obtained by setting
\begin{equation}
  \label{PRODEC}
  \rho_V=\rho_U\rho_W,\qquad W\DEF U^{-1}V.  
\end{equation}
There is an ambiguity in the definition of $\rho_U$ which cannot be
resolved. It amounts to the freedom in the choice of an \ONB in the
``minimized'' space $\mathfrak{h}_V$ associated with the \FS
$\omega_V$.

The fact that any implementable quasi--free endomorphism can be
written as a product of two very simple factors is then used to
determine the connected components of the semigroup of implementable
endomorphisms. It is well--known that the restricted orthogonal group
(the group of implementable automorphisms) has two components which
are distinguished by the Araki--Evans index, i.e.\ by the parity of
the dimension of the space $\mathfrak{h}_V$. We find that the
Araki--Evans index is not an invariant of genuine endomorphisms, and
that each set of endomorphisms $\rho_V$ with fixed nonzero value of
$\IND V$ is connected.  The product decomposition of endomorphisms
will further be used in Section~\ref{sec:FORM} to reduce the proof of
the completeness relation \eqref{CUNTZb} for the charged fields to the
simpler case of endomorphisms which leave $\omega$ invariant.
\medskip

\emph{Section~\ref{sec:FORM}.} We show in this section how to
construct an \ONB in the \HSP of isometries $H(\rho_V)$ which
implements a given quasi--free endomorphism $\rho_V$ on Fock space.
The implementers can be expressed in terms of annihilation and
creation operators.  This will make it necessary to employ special
Fock space techniques, which have been avoided in the previous
sections.  The strategy underlying this construction is the following:
We first generalize the known methods of constructing unitary
implementers for automorphisms to the case of genuine
endomorphisms\footnote{It is well--known that an automorphism $\rho_U$
  is implementable \IFF there exists a unit vector $\Omega_U$ in Fock
  space which lies in the joint kernel of all transformed annihilation
  operators. This vector $\Omega_U$ induces the state
  $\omega\0{\rho_U}^{-1}$.  Once $\Omega_U$ is known, the unitary
  implementer \PU for $\rho_U$ can be constructed essentially by
  setting $\PU\pi(a)\Omega=\pi(\rho_U(a))\Omega_U$, where $\pi$
  denotes the Fock representation and $\Omega$ the original Fock
  vacuum vector. In the case of endomorphisms, the state $\omega_V$
  plays the r\^ole of $\omega\0{\rho_U}^{-1}$. But note that the
  cyclic vector $\Omega_V$ associated with $\omega_V$ is no longer
  uniquely determined by the requirement that it be destroyed by the
  transformed annihilation operators. In fact, any element of the \HSP
  $H(\rho_V)\Omega$ has the latter property. This corresponds to the
  fact that $\omega_V$ is not the only possible extension of the
  partial \FS mentioned above.  What remains true is that any
  implementer is uniquely determined by its value on $\Omega$. Also
  note that, given $\omega_V$, it is still not obvious that the
  construction of \PV{0} goes through in the known way: In the case of
  automorphisms, both relations $U^*U=\1$ and $UU^*=\1$ usually enter
  the construction on the same footing;
  whereas in the case of endomorphisms, the second relation is lost,
  and one has to take this into account with care.}. This
generalization makes essential use of the \FS $\omega_V$ introduced in
Section~\ref{sec:IP} and permits us to define \emph{one} isometric
implementer \PV{0} for $\rho_V$. \PV{0} is characterized by the
property that its value on the Fock vacuum $\Omega$ reproduces the
cyclic vector inducing the state $\omega_V$. We then construct a
complete set of implementers by multiplying \PV{0} with suitable
partial isometries from the left.  This approach is suggested by the
observation that, if one has an \ONB \PV{\alpha} in $H(\rho_V)$, then
the operators $\PV{\alpha}\PV{0}^*$ are partial isometries in the
commutant of the range of $\rho_V$, and the \PV{\alpha} can be
reconstructed from these operators (together with \PV{0}) by setting
\begin{equation}
  \label{PARIS}
  \PV{\alpha}\DEF\Bigl(\PV{\alpha}\PV{0}^*\Bigr)\PV{0}.  
\end{equation}

The first step in this program is to give an appropriate definition
of ``\BHs''. This will be achieved by combining the algebraic approach
of Araki with the more analytic approach of Ruijsenaars. As is
well--known, the \CAR contains the spin group, the universal covering
group of the connected component of the identity of the group of
Bogoliubov operators which induce inner quasi--free automorphisms,
together with its Lie algebra. The elements of this Lie algebra are
called ``\BHs'' because they are bilinear expressions in the
generators of the \CAR, and can be identified with certain trace class
operators $H$ on the underlying \HSP. But we need a more general
definition of \BHs if we want to cover the case of general
implementable transformations. A natural extension of this Lie algebra
is the current algebra, the Lie algebra of skew--adjoint operators on
Fock space whose exponentials implement one--parameter groups of
quasi--free automorphisms. The current algebra can be identified with
a larger class of operators $H$, namely with the Lie algebra of the
restricted orthogonal group (if one allows for the occurrence of
Schwinger terms). But even in the case of quasi--free automorphisms
implementers can in general not be obtained as exponentials of these
currents, so that this class is still too narrow.

The way out is to consider \emph{Wick ordered} exponentials. Such Wick
ordered exponentials can be defined, a priori as \QFs on a dense
domain in Fock space, for Wick ordered \BHs induced by arbitrary
bounded operators $H$. Under a certain Hilbert--Schmidt condition on
$H$, these Wick ordered exponentials are the \QFs of densely defined,
in general unbounded, operators. The \CRs of these operators with
creation and annihilation operators can be explicitly computed, and
are used to determine all operators $H$ such that the Wick ordered
exponential of the \BH induced by $H$ has the ``correct'' intertwining
properties relative to a given quasi--free endomorphism $\rho_V$.  The
operators $H$ with this property are in one--to--one correspondence
with certain operators $T$ obtained during the study of the state
extension problem in Section~\ref{sec:IP}. The minimal choice $T_V$
made in that section leads then to a unique operator $H_V$ associated
with $\rho_V$, and the isometric implementer \PV{0} is obtained as a
finite sum of terms, each involving the Wick ordered exponential of
the \BH induced by $H_V$ plus some additional operators which
essentially fill up the ``Dirac sea'' corresponding to the finite
dimensional subspace $\mathfrak{h}_V$.

A complete \ONB in $H(\rho_V)$ is then constructed from \PV{0} by the
method outlined above. Here it is used that the commutant of the range
of $\rho_V$ can be easily described with the help of the twisted Fock
representation. It is also straightforward to obtain partial
isometries in this commutant since any Fermionic creation or
annihilation operator is already (a multiple of) a partial isometry,
by Pauli's principle. It is less obvious how to characterize partial
isometries which contain the range of \PV{0} in their initial spaces,
as is required by \eqref{PARIS}. This can be done with the help of the
\FS $\omega_V$ and its associated ``parameter'' $T_V$.

We finally arrive at the following scenario. There is a certain
subspace $\mathfrak{k}_V$ of the kernel of $V^*$ which has dimension
$-\2\IND V$. We choose an \ONB in $\mathfrak{k}_V$. The representors
of these basis vectors in the twisted Fock representation then behave
like creation operators relative to the implementer \PV{0}. That is,
we obtain an \ONB in $H(\rho_V)$ by multiplying \PV{0} from the left
with all possible ordered monomials in these operators. It follows
that $H(\rho_V)$ is isomorphic to the antisymmetric Fock space over
$\mathfrak{k}_V$, so that 
$$\dim H(\rho_V)=d_V$$
(cf.\ \eqref{INDFOR}). The basis of
implementers is chosen in such a way that the value of any one of them
on the Fock vacuum $\Omega$ is the state vector of some \FS. Moreover,
the choice of implementers is compatible with the product
decomposition \eqref{PRODEC}. Roughly speaking, the factor $\rho_U$
carries the exponential term plus the operators corresponding to
$\mathfrak{h}_V$, whereas $\rho_W$ is responsible for the additional
partial isometries.
The completeness of implementers \eqref{CUNTZb} is proved by showing
that the ranges of the implementers of the factor $\rho_W$ are equal
to the invariant subspaces which appeared in the decomposition of the
representation $\pi\0\rho_W$ in Section~\ref{sec:REP}.

Let us finally remark that the formulas given for the implementers are
not in ``normal form'' in the strict sense, i.e.\ they are not
completely Wick ordered. There are two reasons for this. The first is
the use of the twisted Fock representation, which involves the second
quantization of $-\1$ as a factor. These factors could be avoided by
incorporating them into the \BH of $H_V$, but the formulas would
become less transparent then, and the combinatorics would be more
complicated. The second reason is that the additional partial
isometries with which \PV{0} is multiplied contain an annihilation
part which should be moved to the right, with the help of the \CRs, in
order to get a completely Wick ordered expression. We have not done so
because the Fock space structure of $H(\rho_V)$ would then no longer
be visible.
\medskip

\emph{Section~\ref{sec:BOSSTAT}.} Here we derive formulas for the
Bosonized statistics operators of quasi--free endomorphisms with
finite statistics.  The basic observation is that partial isometries
of the form $\PV{\alpha}\PV{\beta}^*$ have an explicit representation
as monomials in the basis vectors of $\mathfrak{k}_V$. Special
examples are the operators $\PV{\alpha}\PV{0}^*$ appearing in
\eqref{PARIS}, and the operators $\PV{\alpha}\PV{\alpha}^*$, the
projections onto the ranges of the \PV{\alpha}.

Recall from \eqref{BSO0} that the Bosonized statistics operator is a
certain polynomial in the implementers of the endomorphism. The
knowledge of the operators $\PV{\alpha}\PV{\beta}^*$ and of the
intertwiner properties of the \PV{\alpha} suffices to identify this
polynomial with an element of the even \CAR. As a consistency check,
we compute the ``Bosonized statistics parameter'' by applying the
{\CAL}ic left inverses that were introduced in Section~\ref{sec:QFCAR}
to the Bosonized statistics operators. The result is the inverse of
the statistics dimension, in accordance with \eqref{BSP}.
\bigskip{}

\emph{Section~\ref{sec:CCR}.} This section contains an analogous
analysis of quasi--free endomorphisms of the CCR (or Weyl) algebra. It
is essentially based on \cite{CB3}. The general remarks made above in
the introduction to Section~\ref{sec:CAR} apply here as well. The
methods are now borrowed from the representation theory of the \CCR,
and the theory of Fredholm operators will again be of importance. The
survey of Section~\ref{sec:CCR} will be comparatively short, and we
will try to emphasize the differences between the CAR and CCR cases.
\medskip

\emph{Section~\ref{sec:SDCCR}.} The basic notions are established.  We
find it again convenient to use Araki's ``selfdual'' formulation. In
the CCR case it is necessary to start with a distinguished
``reference'' \FS $\omega$ in order to get a \HSP topology on the
underlying symplectic space. A certain dichotomy will arise in the
following from the fact that the algebraic relations are dictated by
the symplectic form, whereas the analytic aspects refer to the \HSP
inner product. Relevant topics discussed here are Araki's duality
relation and a statement about the affiliation of sums of creation and
annihilation operators to ``local'' Weyl algebras.
\medskip

\emph{Section~\ref{sec:IMP}.} Quasi--free endomorphisms are
introduced. They are given by Bogoliubov operators $V$ acting on the
symplectic space. Bogoliubov operators preserve the symplectic form
and have well--defined Fredholm indices. In contrast to the CAR case,
$\IND V$ is always even, irrespective of implementability.

As in the CAR case, we reduce the question of implementability to the
question of quasi--equivalence of quasi--free states, by showing that
the representation induced by a quasi--free endomorphism \rho in the
given Fock representation is a multiple of the GNS representation of
the state $\omega\0\rho$. The multiplicity is either $1$ or $\infty$.
Invariant subspaces are explicitly described.

The derivation of the necessary and sufficient condition for
implementability is based on the criterion for quasi--equivalence of
quasi--free states in the form given by Araki and Yamagami. Some work
has to be done to get rid of the square roots appearing in this
criterion. We do this with the help of an inequality of Araki and
Yamagami; this inequality enables us to reduce the problem to the CAR
case via polar decomposition of $V$, because the isometric part of $V$
is a CAR Bogoliubov operator. The resulting condition is a
generalization of the well--known condition of Shale which covers the
case of automorphisms.  The two conditions do \emph{not} have the same
form, in contrast to the CAR case. The statistics dimensions are now
given by
\begin{equation}
  \label{ccr:DV}
  d_V=
  \begin{cases}
    1, & \IND V=0,\\
    \infty, & \IND V\neq 0.
  \end{cases}
\end{equation}
The Fredholm index is therefore a finer invariant than the algebraic
index.
\medskip

\emph{Section~\ref{sec:SG}.} We study the semigroup of implementable
quasi--free endomorphisms, an extension of the restricted symplectic
group. Again we aim at showing that this semigroup can be written as a
product of a subgroup of automorphisms close to the identity and a
sub--semigroup of endomorphisms which leave the \FS $\omega$
unchanged. To this end we consider the set of \FSs which are
equivalent to $\omega$. We can parameterize this set similar to the
CAR case. There is however no counterpart of the spaces\footnote{The
  canonical anti\CRs are symmetric in creation and annihilation
  operators, but the canonical \CRs are not. Thus there exist
  endomorphisms of the \CAR which interchange creation and
  annihilation operators, and this is the origin of the spaces
  $\mathfrak{h}$. This possibility is absent in the CCR case.}
$\mathfrak{h}$, and the only parameter that is needed is an element
$Z$ of the infinite dimensional open unit disk. This parameter $Z$
characterizes the cyclic vector in Fock space which induces the
corresponding state, and has properties similar to the operator $T$
occurring in the Fermionic case. For any \FS equivalent to $\omega$
there is a canonical choice of an automorphism in the small subgroup
which transforms this state into $\omega$. Note that, with
$\mathfrak{h}$, also the ambiguity in the choice of this automorphism
has disappeared.

To obtain the product decomposition of a quasi--free endomorphism
$\rho_V$ as in \eqref{PRODEC}, we must again solve the problem of
extending a certain ``partial'' \FS associated with $\rho_V$ to a
proper \FS $\omega_V$. Recall that this problem was solved in the CAR
case by ``minimizing'' the parameters $T$ and $\mathfrak{h}$ in some
sense. In particular, it is possible to define the operator $T_V$ as a
function of (the components of) $V$.  However, 
we could not find a similar prescription for the operators $Z$ (and
presumably such a prescription does not exist). This complication is
caused by the fact that $Z$ has to fulfill an additional requirement
related to the positivity of the state, viz.\ its norm has to be
smaller than one. (There is no such restriction on the operators $T$;
the admissible $T$ form in fact a \HSP.) Instead we discovered a
canonical method, based on spectral theory, how to extend the partial
state \emph{directly,} i.e.\ without having recourse to the parameter
$Z$. This state extension $\omega_V$ is then used to \emph{define} the
parameter $Z_V$; remember that $Z_V$ will be needed later for the
construction of implementers. Having assigned a \FS $\omega_V$ to
$\rho_V$ so that $\omega_V\0\rho_V=\omega$ holds, there is then an
unambiguous choice of an automorphism $\rho_U$ in the small subgroup
such that $\omega_V\0\rho_U=\omega$, and the desired product
decomposition is finally obtained as in \eqref{PRODEC}.

As a corollary, we determine the connected components of the
semigroup. It turns out that any subset of endomorphisms $\rho_V$ with
$\IND V$ constant is connected.
\medskip

\emph{Section~\ref{sec:CON}.} The construction of a complete set of
implementers for a given endomorphism $\rho_V$ is performed. We start
by defining Wick ordered Bosonic \BHs, and Wick ordered exponentials
thereof, on Fock space. These are in general \QFs, but determine
densely defined operators under some conditions on the associated
operators $H$. One can again compute \CRs of Wick ordered exponentials
with creation and annihilation operators, and select the operators $H$
with the property that the corresponding Wick ordered exponential
fulfills appropriate intertwiner relations \WRT $\rho_V$. These
operators $H$ are in one--to--one correspondence with the operators
$Z$ parameterizing the \FSs which solve the extension problem from
Section~\ref{sec:SG}.

One implementer \PV{0} is then obtained as the (normalized) Wick
ordered exponential of the \BH induced by the operator $H_V$ which
corresponds to $Z_V$. The value of \PV{0} on the Fock vacuum $\Omega$
is the cyclic vector associated with the state $\omega_V$. To get a
complete set of implementers, we choose a certain basis in a subspace
$\mathfrak{k}_V$ of the kernel of the symplectic adjoint of $V$. The
dimension of this subspace is $-\2\IND V$. Polar decomposition of the
representors of these basis elements yields a set of isometries which
commute with each other and with the elements of the range of
$\rho_V$. One can then show that the operators \PV{\alpha} obtained by
multiplying \PV{0} from the left with all possible ordered monomials
in these isometries satisfy the relations of (an essential
representation of) the \CA $\mathcal{O}_\infty$. Since these
isometries behave like (isometric parts of) creation operators \WRT
the ``vacuum'' \PV{0}, one finds that the \HSP $H(\rho_V)$ is
canonically isomorphic to the symmetric Fock space over
$\mathfrak{k}_V$. 

Our choice of implementers is again compatible with the product
decomposition $\rho_V=\rho_U\rho_W$. Roughly, the factor $\rho_U$ is
responsible for the Wick ordered exponential, and the factor $\rho_W$
for the additional isometries. The completeness of implementers
follows from the fact that the ranges of the implementers of the
factor $\rho_W$ coincide with the invariant subspaces for the
representation $\pi\0\rho_W$ described in Section~\ref{sec:IMP}. 
The expressions for the implementers are again not completely Wick
ordered, because the additional isometries will in general contain an
annihilation part. But strict Wick ordering would hide the inherent
Fock space structure of $H(\rho_V)$.  \bigskip{}

\emph{Section~\ref{sec:SECTORS}.} The general theory of the
implementation of quasi--free endomorphisms of the CAR and \CCRs has
been completely developed in Sections~\ref{sec:CAR} and \ref{sec:CCR}.
The detailed knowledge of the structure of the implementing \HSPs will
now be used to gain insight into the charge structure of quasi--free
endomorphisms and to determine the possible charge quantum numbers.

The setting will be tailored to the situation implied by the
Doplicher--Roberts theory. That is, we will consider the CAR and \CCRs
as \FAs which contain the observables as fixed points under the action
of a given global \GG $G$. The \GG will be assumed to consist of
quasi--free automorphisms which leave a fixed \FS, the vacuum state of
the \FA, invariant. Therefore $G$ acts by usual second quantization on
Fock space. 

As a consequence, quasi--free endomorphisms are a priori endomorphisms
of the \FA. (Note that by the results of \DR, localized endomorphisms
of the observable algebra have a natural extension to the \FA in terms
of their implementers, by replacing the observable $a$ in \eqref{IMP}
by elements of the \FA.) Thus we have to single out a subset of
quasi--free endomorphisms which restrict to endomorphisms of the gauge
invariant subalgebra.

The relevant subset is the semigroup of \emph{gauge invariant
  endomorphisms,} i.e.\ of endomorphisms \rho commuting with $G$,
because these are precisely the endomorphisms whose implementing \HSPs
$H(\rho)$ carry a representation of $G$. By the discussion following
Eq.~\eqref{HRHO}, the determination of the charge quantum numbers of
\rho is equivalent to the determination of the representation of $G$ on
$H(\rho)$. This representation is further equivalent to the
representation of $G$ on the \HSP $H(\rho)\Omega$, which is easier to
handle than $H(\rho)$ itself (here $\Omega$ denotes the Fock vacuum
vector).

It should be noted that our assumptions are satisfied in models like
the $N$--component Dirac field with \GG $U(N)$.
\medskip

\emph{Section~\ref{sec:CARCH}.} We compute the charge quantum numbers
of gauge invariant quasi--free endomorphisms $\rho_V$ of the \CAR. It
turns out that they are essentially determined by the subspaces
$\mathfrak{h}_V$ and $\mathfrak{k}_V$ introduced in
Section~\ref{sec:CAR}.

We have to study the behaviour of the implementers under gauge
transformations. The values of the implementers of $\rho_V$ on the
vacuum vector $\Omega$ have the following structure: To the left
stands a product of partial isometries associated with the subspace
$\mathfrak{k}_V$, followed by the ``filled Dirac sea'' corresponding
to the finite dimensional subspace $\mathfrak{h}_V$, and finally the
``pure creation part'' of the Wick ordered exponential of the \BH of
$H_V$, applied to $\Omega$. We show that the subspaces
$\mathfrak{h}_V$ and $\mathfrak{k}_V$ are representation spaces of
$G$, and that the operators related to these subspaces transform
linearly under $G$. The exponential term on the other hand is
invariant. This follows from the fact that the operators $T_V$ can be
expressed as a function of the components of $V$, and confirms that
the ``minimal'' choice of $T_V$ made in Section~\ref{sec:IP} is a
reasonable one.

The transformation law of implementers implies that the representation
$\mathcal{U}_V$ of $G$ on $H(\rho_V)$ has the form 
\begin{equation}
  \label{REPGG}
  \mathcal{U}_V\simeq{\det}_{\mathfrak{h}_V}\otimes
  \Lambda_{\mathfrak{k}_V}.
\end{equation}
Here ${\det}_{\mathfrak{h}_V}$ is the
one--dimensional representation obtained by taking the determinant on
$\mathfrak{h}_V$ of the Bogoliubov operators in $G$, and
$\Lambda_{\mathfrak{k}_V}$ is the $d_V$--dimensional representation of
$G$ on the antisymmetric Fock space over $\mathfrak{k}_V$. (Recall
that $H(\rho_V)$  is isomorphic to this Fock space.)

It follows that $\mathcal{U}_V$, and hence $\rho_V$, is reducible if
$\IND V\neq0$, because $\Lambda_{\mathfrak{k}_V}$ is reducible. The
representation $\Lambda_{\mathfrak{k}_V}$ contains together with the
representation on $\mathfrak{k}_V$ all the higher antisymmetric tensor
powers thereof. The ``least reducible'' case is obtained if
$\mathfrak{k}_V$ is irreducible.

A special case worth mentioning is the case $G=\TT$ and $\IND V=0$,
i.e.\ the case of the \emph{restricted unitary group.} It is
well--known from the work on the external field problem that the
charge of elements of the restricted unitary group is given by a
certain Fredholm index $\IND V_{++}$ (which has nothing to do with the
index of $V$). This fact can be easily derived from our much more
general result: The factor $\Lambda_{\mathfrak{k}_V}$ in \eqref{REPGG}
becomes trivial, and the factor ${\det}_{\mathfrak{h}_V}$ yields the
index of $V_{++}$.

We study the question which representations of $G$ can possibly occur
on the subspaces $\mathfrak{h}_V$ and $\mathfrak{k}_V$. In typical
cases, any representation of $G$ that is realized on Fock space
appears as a subrepresentation on some $H(\rho_V)$.  Then we compare
our findings with the generic superselection structure of \QFT. The
semigroup of gauge invariant endomorphisms is not closed under taking
subobjects or direct sums. It is closed under taking conjugates if one
makes natural assumptions on the action of $G$.  Under these
assumptions, one can assign to each gauge invariant Bogoliubov
operator $V$ another such operator $V^c$ such that
$\mathfrak{h}_{V^c}$ and $\mathfrak{k}_{V^c}$ are antiunitarily
equivalent to $\mathfrak{h}_V$ and $\mathfrak{k}_V$, so that the
representation $\mathcal{U}_{V^c}$ is unitarily equivalent to the
complex conjugate representation of $\mathcal{U}_V$.

Finally we give an explicit example of a localized implementable gauge
invariant endomorphism with statistics dimension $2^N$ of the free
massless $N$--component Dirac field in two dimensions. The
construction rests on the use of ``local'' Fourier bases for the
chiral components, and is in this respect similar to the known
examples of localized endomorphisms in conformal field theory.
\medskip

\emph{Section~\ref{sec:CCRCH}.}  Here we compute the charge quantum
numbers of gauge invariant quasi--free endomorphisms $\rho_V$ of the
\CCR. Since the spaces $\mathfrak{h}_V$ are absent in the CCR case,
the charge quantum numbers are entirely determined by the
representation of $G$ on $\mathfrak{k}_V$.

We derive the transformation law of the implementers of $\rho_V$.
Recall that the implementers are obtained by multiplying the
distinguished implementer \PV{0} from the left with certain isometries
which are associated with the space $\mathfrak{k}_V$. \PV{0} itself is
a Wick ordered exponential of a \BH which is characterized by the
operator $Z_V$. Similar to the CAR case, this Wick ordered exponential
is gauge invariant, and we take this as a confirmation that we gave
the ``correct'' definition of $Z_V$ in Section~\ref{sec:SG}. Though
there is no explicit formula for $Z_V$ in terms of $V$, it is still
true that $Z_V$ is in some sense a function of $V$, and this suffices
to prove the invariance of \PV{0}.

The subspace $\mathfrak{k}_V$ is again $G$--invariant. But in contrast
to the CAR case, the isometries associated with $\mathfrak{k}_V$ do
\emph{not} transform linearly under $G$. One can however show that
they obey a linear transform law when restricted to the range of
\PV{0}, and that is essentially all we need.

We conclude that the representation $\mathcal{U}_V$ of the \GG $G$ on
the \HSP $H(\rho_V)$ is unitarily equivalent to the representation on
the symmetric Fock space over $\mathfrak{k}_V$ induced by the
representation on $\mathfrak{k}_V$. Thus one can say that the
endomorphisms in the CCR case are even ``more'' reducible than in the
CAR case, because $\mathcal{U}_V$ (and hence $\rho_V$) always splits
into an \emph{infinite} direct sum of irreducibles if $\IND V\neq0$.
Another consequence is that automorphisms ($\IND V=0$) carry no
charge.

Finally, we investigate which representations can be realized on
$\mathfrak{k}_V$, and under which conditions charge conjugation is
ensured.  The remarks made in this context in the Fermionic case apply
here as well.


%% file: cu.tex
\chapter*{Endomorphisms of CAR and CCR Algebras}
\label{cha:CAR}
For the sake of completeness, we would like to include a brief survey
on the \CAs before discussing the main subject of our thesis. The
content of Section~\ref{sec:CUN} is in no way essential for an
understanding of Sections~\ref{sec:CAR} -- \ref{sec:SECTORS}, so that
this section may be skipped at a first reading of the text.

We would further like to remark that proofs of statements in
Sections~\ref{sec:CAR} and \ref{sec:CCR} that already appeared in our
papers \cite{CB1,CB2,CB3} will occasionally only be sketched here.
\specialsection{The Cuntz Algebras $\mathcal{O}(H)$}
\label{sec:CUN}
The universal \CALs $\mathcal{O}(H)$ generated by separable complex
\HSPs $H$ have been introduced by Cuntz in 1977 \cite{C77}. (\HSPs
inside operator algebras had before been considered by \DR
\cite{DR72,R76}). They provided new examples of \CALs with unexpected
properties, but also play an important role in the general structure
theory of \CALs.  See \cite{EK} for an account of $\mathcal{O}(H)$ in
a textbook.

Let $H$ be a separable complex \HSP of dimension\footnote{If
  $\dim H=1$, then the \CAL defined by these relations is the Abelian
  algebra generated by a single unitary operator $s$, isomorphic to
  the algebra of continuous functions on the spectrum of $s$. The
  Cuntz--Toeplitz algebra $\mathcal{T}(H)$ introduced below then
  reduces to the \CAL generated by the unilateral shift, the Toeplitz
  extension of $C(\TT)$ by the compacts \cite{Co}.} $\geq2$, with \ONB
$\{s_j\}_{j\in J}$. The unital *--algebra generated by the elements of
$H$, with relations (cf.~\eqref{CUNTZ})
$$s^*t=\langle s,t\rangle\1,\quad s,t\in H$$
and
$$\sum_{j\in J}s_js_j^*=\1$$
(resp.\ $\sum_{j\in J_0}s_js_j^*<\1$ for
any finite subset $J_0\subset J$, should $H$ be infinite dimensional),
possesses a unique $C^*$--norm. Its completion in this norm, the {\em
  \CA} $\mathcal{O}(H)$, is a separable, simple, nuclear
\CAL \cite{C77}. It can be constructed from the full Fock space
$\mathcal{F}(H)$ over $H$ as follows \cite{DEE80}. Let $\tau(s)$ be
the operator of norm $\NORM{s}$ which acts on $\mathcal{F}(H)$ by
tensoring on the left with $s\in H$. Then
$$\tau(s)^*\tau(t)=\langle s,t\rangle\1,\qquad \sum_{j\in
  J}\tau(s_j)\tau(s_j)^*=\1-p_\Omega\quad\text{(strongly)}$$
where
$p_\Omega$ denotes the rank--one projection onto the Fock vacuum
$\Omega$. The \CAL $\mathcal{T}(H)$ generated by all $\tau(s),\ s\in
H$, is the \emph{Cuntz--Toeplitz algebra} over $H$. If $H$ is finite
dimensional, then $\mathcal{T}(H)$ contains the ideal generated by
$p_\Omega$, the compact operators on $\mathcal{F}(H)$. In this case
$\mathcal{O}(H)$ is isomorphic to the quotient of $\mathcal{T}(H)$ by
this ideal. If $H$ is infinite dimensional, then $\mathcal{T}(H)$ is
simple, and in fact isomorphic to $\mathcal{O}(H)$.

The isomorphism class of $\mathcal{O}(H)$ depends only on the
dimension of $H$, thus it is also customary to write $\mathcal{O}_n$
instead of $\mathcal{O}(H)$ ($n\DEF\dim H$). One has the following
inclusions, since e.g.\ the elements
$s_1^2,s_1s_2,\dots,s_1s_n,s_2,\dots,s_n$ generate a copy of
$\mathcal{O}_{2n-1}$ in $\mathcal{O}_n$
$$\mathcal{O}_n\supset\mathcal{O}_{n+m(n-1)},\qquad m\leq\infty.$$
In
particular, $\mathcal{O}_2$ contains copies of all other
$\mathcal{O}_n$.  Algebras $\mathcal{O}_m$, $\mathcal{O}_n$ with
$m\neq n$ cannot be isomorphic as their $K_0$--groups differ
$$K_0(\mathcal{O}_n)=\ZZ/(n-1),\qquad K_0(\mathcal{O}_\infty)=\ZZ.$$
Any orthogonal projection in $\mathcal{O}(H)$ is equivalent to a
projection of the form $\sum_{\text{finite}}s_js_j^*$ or
$\1-\sum_{\text{finite}}s_js_j^*$. The unitary group
$\mathcal{U}(\mathcal{O}(H))$ of $\mathcal{O}(H)$ is connected, thus
one has $K_1(\mathcal{O}(H))=\{0\}$ \cite{C81}. If $\dim H$ is finite,
$\mathcal{U}(\mathcal{O}(H))$ is homeomorphic to the semigroup of all
unital *--endomorphisms of $\mathcal{O}(H)$: given
$u\in\mathcal{U}(\mathcal{O}(H))$, there is (by the universality of
$\mathcal{O}(H)$) a unique endomorphism $\rho$ such that $\rho(s)=us,\ 
s\in H\subset\mathcal{O}(H)$; conversely, an endomorphism $\rho$
determines the unitary $u=\sum_{j\in
  J}\rho(s_j)s_j^*\in\mathcal{O}(H)$ \cite{C80}.

Special examples of endomorphisms of $\mathcal{O}(H)$ have been given
by Izumi \cite{Iz}, using fusion rules of sectors (cf.\ the definition
\eqref{DIRSUM} of direct sums of endomorphisms). Another class of
endomorphisms of $\mathcal{O}(H)$ are the \emph{quasi--free
  endomorphisms,} namely those which leave $H$ invariant.
Non--surjective quasi--free endomorphisms exist only in the case $\dim
H=\infty$.  Quasi--free \emph{auto}morphisms are extensions of unitary
operators on $H$.  As a first step towards their duality theory for
compact groups, \DR studied the quasi--free action of a closed
subgroup $G$ of the special unitary group of $H$ ($\dim H<\infty$) on
$\mathcal{O}(H)$ \cite{DR87}, and in particular the relation between
$\mathcal{O}(H)$ and its fixed point subalgebra $\mathcal{O}(H)^G$
(serving as a prototype for the relation between field algebra and
observable algebra). They showed that $\mathcal{O}(H)^G$ is simple,
with trivial relative commutant in $\mathcal{O}(H)$, and that any
automorphism of $\mathcal{O}(H)$ which acts trivially on
$\mathcal{O}(H)^G$ is given by an element of $G$. Moreover, if
$\rho_H(a)=\sum s_jas_j^*$ is the endomorphism of $\mathcal{O}(H)$
induced by $H$, then the $G$--invariant intertwiners between powers of
${\rho_H}|_{\mathcal{O}(H)^G}$ are the same as the intertwiners
between tensor products of the representation of $G$ (the tensor
powers of $H$ are canonically embedded in $\mathcal{O}(H)$), and these
intertwiners generate $\mathcal{O}(H)^G$. These results have been
extended to the case of Hopf algebra actions by Cuntz \cite{C91}.

An important tool in the study of $\mathcal{O}(H)$ is to consider
$\mathcal{O}(H)^\TT$, the fixed point algebra under the quasi--free
action of the circle group \TT. As a Banach space,
$\mathcal{O}(H)^\TT$ is generated by monomials $s_{i_1}\dotsm
s_{i_m}s_{j_m}^*\dotsm s_{j_1}^*$ (same number of $s$ and $s^*$). If
$n\DEF\dim H$ is finite, then the monomials of the above form, with
$m$ fixed, constitute a system of $n^m\times n^m$ matrix units, hence
span an algebra isomorphic to $M(n^m,\CC)$. Since the embedding
$s_{i_1}\dotsm s_{i_m}s_{j_m}^*\dotsm s_{j_1}^*\mapsto
\sum_js_{i_1}\dotsm s_{i_m}s_js_j^*s_{j_m}^*\dotsm s_{j_1}^*$
corresponds to the usual embedding $M(n^m,\CC)\to M(n^{m+1},\CC),\ 
A\mapsto A\otimes \1_n$, it follows that $\mathcal{O}(H)^\TT$ is an
UHF algebra of type $n^\infty$, canonically isomorphic to the infinite
tensor product of copies of $M(n,\CC)$. In particular,
${\mathcal{O}_2}^\TT$ is isomorphic to the CAR algebra. If $\dim H$ is
infinite, then the monomials with fixed length $m\geq1$ generate an
algebra isomorphic to the compact operators $\mathcal{J}_\infty(H)$ on
$H$, and $\mathcal{O}(H)^\TT$ is isomorphic to a non--simple
AF--subalgebra of the infinite tensor product of copies (with unit
adjoined) of $\mathcal{J}_\infty(H)$ \cite{C77}.

The representation of $\mathcal{O}(H)^\TT$ as (a subalgebra of) an
infinite tensor product allows to define \emph{quasi--free states} over
$\mathcal{O}(H)$ as gauge invariant extensions of product states to
$\mathcal{O}(H)$ \cite{DEE80}, by utilizing the canonical \CE
$\mathcal{O}(H)\to\mathcal{O}(H)^\TT$. Specifically, any sequence of
positive trace class operators $\{K_j\}$ on $H$ with $\TR K_j=1$
(resp.\ $\TR K_j\leq1$ if $\dim H=\infty$) yields a quasi--free state
$\omega_{\{K_j\}}$
$$\omega_{\{K_j\}}(f_1\dotsm f_kg_l^*\dotsm g_1^*)=
\delta_{kl}\sum_{j=1}^k\langle f_j,K_jg_j\rangle,\quad f_j,g_j\in H.$$
Quasi--free automorphisms restrict to product automorphisms on
$\mathcal{O}(H)^\TT$, and $\rho_H$ restricts to the unilateral shift.

Quasi--free automorphisms and quasi--free states of $\mathcal{O}(H)$
have been studied by Evans et al.\ \cite{DEE80,ACE}. Among their
results are conditions for existence and uniqueness of KMS and ground
states for one--parameter groups of quasi--free automorphisms,
characterizations of primary quasi--free states, and criteria for
implementability of quasi--free automorphisms in quasi--free states.
For instance, they showed that $\mathcal{O}(H)$ has no inner
quasi--free automorphisms besides the identity. Equivalence of
quasi--free states over $\mathcal{O}_\infty$ has been studied by Laca
\cite{La93b}. 

Representations of $\mathcal{O}(H)$ are closely related to
endomorphisms of $\BB(\HH)$. A nondegenerate representation of
$\mathcal{O}(H)$ on a separable \HSP \HH gives rise, via the
endomorphism $\rho_H$ induced by $H$, to a unital *--endomorphism of
$\BB(\HH)$ (if $\dim H=\infty$, the representation $\pi$ of
$\mathcal{O}(H)$ has to be \emph{essential} \cite{La93a}, i.e.\ 
$\sum\pi(s_j)\pi(s_j)^*=\1$ in the strong topology). Conversely, since
$\BB(\HH)$ has only one normal representation up to
quasi--equivalence, each unital *--endomorphism of $\BB(\HH)$ is
inner. The representation of $\mathcal{O}(H)$ corresponding to an
endomorphism of $\BB(\HH)$ is, however, only determined modulo
quasi--free automorphisms. To illustrate how properties of
representations of $\mathcal{O}(H)$ are linked to properties of
endomorphisms of $\BB(\HH)$, one has e.g.\ that the commutant of an
essential representation $\pi$ of $\mathcal{O}(H)$ is equal to the
algebra of fixed points of the corresponding endomorphism $\rho_\pi$,
and that the commutant of the restriction of $\pi$ to
$\mathcal{O}(H)^\TT$ is equal to the intersection of the ranges of all
powers of $\rho_\pi$. Thus $\pi$ is irreducible \IFF $\rho_\pi$ is
ergodic, and $\pi|_{\mathcal{O}(H)^\TT}$ is irreducible \IFF
$\rho_\pi$ is a shift \cite{La93a,BJP}. A classification of certain
ergodic endomorphisms of $\BB(\HH)$ up to conjugacy has been achieved
by Laca and Fowler, by describing all extensions of pure states over
$\mathcal{O}(H)^\TT$ to $\mathcal{O}(H)$ \cite{FL}.

Interest in the connection between representations of $\mathcal{O}(H)$
and endomorphisms of $\BB(\HH)$ arose from the theory of
$E_0$--semigroups of endomorphisms which was initiated by Powers
\cite{P88}, with contributions by Arveson, Bratteli, J\o rgensen,
Laca, Price, Robinson and others (see the review \cite{A94} \ART). One
has an index theory for $E_0$--semigroups which gives partial results
for the classification of $E_0$--semigroups up to outer conjugacy. The
basic examples of $E_0$--semigroups are semigroups of Bogoliubov
endomorphisms of CAR and CCR algebras.

Representations of $\mathcal{O}(H)$ related to wavelet theory have
been studied by Bratteli et.~al.\ (see \cite{BJ96,BJKW} \ART). They
obtained the decomposition of a special class of such representations
(and of their restrictions to $\mathcal{O}(H)^\TT$) into irreducibles
via number theory.

The \CAs can be regarded as elementary building blocks of infinite
\CALs (algebras containing non--unitary isometries). Any unital simple
infinite \CAL contains copies of all $\mathcal{O}_n$ as subquotients
(quotients of subalgebras) \cite{C77}. By Kirchberg's results
\cite{Ki,Ki94,Ki95,KP}, any separable unital exact \CAL is isomorphic
to a subalgebra of $\mathcal{O}_2$ (and to a subquotient of the \CAR);
any separable unital nuclear (s.u.n.) \CAL is isomorphic to the range
of a \CE of $\mathcal{O}_2$; and any simple s.u.n.\ \CAL which
contains a central sequence of copies of $\mathcal{O}_2$ is itself
isomorphic to $\mathcal{O}_2$ (cf.\ footnote (\ref{foo:B}) on
page~\pageref{foo:B}).  An important open question concerning
Elliott's classification program for nuclear \CALs \cite{Ell,EK} is
whether any simple s.u.n.\ \CAL with vanishing $K_0$ and $K_1$--groups
is already isomorphic to $\mathcal{O}_2$. Inductive limits of matrix
algebras over $\mathcal{O}_n$ have been classified by R{\o}rdam
\cite{R93}.


%% file: ca.tex
\specialsection{Quasi--free Endomorphisms of the CAR Algebra}
\label{sec:CAR}
Basics on the \CAR can be found in the textbooks \cite{BR,EK,PR} and
will be cited here only as far as necessary. Because of its
convenience for handling Bogoliubov transformations, Araki's formalism
of selfdual \CARs \cite{A68,A70,A87} will be used throughout. A
selfdual \CAR is simply the complexification of a real Clifford
algebra.

Most of the results in this section are contained in \cite{CB1}.
However, it has become clear in the meantime that a different
arrangement of the material would be desirable. Thus the present setup
deviates from the one in \cite{CB1}, most significantly in
Section~\ref{sec:IP}. Minor improvements can be found throughout the
text.
\subsection{Quasi--free endomorphisms and quasi--free states}
\label{sec:QFCAR}
Let \KK be an infinite dimensional separable complex Hilbert space,
equipped with a complex conjugation $f\mapsto f^*$. 
The {\em(selfdual) \CAR} \CK over \KK is the unique (simple) \CAL
generated by \1 and the elements of \KK, subject to the
anticommutation relation
$$\{f^*,g\}\DEF f^*g+gf^*=\langle f,g\rangle\1,\qquad f,g\in\KK.$$
\CK is the complexified (or $C^*$--) Clifford algebra over the real \HSP
$$\RE\KK\DEF\SET{f\in\KK}{f^*=f},$$
a UHF algebra of type $2^\infty$.
The $C^*$--norm on \CK extends the norm on $\RE\KK$ (but not the norm
on \KK) up to a factor $\sqrt{2}$. \KK will henceforth be viewed as a
subspace of \CK.

\emph{Quasi--free endomorphisms} (or \emph{Bogoliubov transformations})
are precisely the unital *--endomorphisms of \CK that leave \KK
invariant. Put differently, every isometry $V$ on \KK that commutes
with complex conjugation (and therefore restricts to a real--linear
isometry of $\RE\KK$) extends to a unital isometric *--endomorphism
$\rho_V$ of \CK:
$$\rho_V(f)=Vf,\qquad f\in\KK.$$
Such isometries $V$ are called {\em
  Bogoliubov operators}, and the semigroup of Bogoliubov operators is
denoted by
\begin{equation}
  \label{IK}
  {\II}(\KK)\DEF\SET{V\in\BB(\KK)}{V^*V=\1,\ \4{V}=V}
\end{equation}
where the bar indicates complex conjugation
\begin{equation}
  \label{CONJ}
  \4{A}f\DEF A(f^*)^*,\qquad f\in\KK
\end{equation}
for bounded linear operators $A\in\BB(\KK)$. The map $V\mapsto\rho_V$
is a unital isomorphism from ${\II}(\KK)$ onto the semigroup of
quasi--free endomorphisms; for fixed $a\in\CK$, the map
$V\mapsto\rho_V(a)$ is continuous \WRT the strong topology on
${\II}(\KK)$ and the norm topology on \CK.

Let $V\in{\II}(\KK)$. Since $\RAN V$ is closed and $\ker
V=\{0\}$, $V$ and $V^*$ are semi--Fredholm operators \cite{K} and have
well--defined Fredholm indices. The map
$${\II}(\KK)\to\NN\cup\{\infty\},\quad V\mapsto-\IND
V\DEF\dim\ker V^*$$
is a surjective homomorphism of semigroups
($0\in\NN$ by convention). Additivity of the Fredholm index follows in
this special case simply from $\ker(VW)^*=\ker V^*\oplus V(\ker W^*)$.
The semigroup ${\II}(\KK)$ is the disjoint union of subsets
\begin{equation}
  \label{INK}
  {\II}(\KK)=\bigcup_{n\in\NN\cup\{\infty\}}{\II}^n(\KK),\qquad
  {\II}^n(\KK)\DEF\SET{V\in{\II}(\KK)}{\IND V=-n}. 
\end{equation}
Note that $\rho_V$ is an automorphism \IFF $V$ belongs to
${\II}^0(\KK)$, the group of unitary Bogoliubov operators.
${\II}^0(\KK)$ acts on ${\II}(\KK)$ by left multiplication, the orbits
of this action are the sets ${\II}^n(\KK)$, and the stabilizer of
$V\in{\II}^n(\KK)$ is isomorphic to $O(n)$ (the orthogonal group of an
$n$--dimensional real \HSP).

The quasi--free automorphism $\rho_{-\1}$ induces a $\ZZ_2$--grading on
\CK:
$$\begin{array}{l}
  \CK=\CK_0\oplus\CK_1,\\ 
  \CK_g\cdot\CK_{g'}\subset\CK_{g+g'},\\
  \CK_g\DEF\SET{a}{\rho_{-\1}(a)=(-1)^ga},\qquad g,g'\in\ZZ_2=\{0,1\}.
\end{array}$$
Doplicher and Powers proved that the \emph{even subalgebra} $\CK_0$ is
a simple \CAL \cite{DP}. St{\o}rmer sharpened this result by showing
that $\CK_0$ is UHF of type $2^\infty$, hence *-isomorphic to \CK
itself \cite{St}. We found that any $V\in{\II}^1(\KK)$ gives rise to
an isomorphism from \CK onto $\CK_0$ in the following way \cite{CB2}.
Let $f_V$ be the unique (up to a sign) unitary skew--adjoint element
in $\ker V^*\subset\CK$, and let $u_V\DEF\frac{1}{\sqrt{2}}(\1+f_V)$.
Then $u_V$ is unitary, $u_V^2=f_V$, and the map
\begin{equation}
  \label{EVEN}
  \sigma_V:a\mapsto u_V\rho_V(a){u_V}^*
\end{equation}
defines a unital *-isomorphism from \CK onto $\CK_0$. $\sigma_V$ acts
on even elements like $\rho_V$, and on odd elements like $\rho_V$
followed by left multiplication with $f_V$. A similar construction has
been given independently by \cite{Ro}.
\medskip

Next we describe the set of states we are interested in. A state
$\omega$ over \CK is called \emph{quasi--free} if it is invariant under
$\rho_{-\1}$ (i.e.\ it vanishes on $\CK_1$), and if its even $n$--point
functions have the form \cite{A70}
$$\omega(f_1\cdots f_{2m})=(-1)^\frac{m(m-1)}{2}\sum_\sigma
\SGN\sigma\cdot\omega\Big(f_{\sigma(1)}f_{\sigma(m+1)}\Big)
\cdots\omega\Big(f_{\sigma (m)}f_{\sigma (2m)}\Big)$$
where the sum runs over all permutations $\sigma$ satisfying
$\sigma(1)<\ldots<\sigma(m)$ and $\sigma(j)<\sigma(j+m),\ 
j=1,\dots,m$. Therefore quasi--free states are completely determined
by their two--point functions, and one has a bijection
$$S\mapsto\omega_S,\quad\omega_S(f^*g)=\langle f,Sg\rangle$$
between the convex set
$$\QK\DEF\SET{S\in\BB(\KK)}{0\le S\le\1,\ \4{S}=\1-S}$$
and the set of
quasi--free states. A (non--trivial) convex combination of two
distinct quasi--free states $\omega_S,\omega_{S'}$ is quasi--free \IFF
$S-S'$ has rank two \cite{W75}. Quasi--free endomorphisms act from
the right on quasi--free states according to
\begin{equation}
  \label{OP}
  \omega_S\0\rho_V=\omega_{V^*SV}.
\end{equation}
Any *--automorphism which maps the set of quasi--free states onto
itself is known to be quasi--free \cite{W75}.

Projections in $\QK$ are called {\em\BPs}, and the corresponding
states are called {\em\FSs;} the latter are precisely the \emph{pure}
quasi--free states \cite{MRT}. The group of quasi--free automorphisms
acts transitively on the set of \FSs, because
${\II}^0(\KK)$ acts transitively on the set of \BPs. Note that
for a \BP $P$, the complementary (basis) projection is simply given by
$\4{P}$. Since $\omega_P(f^*f)=0$ if $f\in\4{P}(\KK)$, the elements of
$\4{P}(\KK)$ (resp.\ $P(\KK)$) correspond to annihilation (resp.\ 
creation) operators in the state $\omega_P$. A (faithful and
irreducible) GNS representation $\pi_P$ for $\omega_P$ is given by
\begin{equation}
  \label{PI}
  \pi_P(f)\DEF a(Pf)^*+a(P(f^*))
\end{equation}
on the antisymmetric Fock space
$\mathcal{F}_a(P(\KK))$ over $P(\KK)$, with the usual Fock vacuum
$\Omega_P$ as cyclic vector and annihilation operators $a(g),\ g\in
P(\KK)$. In a Fock representation $\pi_P$, a quasi--free endomorphism
$\rho_V$ induces the transformation
\begin{equation*}
  a(g)\mapsto a(PVPg)+a(PV\4{P}(g^*))^*,\quad g\in P(\KK),
\end{equation*}
which shows the connection to the (state--dependent) description of
Bogoliubov transformations by pairs of operators as preferred by some
authors. GNS representations of arbitrary quasi--free states can be
obtained, by the ``doubling procedure'' of Powers and St{\o}rmer
\cite{PS,A70}, as restrictions of Fock representations of an enlarged
\CAR. 
\medskip

The grading automorphism $\rho_{-\1}$ is not inner in \CK. However,
since every quasi--free state is invariant under this automorphism,
$\rho_{-\1}$ is canonically implemented in any quasi--free state. Let
$P$ be a \BP, and let $\Psi_P(-\1)$ be the self--adjoint unitary on
Fock space given by 
\begin{equation}
  \label{PME}
  \Psi_P(-\1)\pi_P(a)\Omega_P=\pi_P(\rho_{-\1}(a))\Omega_P,
  \qquad a\in\CK.
\end{equation}
Then $\Theta_P(-\1)\DEF\tfrac{1}{\sqrt{2}}\big(\1-i\Psi_P(-\1)\big)$
is also unitary. Define a new representation $\psi_P$ of \CK on
$\mathcal{F}_a(P(\KK))$, which is unitarily equivalent to $\pi_P$
\begin{equation}
  \label{PSI}
  \psi_P(a)\DEF\Theta_P(-\1)\pi_P(a)\Theta_P(-\1)^*.
\end{equation}
We will call $\psi_P$ the \emph{twisted Fock representation} induced by
$P$. Then one has
\begin{align}
  \psi_P|_{\CK_0}&=\pi_P|_{\CK_0},\label{PSI0}\\
  \psi_P(a)&=i\pi_P(a)\Psi_P(-\1),\qquad a\in\CK_1,\label{PSI1}\\
  [\pi_P(f)^*,\psi_P(g)]&=i\langle f,g\rangle\Psi_P(-\1),
  \qquad f,g\in\KK,\label{PSI2}
\end{align}
and, as shown by Foit \cite{Fo}, \emph{twisted duality} holds for all
*--invariant subspaces $\HH\subset\KK$, an adaptation of Haag duality
(see p.~\pageref{HD}) to Fermi fields:
\begin{equation}
  \label{TD}
  \pi_P(\mathfrak{C}(\HH))'=\psi_P(\mathfrak{C}(\HH^\bot))''.
\end{equation}
Here $\mathfrak{C}(\HH)$ is the $C^*$--subalgebra of \CK generated by
the elements of $\HH$, and similar for $\mathfrak{C}(\HH^\bot)$.

Given a \BP $P$, a state over \CK is said to be (\TT--) \emph{gauge
  invariant} if it is invariant under the one--parameter group of
quasi--free automorphisms $(\rho_{U_\lambda })_{\lambda\in\RR}$ with
\begin{equation}
  \label{U1}
  U_\lambda\DEF e^{i\lambda}P+e^{-i\lambda}\4{P}\in{\II}^0(\KK). 
\end{equation}
A quasi--free state $\omega_S$ is gauge
invariant \IFF $[P,S]=0$.

The so--called \emph{central state} $\omega_{1/2}$ \cite{SS64,Man,A70}
is the unique tracial state over \CK. By uniqueness, $\omega_{1/2}$ is
invariant under all unital *--endomorphisms of \CK. It can be used to
define \CEs on \CK. If $V$ is a Bogoliubov operator with $-\IND
V<\infty$, then there is a unique minimal (faithful) \CE $E_V$ from
\CK onto $\rho_V(\CK)$, determined by $E_V(ab)=a\omega_{1/2}(b)$ if
$a\in\mathfrak{C}(\RAN V)=\rho_V(\CK)$, $b\in\mathfrak{C}(\ker V^*)$
(see \cite{CB1}). Using an explicit ``quasi--basis'' for $E_V$, we
computed the Watatani index \cite{Wat} of $E_V$ in \cite{CB1}
$$\IND E_V=2^{-\IND V}.$$
Thus we found the fundamental index formula
\begin{equation}
  \label{DV}
  d_V\DEF [\CK:\rho_V(\CK)]^\2=2^{M_V},\qquad 
  M_V\DEF-\2\IND V
\end{equation}
which relates the Fredholm index of $V$ to the Watatani index
$[\CK:\rho_V(\CK)]$ of $\rho_V$. We shall take \eqref{DV} as the
definition of the numbers $d_V$ and $M_V$ also in the case $-\IND
V=\infty$. $d_V$ may be regarded as the \emph{statistics dimension} of
the quasi--free endomorphism $\rho_V$, cf.\ p.~\pageref{page:D}. One
obviously has $M_{VV'}=M_V+M_{V'}$ and $d_{VV'}=d_Vd_{V'}$. Also note
that the \CEs $E_V$ allow to define \emph{left inverses} (see \cite{H})
$\rho_V^{-1}\0E_V$ for quasi--free endomorphisms. Explicitly, for a
quasi--free endomorphism $\rho_V$, a left inverse $\phi_V$ is given by
\begin{equation}
  \label{LI}
  \phi_V(ab)\DEF\rho_V^{-1}(a)\omega_{1/2}(b)\qquad\text{if }
  a\in\mathfrak{C}(\RAN V),\ b\in\mathfrak{C}(\ker V^*).
\end{equation}
Not
surprisingly, the \CEs from \CK onto $\CK_0$ that are obtained in this
way from the isomorphisms $\sigma_V$ of \eqref{EVEN} ($d_V=\sqrt{2}$)
are equal to the mean over the action of $\ZZ_2$.

The \VN algebra generated by \CK in the central state $\omega_{1/2}$
is the hyperfinite ${\rm II}_1$ factor. In general, the types of
quasi--free factor states can be computed from spectral properties of
the associated operators $S\in\QK$. See \cite{MY} for a complete
classification (extending earlier results in \cite{DA,Ri,PS}),
including the fine classification of type III.
\medskip

Of uppermost importance for our study of implementable quasi--free
endomorphisms are the criteria for quasi--equivalence of quasi--free
states. First results in this direction were obtained by Shale and
Stinespring \cite{SS65}. These authors showed that a quasi--free
automorphism $\rho_U,\ U\in{\II}^0(\KK)$, is unitarily implementable
in a Fock representation $\pi_P$ \IFF
\begin{equation}
  \label{SS}
  [P,U] \text{ is Hilbert--Schmidt.}
\end{equation}
Equivalently, two \FSs $\omega_P,\omega_{P'}$ (i.e.\ their GNS
representations) are unitarily equivalent \IFF $P-P'$ is
Hilbert--Schmidt. A sufficient condition for quasi--equivalence of
gauge invariant quasi--free states followed from the work of
Dell'Antonio \cite{DA} and Rideau \cite{Ri}. Powers and St{\o}rmer
proved this condition also to be necessary \cite{PS}, and Araki
extended the result to arbitrary quasi--free states \cite{A70}. One
has
\begin{equation}
  \label{QEQ}
  \omega_S\approx\omega_{S'}\iff S^{1/2}-{S'}^{1/2} 
  \text{ is Hilbert--Schmidt}.
\end{equation}
Here ``$\approx$'' means ``quasi--equivalent''. It has been observed
by Powers \cite{P87} that this criterion can be simplified if one of
the operators $S,S'$ is a projection. Namely, if $P$ is a \BP, then
\begin{equation}
  \label{PSC}
  \omega_P\approx\omega_S\iff\TR\4{P}S\4{P}<\infty.
\end{equation}
Quasi--equivalence of the restrictions of gauge invariant quasi--free
states to gauge invariant \CARs $\CK^G$ (now \WRT the quasi--free
action of an arbitrary compact group $G$) has been investigated by
Matsui \cite{M87a}, extending results of Araki and Evans for the case
$G=\ZZ_2$ \cite{AE}, and of Baker and Powers for the groups $\ZZ_2,\ 
\TT$ and $SU(2)$ \cite{BPa,BPb}. If $P,P'$ are \BPs commuting with
the action of $G$, then one has \cite{M87a}
\begin{align}
  \omega_P|_{\CK^G}\simeq\omega_{P'}|_{\CK^G}\iff &P-P'\text{ is
    Hilbert--Schmidt,}\label{MATSUI1}\\
  &{\det}_{P(\KK)\cap\4{P'}(\KK)}(g)=1\text{ for all }g\in G.
    \label{MATSUI2}
\end{align}
Here ``$\simeq$'' means ``unitarily equivalent'', and
$P(\KK)\cap\4{P'}(\KK)$ is a finite dimensional $G$--invariant
subspace if $P-P'$ is Hilbert--Schmidt. The condition on the
determinant is a generalization of the $\ZZ_2$--index of Araki and
Evans \cite{AE,A87,EK}. If the GNS representations of $\omega_P$ and
$\omega_{P'}$ are both realized on $\mathcal{F}_a(P(\KK))$ (this is
possible under the Hilbert--Schmidt condition), then the
transformation law of the cyclic vector $\Omega_{P'}$ is exactly given
by the character $\det_{P(\KK)\cap\4{P'}(\KK)}(g)$.
We will rediscover this character in Section~\ref{sec:CARCH}. On the
other hand, if $S,S'\in\QK$ commute with the action of $G$ and have
trivial kernels, then one has \cite{M87a}
$$\omega_S|_{\CK^G}\approx\omega_{S'}|_{\CK^G}\iff S^\2-{S'}^\2\text{
  is Hilbert--Schmidt}.$$
\subsection[Representations of the form $\pi\0\rho$]{Representations
  of the form $\boldsymbol{\pi\0\rho}$} 
\label{sec:REP}
As mentioned in the introduction on p.~\pageref{page:LOCEND}, the
representations describing \SSSs in the algebraic approach have the
form $\pi_0\0\rho$ where $\pi_0$ is a vacuum representation and \rho
is some localized endomorphism of the observable algebra. Here we
study representations $\pi_P\0\rho_V$ of \CK where $\pi_P$ is a Fock
representation and $\rho_V$ a quasi--free endomorphism, and in
particular the decomposition of such representations into cyclic
and irreducible subrepresentations. Among the results are a
necessary and sufficient condition for implementability of quasi--free
endomorphisms (a generalization of the Shale--Stinespring condition),
and alternative proofs of results of B\"ockenhauer \cite{JMB2}. 

Let us first repeat what is meant by ``implementability of
endomorphisms''.
\begin{Def}
  \label{def:IMP}\hspace*{\fill}\\
  A *--endomorphism \rho of a \CAL \AA is \emph{implementable}
  in a representation $(\pi,\HH)$ if there exists a (possibly finite)
  sequence $(\Psi_n)_{n\in I}$ in $\BB(\HH)$ with relations
  \begin{align}
    &\Psi_m^*\Psi_n=\delta_{mn}\1,\qquad 
      \sum\limits_{n\in I}\Psi_n\Psi_n^*=\1,\footnotemark\label{CUNTZ1}
      \addtocounter{footnote}{-1}\\
    \intertext{which implements \rho via}
    &\pi(\rho(a))=\sum_{n\in I}\Psi_n\pi(a)\Psi_n^*,\footnotemark
      \qquad a\in\AA.\label{CUNTZ2}
  \end{align}
  \footnotetext{in the strong topology if $I$ is infinite. In the
    terminology of Laca \cite{La93a}, we only consider essential
    representations of \CAs. See Section~\ref{sec:CUN}.}
\end{Def} 
$\HH$ then decomposes into the orthogonal direct sum of the ranges of
the isometries $\Psi_n$, and $\pi\0\rho$ decomposes into
subrepresentations $\pi\0\rho|_{\RAN\Psi_n}$, each of them unitarily
equivalent to $\pi$.  But the converse is also true, so $\rho$ is
implementable in $\pi$ \IFF $\pi\0\rho$ is equivalent to a multiple of
$\pi$. For irreducible $\pi$ this means
\begin{equation}
  \label{IMPQEQ}
  \rho\text{ is implementable in }\pi\iff\pi\0\rho\approx\pi.
\end{equation}
The isometries $(\Psi_n)_{n\in I}$ constitute an \ONB of the \HSP
$H\DEF\4{\SPAN}(\Psi_n)$ in $\BB(\HH)$, with scalar product given by
$\Psi^*\Psi'=\langle \Psi,\Psi'\rangle\1$.
Every element $\Psi$ of $H$ is an intertwiner from $\pi$ to
$\pi\0\rho$:
\begin{equation}
  \label{INTER}
  \Psi\pi(a)=\pi(\rho(a))\Psi,\quad a\in{\AA}.
\end{equation}
$H$ coincides with the space of intertwiners from $\pi$ to $\pi\0\rho$
\IFF $\pi$ is irreducible. If $\pi$ is reducible, there may exist
several \HSPs implementing \rho, mutually related by unitaries in
$\pi(\rho(\AA))'$. More precisely, if $(\Psi_n)_{n\in I}$ and
$(\Psi_n')_{n\in I}$ both implement \rho in $\pi$, then
$\Psi\DEF\sum_n\Psi_n'\Psi_n^*$ is a unitary in $\pi(\rho(\AA))'$, and
$\Psi_n'=\Psi\Psi_n$. Conversely, given $(\Psi_n)_{n\in I}$ and a
unitary $\Psi\in\pi(\rho(\AA))'$, $(\Psi\Psi_n)_{n\in I}$ is a set of
implementing isometries.

An implementable endomorphism \rho gives rise to normal
*--endomorphisms $\rho_H(a)\DEF\sum_{n\in I}\Psi_na\Psi_n^*$ of
$\BB(\HH)$, with index \cite{L}
$$[\BB(\HH):\rho_H(\BB(\HH))]=(\dim H)^2,$$
where $\dim H$ does not
depend on the choice of $H=\4{\SPAN}(\Psi_n)$. Let us outline the
computation of the index in the algebraic setting of Watatani
\cite{Wat}, for the case $\dim H<\infty$. $\phi_H(a)\DEF(\dim
H)^{-1}\sum_n\Psi_n^*a\Psi_n$ is a left inverse for $\rho_H$ (cf.\ 
\eqref{LI0}), yielding the minimal \CE $E_H\DEF\rho_H\0\phi_H$ from
$\BB(\HH)$ onto $\rho_H(\BB(\HH))$. $(\sqrt{\dim
  H}\cdot\Psi_n^*)_{n=1,\ldots,\dim H}$ is a quasi--basis for $E_H$,
hence $\IND E_H=\dim H\cdot\sum_n\Psi_n^*\Psi_n=(\dim H)^2$. Of
course, we will see that $\dim H=d_V$ (defined by \eqref{DV}) if
$\rho_V$ is a quasi--free endomorphism which is implementable in some
Fock representation.
 
Let us add a last remark on the general situation. Suppose we are
given a set of implementers $(\Psi_n)_{n\in I}$. Then for $m,n\in I$,
$\Psi_m\Psi_n^*\in\pi(\rho(\AA))'$ is a partial isometry containing
$\RAN\Psi_n$ in its initial space, and $\Psi_m=(\Psi_m\Psi_n^*)
\Psi_n$. This suggests to construct a complete set of implementing
isometries by multiplying one isometry $\Psi$ fulfilling \eqref{INTER}
with certain partial isometries in $\pi(\rho(\AA))'$. We shall employ
this idea later in Section~\ref{sec:FORM}.
\medskip

After this digression we concentrate on Bogoliubov transformations
again. Inspection of \eqref{IMPQEQ} leads one to study the
representations $\pi_P\0\rho_V$; as will turn out, they are
quasi--equivalent to GNS representations associated with the states
$\omega_P\0\rho_V$ (a similar observation has been made, in a
different setting, by Rideau \cite{Ri}). Thus the question of
quasi--equivalence of such representations can be traced back to the
question of quasi--equivalence of the corresponding states. 

Let $P$ be a \BP, let $V\in{\II}(\KK)$, and regard
\begin{equation}
  \label{V}
  v\DEF PVV^*P
\end{equation} 
as an operator on $P(\KK)$. The direct sum decomposition $P(\KK)=\ker
v\,\oplus \,\4{\RAN v}$ induces a tensor product decomposition of Fock
space: $\mathcal{F}_a(P(\KK))\cong\mathcal{F}_a(\ker v)\otimes
\mathcal{F}_a(\4{\RAN v})$. Choose an \ONB $(f_j)_{j=1,\ldots,N_V}$
for $\ker v=P(\KK)\cap\ker V^*$, where
\begin{equation}
  \label{NV}
  N_V\DEF\dim\ker v\leq M_V
\end{equation}
(the inequality follows from $\ker v\oplus\ker\4{v}\subset\ker V^*$).
Let $\psi_P$ be the twisted Fock representation defined in
\eqref{PSI}, and let $I_{N_V}$ be the set of $2^{N_V}$ multi--indices
$\alpha=(\alpha_1,\dots,\alpha_l),\ \alpha_j\in\NN,\ l<\infty$,
obeying
\begin{equation} 
  \label{I}
  0\leq l\leq N_V,\quad1\leq\alpha_1<\dots<\alpha_l\leq N_V\quad
  (\alpha\DEF0\text{ if }l=0).
\end{equation}
For $\alpha\in I_{N_V}$, set
\begin{equation}
  \label{DEF}
  \begin{split} 
    \psi_\alpha&\DEF\psi_P(f_{\alpha_1}\dotsm f_{\alpha_l})
      \quad(\psi_0\DEF\1),\\ 
    \phi_\alpha&\DEF\psi_\alpha\Omega_P,\\
    \mathcal{F}_\alpha&\DEF\4{\pi_P(\rho_V(\CK))\phi_\alpha},\\
    \pi_\alpha&\DEF\pi_P\0\rho_V|_{\mathcal{F}_\alpha}.
  \end{split}
\end{equation}
Note that, by the CAR and \eqref{TD}, the $\psi_\alpha$ are partial
isometries in $\pi_P(\rho_V(\CK))'$. 
\begin{Prop}
  \label{prop:DECCAR}\hspace*{\fill}\\
  Each of the $2^{N_V}$ cyclic subrepresentations
  $(\pi_\alpha,\mathcal{F}_\alpha,\phi_\alpha)$ induces the state
  $\omega_P\0\rho_V$, and $\pi_P\0\rho_V$ splits into their direct
  sum: 
  $$\quad\pi_P\0\rho_V=\bigoplus_{\alpha\in I_{N_V}}\pi_\alpha.$$
\end{Prop}
\begin{proof}
  It is clear by definition that $\mathcal{F}_\alpha$ is an invariant
  subspace for $\pi_P\0\rho_V$ with cyclic vector $\phi_\alpha$. Since
  $\psi_\alpha\in\pi_P(\rho_V(\CK))'$ and
  $\psi^*_\alpha\psi_\alpha\Omega_P=\Omega_P$, we have $\langle
  \phi_\alpha,\pi_\alpha(a)\phi_\alpha\rangle=\langle\Omega_P,
  \pi_P(\rho_V(a))\Omega_P\rangle=\omega_P(\rho_V(a)),\ a\in\CK$. Thus
  $(\pi_\alpha,\mathcal{F}_\alpha,\phi_\alpha)$ is a GNS
  representation for $\omega_P\0\rho_V$ (and the representations
  $\pi_\alpha$ are mutually unitarily equivalent).
  
  Next we show $\mathcal{F}_\alpha\bot\mathcal{F}_\beta$ for
  $\alpha\neq\beta$. Since at least one of the vectors
  $\psi_\alpha^*\psi_\beta\Omega_P,\ \psi_\beta^*\psi_\alpha\Omega_P$
  vanishes if $\alpha\neq\beta$, we have for $a,b\in\CK$
  $$\langle\pi_P(\rho_V(a))\phi_\alpha,\pi_P(\rho_V(b))
    \phi_\beta\rangle=\langle\psi_\alpha\Omega_P,
    \pi_P(\rho_V(a^*b))\psi_\beta\Omega_P\rangle=0,$$
  implying
  orthogonality of $\mathcal{F}_\alpha$ and $\mathcal{F}_\beta$.
  
  Finally we have to prove
  $\mathcal{F}_a(P(\KK))=\oplus_\alpha\mathcal{F}_\alpha$. Using
  $\pi_P(\rho_V(f))=a(PVf)^*+a(PVf^*),\ f\in\KK$, one can show by
  induction on the particle number 
  $$\mathcal{F}_0=\4{\pi_P(\rho_V(\CK))\Omega_P}=
    \mathcal{F}_a(\4{\RAN PV}) =\mathcal{F}_a(\4{\RAN v}).$$
  Since the
  $\phi_\alpha$ form an \ONB for $\mathcal{F}_a(\ker v)$, the
  assertion follows.
\end{proof}
The decomposition of these cyclic representations into irreducibles
will be examined after stating the implementability condition.
Remember that $\4{P}=\1-P$.
\begin{Th}
  \label{th:IMPCAR}\hspace*{\fill}\\
  A quasi--free endomorphism $\rho_V$ is implementable in a Fock
  representation $\pi_P$ \IFF $PV\4{P}$ is a Hilbert--Schmidt
  operator.
\end{Th}
\begin{proof}
  In view of (\ref{IMPQEQ}) and Proposition~\ref{prop:DECCAR}, $\rho_V$
  is implementable in $\pi_P$ \IFF $\omega_P\0\rho_V\approx\omega_P$.
  Since $\omega_P\0\rho_V=\omega_{V^*PV}$ by \eqref{OP}, the
  Powers--St\o rmer--Araki criterion in the form~\eqref{PSC} implies
  that $\omega_P\0\rho_V\approx\omega_P$ \IFF
  $\TR\4{P}V^*PV\4{P}<\infty$. The latter condition is clearly
  equivalent to $PV\4{P}$ being Hilbert--Schmidt.
\end{proof}
Note that $PV\4{P}$ is Hilbert--Schmidt \IFF $[P,V]=PV\4{P}-\4{P}VP$
is, so the Shale--Stinespring condition (\ref{SS}) remains valid. We
denote the semigroup of Bogoliubov operators fulfilling this condition
by
$${\II}_P(\KK)\DEF\SET{V\in{\II}(\KK)}{PV\4{P}\text{ is
  Hilbert--Schmidt}}.$$
Since $PV\4{P}$ and $\4{P}VP$ are compact for
$V\in{\II}_P(\KK)$, $PVP+\4{P}V\4{P}=V-PV\4{P}-\4{P}VP$ is
semi--Fredholm, and
$$M_V=-\IND PVP\in\NN\cup\{\infty\}.$$
Thus we have a decomposition  (cf.~\eqref{INK})
$${\II}_P(\KK)=\bigcup_{m\in\NN\cup\{\infty\}}
{\II}_P^{2m}(\KK),\qquad
{\II}_P^{2m}(\KK)\DEF\SET{V\in{\II}_P(\KK)}{M_V=m}.$$
The group
$\II_P^0(\KK)$ is usually called the \emph{restricted orthogonal group}
\cite{PrSe}. Note that the ``statistics dimension'' $d_V$ defined by
\eqref{DV} is contained in $\NN\cup\{\infty\}$ if $V\in{\II}_P(\KK)$.
Note also that non--surjective quasi--free endomorphisms cannot be
inner in \CK since the \CAR, being AF and thus finite, does not
contain non--unitary isometries.
\medskip

In the course of constructing localized endomorphisms for the
conformal WZW models, J.~B\"ockenhauer described the decomposition of
representations $\pi_P\0\rho_V$ and of their restrictions to the even
subalgebra $\CK_0$ into irreducibles \cite{JMB2} (see also \cite{Sz}).
His methods work only for Bogoliubov operators with finite index,
i.e.\ those belonging to the sub--semigroup
$$\IFIN(\KK)\DEF\SET{V\in{\II}(\KK)}{M_V<\infty}.$$
We shall now
present alternative proofs of his results which have the merit of
being completely basis--independent.

For $V\in{\II}(\KK)$, let $Q_V$ be the orthogonal projection
onto $\ker V^*$, and let $S_V$ be the operator characterizing the
quasi--free state $\omega_P\0\rho_V$
$$Q_V\DEF[V^*,V]=\1-VV^*,\qquad S_V\DEF V^*PV\in\QK.$$
The operators $Q_V$ and $S_V\4{S_V}$ have finite rank if
$V\in\IFIN(\KK)$.

Let us first determine when two representations of the form
$\pi_P\0\rho_V$ ($P$ fixed) are unitarily equivalent.
\begin{Lem}
  \label{lem:EQ}\hspace*{\fill}\\
  Let $V,V'\in\IFIN(\KK)$. Then the following
  conditions are equivalent:\\
  \indent a) $\pi_P\0\rho_V$ and $\pi_P\0\rho_{V'}$ are unitarily
  equivalent;\\ 
  \indent b) there exists $U\in{\II}^0_P(\KK)$ with $V'=UV$;\\
  \indent c) $\IND V=\IND V'$, and $S_V-S_{V'}$ is Hilbert--Schmidt.
\end{Lem}
\begin{proof}
  We first show a) $\Rightarrow$ c). By Proposition~\ref{prop:DECCAR},
  $\pi_P\0\rho_V\simeq\pi_P\0\rho_{V'}$ implies $\omega_{S_V}=
  \omega_P\0\rho_V\approx\omega_P\0\rho_{V'}=\omega_{S_{V'}}$. Hence
  by \eqref{QEQ}, $S_V^\2-S_{V'}^\2$ is Hilbert--Schmidt (HS) which
  is, for $V,V'\in\IFIN(\KK)$, equivalent to $S_V-S_{V'}$ being HS
  (see \cite{JMB2}).  Moreover, equivalent representations have
  isomorphic commutants. We have by \eqref{TD}
  $\pi_P(\rho_V(\CK))'=\psi_P(\mathfrak{C}(\ker V^*))''\simeq
  \pi_P(\mathfrak{C}(\ker V^*))''$. Hence the commutants have
  dimensions $d_V^2=2^{-\IND V}$ resp.\ $d_{V'}^2=2^{-\IND V'}$, and
  the indices of $V$ and $V'$ must be equal.
  
  Next we show c) $\Rightarrow$ b). Let $u$ be a partial isometry with
  initial space $\ker V^*$, final space $\ker{V'}^*$, and $u=\4{u}$
  (such $u$ exists due to *--invariance and equality of dimensions of
  the kernels). Then $U\DEF V'V^*+u$ is an element of
  ${\II}^0(\KK)$ and fulfills $V'=UV$. We have to prove that
  $PU\4{P}$ is HS\@. But $u$ has finite rank, so it suffices to show
  that $A\DEF\4{P}VS_{V'}V^*\4{P}$ is of trace class. Since
  $S_V\4{S_V}$ and $S_{V'}\4{S_{V'}}$ have finite rank,
  $S_{V'}\4{S_{V}}+S_{V}\4{S_{V'}}=(S_{V'}-S_V)(\4{S_V}-\4{S_{V'}})
  +S_V\4{S_V}+S_{V'}\4{S_{V'}}$ is trace class. So the same is true
  for $A=AQ_V+AVV^*=AQ_V+\4{P}V(S_{V'}\4{S_{V}}+S_{V}\4{S_{V'}})V^*+
  \4{P}Q_VPV\4{S_{V'}}V^*$.

  b) $\Rightarrow$ a) is obvious.
\end{proof}
In order to apply part c) of the lemma, we need information about the
operators $S_V$. An orthogonal projection $E$ on $\KK$ is called a
{\em partial \BP} \label{page:PBP}\cite{A70} if $E\4{E}=0$. By
definition, the \emph{codimension} of $E$ is the dimension of
$\ker(E+\4{E})$. For instance, $VPV^*$ is a partial \BP with
codimension $2M_V=-\IND V$ for any $V\in{\II}(\KK)$. The following
lemma holds for arbitrary $S\in\QK$ (except for the formula for the
codimension, of course) as long as $S\4{S}$ has finite rank.
\begin{Lem}
  \label{lem:SV}\hspace*{\fill}\\
  Let $V\in\IFIN(\KK)$, and let $E_V$ denote the orthogonal projection
  onto $\ker S_V\4{S_V}$. Then $S_VE_V=E_VS_V$ is a partial \BP with
  finite codimension $2(M_V-N_V)$. Moreover, there exist
  $\lambda_1,\ldots,\lambda_r\in(0,\2)$, $r\leq M_V-N_V$, partial \BPs
  $E_1,\ldots,E_r$, and an orthogonal projection $E_\2=\4{E_\2}$ such
  that
  \begin{gather}
    E_V+E_\2+\sum_{j=1}^r(E_j+\4{E_j})=\1,\nonumber\\
    S_V=S_VE_V+\2E_\2+
      \sum_{j=1}^r\Big(\lambda_jE_j+(1-\lambda_j)\4{E_j}\Big).\label{SV}
  \end{gather}
\end{Lem}
\begin{proof}
  Since $S_V\4{S_V}=S_V-S_V^2$, $S_V$ commutes with $E_V$ and fulfills
  $S_VE_V=S_V^2E_V$ and $(S_VE_V)(\4{S_VE_V})=S_V\4{S_V}E_V=0$.  Hence
  $S_VE_V$ is a partial \BP. The dimension of $\ker(S_VE_V+\4{S_VE_V})
  =\ker E_V$ (the codimension of $S_VE_V$) equals the rank of
  $S_V\4{S_V}$. By $S_V\4{S_V}=V^*PQ_VPV$, the rank of $S_V\4{S_V}$ is
  equal to $\dim V^*P(\ker V^*)$. Now consider the decomposition
  $$\ker V^*= \ker v\oplus \ker\4{v}\oplus\Big(\ker V^*\ominus(\ker
    v\oplus\ker\4{v})\Big)$$
  with $v$ given by \eqref{V}. $V^*P$
  vanishes on $\ker v\oplus\ker\4{v}$, but the restriction of $V^*P$
  to $\ker V^*\ominus(\ker v\oplus\ker\4{v})$ is one to one since
  $V^*Pk=0=V^*k$ implies $V^*\4{P}k=0$, i.e.\ $k\in\ker
  v\oplus\ker\4{v}$. Hence the codimension of $S_VE_V$ equals
  $\dim(\ker V^*\ominus(\ker v\oplus\ker\4{v}))=-\IND V-2N_V$.
  
  Let $s_V$ denote the restriction of $S_V$ to $\RAN S_V\4{S_V}$.
  $s_V$ is a positive operator on a finite dimensional \HSP and has a
  complete set of eigenvectors with eigenvalues in $(0,1)$. If
  $\lambda$ is an eigenvalue of $s_V$, then $1-\lambda$ is also an
  eigenvalue (with the same multiplicity) due to $s_V+\4{s_V}=\1-E_V$.
  Thus there exist $\lambda_1,\ldots,\lambda_r\in(0,\frac{1}{2})$ and
  spectral projections $E_\frac{1}{2},E_1,\ldots,E_r$ with
  $\4{E_\frac{1}{2}}=E_\frac{1}{2},\ E_j\4{E_j}=0$ such that
  $E_\frac{1}{2}+\sum_{j=1}^r(E_j+\4{E_j})=\1-E_V$ and
  $s_V=\frac{1}{2}E_\frac{1}{2}+\sum_{j=1}^r(\lambda_jE_j+
  (1-\lambda_j)\4{E_j})$.  
\end{proof}
As a consequence, operators $S_V$ with $M_V=\2$ necessarily have
the form $S_V=S_VE_V+\2E_\2$ where $E_\2=\1-E_V$ has rank one. By
taking direct sums of $V\in{\II}^1(\KK)$ with operators $V(\varphi)$
from Example~\ref{ex:EX1} below, we see that any combination of
eigenvalues and multiplicities that is allowed by Lemma~\ref{lem:SV}
actually occurs for some $S_{V'}$. Therefore any $S\in\QK$ such that
$S\4{S}$ has finite rank is of the form $S=S_V$ for some
$V\in\IFIN(\KK)$.

We further remark that a quasi--free state $\omega_S$ with $S$ of the
form \eqref{SV} is a product state\footnote{A state $\omega$ is a {\em
    product state} \WRT a decomposition $\KK=\oplus_j\KK_j$ of $\KK$
  into closed, *--invariant subspaces if
  $\omega(ab)=\omega(a)\omega(b)$ whenever $a\in\mathfrak{C}(\KK_j),\ 
  b\in\mathfrak{C}(\KK_j^\bot)$.} as defined by Powers \cite{P67} (see
also \cite{MRT,Man}) \WRT the decomposition $\KK=\ker S\4{S}\oplus\RAN
E_\2\oplus\bigoplus_j\RAN(E_j+\4{E_j})$.  Clearly, the restriction of
$\omega_S$ to $\mathfrak{C}(\ker S\4{S})$ is a \FS, the
restriction to $\mathfrak{C}(\RAN E_\2)$ the central state.
\begin{Ex}
  \label{ex:EX1}\hspace*{\fill}\\
  Let $(f_n)_{n\in\NN}$ be an \ONB for $P(\KK)$, and set $E_n\DEF
  f_n\langle f_n,.\,\rangle$, with $f_n^+\DEF(f_n+f_n^*)/\sqrt{2},\ 
  f_n^-\DEF i(f_n-f_n^*)/\sqrt{2}$. Then $(f_n^s)_{s=\pm,\,n\in\NN}$
  is an \ONB for $\KK$ consisting of *--invariant vectors. For
  $\varphi\in\RR$, define a Bogoliubov operator
  \begin{equation*}
    \begin{split}
      V(\varphi)\DEF\ &(f_0^+\cos\varphi+f_1^-\sin\varphi)\langle
        f_0^+,.\,\rangle+(f_0^-\sin\varphi-f_1^+\cos\varphi)\langle
        f_0^-,.\,\rangle\\
      &+\sum_{s=\pm,\,n\geq1}f_{n+1}^s\langle f_n^s,.\,\rangle.
    \end{split}
  \end{equation*}
  Then $V(\varphi)\in{\II}^2(\KK)$, and the eigenvalue
  $\lambda_\varphi=\2(1+\sin2\varphi)$ of
  $S_{V(\varphi)}=\lambda_\varphi E_0+(1-\lambda_\varphi)\4{E_0}
  +\sum_{n\geq1}E_n$ assumes any value in $[0,1]$ as $\varphi$ varies
  over $[-\pi/4,\pi/4]$.
\end{Ex}

Next we characterize the Bogoliubov operators $V$ for which $S_V$
takes a particularly simple form. A distinction arises between the
cases of even and odd Fredholm index. Note that the parity of $-\IND
V$ is equal to the parity of $\dim\ker(S_V-\2)$ (it is well--known
that quasi--free states $\omega_S$ with $\dim\ker(S-\2)$ even resp.\ 
odd behave differently (cf.\ \cite{MY})).
\begin{Lem}
  \label{lem:PURE}\hspace*{\fill}\\
  \indent a) Let $W\in{\II}(\KK)$. Then the following conditions are
  equivalent:\\ 
  \indent\indent(i) $\omega_P\0\rho_W$ is a pure state;\\
  \indent\indent(ii) $S_W$ is a \BP;\\
  \indent\indent(iii) $[P,WW^*]=0$.\\
  \indent If any of these conditions is fulfilled, then $M_W=N_W$ and
  $\pi_P\0\rho_W\simeq d_W\cdot\pi_{S_W}$.

  b) For any \BP $P'$ and $m\in\NN\cup\{\infty\}$, there exists
  $W\in{\II}^{2m}(\KK)$ with\\
  \indent\phantom{b) }$S_W=P'$.

  c) Let $W\in\IFIN(\KK)$. Then the following
  conditions are equivalent:\\ 
  \indent\indent(i) $\omega_P\0\rho_W$ is a mixture of two disjoint pure
    states;\\  
  \indent\indent(ii) $S_WE_W$ is a partial \BP with codimension $1$;\\
  \indent\indent(iii) $[P,WW^*]$ has rank $2$;\\
  \indent\indent(iv) $M_W=N_W+\2$.

  d) For any partial \BP $P'$ with codimension $1$ and
  $m\in\NN\cup\{\infty\}$,\\ 
  \indent\phantom{d) }there exists $W\in{\II}^{2m}(\KK)$ with
  $S_WE_W=P'$.
\end{Lem}
\begin{proof}
  a) We know from Section~\ref{sec:QFCAR} that $\omega_P\0\rho_W$ is
  pure \IFF $S_W$ is a projection. We have
  $$S_W^2=S_W\iff W^*PQ_WPW=0\iff Q_WPWW^*=0\iff[P,WW^*]=0.$$
  If this
  is fulfilled, then $\ker WW^*=\ker(PWW^*P)\oplus
  \ker(\4{P}WW^*\4{P})$ has dimension $2M_W=2N_W$. By
  Proposition~\ref{prop:DECCAR}, $\pi_P\0\rho_W$ is the direct sum of
  $2^{N_W}=d_W$ irreducible subrepresentations, each equivalent to the
  Fock representation $\pi_{S_W}$.
  
  b) Let $m$ and $P'$ be given. There clearly exists
  $W'\in\II^{2m}(\KK)$ with $[P,W']=0$. Since $\II^{0}(\KK)$ acts
  transitively on the set of \BPs, we may choose $U\in\II^{0}(\KK)$
  with $U^*PU=P'$. Then $W\DEF W'U$ has the desired properties.
  
  c) (ii) $\Leftrightarrow$ (iii) follows from the facts that the
  codimension of $S_WE_W$ equals the rank of $WW^*PQ_W$ (cf.\ the
  proof of Lemma~\ref{lem:SV}) and that $[P,WW^*]=Q_WPWW^*-WW^*PQ_W$.
  (ii) is equivalent to (iv) by virtue of Lemma~\ref{lem:SV}. (ii)
  $\Rightarrow$ (i) has been shown by Araki~\cite{A70}. To prove (i)
  $\Rightarrow$ (iv), assume that $M_W>N_W+\2$ (if $M_W=N_W$, then
  $S_W$ is a \BP and $\omega_P\0\rho_W$ pure). Then one can show that
  $\omega_P\0\rho_W$ is a mixture of two quasi--equivalent orthogonal
  states (see \cite{CB1}), hence cannot be a mixture of two disjoint
  pure states. This proves part c).

  d) Let $(f_n)_{n\in\NN}$ be an \ONB for $P(\KK)$, $(g_n)_{n\geq1}$
  an \ONB for $P'(\KK)$, and $g_0$ a unit vector in $\ker(P'+\4{P'})$.
  Set 
  $$V\DEF f_0^+\langle g_0,.\,\rangle+\sum_{s=\pm,\,n\geq1}
  f_n^s\langle g_n^s,.\,\rangle$$
  (we used the notation of
  Example~\ref{ex:EX1}). Then $V\in\II^{1}(\KK)$ and $S_V=\2g_0\langle
  g_0,.\,\rangle+P'$. This implies $S_VE_V=P'$, and if we choose $W'$
  as in the proof of b), then $W\DEF W'V$ has the desired properties.
\end{proof}
Now we are in a position to discuss the decomposition of
representations $\pi_P\0\rho_V$ with $V\in\IFIN(\KK)$.  If $-\IND V$
is even (resp.\ odd), then $S_VE_V$ is a partial \BP with even (odd)
codimension by Lemma~\ref{lem:SV}, and there exists a \BP (partial \BP
with codimension $1$) $P'$ such that $P'-S_V$ is Hilbert--Schmidt (we
may choose $P'$ to coincide with $S_VE_V$ on $\ker S_V\4{S_V}$; then
$P'-S_V$ has finite rank). By Lemma~\ref{lem:PURE}, there exists $W$
with $\IND W=\IND V$ and $S_WE_W=P'$, and Lemma~\ref{lem:EQ} implies
$\pi_P\0\rho_V\simeq\pi_P\0\rho_W$. The latter representation splits
into $2^{N_W}$ copies of the GNS representation $\pi_{S_W}$ for the
state $\omega_P\0\rho_W$ by Proposition~\ref{prop:DECCAR}. If $-\IND
V$ is even, then $\pi_{S_W}=\pi_{P'}$ and $2^{N_W}=d_V$. If $-\IND V$
is odd, then $\pi_{S_W}=\pi^+\oplus\pi^-$ where $\pi^\pm$ are mutually
inequivalent, irreducible, so--called \emph{pseudo Fock
  representations} by virtue of a lemma of Araki (see \cite{A70} for
details), and $2^{N_W}=2^{-\2}d_V$.

Summarizing, we rediscover B\"ockenhauer's first result on
representations of the \CAR \cite{JMB2}:
\begin{Th}
  \label{th:DEC}\hspace*{\fill}\\
  Let $P$ be a \BP and $V\in\IFIN(\KK)$. If $-\IND V$ is even, then
  there exist \BPs $P'$ such that $P'-S_V$ is Hilbert--Schmidt, and
  for each such $P'$
  $$\pi_P\0\rho_V\simeq d_V\cdot\pi_{P'}.$$
  If $-\IND V$ is odd, there
  exist partial \BPs $P'$ with codimension 1 such that $P'-S_V$ is
  Hilbert--Schmidt. For each such $P'$,
  $$\pi_P\0\rho_V\simeq 2^{-\2}d_V\cdot(\pi_{P'}^+\oplus\pi_{P'}^-)$$
  where $\pi_{P'}^\pm$ are the \textup{(}inequivalent,
  irreducible\textup{)} pseudo Fock representations induced by $P'$.
\end{Th}
Using the isomorphisms $\sigma_{V_1}:\CK\to\CK_0,\ 
V_1\in{\II}^1(\KK)$, from \eqref{EVEN}, we can immediately see
how the restriction of $\pi_P\0\rho_V$ to $\CK_0$ decomposes. From
$$\pi_P\0\rho_V(\CK_0)=\pi_P\0\rho_V\0\sigma_{V_1}(\CK)\simeq
\pi_P\0\rho_{VV_1}(\CK)$$
and $M_{VV_1}=M_V+\2,\ d_{VV_1}=\sqrt{2}d_V$, we infer
B\"ockenhauer's second result \cite{JMB2}:
\begin{Cor}
  \label{cor:DEC}\hspace*{\fill}\\
  Let $P$ be a \BP and $V\in\IFIN(\KK)$. If $-\IND V$
  is even, then $\pi_P\0\rho_V|_{\CK_0}$ is equivalent to a multiple
  of the direct sum of two inequivalent irreducible representations,
  with multiplicity $d_V$. If $-\IND V$ is odd, then
  $\pi_P\0\rho_V|_{\CK_0}$ is equivalent to a multiple of an
  irreducible representation, with multiplicity $\sqrt{2}d_V$.
\end{Cor}
\subsection{The semigroup of implementable endomorphisms}
\label{sec:IP}
From now on we choose a fixed \BP $P_1$, and we set
$P_2\DEF\1-P_1=\4{P_1}$. Let us introduce the following notation.
The components of an operator $A\in\BB(\KK)$ are denoted by
$$A_{mn}\DEF P_mAP_n,\quad m,n=1,2$$
and are regarded as operators
from $\KK_n\DEF P_n(\KK)$ to $\KK_m$. Thus $\ker A_{mn}$, $(\ker
A_{mn})^\bot$, $(\RAN A_{nm})^\bot$ etc.\ are viewed as subspaces of
$\KK_n$, and relations like
$${A_{mn}}^*={A^*}_{nm},\quad\4{A_{11}}=\4{A}_{22},\quad\text{etc.}$$
will frequently be used.  We also use matrix notation
$A=\MAT{A_{11}}{A_{12}}{A_{21}}{A_{22}}$ \WRT the decomposition
$\KK=\KK_1\oplus\KK_2$. $A$ is called \emph{antisymmetric} if
\begin{equation}
  \label{AS}
  A^\tau\DEF\4{A^*}=-A.
\end{equation}
If \HH is a subspace of \KK, then $\HH^*$ will denote the complex
conjugate space (and not the dual space)
$$\HH^*\DEF\SET{f^*}{f\in\HH}.$$
Thus one has e.g.\ ${\KK_1}^*=\KK_2$.
The reader is kindly asked to pay attention to the various meanings of
the star ``*''. The correct one should always be clear from the
context.

Let $V\in\II_{P_1}(\KK)=\SET{V\in{\II}(\KK)}{V_{12}\text{ is
    Hilbert--Schmidt}}$. The relation $V^*V=\1$ reads in components
\begin{subequations}
  \begin{align}
    {V_{11}}^*V_{11}+{V_{21}}^*V_{21}&=P_1,\label{REL1}\\
    {V_{22}}^*V_{22}+{V_{12}}^*V_{12}&=P_2,\label{REL2}\\
    {V_{11}}^*V_{12}+{V_{21}}^*V_{22}&=0,\label{REL3}\\
    {V_{22}}^*V_{21}+{V_{12}}^*V_{11}&=0,\label{REL4}
  \end{align}
\end{subequations}
and the relation $\4{V}=V$ entails that
$$\4{V_{11}}=V_{22},\quad\4{V_{12}}=V_{21}.$$
${V_{11}}^*V_{12},\ {V_{21}}^*V_{22},\ {V_{22}}^*V_{21}$ and
${V_{12}}^*V_{11}$ are antisymmetric by \eqref{REL3} and \eqref{REL4}.

Since $V_{12}$ is a Hilbert--Schmidt operator by
Theorem~\ref{th:IMPCAR}, ${V_{22}}^*V_{22}$ is Fredholm (with
vanishing index) by \eqref{REL2}. This implies
\begin{equation}
  \label{LVFIN}
  L_V\DEF\dim\ker V_{22}=\dim\ker{V_{22}}^*V_{22}<\infty.
\end{equation}
Note that $V_{12}|_{\ker V_{22}}$ is isometric and, by \eqref{REL3}, 
\begin{equation}
  \label{KER}
  V_{12}(\ker V_{22})\subset\ker{V_{11}}^*\cap\RAN V.
\end{equation}
As mentioned after Theorem~\ref{th:IMPCAR}, ${V_{11}}$ is
semi--Fredholm with $-\IND{V_{11}}=M_V$. These observations imply
\begin{equation}
  \label{DIMKER}
  \dim(\ker{V_{11}}^*\ominus V_{12}(\ker V_{22}))=M_V.
\end{equation}

The semigroup of implementable quasi--free endomorphisms is isomorphic
to the semigroup ${\II}_{P_1}(\KK)$. The latter is a topological
semigroup relative to the metric (cf.~\cite{A87})
$$\delta_{P_1}(V,V')\DEF\NORM{V-V'}+\HSNORM{V_{12}-V'_{12}},$$
where $\HSNORM{\ }$ denotes the Hilbert--Schmidt
norm. ${\II}_{P_1}(\KK)$ contains the closed sub--semigroup of
diagonal Bogoliubov operators (which are also called \emph{gauge
  invariant,} because they commute with the \TT--action \eqref{U1}):
\begin{equation}
  \label{IDIAG}
  \IDIAG(\KK)\DEF\SET{W\in\II(\KK)}{[P_1,W]=0}.
\end{equation}
$\IDIAG(\KK)$ is isomorphic to the semigroup of isometries
of $\KK_1$, via the map $V\mapsto V_{11}$. 

The restricted orthogonal group $\II_{P_1}^0(\KK)$ has a natural
normal subgroup
$$\IHS(\KK)\DEF\SET{U\in\II(\KK)}{U-\1\text{ is Hilbert--Schmidt}}.$$
An automorphism $\rho_U,\ U\in\II_{P_1}^0(\KK)$, is implementable in
\emph{all} Fock representations \IFF $U\in\IHS(\KK)$ or
$-U\in\IHS(\KK)$ \cite{SS65,A70}. The group $\IHS(\KK)$ is also
related to the group of quasi--free automorphisms which are weakly
inner in the representation associated with the central state
$\omega_\2$ \cite{Blat}.

We shall prove that each $V\in\II_{P_1}(\KK)$ can be written as a
product 
\begin{equation}
  \label{UVW}
  V=UW,\quad\text{with }U\in\IHS(\KK),\ W\in\IDIAG(\KK)
\end{equation}
(it is known that each $V\in\II_{P_1}^0(\KK)$ has this form
\cite{COB}). Suppose for the moment that such $U$ and $W$ exist. Then
$P\DEF UP_1U^*$ is a \BP which extends the partial \BP $VP_1V^*$ such
that
\begin{equation}
  \label{PVCAR}
  P_1-P\text{ is Hilbert--Schmidt},\qquad V^*PV=P_1.
\end{equation}
The corresponding \FS $\omega_{P}$ is unitarily equivalent to
$\omega_{P_1}$ and fulfills $\omega_{P}\0\rho_V=\omega_{P_1}$. The
proof of the product decomposition $V=UW$ involves the construction of
such \BPs $P$. So let us start with a
parameterization\footnote{Similar ideas can be found in the book of
  Pressley and Segal \cite{PrSe}.} of the set
\begin{equation}
  \label{PP1DEF}
  \PP_{P_1}\DEF\SET{P\in\QK}{P^2=P,\ P_1-P\text{ is
      Hilbert--Schmidt}}
\end{equation}
of \BPs of \KK which differ from $P_1$ only by a Hilbert--Schmidt
operator. $\PP_{P_1}$ is isomorphic to the set of all \FSs which are
equivalent to $\omega_{P_1}$. Let $\mathfrak{H}_{P_1}$ be the \HSP of
all antisymmetric (see \eqref{AS}) Hilbert--Schmidt operators from
$\KK_1$ to $\KK_2$
\begin{equation}
  \label{HP1}
  \mathfrak{H}_{P_1}\DEF\SET{T\in\BB(\KK_1,\KK_2)}{T^\tau=-T,\
    T\text{ is Hilbert--Schmidt}},
\end{equation}
let $\mathfrak{F}_{P_1}$ be the set of all finite dimensional
subspaces of $\KK_1$
$$\mathfrak{F}_{P_1}\DEF\SET{\mathfrak{h}\subset\KK_1}{\mathfrak{h}\text{
    is a finite dimensional subspace}},$$
and let
\begin{equation}
  \label{PP1}
  \tilde\PP_{P_1}\DEF\SET{(T,\mathfrak{h})\in\mathfrak{H}_{P_1}
    \times\mathfrak{F}_{P_1}}{\mathfrak{h}\subset\ker T}.
\end{equation}
Then the following holds.
\begin{Prop}
  \label{prop:PP1}\hspace*{\fill}\\
  The map 
  \begin{align}
    &P\mapsto(P_{21}{P_{11}}^{-1},\ker P_{11})\nonumber\\
    \intertext{is a bijection from $\PP_{P_1}$ onto
      $\tilde\PP_{P_1}$, with inverse given by} 
    &(T,\mathfrak{h})\mapsto P_{(T,\mathfrak{h})}\DEF P_T-p_\mathfrak{h}
      +\4{p_\mathfrak{h}}.\label{PTH}
  \end{align}
  Here ${P_{11}}^{-1}\in\BB(\KK_1)$ is defined as the inverse of
  $P_{11}$ on the closed subspace $\RAN P_{11}$ and as zero on the
  finite dimensional subspace $\ker P_{11}$. $P_T$ is the
  \textup{(}basis\textup{)} projection onto $\RAN(P_1+T)$
  \begin{equation} 
    \label{PT}
    P_T\DEF (P_1+T)(P_1+T^*T)^{-1}(P_1+T^*),
  \end{equation}
  and $p_\mathfrak{h}$ is the \textup{(}partial basis\textup{)}
  projection onto $\mathfrak{h}\subset\KK_1$.
\end{Prop}
\begin{proof}
  Let us rewrite the conditions on $P$ to be a \BP
  $$P=P^*=\1-\4{P}=P^2$$
  in components:
  \begin{subequations}
    \begin{align}
      P_{11}={P_{11}}^*&=P_1-\4{P_{22}},\label{BP1}\\
      P_{22}={P_{22}}^*&=P_2-\4{P_{11}},\label{BP2}\\
      P_{21}={P_{12}}^*&=-{P_{21}}^\tau,\label{BP3}\\
      P_{11}-{P_{11}}^2&={P_{21}}^*P_{21},\label{BP4}\\
      P_{22}-{P_{22}}^2&={P_{12}}^*P_{12},\label{BP5}\\
      (P_1-P_{11})P_{12}&=P_{12}P_{22},\label{BP6}\\
      (P_2-P_{22})P_{21}&=P_{21}P_{11},\label{BP7}
    \end{align}
  \end{subequations}
  and note that, because $P_1-P=\4{P_{22}}-P_{21}-{P_{21}}^*-P_{22}$,
  $$P_1-P\text{ is Hilbert--Schmidt (HS)}\iff P_2P\text{ is HS}.$$
  
  Now let $P\in\PP_{P_1}$. Then $P_{21}$ is HS, $P_{22}$ is of
  trace class, and $P_{11}$ is Fredholm with index zero by
  \eqref{BP1}.  Therefore $\RAN P_{11}$ is closed and $\ker P_{11}$
  has finite dimension. It follows that ${P_{11}}^{-1}$ is a
  well--defined bounded operator such that
  $P_{11}{P_{11}}^{-1}={P_{11}}^{-1}P_{11}$ is the projection onto
  $\RAN P_{11}$. We have $\ker P_{11}\subset\ker P_{21}$ by
  \eqref{BP4}, hence by \eqref{BP1}--\eqref{BP3} and \eqref{BP7}
  \begin{equation*}
    \begin{split}
      P_{21}{P_{11}}^{-1}+(P_{21}{P_{11}}^{-1})^\tau
        &=P_{21}{P_{11}}^{-1}-\4{{P_{11}}}^{\,-1}P_{21}\\ 
      &=\4{{P_{11}}}^{\,-1}\bigl(\4{P_{11}}P_{21}-P_{21}P_{11}\bigr)
        {P_{11}}^{-1}\\
      &=\4{{P_{11}}}^{\,-1}\bigl((P_2-P_{22})P_{21}
        -P_{21}P_{11}\bigr){P_{11}}^{-1}\\
      &=0,
    \end{split} 
  \end{equation*}
  so $P_{21}{P_{11}}^{-1}$ is antisymmetric in the sense of
  \eqref{AS}. This proves
  $P_{21}{P_{11}}^{-1}\in\mathfrak{H}_{P_1}$. Since by definition
  $\ker{P_{11}}^{-1}=\ker P_{11}$, it follows that
  $(P_{21}{P_{11}}^{-1},\ker P_{11})\in\tilde\PP_{P_1}$.
  
  Conversely, let $(T,\mathfrak{h})\in\tilde\PP_{P_1}$ be given. We
  associate with $T$ the bounded operator
  \begin{equation}
    \label{X}
    X\DEF(P_T)_{11}=(P_1+T^*T)^{-1}    
  \end{equation}
  so that, by \eqref{PT}
  $$P_T=(P_1+T)X(P_1+T^*)=
  \begin{pmatrix}
    X&XT^*\\
    TX&TXT^*
  \end{pmatrix}.$$
  $P_T$ is a projection because $(P_1+T^*)(P_1+T)=X^{-1}$. To prove
  that $P_T+\4{P_T}=\1$ holds, note that $TX^{-1}=\4{X}^{\,-1}T$ and
  therefore $\4{X}T=TX$, $XT^*=T^*\4{X}$. It follows that
  \begin{equation*}
    \begin{split}
      P_T+\4{P_T}&=(P_1+T)X(P_1+T^*)+(P_2-T^*)\4{X}(P_2-T)\\
      &=X+TX+XT^*+TT^*\4{X}+\4{X}-XT^*-TX+T^*TX\\
      &=(P_1+T^*T)X+(P_2+TT^*)\4{X}\\
      &=P_1+P_2\\
      &=\1. 
    \end{split} 
  \end{equation*}
  Thus $P_T$ is a \BP. It is obvious that $\RAN P_T=\RAN(P_1+T)$.
  Since $P_TP_2$ is HS, $P_T$ belongs to $\PP_{P_1}$. The
  condition $T\mathfrak{h}=0$ ensures that $p_\mathfrak{h}\subset
  P_T$.  Therefore
  $P_{(T,\mathfrak{h})}=P_T-p_\mathfrak{h}+\4{p_\mathfrak{h}}$ is a
  \BP. It is also contained in $\PP_{P_1}$ because
  $p_\mathfrak{h}$ has finite rank.
  
  To show that the two maps $(T,\mathfrak{h})\mapsto
  P_{(T,\mathfrak{h})}$ and $P\mapsto(P_{21}{P_{11}}^{-1},\ker
  P_{11})$ are each other's inverse, let first
  $(T,\mathfrak{h})\in\tilde\PP_{P_1}$ be given and set $P\DEF
  P_{(T,\mathfrak{h})}$. Since $X=(P_T)_{11}$ is bijective and since
  $T\mathfrak{h}=0$, we have
  \begin{equation*}
    \begin{split}
      \ker P_{11}&=\ker X^{-1}P_{11} =\ker X^{-1}(X-p_\mathfrak{h})
      =\ker\big(P_1-(P_1+T^*T)p_\mathfrak{h}\big)\\
      &=\ker(P_1-p_\mathfrak{h})=\mathfrak{h}.
    \end{split}
  \end{equation*}
  By $T\mathfrak{h}=0$, and because $(X-p_\mathfrak{h})
  (X-p_\mathfrak{h})^{-1}$ is the projection onto $\RAN
  P_{11}=\mathfrak{h}^\bot$,
  $$P_{21}{P_{11}}^{-1}=(P_T)_{21}\big((P_T)_{11}-p_\mathfrak{h}\big)^{-1}
  =TX(X-p_\mathfrak{h})^{-1}=T(X-p_\mathfrak{h})(X-p_\mathfrak{h})^{-1}=T.$$
  
  Conversely, let $P\in\PP_{P_1}$ be given and set $T\DEF
  P_{21}{P_{11}}^{-1},\ \mathfrak{h}\DEF\ker P_{11}$. Then we have,
  using \eqref{BP4}, 
  $$P_1+T^*T={P_{11}}^{-1}+P_1-{P_{11}}^{-1}P_{11}={P_{11}}^{-1}
  +p_\mathfrak{h}=(P_{11}+p_\mathfrak{h})^{-1}.$$ 
  Together with $T\mathfrak{h}=0$, \eqref{BP3} and \eqref{BP2} we get
  \begin{equation*}
    \begin{split}
      P_{(T,\mathfrak{h})}&=(P_1+T)(P_1+T^*T)^{-1} 
        (P_1+T^*)-p_\mathfrak{h}+\4{p_\mathfrak{h}}\\
      &=(P_1+T)(P_{11}+p_\mathfrak{h})(P_1+T^*)-p_\mathfrak{h}
        +\4{p_\mathfrak{h}}\\
      &=P_{11}+TP_{11}+P_{11}T^*+TP_{11}T^*+p_\mathfrak{h}-p_\mathfrak{h}
        +\4{p_\mathfrak{h}}\\
      &=P_{11}+P_{21}+P_{12}+TT^*\4{P_{11}}+\4{p_\mathfrak{h}}\\
      &=P-P_{22}+TT^*\4{P_{11}}+\4{p_\mathfrak{h}}\\
      &=P-(P_2-\4{p_\mathfrak{h}})+\4{P_{11}}^{-1}\4{P_{11}}\\
      &=P.
    \end{split} 
  \end{equation*}
  This completes the proof.
\end{proof}
\begin{Rem}
  \label{page:KERP}
  Note that 
  \begin{equation}
    \label{P11}
    \ker P_{11}=\KK_1\cap\ker P.
  \end{equation}
  The \BPs of the form
  $P_{(0,\mathfrak{h})},\ \mathfrak{h}\in\mathfrak{F}_{P_1}$, are
  precisely the elements of $\PP_{P_1}$ which commute with $P_1$,
  and the projections of the form $P_T=P_{(T,\{0\})},\ 
  T\in\mathfrak{H}_{P_1}$, are precisely the elements
  $P\in\PP_{P_1}$ with $\ker P_{11}=\{0\}$. For a general element
  $P_{(T,\mathfrak{h})}\in\PP_{P_1}$, it is well known that the
  unique (up to a phase) cyclic vector in $\mathcal{F}_a(\KK_1)$ which
  induces the state $\omega_{P_{(T,\mathfrak{h})}}$ is proportional to
  \begin{equation}
    \label{OPT}
    a(e_1)^*\dotsm a(e_L)^*\exp(\3\4{T}a^*a^*)\Omega_{P_1}
  \end{equation}
  where the $e_j$ form an \ONB in $\mathfrak{h}$, and the exponential
  term will be explained later on.
\end{Rem}
\begin{Lem}
  \label{lem:UTH}\hspace*{\fill}\\
  For $T\in\mathfrak{H}_{P_1}$, define a Bogoliubov operator
  \begin{equation*}
    \begin{split}
      U_T&\DEF
      (P_1+T)(P_1+T^*T)^{-\2}+(P_2-T^*)(P_2+TT^*)^{-\2}\in\IHS(\KK).
    \end{split}
  \end{equation*}
  For $\mathfrak{h}\in\mathfrak{F}_{P_1}$, choose an \ONB
  $\{e_1,\dots,e_L\}$ in $\mathfrak{h}$, define the partial isometry
  $u_\mathfrak{h}\DEF\sum_{r=1}^L e_r^*\langle e_r,.\,\rangle$ with
  initial space $\mathfrak{h}$ and final space $\mathfrak{h}^*$, and
  define a self--adjoint Bogoliubov operator
  \begin{equation}
    \label{UH}
    U_\mathfrak{h}\DEF\1-\begin{pmatrix}
      p_\mathfrak{h}&\4{u_\mathfrak{h}}\\
      u_\mathfrak{h}&\4{p_\mathfrak{h}}
    \end{pmatrix}\in\IHS(\KK).
  \end{equation}
  Then one has, in the notation of Proposition~\ref{prop:PP1},
  $$U_TP_1U_T^*=P_T,\qquad U_\mathfrak{h}P_1U_\mathfrak{h}^*=
  P_{(0,\mathfrak{h})}.$$
  If $(T,\mathfrak{h})\in\tilde\PP_{P_1}$, one has in addition
  $$[U_T,U_\mathfrak{h}]=0,\qquad
  U_TU_\mathfrak{h}P_1U_\mathfrak{h}^*U_T^*=P_{(T,\mathfrak{h})}.$$
  It
  follows that the action $P\mapsto UPU^*$ of the restricted
  orthogonal group $\II^0_{P_1}(\KK)$ on $\PP_{P_1}$ restricts to a
  transitive action of $\IHS(\KK)$.
\end{Lem}
\begin{proof}
  Let $T\in\mathfrak{H}_{P_1}$ be given. The unitary $U_T$ results
  from polar decomposition of $\1+T+\4{T}=U_T\ABS{\1+T+\4{T}}$. With
  $X$ defined by \eqref{X}, $U_T$ can be written as
  $$U_T=
    \begin{pmatrix}
      X^\2 & \4{T}\,\4{X}^\2\\
      TX^\2 & \4{X}^\2
    \end{pmatrix}.$$
  It is straightforward to see that $U_T$ is a Bogoliubov operator
  which transforms $P_1$ into $P_T$. To prove that $U_T-\1$ is HS,
  note that
  $$X^\2-P_1=X^\2(P_1-X^{-1})(P_1+X^{-\2})^{-1}
  =-X^\2T^*T(P_1+X^{-\2})^{-1}$$
  is of trace class. Therefore
  $(U_T-\1)P_1=(P_1+T)X^\2-P_1=X^\2-P_1+TX^\2$ is HS, which implies that
  $U_T-\1=(U_T-\1)P_1+\4{(U_T-\1)P_1}$ is HS.
  
  Now let $\mathfrak{h}\in\mathfrak{F}_{P_1}$. Then $U_\mathfrak{h}$
  is clearly a Bogoliubov operator with $U_\mathfrak{h}-\1$ of finite
  rank and with $U_\mathfrak{h}P_1U_\mathfrak{h}^*=
  P_{(0,\mathfrak{h})}$. $U_\mathfrak{h}$ is therefore contained in
  $\IHS(\KK)$. It is self--adjoint because $u_\mathfrak{h}$ is
  symmetric: $u_\mathfrak{h}^\tau=u_\mathfrak{h}$. (Actually,
  self--adjointness of $U_\mathfrak{h}$ or symmetry of
  $u_\mathfrak{h}$ will not be needed in the sequel. The above
  definitions aim at reducing the ambiguity in the choice of
  $U_\mathfrak{h}$. $u_\mathfrak{h}$ is now determined up to the
  action of the unitary group of $\mathfrak{h}$.)
  
  If $(T,\mathfrak{h})\in\tilde\PP_{P_1}$, then one has
  $Tp_\mathfrak{h}=0$ and, by functional calculus,
  $X^\2p_\mathfrak{h}=p_\mathfrak{h}X^\2=p_\mathfrak{h}$. Using this,
  one gets by straightforward computations
  \begin{align}
    &U_TU_\mathfrak{h}=U_\mathfrak{h}U_T=U_T+U_\mathfrak{h}-\1=    
      \begin{pmatrix}
        X^\2p_{\mathfrak{h}^\bot}
          &-T^*\4{X}^\2\4{p_{\mathfrak{h}^\bot}}-\4{u_\mathfrak{h}}\\
        TX^\2 p_{\mathfrak{h}^\bot}-u_\mathfrak{h} 
          &\4{X}^\2\4{p_{\mathfrak{h}^\bot}}
      \end{pmatrix},\label{UTH1}\\
    &U_TU_\mathfrak{h}P_1U_\mathfrak{h}^*U_T^*
      =U_TP_{(0,\mathfrak{h})}U_T^*=U_\mathfrak{h}P_T
      U_\mathfrak{h}^*=P_{(T,\mathfrak{h})}.\label{UTH2}
  \end{align}
  Here $p_{\mathfrak{h}^\bot}$ denotes the orthogonal projection onto
  $\mathfrak{h}^\bot\subset\KK_1$. It is obvious that
  $\II^0_{P_1}(\KK)$ acts transitively from the left on
  $\PP_{P_1}$ as indicated. By \eqref{UTH2} and by
  $U_TU_\mathfrak{h}\in\IHS(\KK)$, $\IHS(\KK)$ acts already
  transitively.
\end{proof}
\begin{Rem}
  $T$ and $\mathfrak{h}$ can be recovered from $U\DEF
  U_TU_\mathfrak{h}$ as $\mathfrak{h}=\ker U_{11},\ 
  T=U_{21}{U_{11}}^{-1}$ (see below for the definition of
  ${U_{11}}^{-1}$).
\end{Rem}

Next we would like to assign to each $V\in\II_{P_1}(\KK)$
distinguished operators $P_V,\ U_V$ and $W_V$ having the properties
stated in \eqref{PVCAR} and \eqref{UVW}. Note that any \BP $P$
fulfilling \eqref{PVCAR} has the form
$$P=VP_1V^*+q$$
where $q$ is a partial \BP (see page~\pageref{page:PBP}) with 
$$q+\4{q}=Q_V\DEF\1-VV^*.$$
It follows from this and from \eqref{KER} that necessarily
\begin{equation}
  \label{KERP}
  V_{12}(\ker V_{22})\subset\ker P_{11}.
\end{equation}
In general, $\dim\ker P_{11}$ can take any value between $L_V$
(defined in \eqref{LVFIN}) and $L_V+M_V$. Similarly, if $U$ is a
Bogoliubov operator such that $[U^*V,P_1]=0$ (cf.~\eqref{UVW}), then
\begin{equation}
  \label{KERU}
  V_{12}(\ker V_{22})\subset\ker {U_{11}}^*
\end{equation}
(this follows from $0=P_1U^*VP_2={U_{11}}^*V_{12}+{U_{21}}^*V_{22}$).
We shall choose $U_V$ and $P_V$ such that equality holds in
\eqref{KERP} and \eqref{KERU}.

Now let ${V_{11}}^{-1}\in\BB(\KK_1)$ be defined as the inverse of
$V_{11}$ on the closed subspace $\RAN{V_{11}}$ and as zero on $\ker
{V_{11}}^*$ (and define $V_{22}$ analogously). For a closed subspace
$\HH\subset\KK$, let $p_\HH$ be the orthogonal projection onto \HH.
Then
\begin{alignat}{2}
  &\RAN{V_{11}}^{-1}=\RAN{V_{11}}^*,
    &\qquad\qquad&V_{11}{V_{11}}^{-1}=p_{\RAN V_{11}},\label{INVREL1}\\
  &\ker{V_{11}}^{-1}=\ker{V_{11}}^*,
    &&{V_{11}}^{-1}V_{11}=p_{{\RAN V_{11}}^*}.\label{INVREL2}
\end{alignat}
Define $(T_V,\mathfrak{h}_V)\in\tilde\PP_{P_1}$ by
\begin{align}
  T_V&\DEF V_{21}{V_{11}}^{-1}-{V_{22}}^{-1*}{V_{12}}^*p_{{\ker
  V_{11}}^*},\label{TV}\\
  \mathfrak{h}_V&\DEF V_{12}(\ker V_{22}).\label{HV}
\end{align}
By Proposition~\ref{prop:PP1}, any \BP $P\in\PP_{P_1}$ with $\ker
P_{11}=V_{12}(\ker V_{22})$ has the form $P=P_{(T,\mathfrak{h}_V)}$
for some $T\in\mathfrak{H}_{P_1}$. The possible choices of $T$ such
that \eqref{PVCAR} holds are determined in
\begin{Lem}
  \label{lem:TV}\hspace*{\fill}\\
  Let $\mathfrak{h}_V$ be defined by \eqref{HV}. A \BP
  $P_{(T,\mathfrak{h}_V)}\in\PP_{P_1}$ satisfies
  \begin{equation}
    \label{PV}
    V^*P_{(T,\mathfrak{h}_V)}V=P_1,
  \end{equation}
  i.e.\ extends the partial \BP $VP_1V^*$, \IFF 
  \begin{equation}
    \label{T'}
    T=T_V+T',
  \end{equation}
  where $T'\in\mathfrak{H}_{P_1}$ is an operator from
  $\ker{V_{11}}^*\ominus\mathfrak{h}_V$ to
  $\ker{V_{22}}^*\ominus(\mathfrak{h}_V)^*$. For such $T$ one has
  \begin{equation}
    \label{THS}
    \HSNORM{T}^2=\HSNORM{T_V}^2+\HSNORM{T'}^2.
  \end{equation}
\end{Lem}
\begin{proof}
  The formula
  $$p_{\mathfrak{h}_V}=V_{12}p_{\ker V_{22}}{V_{12}}^*$$
  entails that
  $V^*p_{\mathfrak{h}_V}V={V_{12}}^*V_{12}p_{\ker
    V_{22}}{V_{12}}^*V_{12}=p_{\ker V_{22}}$ by \eqref{REL2},
  \eqref{KER}. It follows from \eqref{PTH} that
  $$V^*P_{(T,\mathfrak{h}_V)}V=V^*P_TV-p_{\ker V_{22}}+p_{\ker
    V_{11}}.$$
  Therefore \eqref{PV} is equivalent to
  $$V^*P_TV=P_{(0,\ker V_{11})}\DEF p_{\RAN{V_{11}}^*}+p_{\ker
    V_{22}}$$
  or, by Lemma~\ref{lem:PURE}a, to $P_TV=VP_{(0,\ker
    V_{11})}$. This is further equivalent to $0=\4{P_T}VP_{(0,\ker
    V_{11})}$ and, since $\ker\4{P_T}=\ker(P_2-T)$, to
  $0=(P_2-T)VP_{(0,\ker V_{11})}=V_{21}p_{\RAN{V_{11}}^*}-TV_{11}-
  TV_{12}p_{\ker V_{22}}$. Looking at the components, we finally
  obtain the following conditions, which are equivalent to \eqref{PV}
  \begin{equation}
    \label{T}
    \mathfrak{h}_V\subset\ker T,\qquad
    TV_{11}=V_{21}p_{\RAN{V_{11}}^*}.
  \end{equation}
  (Of course, the first relation of \eqref{T} is also necessary for
  having $(T,\mathfrak{h}_V)\in\tilde\PP_{P_1}$, cf.\ \eqref{PP1}.)
  Let us show that $T_V$ is a special solution of this problem. $T_V$
  is Hilbert--Schmidt since $V_{21}$ and $V_{12}$ are. It is
  antisymmetric because
  \begin{equation*}
    \begin{split}
      T_V+T_V^\tau&=V_{21}{V_{11}}^{-1}-{V_{22}}^{-1*} 
        {V_{12}}^*p_{{\ker V_{11}}^*}+{V_{22}}^{-1*}{V_{12}}^*-
        p_{{\ker V_{22}}^*}V_{21}{V_{11}}^{-1}\\
      &=p_{\RAN V_{22}}V_{21}{V_{11}}^{-1}+{V_{22}}^{-1*}{V_{12}}^*
        p_{\RAN V_{11}}\\
      &={V_{22}}^{-1*}({V_{22}}^*V_{21}+{V_{12}}^*V_{11})
        {V_{11}}^{-1}\\
      &=0
    \end{split}
  \end{equation*}
  (we used \eqref{REL4} and \eqref{INVREL1}). Thus $T_V$ belongs to
  $\mathfrak{H}_{P_1}$. One clearly has
  $T_VV_{11}=V_{21}p_{\RAN{V_{11}}^*}$ and $T_V\mathfrak{h}_V=0$ (use
  \eqref{KER}, \eqref{INVREL1} and \eqref{REL2} to see the latter).
  Thus $T_V$ solves \eqref{T}. 
  
  If $T\in\mathfrak{H}_{P_1}$ is another solution of \eqref{T}, then
  $T'\DEF T-T_V$ is contained in $\mathfrak{H}_{P_1}$. \eqref{T} is
  equivalent to $\RAN V_{11}\oplus \mathfrak{h}_V\subset\ker T'$,
  i.e.\ to $(\ker T')^\bot\subset\ker{V_{11}}^*\ominus\mathfrak{h}_V$.
  By antisymmetry of $T'$, this is also equivalent to $\RAN
  T'\subset\ker{V_{22}}^* \ominus(\mathfrak{h}_V)^*$.
  
  Finally, \eqref{THS} holds because $T_V$ and $T'$ are orthogonal as
  elements of the \HSP $\mathfrak{H}_{P_1}$.
\end{proof}
\begin{Rem}
  \label{page:TV}
  $T_V$ takes the simpler form $T_V=V_{21}{V_{11}}^{-1}$, which is
  well--known from the case of automorphisms ($\IND V=0$), whenever
  $[P_1,VV^*]=0$, i.e.\ whenever the state $\omega_{P_1}\0\rho_V$ is
  pure (cf.\ Lemma~\ref{lem:PURE}a). One can show that $T_V=0$ \IFF
  $V_{21}{V_{11}}^*=0$ \IFF one (and hence, in view of \eqref{REL1},
  \eqref{REL2}, all) of the operators $V_{11},V_{21},V_{12},V_{22}$ is
  a partial isometry.
\end{Rem}
Having specified $(T_V,\mathfrak{h}_V)$ in \eqref{TV}, \eqref{HV}, let
us now associate with $V\in\II_{P_1}(\KK)$ the following operators:
\begin{align}
  P_V&\DEF P_{(T_V,\mathfrak{h}_V)},\label{PV1}\\
  U_V&\DEF U_{T_V}U_{\mathfrak{h}_V},\label{UV}\\
  W_V&\DEF {U_V}^*V,\label{WV}
\end{align}
and let us collect their properties.
\begin{Prop}
  \label{prop:VUW}\hspace*{\fill}\\
  $P_V$ belongs to $\PP_{P_1}$ and satisfies $V^*P_VV=P_1$. It is
  chosen such that $\ker(P_V)_{11}$ and
  $\HSNORM{(P_V)_{21}{(P_V)_{11}}^{-1}}$ are minimal \textup{(}cf.\ 
  Proposition~\ref{prop:PP1} and \eqref{KERP}, \eqref{THS}\textup{)}.
  $U_V$ belongs to $\IHS(\KK)$ and $W_V$ to $\IDIAG(\KK)$; their
  definition depends on the choice of an \ONB in $\mathfrak{h}_V$.
  The operators $U_V$ and $W_V$ fulfill
  \begin{align}
    M_{U_V}&=0, & T_{U_V}&=T_V,
      & \mathfrak{h}_{U_V}&=\mathfrak{h}_V,& P_{U_V}&=P_V;\label{UV1}\\
    M_{W_V}&=M_V, & T_{W_V}&=0,
      & \mathfrak{h}_{W_V}&=\{0\}, & P_{W_V}&=P_1,\label{WV1}
  \end{align}
  \begin{equation}
    \label{WV2}
    W_V=\begin{pmatrix}
      (P_1+{T_V}^*T_V)^\2V_{11}-\4{u_{\mathfrak{h}_V}}V_{21}&0\\
      0&(P_2+T_V{T_V}^*)^\2V_{22}-u_{\mathfrak{h}_V}V_{12}
    \end{pmatrix},
  \end{equation}
  and 
  \begin{equation}
    \label{VUW}
    V=U_VW_V.
  \end{equation}
  If $M_V=0$, then the formulas reduce to
  \begin{equation}
    \label{UVWV}
    U_V=\begin{pmatrix}
      \ABS{{V_{11}}^*} & V_{12}{v_{22}}^*-\4{u_{\mathfrak{h}_V}} \\ 
      V_{21}{v_{11}}^*-u_{\mathfrak{h}_V} & \ABS{{V_{22}}^*}
    \end{pmatrix},\quad
    W_V=\begin{pmatrix}
      v_{11}-\4{u_{\mathfrak{h}_V}}V_{21}&0\\
      0&v_{22}-u_{\mathfrak{h}_V}V_{12}
    \end{pmatrix}
  \end{equation}
  where $v_{11}$ and $v_{22}$ are the partial isometries appearing in
  the polar decomposition of $V_{11}=v_{11}\ABS{V_{11}}$ and
  $V_{22}=v_{22}\ABS{V_{22}}$. On the other hand, if
  $V\in\IDIAG(\KK)$, then $U_V=\1$ and $W_V=V$.
\end{Prop}
\begin{Rem}
  Equipped with the metric induced by $\HSNORM{\ }$, $\PP_{P_1}$
  becomes a topological space. It consists of two connected components
  which are distinguished by the parity of $\dim\ker P_{11}$. It
  follows that the map $V\mapsto P_V$ is \emph{not} continuous on
  $\II_{P_1}(\KK)$, because $\ker(P_V)_{11}=\mathfrak{h}_V$, and
  $(-1)^{\dim\mathfrak{h}_V}$ is not constant on the connected
  components of $\II_{P_1}(\KK)$ (see Corollary~\ref{cor:CONN} and
  Example~\ref{ex:EX2} below).
\end{Rem}
\begin{proof}
  The assertions made about $P_V$ have been proved in
  Proposition~\ref{prop:PP1} and Lemma~\ref{lem:TV}. The claims
  concerning $U_V$ are implied by Lemma~\ref{lem:UTH} and the remark
  following that lemma. $W_V$ is diagonal because
  $P_1W_V=P_1{U_V}^*V={U_V}^*P_VV={U_V}^*VP_1=W_VP_1$, where we used
  Lemma~\ref{lem:UTH} and an argument from the proof of
  Lemma~\ref{lem:TV}. The formula \eqref{WV2} can be derived from
  \eqref{UTH1} together with the relation
  $V_{11}+{T_V}^*V_{21}=X_V^{-1}V_{11},\ X_V\DEF(P_1+{T_V}^*T_V)^{-1}$.
  $V=U_VW_V$ holds by definition of $W_V$ and entails that $\IND
  W_V=\IND V$. The remaining statements on $W_V$ clearly follow from
  $W_V\in\IDIAG(\KK)$.

  If $\IND V=0$, then one has $\mathfrak{h}_V=\ker {V_{11}}^*$ (cf.\
  \eqref{DIMKER}) and $X_V=p_{\ker {V_{11}}^*}+\ABS{{V_{11}}^*}^2$ so
  that $X_V^\2p_{\mathfrak{h}_V^\bot}=\ABS{{V_{11}}^*}$. \eqref{UVWV}
  then follows by straightforward computation from \eqref{UTH1}.

  If $V$ is diagonal, then $T_V=0$ and $\mathfrak{h}_V=\{0\}$, so that
  obviously $U_V=\1$ and $W_V=V$.
\end{proof}
\begin{Rem}
  A different product decomposition $V=\tilde U\tilde W$ was
  established in \cite{CB1}, essentially by polar decomposition of
  $V_{11}$:
  $$\tilde W=\begin{pmatrix}
    v_{11}&V_{12}p_{\ker V_{22}}\\
    V_{21}p_{\ker V_{11}}&v_{22}
  \end{pmatrix},\qquad\tilde U=V\tilde W^*+\dotsb$$
  $\tilde W$ is obviously not diagonal (unless $\ker V_{11}=\{0\}$),
  but has off--diagonal components of finite rank. (However, using the
  operators $U_\mathfrak{h}$ from \eqref{UH}, the definition of
  $\tilde W$ could easily be modified to make $\tilde W$ diagonal.)
  $\tilde W\tilde W^*$ commutes with $P_1$. $\tilde U$ is in general
  not contained in $\IHS(\KK)$, but only in $\II^0_{P_1}(\KK)$. These
  operators have nevertheless some useful properties, similar to $U_V$
  and $W_V$:
  \begin{align*}
    M_{\tilde U}&=0, & T_{\tilde U}&=T_V,
      & \mathfrak{h}_{\tilde U}&=\{0\},& P_{\tilde U}&=P_{T_V};\\  
    M_{\tilde W}&=M_V, & T_{\tilde W}&=0,
      & \mathfrak{h}_{\tilde W}&=\mathfrak{h}_V, & P_{\tilde
      W}&=P_{(0,\mathfrak{h}_V)}. 
  \end{align*}
\end{Rem}
\begin{Cor}
  \label{cor:IP}\hspace*{\fill}\\
  $\II_{P_1}(\KK)=\IHS(\KK)\cdot\IDIAG(\KK)$. The
  $\II^0_{P_1}(\KK)$--orbits in $\II_{P_1}(\KK)$ with respect to left
  multiplication are precisely the sets $\II^{2m}_{P_1}(\KK),\ 
  m\in\NN\cup\{\infty\}$.
\end{Cor}
\begin{proof}
  The product decomposition of $\II_{P_1}(\KK)$ has been obtained
  above. To show that $\II^0_{P_1}(\KK)$ acts transitively on each
  $\II^{2m}_{P_1}(\KK)$, let $V,V'\in\II^{2m}_{P_1}(\KK)$ be given,
  with decompositions $V=UW,\ V'=U'W'$ as in \eqref{VUW}. Since $P_1$
  leaves $\ker{W'}^*$ and $\ker W^*$ invariant, we can choose a
  partial isometry $\hat u$ with initial space $\ker{W'}^*$ and final
  space $\ker W^*$ such that $\4{\hat u}=\hat u$ and $[P_1,\hat u]=0$.
  Then $\hat U\DEF W{W'}^*+\hat u\in\II^0_{P_1}(\KK)$ fulfills $\hat
  UW'=W$.  This implies that $(U\hat U{U'}^*)V'=V$.
\end{proof}

These results can be used to determine the connected components of the
semigroup $\II_{P_1}(\KK)$. It is known that the restricted orthogonal
group $\II^0_{P_1}(\KK)\subset \II_{P_1}(\KK)$ has two connected (and
simply connected) components $\II^0_{P_1}(\KK)^\pm$
\cite{Ca84,PrSe,A87}. Namely,
$$\chi(U)\DEF(-1)^{\dim\ker U_{11}}=(-1)^{\dim\mathfrak{h}_U}$$
defines a continuous character $\chi$ on $\II^0_{P_1}(\KK)$, and
$\chi|_{\II^0_{P_1}(\KK)^\pm}=\pm1$. (This character is equal to the
Araki--Evans index of the pair of \BPs $(P_1,UP_1U^*)$ \cite{AE,A87}.)
We shall see that $\II^{2m}_{P_1}(\KK)$ is connected if $m>0$, and
that the map $\chi:V\mapsto(-1)^{\dim\ker V_{11}}=
(-1)^{\dim\mathfrak{h}_V}$ remains neither multiplicative nor
continuous when extended to the whole semigroup $\II_{P_1}(\KK)$.  We
need the following preparatory result.
\begin{Lem}
  \label{lem:CONN}\hspace*{\fill}\\
  The set of all isometries with a given fixed index on an
  infinite dimensional complex \HSP is arcwise connected in the
  uniform topology.
\end{Lem}
\begin{proof}
  Let $V,V'$ be two isometries with $\IND V=\IND V'$. Since $\dim\ker
  V^*=\dim\ker{V'}^*$, there exists a unitary operator $U$ with
  $V'=UV$ (choose a partial isometry $u$ with initial space $\ker V^*$
  and final space $\ker{V'}^*$, and set $U\DEF V'V^*+u$). Since the
  unitary group is arcwise connected, there exists a continuous curve
  $U(t)$ of unitary operators with $U(0)=\1$ and $U(1)=U$. Then
  $V(t)\DEF U(t)V$ is a continuous curve of isometries with $V(0)=V$
  and $V(1)=V'$.
\end{proof}
\begin{Cor}
  \label{cor:CONN}\hspace*{\fill}\\
  The connected components of $\II_{P_1}(\KK)$ are precisely the subsets
  $\II^0_{P_1}(\KK)^\pm$ and $\II^{2m}_{P_1}(\KK),\ 1\leq m\leq\infty$.
\end{Cor}
\begin{proof}
  Let $V,V'\in\II^{2m}_{P_1}(\KK)$ with $\chi(V)=\chi(V')$ be given,
  and let
  $$V=UW,\quad V'=U'W'$$
  be decompositions as in \eqref{VUW}. It follows from \eqref{UV1}
  that
  $$\chi(U)=\chi(V)=\chi(V')=\chi(U')$$
  so that there exists a continuous curve in $\II^0_{P_1}(\KK)$
  connecting $U$ to $U'$. Since $W$ and $W'$ are both diagonal and
  have index $-2m$, there exists a continuous curve in
  $\II^{2m}_{P_1}(\KK)$ connecting $W$ to $W'$ by
  Lemma~\ref{lem:CONN}. It follows that $UW=V$ and $U'W'=V'$ can be
  connected by a continuous curve in $\II^{2m}_{P_1}(\KK)$. Therefore
  either of the two subsets
  $$\II^{2m}_{P_1}(\KK)^\pm\DEF
    \SET{V\in\II^{2m}_{P_1}(\KK)}{\chi(V)=\pm1}$$
  is arcwise connected.
  Below, we give an example of a continuous curve in
  $\II^{2m}_{P_1}(\KK)$ which connects $\II^{2m}_{P_1}(\KK)^+$ to
  $\II^{2m}_{P_1}(\KK)^-$. Hence $\II^{2m}_{P_1}(\KK)$ itself is
  connected.
\end{proof}
\begin{Ex}
  \label{ex:EX2}\hspace*{\fill}\\
  Let $V(\varphi)$ be the Bogoliubov operator introduced in
  Example~\ref{ex:EX1} in Section~\ref{sec:REP} (with $P=P_1$). Then
  $V(\varphi)\in\II^2_{P_1}(\KK)$ since
  ${V(\varphi)_{12}}^*V(\varphi)_{12}=(1-\lambda_\varphi)\4{E_0}$ has
  finite rank, and $\varphi\mapsto V(\varphi)$ is a continuous curve
  in $\II^2_{P_1}(\KK)$. We have $\ker V(\varphi)_{11}=
  \ker(\lambda_\varphi E_0+\sum_{n\geq1}E_n)$, hence
  $$\chi(V(\varphi))=\begin{cases}
    1, & \varphi\notin(4\ZZ+3)\frac{\pi}{4}\\
    -1, & \varphi\in(4\ZZ+3)\frac{\pi}{4}
  \end{cases}.$$
  Now let $V\in\II^{2m-2}_{P_1}(\KK)$ with $[P_1,V]=0$. Then
  $\chi(VV(\varphi))=\chi(V(\varphi))$ since $V_{11}$ is isometric, so
  $\varphi\mapsto VV(\varphi)$ is a continuous curve in
  $\II^{2m}_{P_1}(\KK)$ which connects $\II^{2m}_{P_1}(\KK)^+$ to
  $\II^{2m}_{P_1}(\KK)^-$. This completes the proof of
  Corollary~\ref{cor:CONN}.
  
  $V(\varphi)$ may also serve to illustrate that $\chi$ is not
  multiplicative on $\II_{P_1}(\KK)$. Define a Bogoliubov operator
  \begin{equation*}
    \begin{split} 
      U\DEF\ &\tfrac{1}{\sqrt{2}}f_0^+\langle f_0^++f_1^-,.\,\rangle
        -\tfrac{1}{\sqrt{2}} f_1^+ \langle f_0^--f_1^+,.\,\rangle
        +\tfrac{1}{\sqrt{2}}f_0^-\langle f_0^-+f_1^+,.\,\rangle\\
      &-\tfrac{1}{\sqrt{2}}f_1^-\langle f_0^+-f_1^-,.\,\rangle
        +\sum_{n\geq2}(E_n+\4{E_n}).
      \end{split}
  \end{equation*}
  Then $U\in\II^0_{P_1}(\KK)$, and a calculation shows that
  $U_{11}=\frac{1}{\sqrt{2}}(E_0+E_1)+\sum_{n\geq2}E_n$ and
  $UV(\frac{3\pi}{4})=V(\frac{\pi}{2})$. This entails
  $$1=\chi\Big(UV(\mbox{$\frac{3\pi}{4}$})\Big)\neq
  \chi(U)\chi\Big(V(\mbox{$\frac{3\pi}{4}$})\Big)=-1$$
  since $\ker U_{11}=\ker(V(\frac{\pi}{2})_{11})=\{0\}$, but $\ker
  (V(\frac{3\pi}{4})_{11})=\CC f_0$. We finally note that the
  eigenvalues $\pm(1-\lambda_\varphi)$ of
  $P_1-S_{V(\varphi)}=(1-\lambda_\varphi)(E_0-\4{E_0})$ have
  multiplicity one if $\lambda_\varphi\neq1$, in contrast to the
  unitary case where the multiplicities of eigenvalues in $(0,1)$ are
  always even~\cite{AE,A87,EK}.
\end{Ex}
\subsection{Normal form of implementers}
\label{sec:FORM}
Unitary operators which implement quasi--free automorphisms of the
\CAR have been constructed by several authors, notably by Friedrichs
\cite{Fri}, Berezin \cite{Ber}, Schroer, Seiler and Swieca \cite{SSS},
Labont\'e \cite{Lab}, Fredenhagen \cite{F77}, Klaus and Scharf
\cite{KS}, Ruijsenaars \cite{R77a,R78}. Our construction of isometric
implementers for quasi--free endomorphisms follows Ruijsenaars'
approach in \cite{R78} which is to our knowledge the most complete
treatment of the implementation problem for quasi--free automorphisms.
Another advantage of \cite{R78} for our purposes is the (implicit) use
of Araki's ``selfdual'' formalism.

Let us begin with a generalization of the definition of ``\BHs'' from
the case of trace class operators to the case of bounded operators.
Bilinear Hamiltonians have been introduced by Araki \cite{A68} as
infinitesimal generators of one--parameter groups of inner Bogoliubov
automorphisms. More specifically, one may assign to a finite rank
operator $H=\sum_jf_j\langle g_j,.\,\rangle$ on \KK the \BH
$$b(H)\DEF\sum_jf_jg_j^*$$
and extend $b$ by continuity to a linear
map from the ideal of all trace class operators on \KK to \CK. If a
trace class operator $H$ satisfies $\4{H}=H$ and $H^\tau=-H$, then
$\2b(H)$ is the generator of the one--parameter group
$(\rho_{\exp(tH)})_{t\in\RR}$:
$$\rho_{\exp(tH)}(a)=\exp\big(\3tb(H)\big)a\exp\big(-\3tb(H)\big),
\qquad a\in\CK.$$
The map $H\mapsto\2b(H)$ is an isomorphism from the Lie
algebra formed by all such $H$ onto the Lie algebra of the spin group.
See \cite{A70,A87} for details.

Since the elements $f\in\KK_1$ correspond to creation operators in the
Fock representation $\pi_{P_1}$, we may write
$$\pi_{P_1}\big(b(H)\big)=H_{11}a^*a+H_{12}a^*a^*+H_{21}aa+H_{22}aa^*$$
where the terms on the right are defined by
$H_{11}a^*a\DEF\pi_{P_1}(b(H_{11}))$ etc. Introducing \emph{Wick
  ordering} by $\WO{a(f)a(g)^*}=-a(g)^*a(f)$, we get
\begin{align}
  \WO{H_{22}aa^*}&=-{H_{22}}^\tau a^*a=H_{22}aa^*-(\TR H_{22})\1,\\
  \WO{\pi_{P_1}\big(b(H)\big)}&=(H_{11}-{H_{22}}^\tau)a^*a+H_{12}a^*a^*
    +H_{21}aa.\label{WICK1}
\end{align}
According to \cite{R78,CR}, one can define such Wick ordered
expressions for \emph{bounded} $H$ as follows. Assume from now on
that (without loss of generality) 
$$\KK_1=L^2(\RR^d),$$
and let $\SS\subset\mathcal{F}_a(\KK_1)$ be the
dense subspace consisting of finite particle vectors $\phi$ with
$n$--particle wave functions $\phi^{(n)}$ in the Schwartz space
$\SS(\RR^{dn})$. For $p\in\RR^d$, the unsmeared annihilation operator
$a(p)$ with (invariant) domain $\SS$ is defined by
$$\big(a(p)\phi\big)^{(n)}(p_1,\dots,p_n)\DEF\sqrt{n+1}\phi^{(n+1)}
(p,p_1,\dots,p_n).$$
Since $a(p)$ is not closable, one defines
$a(p)^*$ as the \QF adjoint of $a(p)$ on $\SS\times\SS$.
Then Wick ordered monomials $a(q_m)^*\dotsm a(q_1)^*a(p_1)\dotsm
a(p_n)$ are well--defined \QFs on $\SS\times\SS$, and for
$\phi,\phi'\in\SS$,
$$\langle\phi,a(q_m)^*\dotsm a(q_1)^*a(p_1)\dotsm a(p_n)\phi'\rangle
\DEF\langle a(q_1)\dotsm a(q_m)\phi,a(p_1)\dotsm a(p_n)\phi'\rangle$$
is a function in $\SS(\RR^{d(m+n)})$ to which tempered distributions
may be applied. For example, one has in the \QF sense
\begin{align*}
  a(f)&=\int\4{f(p)}a(p)\,dp,\\
  a(f)^*&=\int f(p)a(p)^*\,dp,\quad f\in\KK_1.
\end{align*}

Now let $H$ be a bounded operator on \KK which is antisymmetric in the
sense of \eqref{AS}\footnote{The \BH corresponding to a symmetric
  operator vanishes.}. Then there exist tempered distributions
$H_{mn}(p,q),\ m,n=1,2$, given by
\begin{align*}
  \langle f,H_{11}g\rangle &=\int\4{f(p)}H_{11}(p,q)g(q)\,dp\,dq,\\
  \langle f,H_{12}g^*\rangle &=\int\4{f(p)}H_{12}(p,q)\4{g(q)}\,dp\,dq,\\
  \langle f^*,H_{21}g\rangle &=\int f(p)H_{21}(p,q)g(q)\,dp\,dq,\\
  \langle f^*,H_{22}g^*\rangle &=\int f(p)H_{22}(p,q)\4{g(q)}\,dp\,dq, 
    \quad f,g\in\SS(\RR^d)\subset\KK_1.
\end{align*}
Hence the following expressions are \QFs on $\SS\times\SS$
\begin{align*}
  H_{11}a^*a &\DEF\int a(p)^*H_{11}(p,q)a(q)\,dp\,dq,\\
  H_{12}a^*a^* &\DEF\int a(p)^*H_{12}(p,q)a(q)^*\,dp\,dq,\\
  H_{21}aa &\DEF\int a(p)H_{21}(p,q)a(q)\,dp\,dq,\\
  \WO{H_{22}aa^*} &\DEF-\int a(q)^*H_{22}(p,q)a(p)\,dp\,dq
    \ =\ H_{11}a^*a.
\end{align*}
Wick ordering of $H_{22}aa^*$ is necessary to make this expression
well--defined. The last equality follows from antisymmetry of $H$:
$$H_{11}(p,q)=-H_{22}(q,p),\quad H_{12}(p,q)=-H_{12}(q,p),\quad
H_{21}(p,q)=-H_{21}(q,p).$$
The Wick ordered \BH induced by $H$ is then defined in analogy to
\eqref{WICK1} as 
$$\WO{b(H)}\DEF H_{12}a^*a^*+2H_{11}a^*a+H_{21}aa;$$
it is linear in $H$. We define its Wick ordered powers as
\begin{equation} 
  \label{WICK4}
  \WO{b(H)^l}\DEF l!\sum_{\substack{l_1,l_2,l_3=0\\ l_1+l_2+l_3=l}}^l
  \frac{1}{l_1!l_2!l_3!}
  (H_{12})^{l_1}(2H_{11})^{l_2}(H_{21})^{l_3}
  a^{*2l_1+l_2}a^{l_2+2l_3}
\end{equation}
where the terms on the right hand side are \QFs on
$\SS\times\SS$ (cf.~\cite{R78})
\begin{equation}
  \label{MONOM}
  \begin{split}
    (H_{12})^{l_1}(&H_{11})^{l_2}(H_{21})^{l_3}
      a^{*2l_1+l_2}a^{l_2+2l_3}\\ 
    \DEF\int &H_{12}(p_1,q_1)\dotsm H_{12}(p_{l_1},q_{l_1})
      H_{11}(p_1',q_1')\dotsm H_{11}(p_{l_2}',q_{l_2}')\\
    \cdot\, &H_{21}(p_1'',q_1'')\dotsm H_{21}(p_{l_3}'',
      q_{l_3}'')a(p_1)^*\dotsm a(p_{l_1})^*a(q_{l_1})^*\dotsm 
      a(q_1)^*\\
    \cdot\, &a(p_1')^*\dotsm a(p_{l_2}')^*a(q_{l_2}')\dotsm 
      a(q_1')a(p_1'')\dotsm a(p_{l_3}'')a(q_{l_3}'')\dotsm a(q_1'')\\
    \cdot\, &dp_1\,dq_1\ldots dp_{l_1}\,dq_{l_1}\,dp_1'\,dq_1'\ldots 
      dp_{l_2}'\,dq_{l_2}'\,dp_1''\,dq_1''\ldots dp_{l_3}''\,dq_{l_3}''.
  \end{split}
\end{equation}
Finally, we define the Wick ordered exponential
\begin{equation}
  \label{WICK5}
  \EH\DEF\sum_{l=0}^\infty
  \frac{1}{l!2^l}\WO{b(H)^l}
\end{equation}
which is also a well--defined \QF on $\SS\times\SS$, since
the sum in \eqref{WICK5} is finite when applied to vectors
$\phi,\phi'\in\SS$.

By Ruijsenaars' result \cite{R78}, \EH is the \QF of a unique linear
operator, defined on the dense subspace \DD of algebraic tensors in
$\mathcal{F}_a(\KK_1)$, provided that $H_{12}$ is Hilbert--Schmidt.
(This is equivalent to $\4{H_{12}}\in\mathfrak{H}_{P_1}$.) In this
case, the series \eqref{WICK5} converges strongly on \DD, \EH (viewed
as an operator) maps \DD into the dense subspace of
$C^\infty$--vectors for the number operator, and
\begin{align}
  \EH\Omega_{P_1}&=\exp(\3H_{12}a^*a^*)\Omega_{P_1},\label{EOM}\\
  \bigNORM{\exp(\3H_{12}a^*a^*)\Omega_{P_1}}&=
    \big({\det}_{\KK_1}(P_1+H_{12}{H_{12}}^*)\big)^{1/4}.
    \label{OMNORM}
\end{align}
Let us compute the commutation relations of the operators \EH with
creation and annihilation operators.
\begin{Lem}
  \label{lem:REL}\hspace*{\fill}\\
  Let $H\in\BB(\KK)$ be antisymmetric with $H_{12}$ Hilbert--Schmidt.
  For $f,g\in\KK_1$, the following relations hold on \DD
  \begin{align*}
    [\EH,a(f)^*]\ &=\ a(H_{11}f)^*\EH+\EH a\bigl((H_{21}f)^*\bigr),\\
    [\EH,a(g)]\ &=\ a(H_{12}g^*)^*\EH-\EH a({H_{11}}^*g).
  \end{align*}
\end{Lem}
\begin{proof}
  It is a lengthy but straightforward exercise in anticommutation
  relations to calculate the commutation relations of Wick monomials
  of the form \eqref{MONOM} with creation and annihilation operators:
  \begin{align*}
    [H_{l_1,l_2,l_3},a(f)^*]\ &=\ l_2a(H_{11}f)^*H_{l_1,l_2-1,l_3}
      +2l_3H_{l_1,l_2,l_3-1}a\bigl((H_{21}f)^*\bigr),\\
    [H_{l_1,l_2,l_3},a(g)]\ &=\ 2l_1a(H_{12}g^*)^*H_{l_1-1,l_2,l_3}-
      l_2H_{l_1,l_2-1,l_3} a({H_{11}}^*g),
  \end{align*}
  where $H_{l_1,l_2,l_3}\DEF(H_{12})^{l_1}(H_{11})^{l_2}(H_{21})^{l_3}
  a^{*2l_1+l_2}a^{l_2+2l_3}$. From this one obtains
  \begin{equation*}
    \begin{split}
    [\EH&,a(f)^*]=\sum_{l=0}^\infty 2^{-l}\sum_{l_1+l_2+l_3=l}
      \frac{2^{l_2}}{l_1!l_2!l_3!}[H_{l_1,l_2,l_3},a(f)^*]\\
    =&\ a(H_{11}f)^* \sum_{l=1}^\infty 2^{-(l-1)}\sum_{l_1+l_2+l_3=l}
      \frac{2^{l_2-1}}{l_1!(l_2-1)!l_3!} H_{l_1,l_2-1,l_3}\\
    &\ +\sum_{l=1}^\infty 2^{-(l-1)}\sum_{l_1+l_2+l_3=l}
      \frac{2^{l_2}}{l_1!l_2!(l_3-1)!}H_{l_1,l_2,l_3-1}
      a\bigl((H_{21}f)^*\bigr)\\
    =&\ a(H_{11}f)^*\EH+\EH a\bigl((H_{21}f)^*\bigr)\\
    \intertext{and}
    [\EH&,a(g)]=\sum_{l=0}^\infty 2^{-l}\sum_{l_1+l_2+l_3=l}
      \frac{2^{l_2}}{l_1!l_2!l_3!}[H_{l_1,l_2,l_3},a(g)]\\
    =&\ a(H_{12}g^*)^*\sum_{l=1}^\infty 2^{-(l-1)}\sum_{l_1+l_2+l_3=l}
      \frac{2^{l_2}}{(l_1-1)!l_2!l_3!}H_{l_1-1,l_2,l_3}\\
    &\ -\sum_{l=1}^\infty 2^{-(l-1)}\sum_{l_1+l_2+l_3=l}
      \frac{2^{l_2-1}}{l_1!(l_2-1)!l_3!}H_{l_1,l_2-1,l_3}
      a({H_{11}}^*g)\\
    =&\ a(H_{12}g^*)^*\EH-\EH a({H_{11}}^*g).
    \end{split}
  \end{equation*}
\end{proof}
From now on let a fixed $V\in\II_{P_1}(\KK)$ be given. To construct
implementers for $\rho_V$, we have to look for antisymmetric operators
$H$ with $\4{H_{12}}\in\mathfrak{H}_{P_1}$ and (cf.~\eqref{INTER})
\begin{equation}
  \label{INT}
  \EH\pi_{P_1}(f)=\pi_{P_1}(Vf)\EH,\qquad f\in\KK_1\oplus\RAN{V_{22}}^*
\end{equation}
on \DD. Note that \eqref{INT} cannot be fulfilled for nonzero
$f\in\ker V_{22}$ since for such $f$, $\pi_{P_1}(f)=a(f^*)$ is an
annihilation operator and $\pi_{P_1}(Vf)=a(V_{12}f)^*$ a creation
operator, so that the left hand side (but not the right hand side) of
\eqref{INT} vanishes on $\Omega_{P_1}$. This defect can be cured by
``filling up the Dirac sea'' corresponding to $\mathfrak{h}_V$
(cf.~\eqref{OPT} and \eqref{PAV}) if we impose the following relation
for vectors in $\ker V_{22}$:
\begin{equation}
  \label{INTER3}
  \EH\pi_{P_1}(g^*)=0,\qquad g\in\ker V_{22}.
\end{equation}
It turns out that the solutions $H$ of \eqref{INT} and \eqref{INTER3}
are in one--to--one correspondence with the operators $T$ described in
Lemma~\ref{lem:TV}:
\begin{Lem}
  \label{lem:HT}\hspace*{\fill}\\
  The antisymmetric solutions $H$ of \eqref{INT} and \eqref{INTER3}
  with $H_{12}$ Hilbert--Schmidt are precisely the operators of the
  form
  \begin{equation}
    \label{H}
    H=\begin{pmatrix}
      V_{11}-P_1+T^*V_{21} & -T^* \\
      ({V_{22}}^*-{V_{12}}^*T^*)V_{21} & P_2-{V_{22}}^*+{V_{12}}^*T^*
    \end{pmatrix},
  \end{equation}
  where $T\in\mathfrak{H}_{P_1}$ fulfills \eqref{T}, i.e.\ is of the
  form \eqref{T'}.
\end{Lem}
\begin{proof}
  First note that a Wick ordered expression of the form $a(f)^*\EH+\EH
  a(g)$ vanishes \IFF $f$ and $g$ both vanish. In fact, application to
  the vacuum gives $a(f)^*\exp(\2H_{12}a^*a^*)\Omega_{P_1}$ which is
  zero \IFF $f=0$ (to see this, look for instance at the one--particle
  component).  Similarly, $\EH a(g)a(g)^*\Omega_{P_1}=\NORM{g}^2_\KK
  \exp(\2H_{12}a^*a^*)\Omega_{P_1}$ vanishes \IFF $g=0$.
  
  Hence we get all solutions of \eqref{INT} and \eqref{INTER3} if we
  write these equations in Wick ordered form and then compare term by
  term. We have by Lemma~\ref{lem:REL} and by the definition
  \eqref{PI} of $\pi_{P_1}$, with the shorthand $\eta_H\DEF\EH$,
  \begin{align*}
    \eta_H\pi_{P_1}(f) &=a\big((P_1+H_{11})f\big)^*\eta_H
      +\eta_H a\big((P_1+\4{H_{21}})f^*\big),\\
    \pi_{P_1}(Vf)\eta_H &=a\big((P_1-H_{12})Vf\big)^*\eta_H
      +\eta_H a\big((P_1-\4{H_{22}})Vf^*\big),\\
    \eta_H\pi(f)^* &=a\big((P_1+H_{11})f^*\big)^*\eta_H
      +\eta_H a\big(\4{H_{21}}f\big),\qquad f\in\KK.
  \end{align*}
  Thus \eqref{INT} is equivalent to
  \begin{subequations}
    \begin{align}
      P_1+H_{11}+(H_{12}-P_1)V(P_1+p_{\RAN{V_{22}}^*}) &=0,
      \label{H1}\\
      p_{\RAN{V_{22}}^*}+H_{21}+(H_{22}-P_2)V(P_1+p_{\RAN{V_{22}}^*}) 
      &=0,\label{H2}
    \end{align}
    and \eqref{INTER3} is equivalent to
    \begin{align}
      (P_1+H_{11})p_{\ker V_{11}} &=0,\label{H3}\\
      H_{21}p_{\ker V_{11}} &=0,\label{H4}
    \end{align}
  \end{subequations}
  where complex conjugation was occasionally applied. Note that
  \eqref{H1}--\eqref{H4} are equivalent to the single equation
  \begin{equation}
    \label{H0}
    \1-V+H(P_1+P_2V)+(V_{12}+\4{H})p_{\ker V_{22}}=0,
  \end{equation}
  which is a generalization of Eq.~(4.4) in \cite{R78}.

  Let us show that each solution $H$ of \eqref{H0} is completely
  determined by its component $H_{12}$. Given $H_{12}$, $H_{11}$ is
  fixed by \eqref{H1}:
  \begin{equation}
    \label{H11}
    H_{11}=V_{11}-P_1-H_{12}V_{21},
  \end{equation}
  $H_{22}$ is fixed by antisymmetry:
  \begin{equation}
    \label{H22}
    H_{22}=-{H_{11}}^\tau=P_2-{V_{22}}^*-{V_{12}}^*H_{12},
  \end{equation}
  and $H_{21}$ is determined by \eqref{H2}:
  \begin{equation}
    \label{H21}
    H_{21}=(P_2-H_{22})V_{21}=({V_{22}}^*+{V_{12}}^*H_{12})V_{21}.
  \end{equation}
  Therefore $H$ must have the form \eqref{H}, with $\4{T}\DEF H_{12}$.
  
  It remains to determine the admissible components $H_{12}$.
  \eqref{H1} implies that
  $$H_{12}V_{22}=V_{12}p_{\RAN{V_{22}}^*}.$$
  Inserting \eqref{H11}, \eqref{H3} is equivalent to 
  $$V_{21}(\ker V_{11})\subset\ker H_{12},$$
  and under this condition,
  \eqref{H4} holds automatically. Thus $T\DEF \4{H_{12}}$ has to
  fulfill \eqref{T}. Conversely, it is straightforward to verify that
  via \eqref{H}, any $T\in\mathfrak{H}_{P_1}$ obeying \eqref{T} gives
  rise to a solution $H$ of \eqref{H0}. In fact, one only has to check
  that $H_{21}$ defined by \eqref{H21} is antisymmetric and that
  $p_{\RAN{V_{22}}^*}+(H_{22}-P_2)V_{22}=0$ (the rest is clear by
  construction). By antisymmetry of $H_{12}$ and by \eqref{REL4}
  $$H_{21}+{H_{21}}^\tau=({V_{22}}^*+{V_{12}}^*H_{12})V_{21}+
  {V_{12}}^*(V_{11}-H_{12}V_{21})=0,$$
  so $H_{21}$ is antisymmetric. By \eqref{H22},
  \eqref{T} and \eqref{REL1},
  \begin{equation*}
    \begin{split}
      p_{\RAN{V_{22}}^*}+(H_{22}-P_2)V_{22}
        &=p_{\RAN{V_{22}}^*}-{V_{22}}^*V_{22}-{V_{12}}^*H_{12}V_{22}\\
      &=(P_1-{V_{22}}^*V_{22}-{V_{12}}^*V_{12})p_{\RAN{V_{22}}^*}\\
      &=0.
    \end{split}
  \end{equation*}
\end{proof}

Now we can proceed to exhibit the normal form of a complete set of
implementers for $\rho_V$. Let $H_V$ be defined by \eqref{H}, with
$T=T_V$ (see \eqref{TV}). Using
$p_{\ker{V_{22}}^*}=P_2-V_{22}{V_{22}}^{-1}$, one computes that
\begin{align*}
  (H_V)_{11}&={V_{11}}^{-1*}-P_1-p_{{\ker V_{11}}^*}
    V_{12}{V_{22}}^{-1}V_{21},\\
  (H_V)_{12}&=V_{12}{V_{22}}^{-1}-{V_{11}}^{-1*}{V_{21}}^*
    p_{{\ker V_{22}}^*},\\
  (H_V)_{21}&=({V_{22}}^{-1}-{V_{12}}^*{V_{11}}^{-1*}{V_{21}}^*
    p_{{\ker V_{22}}^*})V_{21},\\
  (H_V)_{22}&=P_2-{V_{22}}^{-1}+{V_{12}}^*{V_{11}}^{-1*}{V_{21}}^*
    p_{{\ker V_{22}}^*}.
  \end{align*}
$H_V$ is the analogue of Ruijsenaars' ``associate'' $\Lambda$ \cite{R78}.

Furthermore let $\{e_1,\dots,e_{L_V}\}$ be the \ONB in
$\mathfrak{h}_V=V_{12}(\ker V_{22})$ that was already used to define
$U_V$ in \eqref{UV} (cf.~Lemma~\ref{lem:UTH}), let
$\{e'_1,\dots,e'_{L_V}\}$ be the \ONB in $\ker V_{22}$ given by
$$e'_r\DEF V^*e_r={V_{12}}^*e_r,\qquad r=1,\dots,L_V,$$
and let $\{g_1,\dots,g_{M_V}\},\ M_V\DEF-\2\IND V$, be an \ONB in
\begin{equation}
  \label{KV}
  \mathfrak{k}_V\DEF P_V(\ker V^*).
\end{equation}
(Note that $P_V$ commutes with $VV^*$ and therefore
restricts to a \BP of $\ker V^*$.) Recall that the statistics
dimension $d_V$ \eqref{DV} of $\rho_V$ is given by
$$d_V=2^{M_V}$$
and that the twisted Fock representation $\psi_{P_1}$
defined in \eqref{PSI} fulfills
\begin{equation}
  \label{COMM}
  \pi_{P_1}\Big(\rho_V\big(\CK\big)\Big)'= 
  \psi_{P_1}\big(\mathfrak{C}(\ker V^*)\big)''.
\end{equation}
One has
$$\psi_{P_1}(e_r)=ia(e_r)^*\PME,\qquad\psi_{P_1}(e'_r)=
ia({e'_r}^*)\PME$$
(for notational convenience, we shall drop the
index $P_1$ on implementers like $\PME$ from now on, cf.\ 
\eqref{PME}). Finally define the following operators on \DD, with
range contained in the space of $C^\infty$--vectors for the number
operator
\begin{equation}
  \label{PAV}
  \begin{split}
    \PV{\alpha}&\DEF\big({\det}_{\KK_1}(P_1+{T_V}^*T_V)\big)^{-1/4}
      \psi_{P_1}(g_{\alpha_1}\dotsm g_{\alpha_l})\\
    &\quad\cdot\sum_{(\sigma,s)\in\PLV}(-1)^{L_V-s}\SGN\sigma\,
      \psi_{P_1}(e_{\sigma(1)}\dotsm e_{\sigma(s)})\EHV\\
    &\phantom{{\DEF}\cdot\sum_{(\sigma,s)\in\PLV}}\cdot\psi_{P_1}
      (e'_{\sigma(s+1)}\dotsm e'_{\sigma(L_V)}).
  \end{split}
\end{equation}
Here $\alpha=(\alpha_1,\dots,\alpha_l)\in I_{M_V}$ is a multi--index
as in \eqref{I}, and $\PLV$ is the index set consisting of all pairs
$(\sigma,s)$ with $s\in\{0,\dots,L_V\}$ and $\sigma$ a permutation of
order $L_V$ satisfying $\sigma(1)<\dots<\sigma(s)$ and
$\sigma(s+1)<\dots<\sigma(L_V)$. $\PLV$ is canonically isomorphic to
the power set $\PLVT$ of $\{1,\dots,L_V\}$ through identification of
$(\sigma,s)$ with $\{\sigma(1),\dots,\sigma(s)\}$, hence its
cardinality is $2^{L_V}$. Note that
  \begin{equation}
    \label{PSIAOM}
    \begin{split}
      \PV{\alpha}\Omega_{P_1}&=\big({\det}_{\KK_1}(P_1+{T_V}^*T_V)
        \big)^{-1/4}\psi_{P_1}(g_{\alpha_1}\dotsm g_{\alpha_l})\\
      &\quad\cdot\psi_{P_1}(e_1\dotsm e_{L_V})\exp(\3\4{T_V}a^*a^*)
        \Omega_{P_1}
    \end{split}
  \end{equation}
  by \eqref{EOM} and because the $\psi_{P_1}(e_r')$ annihilate the
  vacuum.
\begin{Th}
  \label{th:HRHO}\hspace*{\fill}\\
  Let $V\in\II_{P_1}(\KK)$. Then the $d_V$ operators $\PV{\alpha},\ 
  \alpha\in I_{M_V}$, have continuous extensions to isometries on
  $\mathcal{F}_a(\KK_1)$ \textup{(}henceforth denoted by the same
  symbols\textup{)} which implement $\rho_V$ in $\pi_{P_1}$ in the
  sense of Definition~\ref{def:IMP}.
\end{Th}
\begin{proof}
  1. We first show that the following intertwiner relation holds on \DD
  \begin{equation} 
    \label{INTER4}
    \PV{\alpha}\pi_{P_1}(f)=\pi_{P_1}(Vf)\PV{\alpha},
    \quad f\in\KK.
  \end{equation}
  Note that it suffices to prove \eqref{INTER4} for $\alpha=0$ because
  \begin{equation}
    \label{PAN}
    \PV{\alpha}=\psi_{P_1}(g_{\alpha_1}\dotsm g_{\alpha_l})\PV{0}
  \end{equation}
  and because the $\psi_{P_1}(g_j)$ belong to $\pi_{P_1}(\rho_V(\CK))'$. 

  Let first $f\in\RAN{V_{11}}^*\oplus\RAN{V_{22}}^*$. Then it follows
  from \eqref{PSI2} that
  $$[\psi_{P_1}(e_r'),\pi_{P_1}(f)]=0=[\psi_{P_1}(e_r),\pi_{P_1}(Vf)]$$
  so that \eqref{INTER4} is a consequence of \eqref{INT}.
  
  To prove \eqref{INTER4} for $f\in\ker V_{11}\oplus\ker V_{22}$, note
  that for fixed $r$, the bijection
  \begin{align*}
    \SET{\mathcal{M}\in\PLVT}{r\in\mathcal{M}}&\to
    \SET{\mathcal{M}'\in\PLVT}{r\notin\mathcal{M}'},\\ 
    \mathcal{M}&\mapsto\mathcal{M}\setminus\{r\}
  \end{align*}
  induces a bijection $(\sigma,s)\mapsto(\sigma',s')$ from
  $\SET{(\sigma,s)\in\PLV}{r\in\{\sigma(1),\dots,\sigma(s)\}}$ onto
  $\SET{(\sigma',s')\in\PLV}{r\notin\{\sigma'(1),\dots,\sigma'(s')\}}$
  with
  \begin{equation}
    \label{RELBIJ}
    s=s'+1,\quad(-1)^s\SGN\sigma=(-1)^r\SGN\sigma',
    \quad\sigma^{-1}(r)+{\sigma'}^{-1}(r)=r+s.
  \end{equation}
  Now let $D_V\DEF({\det}_{\KK_1}(P_1+{T_V}^*T_V))^{-1/4}$, and
  consider the case $f=e_r'\in\ker V_{22}$. We have on \DD, by virtue
  of
  \begin{align*}
    \psi_{P_1}(e_r')\pi_{P_1}(e_r')&=0=\pi_{P_1}(e_r)\psi_{P_1}(e_r),\\
    \{\psi_{P_1}(h),\PME\}&=0=[\EHV,\PME]
  \end{align*}
  and by \eqref{PSI2}, \eqref{PSI1} and \eqref{RELBIJ}, where terms
  under the sign ``$\WT{\quad}$'' are to be omitted
  \begin{equation*}
    \begin{split}
      \PV{0}\pi_{P_1}(e_r') &=D_V\psi_{P_1}(e_r)
        \sum_{\substack{(\sigma,s)\in\PLV\\r\in\{\sigma(1),\dots,\sigma(s)\}}}
        (-1)^{L_V-s+\sigma^{-1}(r)-1}\SGN\sigma\\
      &\quad\cdot\psi_{P_1}(e_{\sigma(1)}\dotsm\WT{e_r}
        \dotsm e_{\sigma(s)})\EHV\\
      &\quad\cdot\psi_{P_1}(e'_{\sigma(s+1)})\dotsm i\PME
        \psi_{P_1}(e_r')\dotsm\psi_{P_1}(e'_{\sigma(L_V)})\\
      &=D_V\pi_{P_1}(e_r)\sum_{\substack{(\sigma',s')\in\PLV\\r\notin
          \{\sigma'(1),\dots,\sigma'(s')\}}}(-1)^{L_V-s'}\SGN\sigma'
        \psi_{P_1}(e_{\sigma'(1)}\dotsm e_{\sigma'(s')})\\
      &\quad\cdot\EHV\psi_{P_1}
        (e'_{\sigma'(s'+1)}\dotsm e'_{\sigma'(L_V)})\\
      &=\pi_{P_1}(Ve_r')\PV{0}.
    \end{split}
  \end{equation*}
  The remaining case $f={e_r'}^*\in\ker V_{11}$ is similarly obtained
  with the help of \eqref{INT} and \eqref{INTER3}
  \begin{equation}
    \label{PEHV}
    \pi_{P_1}(e_r^*)\EHV=\EHV\pi_{P_1}({e_r'}^*)=0,
  \end{equation}
  in connection with
  $$[\pi_{P_1}(e_r^*),\psi_{P_1}(e_s)]=[\pi_{P_1}({e_r'}^*),
  \psi_{P_1}(e_s')]=i\delta_{rs}\PME$$
  (cf.~\eqref{PSI2}):
  \begin{equation*}
    \begin{split}
      \pi_{P_1}(V{e_r'}^*)\PV{0} &=\pi_{P_1}(e_r^*)\PV{0}\\
      &=i\PME D_V\sum_{\substack{(\sigma,s)\in\PLV\\ 
          r\in\{\sigma(1),\dots,\sigma(s)\}}}
        (-1)^{L_V-s+\sigma^{-1}(r)-1} \SGN\sigma\\
      &\quad\cdot\psi_{P_1}(e_{\sigma(1)}\dotsm\WT{e_r}
        \dotsm e_{\sigma(s)})\EHV\\
      &\quad\cdot\psi_{P_1}(e'_{\sigma(s+1)}\dotsm e'_{\sigma(L_V)})\\ 
      &=i\PME D_V\sum_{\substack{(\sigma',s')\in\PLV\\ 
          r\notin\{\sigma'(1),\dots,\sigma'(s')\}}}(-1)^{L_V-s'
        +{\sigma'}^{-1}(r)}\SGN\sigma'\\
      &\quad\cdot\psi_{P_1}(e_{\sigma'(1)}\dotsm
        e_{\sigma'(s')})\EHV\\
      &\quad\cdot\psi_{P_1}(e'_{\sigma'(s'+1)}\dotsm\WT{e_r'}
        \dotsm e'_{\sigma'(L_V)})\\
      &=\PV{0}\pi_{P_1}({e_r'}^*).
    \end{split}
  \end{equation*}
  This completes the proof of \eqref{INTER4}.
  
  2. To show that \PV{\alpha} is isometric, note that one has on \DD
  \begin{equation}
    \label{PSPN}
    \psi_{P_1}(g)^*\PV{0}=0,\qquad g\in\mathfrak{k}_V=P_V(\ker V^*).
  \end{equation}
  To see this, remember that the \BP $P_V$ has the form \eqref{PTH}
  $P_V=P_{T_V}-p_{\mathfrak{h}_V}+\4{p_{\mathfrak{h}_V}}$. Since
  $\mathfrak{h}_V$ is contained in $\RAN V$, we have
  \begin{equation}
    \label{PVQV}
    P_VQ_V=P_{T_V}Q_V
  \end{equation}
  where $Q_V=\1-VV^*$ is the projection onto
  $\ker V^*$. It follows that $\mathfrak{k}_V=\ker V^*\cap\RAN
  P_{T_V}=\ker V^*\cap\RAN(P_1+T_V)$. We thus get from
  $\RAN(P_1+T_V)=\ker(P_2-T_V)$ that $P_2g=T_VP_1g$. By
  Lemma~\ref{lem:REL}, this entails
  $$\pi_{P_1}(g)^*\EHV\Omega_{P_1}=0,\qquad g\in\mathfrak{k}_V,$$
  and we get, using \eqref{INTER4}, \eqref{PSI1} and \eqref{PSI2}, for
  $f_1,\dots,f_n\in\KK$
  \begin{equation*}
    \begin{split}
      \psi_{P_1}(g)^*\PV{0}&\pi_{P_1}(f_1)\dotsm\pi_{P_1}(f_n)
        \Omega_{P_1}\\
      &=-i\pi_{P_1}(Vf_1)\dotsm\pi_{P_1}(Vf_n)\PME\pi_{P_1}(g)^*
        \PV{0}\Omega_{P_1}\\
      &=-iD_V\pi_{P_1}(Vf_1)\dotsm\pi_{P_1}(Vf_n)\PME\\
      &\quad\cdot\psi_{P_1}(e_1\dotsm e_{L_V})\pi_{P_1}(g)^*
        \EHV\Omega_{P_1}\\
      &=0
    \end{split}
  \end{equation*}
  which proves \eqref{PSPN}.
  
  Since the $\psi_{P_1}(g_j)$ are partial isometries whose source and
  range projections sum up to \1 by the CAR
  $$\psi_{P_1}(g_j)^*\psi_{P_1}(g_j)+\psi_{P_1}(g_j)
    \psi_{P_1}(g_j)^*=\1$$
  and because these projections mutually commute for different values
  of $j$, it follows that $\psi_{P_1}(g_{\alpha_1})\dotsm
  \psi_{P_1}(g_{\alpha_l})$ is a partial isometry in
  $\pi_{P_1}(\rho_V(\CK))'$ which contains $\RAN\PV{0}$ in its initial
  space. Therefore \PV{\alpha} will be isometric provided that \PV{0}
  is. We have from \eqref{PSIAOM}, \eqref{PEHV}, \eqref{OMNORM} and
  the CAR
  \begin{equation*}
    \begin{split}
      \NORM{\PV{0}\Omega_{P_1}}^2 &=D_V^2\langle\EHV\Omega_{P_1},
        \psi_{P_1}(e_{L_V}^*\dotsm e_1^*)\\
      &\phantom{{=}D_V^2\langle}\quad\cdot\psi_{P_1}(e_1
        \dotsm e_{L_V})\EHV\Omega_{P_1}\rangle\\
      &=D_V^2\NORM{\EHV\Omega_{P_1}}^2\\
      &=1.
    \end{split}
  \end{equation*}
  Using the CAR and the fact that $\PV{0}\Omega_{P_1}$ serves as a
  vacuum for the transformed annihilation operators, this implies for
  arbitrary $f_1,\dots,f_m,h_1,\dots,h_n\in\KK$
  \begin{equation*}
    \begin{split}
    \langle\PV{0}\pi_{P_1}(f_1\dotsm&f_m)\Omega_{P_1},\PV{0}\pi_{P_1}
      (h_1\dotsm h_n)\Omega_{P_1}\rangle\\
    &=\langle\PV{0}\Omega_{P_1},\pi_{P_1}\big(\rho_V(f_m^*\dotsm
      f_1^*h_1\dotsm h_n)\big)\PV{0}\Omega_{P_1}\rangle\\
    &=\langle\Omega_{P_1},\pi_{P_1}(f_m^*\dotsm f_1^*h_1\dotsm h_n)
      \Omega_{P_1}\rangle.
    \end{split}
  \end{equation*}
  Therefore \PV{0} is isometric on \DD and has a continuous extension
  to an isometry which satisfies \eqref{INTER4} on
  $\mathcal{F}_a(\KK_1)$. By the above, the same holds true for the
  \PV{\alpha}.
  
  3. It remains to show that the \PV{\alpha} fulfill the Cuntz
  relations \eqref{CUNTZ1} (or \eqref{CUNTZ}). Since $\psi_{P_1}$ is a
  representation of the CAR and by \eqref{PAN}, \eqref{PSPN}, the
  \PV{\alpha} are orthonormal:
  \begin{equation*}
    \PV{\alpha}^*\PV{\beta}=\PV{0}^*\psi_{P_1}(g_{\alpha_l}^*\dotsm
    g_{\alpha_1}^*)\psi_{P_1}(g_{\beta_1}\dotsm g_{\beta_m})\PV{0}
    =\delta_{\alpha\beta}\1
  \end{equation*}
  for $\alpha$ as above and $\beta=(\beta_1,\dots,\beta_m)\in I_{M_V}$. 
  
  The proof of the completeness relation
  \begin{equation}
    \label{COMPL}
    \sum_{\alpha\in I_{M_V}}\PV{\alpha}\PV{\alpha}^*=\1
  \end{equation}
  relies on the product decomposition $V=U_VW_V$ from
  Section~\ref{sec:IP}. Recall that the \ONB $\{e_1,\dots,e_{L_V}\}$
  in $\mathfrak{h}_V=\mathfrak{h}_{U_V}$ was used to define $U_V$. By
  Proposition~\ref{prop:VUW}, $U_V$ maps $\ker{W_V}^*$ unitarily onto
  $\ker V^*$ and fulfills $U_VP_1=P_VU_V$. Therefore $U_V$ restricts
  to a unitary isomorphism from $\mathfrak{k}_{W_V}=P_1(\ker{W_V}^*)$
  onto $\mathfrak{k}_V$. We choose
  \begin{equation}
    \label{FJ}
    f_j\DEF(-1)^{L_V}{U_V}^*g_j,\qquad j=1,\dots,M_V\quad(M_V=M_{W_V})
  \end{equation}
  as \ONB in $\mathfrak{k}_{W_V}$. Applying \eqref{PAV} to $W_V$, we
  obtain an orthonormal set of isometries satisfying \eqref{INTER4}
  \WRT $W_V$
  \begin{align*}
    \PW{\alpha}&=\psi_{P_1}(f_{\alpha_1}\dotsm f_{\alpha_l})\EHW,\\
    \text{with } H_{W_V}&=
    \begin{pmatrix}
      (W_V)_{11}-P_1 & 0\\ 0 & P_2-{(W_V)_{22}}^*
    \end{pmatrix}.
  \end{align*}
  Let us show that this set of isometries is complete. We have
  \begin{equation}
    \label{PWO}
    \PW{\alpha}\Omega_{P_1}=\psi_{P_1}(f_{\alpha_1}\dotsm
    f_{\alpha_l})\Omega_{P_1}.
  \end{equation}
  Comparing with \eqref{DEF}, we see that $\RAN\PW{\alpha}$ is exactly
  the cyclic (in fact, irreducible) subspace $\mathcal{F}_\alpha(W_V)$
  for the representation $\pi_{P_1}\0\rho_{W_V}$. Completeness for
  the \PW{\alpha} thus follows from Proposition~\ref{prop:DECCAR}:
  $$\oplus_\alpha\RAN\PW{\alpha}=\oplus_\alpha
  \mathcal{F}_\alpha(W_V)=\mathcal{F}_a(\KK_1).$$
  
  Now let \PUV be the unitary implementer for $\rho_{U_V}$ given by
  \eqref{PAV}. Then the isometries $\PUV\PW{\alpha}$ obviously
  constitute a complete set of implementers for $\rho_V$.  We are
  going to show that actually
  \begin{equation}
    \label{PVUV}
    \PUV\PW{\alpha}=\PV{\alpha}
  \end{equation}
  holds under the above assumptions. Since each implementer is
  completely determined by its value on $\Omega_{P_1}$ (this follows
  from \eqref{INTER}), it suffices to prove \eqref{PVUV} when applied
  to $\Omega_{P_1}$. Note that $\PUV\PME=(-1)^{L_V}\PME\PUV$ so that
  $$\PUV\psi_{P_1}(f)=(-1)^{L_V}\psi_{P_1}(U_Vf)\PUV,\qquad f\in\KK.$$
  Hence we obtain from \eqref{PWO}, \eqref{PSIAOM} and
  Proposition~\ref{prop:VUW}
  \begin{equation*}
    \begin{split}
      \PUV\PW{\alpha}\Omega_{P_1} &=(-1)^{lL_V}\psi_{P_1}
        (U_Vf_{\alpha_1}\dotsm U_Vf_{\alpha_l})\PUV\Omega_{P_1}\\
      &=D_V\psi_{P_1}(g_{\alpha_1}\dotsm g_{\alpha_l})
        \psi_{P_1}(e_1\dotsm e_{L_V})\\
      &\quad\cdot\exp(\3\4{T_V}a^*a^*)\Omega_{P_1}\\
      &=\PV{\alpha}\Omega_{P_1}.
    \end{split}
  \end{equation*}
  Therefore \eqref{PVUV} holds. Since \eqref{CUNTZ2} follows from
  \eqref{CUNTZ1} and \eqref{INTER}, the theorem is proven. 
\end{proof}
\begin{Rem}
  The proof of Theorem~\ref{th:HRHO} shows (cf.\ \eqref{PVUV},
  \eqref{PWO}, \eqref{FJ}, \eqref{OP} and Prop.~\ref{prop:PP1}) that
  the vectors $\PV{\alpha}\Omega_{P_1}$ are cyclic vectors inducing
  certain \FSs, viz.\ the \FSs corresponding to the \BPs
  $$P_V^\alpha\DEF P_V-p_{\mathfrak{k}_V^\alpha}+
  \4{p_{\mathfrak{k}_V^\alpha}}.$$
  Here $\mathfrak{k}_V^\alpha$ is the
  subspace of $\mathfrak{k}_V$ spanned by the vectors
  $g_{\alpha_1},\dots,g_{\alpha_l}$ if $l$ is the length of
  $\alpha$. This is a non--trivial fact because linear combinations of
  such vectors will in general not induce quasi--free states.
\end{Rem}
Since $\psi_{P_1}$ is a representation of the CAR and by \eqref{PSPN},
the \HSP $H(\rho_V)$ generated by the \PV{\alpha} carries a Fock space
structure:
\begin{Cor}
  \label{cor:HRFOCK}\hspace*{\fill}\\
  The map
  $$\PV{\alpha}\mapsto a(g_{\alpha_1})^*\dotsm a(g_{\alpha_l})^*
  \Omega,\qquad\alpha=(\alpha_1,\dots,\alpha_l)\in I_{M_V}$$
  extends to a unitary isomorphism from $H(\rho_V)$ onto the
  antisymmetric Fock space $\mathcal{F}_a(\mathfrak{k}_V)$ over
  $\mathfrak{k}_V$.
\end{Cor}
\noindent Here $a(g)^*$ and $\Omega$ denote the creation operators and the
Fock vacuum in $\mathcal{F}_a(\mathfrak{k}_V)$. We shall see in
Section~\ref{sec:CARCH} that, if $V$ is gauge invariant, then the
isomorphism depicted in Corollary~\ref{cor:HRFOCK} is not only an
isomorphism of (graded) \HSPs but (up to a character of the \GG) also
an isomorphism of $G$--modules.
\subsection{Bosonized statistics}
\label{sec:BOSSTAT}
Though the formula \eqref{PAV} for \PV{\alpha} looks quite
complicated, it is not difficult to write the ``Bosonized statistics
operator'' $\hat\varepsilon_V$ associated with $V\in\II_{P_1}(\KK)$ as
a polynomial in $g_j, g_j^*$ if $\rho_V$ has finite statistics
dimension. Recall from the introduction (see \eqref{BSO0}) that
$\hat\varepsilon_V$ is defined in terms of the implementers as
\begin{equation}
  \label{BSO}
  \hat\varepsilon_V=\sum_{\alpha,\beta\in I_{M_V}}\PV{\alpha}
  \PV{\beta}{\PV{\alpha}}^*{\PV{\beta}}^*.
\end{equation}
Let us first derive a simple formula
for the operators $\PV{\alpha}\PV{\beta}^*$. (These operators are
matrix units for the commutant $\pi_{P_1}(\rho_V(\CK))'$.) We will use
the following notation for multi--indices $\alpha,\beta\in I_{M_V}$,
which is suggested by the identification
$\alpha\DEF(\alpha_1,\dots,\alpha_l)\mapsto\{\alpha_1,\dots,\alpha_l\}$
of $I_{M_V}$ with the power set 
of $\{1,\dots,M_V\}$. $l_\alpha\DEF l$ will denote the length of
$\alpha$, and $\alpha\cap\beta\in I_{M_V}$ will denote the
multi--index whose entries are the elements of the intersection of the
entries of $\alpha$ and $\beta$. $\alpha^c\in I_{M_V}$ will be the
``complementary'' multi--index whose entries are the elements of
$\{1,\dots,M_V\}\setminus\{\alpha_1,\dots,\alpha_{l_\alpha}\}$.
We further set
\begin{align*}
  g_\alpha &\DEF g_{\alpha_1}\dotsm g_{\alpha_{l_\alpha}},\\
  \Gamma_{\alpha\beta} &\DEF g_\alpha({g_{\alpha^c\cap\beta^c}})^* 
    g_{\alpha^c\cap\beta^c}{g_\beta}^*={\Gamma_{\beta\alpha}}^*,
    \qquad\alpha,\beta\in I_{M_V}.
\end{align*}
\begin{Lem}
  \label{lem:PAPB}\hspace*{\fill}\\
  Let $V\in\II_{P_1}(\KK)$ with $-\IND V<\infty$, and let
  $\alpha,\beta\in I_{M_V}$. Then
  $$\PV{\alpha}\PV{\beta}^*=\psi_{P_1}(\Gamma_{\alpha\beta}).$$
\end{Lem}
\begin{proof}
  Let
  $A_\alpha\DEF\psi_{P_1}(g_\alpha{g_{\alpha^c}}^*g_{\alpha^c})$. One
  has, by the CAR and by \eqref{PSPN}, \eqref{PAN}
  $$A_\alpha\PV{0}\PV{0}^*=\psi_{P_1}(g_\alpha)\PV{0}\PV{0}^*
  =\PV{\alpha}\PV{0}^*.$$
  If $\alpha'\neq0$ is another multi--index in $I_{M_V}$, then
  $$A_\alpha\PV{\alpha'}=A_\alpha\psi_{P_1}(g_{\alpha'})\PV{0}=0$$
  because $g_j^2=0,\ j=1,\dots,M_V$. Hence we obtain from
  \eqref{COMPL}
  $$A_\alpha=A_\alpha\sum_{\beta\in I_{M_V}}\PV{\beta}\PV{\beta}^*
  =A_\alpha\PV{0}\PV{0}^*=\PV{\alpha}\PV{0}^*.$$ 
  This entails, for arbitrary $\alpha,\beta\in I_{M_V}$,
  \begin{equation*}
    \begin{split}
      \PV{\alpha}\PV{\beta}^* &=\PV{\alpha}\PV{0}^*
        \big(\PV{\beta}\PV{0}^*\big)^*\\
      &=A_\alpha{A_\beta}^*\\
      &\DEF\psi_{P_1}(g_\alpha{g_{\alpha^c}}^*g_{\alpha^c}{g_{\beta^c}}^*
        g_{\beta^c}{g_\beta}^*).
    \end{split}
  \end{equation*}
  Now, by the CAR, $g_\alpha$ commutes with
  ${g_{\alpha^c}}^*g_{\alpha^c}$, ${g_\beta}^*$ commutes with
  ${g_{\beta^c}}^*g_{\beta^c}$, and one has
  $g_\alpha(g_{\beta^c\cap\alpha})^*g_{\beta^c\cap\alpha} =g_\alpha$
  and $(g_{\alpha^c\cap\beta})^*g_{\alpha^c\cap\beta}
  {g_\beta}^*={g_\beta}^*$. Thus we finally get $\PV{\alpha}
  \PV{\beta}^*=\psi_{P_1}\big(g_\alpha({g_{\alpha^c\cap\beta^c}})^*
  g_{\alpha^c\cap\beta^c}{g_\beta}^*\big)=
  \psi_{P_1}(\Gamma_{\alpha\beta})$.
\end{proof}
\begin{Rem}
  1. As a special case, one obtains the projections onto
  $\RAN\PV{\alpha}$ 
  $$\PV{\alpha}\PV{\alpha}^*=\pi_{P_1}(g_\alpha{g_{\alpha}}^*
  {g_{\alpha^c}}^*g_{\alpha^c}),$$
  (we used \eqref{PSI0}) from which one directly sees that
  \eqref{COMPL} holds
  $$\sum_\alpha\PV{\alpha}\PV{\alpha}^*=\pi_{P_1}\big((g_1{g_1}^*+
  {g_1}^*g_1)\dotsm(g_{M_V}{g_{M_V}}^*+{g_{M_V}}^*g_{M_V})\big)=\1.$$
  2. If $\IND V=-\infty$, then one still has 
  $$\PV{\alpha}\PV{\beta}^*=
  \psi_{P_1}(g_\alpha)\PV{0}\PV{0}^*\psi_{P_1}({g_\beta}^*),$$
  where
  $\PV{0}\PV{0}^*$ can be obtained as a strong limit
  $$\PV{0}\PV{0}^*=\SLIM{n}\pi_{P_1}
  ({g_1}^*g_1\dotsm{g_n}^*g_n).$$
  But this projection is no longer contained in \CK.
\end{Rem}
\begin{Prop}
  \hspace*{\fill}\\
  Let $V$ be as in Lemma~\ref{lem:PAPB}. Then the Bosonized statistics
  operator $\hat\varepsilon_V$ defined by \eqref{BSO} can be written as
  \begin{equation}
    \label{BSO1}
    \begin{split}
      \hat\varepsilon_V &=\pi_{P_1}(\tilde\varepsilon_V),\\
      \tilde\varepsilon_V &=\sum_{\alpha,\beta\in
        I_{M_V}}(-1)^{(l_\alpha+l_\beta)(l_\beta+L_V)}\Gamma_{\alpha\beta}
        \rho_V(\Gamma_{\beta\alpha})\in\CK_0.
    \end{split}
  \end{equation}
\end{Prop}
\begin{proof}
  It follows from \eqref{PAV} that 
  \begin{equation}
    \label{PA-1}
    \PME\PV{\alpha}\PME=(-1)^{(l_\alpha+L_V)}\PV{\alpha}.
  \end{equation}
  Therefore one has for $f\in\KK$, using \eqref{PSI1}
  \begin{equation*}
    \begin{split}
      \PV{\alpha}\psi_{P_1}(f) &=i\PV{\alpha}\pi_{P_1}(f)\PME\\
      &=i(-1)^{(l_\alpha+L_V)}\pi_{P_1}(Vf)\PME\PV{\alpha}\\
      &=(-1)^{(l_\alpha+L_V)}\psi_{P_1}(Vf)\PV{\alpha},
    \end{split}
  \end{equation*}
  and hence
  $$\psi_{P_1}(g_\beta)\PV{\alpha}^*=(-1)^{(l_\alpha+L_V)l_\beta}
  \PV{\alpha}^*\psi_{P_1}(\rho_V(g_\beta)),
  \quad\alpha,\beta\in I_{M_V}.$$
  Thus one gets from Lemma~\ref{lem:PAPB}
  \begin{equation*}
    \begin{split}
      \hat\varepsilon_V &=\sum_{\alpha,\beta}\PV{\alpha}\PV{\beta}
        {\PV{\alpha}}^*{\PV{\beta}}^*\\
      &=\sum_{\alpha,\beta}\PV{\alpha}\psi_{P_1}(\Gamma_{\beta\alpha})
        {\PV{\beta}}^*\\
      &=\sum_{\alpha,\beta}(-1)^{(l_\beta+L_V)(l_\alpha+l_\beta)}
        \PV{\alpha}{\PV{\beta}}^*\psi_{P_1}\big(\rho_V
        (\Gamma_{\beta\alpha})\big)\\
      &=\sum_{\alpha,\beta}(-1)^{(l_\beta+L_V)(l_\alpha+l_\beta)}
        \psi_{P_1}\big(\Gamma_{\alpha\beta}\rho_V
        (\Gamma_{\beta\alpha})\big)\\
      &=\psi_{P_1}(\tilde\varepsilon_V),
    \end{split}
  \end{equation*}
  with $\tilde\varepsilon_V$ as above. But $\tilde\varepsilon_V$ is
  even so that $\psi_{P_1}(\tilde\varepsilon_V)
  =\pi_{P_1}(\tilde\varepsilon_V)$ by \eqref{PSI0}.
\end{proof}
To check the consistency of our constructions, let us finally
calculate the ``Bosonized statistics parameter''
$$\hat\lambda_V\DEF\pi_{P_1}(\phi_V(\tilde\varepsilon_V))$$
which is associated with the Bosonized statistics operator and with
the left inverse $\phi_V$ from Section~\ref{sec:QFCAR} (see
\eqref{LI}). Recall that $\phi_V$ is given by 
$$\phi_V(ab)=\rho_V^{-1}(a)\omega_{1/2}(b)\quad\text{if }
a\in\mathfrak{C}(\RAN V),\ b\in\mathfrak{C}(\ker V^*),$$
where $\omega_\2$ is the trace on \CK.
\begin{Cor}
  \hspace*{\fill}\\
  The ``Bosonized statistics parameter'' $\hat\lambda_V$ of
  $V\in\II_{P_1}(\KK)$ with $-\IND V<\infty$ equals
  $$\hat\lambda_V=\frac{1}{d_V}\1.$$
\end{Cor}
\begin{proof}
  It follows from $\Gamma_{\alpha\beta}\in\mathfrak{C}(\ker V^*)$ that
  $\Gamma_{\alpha\beta}\rho_V(\Gamma_{\beta\alpha})=
  (-1)^{l_\alpha+l_\beta}\rho_V(\Gamma_{\beta\alpha})
  \Gamma_{\alpha\beta}$ and that $\phi_V(\rho_V(\Gamma_{\beta\alpha})
  \Gamma_{\alpha\beta})=\Gamma_{\beta\alpha}\omega_\2
  (\Gamma_{\alpha\beta})$. We claim that
  $$\omega_\2(\Gamma_{\alpha\beta})
  =\frac{1}{d_V}\cdot\delta_{\alpha\beta}.$$
  Consider first a term of
  the form $\omega_\2({g_\alpha}^*g_\alpha)$. Since the $g_j$ are
  mutually orthogonal, one has $\omega_\2({g_\alpha}^*g_\alpha)
  =\omega_\2({g_{\alpha_1}}^*g_{\alpha_1})\dotsm\omega_\2
  ({g_{\alpha_{l_\alpha}}}^*g_{\alpha_{l_\alpha}})= 2^{-{l_\alpha}}$.
  Next assume that $\alpha\neq\beta$. Without loss of generality,
  assume that $\alpha_1$ does not occur in $\beta$. Then there exists
  a quasi--free automorphism \rho which maps $g_{\alpha_1}$ to
  $-g_{\alpha_1}$ and leaves all other $g_j$ unchanged. It follows that
  $\rho({g_\alpha}^*g_\beta)=-{g_\alpha}^*g_\beta$ so that
  $\omega_\2({g_\alpha}^*g_\beta)=0$, because $\omega_\2$ is invariant
  under \rho. This yields
  \begin{equation*}
    \begin{split}
      \omega_\2(\Gamma_{\alpha\beta}) &=\omega_\2({g_\beta}^*g_\alpha
        ({g_{\alpha^c\cap\beta^c}})^*g_{\alpha^c\cap\beta^c})
        =\omega_\2({g_\beta}^*g_\alpha)\omega_\2
        (({g_{\alpha^c\cap\beta^c}})^*g_{\alpha^c\cap\beta^c})\\
      &=2^{-M_V}\delta_{\alpha\beta}=d_V^{-1}\delta_{\alpha\beta}
    \end{split}
  \end{equation*}
  as claimed. Hence we obtain 
  \begin{equation*}
    \begin{split}
      \phi_V(\tilde\varepsilon_V) &=\sum_{\alpha,\beta}
        (-1)^{(l_\alpha+l_\beta)(l_\beta+L_V)}\phi_V
        \big(\Gamma_{\alpha\beta}\rho_V(\Gamma_{\beta\alpha})\big)\\
      &=\sum_{\alpha,\beta}(-1)^{(l_\alpha+l_\beta)(l_\beta+L_V+1)}
        \Gamma_{\beta\alpha}\cdot\omega_\2(\Gamma_{\alpha\beta})\\
      &=\frac{1}{d_V}\sum_\alpha\Gamma_{\alpha\alpha}\\
      &=\frac{1}{d_V}\1.
    \end{split}
  \end{equation*}
\end{proof}
It is clear that one gets the same result by applying the left inverse
$\phi_{H(\rho_V)}$ from Section~\ref{sec:REP} (or from the
Introduction, Eq.~\eqref{LI0}) to $\hat\varepsilon_V$.  Recall that
$\phi_{H(\rho_V)}$ is a left inverse for the normal extension of
$\rho_V$ to $\BB(\mathcal{F}_a(\KK_1))$, given by
$$\phi_{H(\rho_V)}(x)=\frac{1}{d_V}\sum_\alpha\PV{\alpha}^*
x\PV{\alpha},\qquad x\in\BB(\mathcal{F}_a(\KK_1)).$$
One can show, by similar 
computations as above, that $\phi_{H(\rho_V)}$ extends $\phi_V$:
$$\phi_{H(\rho_V)}(\pi_{P_1}(a))=\pi_{P_1}(\phi_V(a)),\qquad
a\in\CK.$$


%% file: cc.tex
\specialsection{Quasi--free Endomorphisms of the CCR Algebra}
\label{sec:CCR}
This section contains an analysis of the semigroup of quasi--free
endomorphisms of the \CCR similar to the analysis done in
Section~\ref{sec:CAR} for the \CAR. Generally speaking, the CCR case
is algebraically simpler, but the analytic aspects are more involved.
For an introduction to the \CCR see the textbooks \cite{BR,Pe}.
\subsection{The selfdual \CCR}
\label{sec:SDCCR}
Let $\KK^0$ be an infinite dimensional complex linear space, equipped
with a nondegenerate hermitian sesquilinear form $\kappa$ and an
antilinear involution $f\mapsto f^*$, such that
$$\kappa(f^*,g^*)=-\kappa(g,f),\qquad f,g\in\KK^0.$$
One should think
of $\KK^0$ as being the complexification of the real linear space
$$\RE\KK^0\DEF\SET{f\in\KK^0}{f^*=f},$$
together with the canonical
conjugation on $\KK^0=\CC\otimes_{\RR}\RE\KK^0$. $-i\kappa$ should be
viewed as the sesquilinear extension of a nondegenerate symplectic
form on $\RE\KK^0$.

The \emph{(selfdual) CCR algebra} $\mathfrak{C}(\KK^0,\kappa)$
\cite{AS71,A71,A82} over $(\KK^0,\kappa)$ is the simple *--algebra
which is generated by \1 and elements $f\in\KK^0$, subject to the
commutation relation
\begin{equation}
  \label{ccr:CCR}
  [f^*,g]=\kappa(f,g)\1,\qquad f,g\in\KK^0.
\end{equation}
We henceforth assume the existence of a distinguished \FS over
$\mathfrak{C}(\KK^0,\kappa)$. As in the CAR case, \FSs correspond to
\BPs. A linear operator $P_1$, defined on the whole of $\KK^0$, is
a {\em\BP} of $(\KK^0,\kappa)$ if it satisfies for $f,g\in\KK^0$
\begin{equation}
  \label{ccr:BP}
  \begin{aligned}
    P_1^2&=P_1, &\qquad \kappa(f,P_1g)&=\kappa(P_1f,g),\\ 
    P_1+\4{P_1}&=\1, &\qquad\kappa(f,P_1f)&>0\quad\text{if }P_1f\not=0.
  \end{aligned}
\end{equation}
Here we used the notation \eqref{CONJ}
$$\4{P_1}f\DEF P_1(f^*)^*$$
for the complex conjugate operator. Let
$$P_2\DEF\1-P_1,\quad C\DEF P_1-P_2,\quad\langle{f,g}\rangle_{P_1}
\DEF\kappa(f,Cg).$$
The positive definite inner product $\langle\ ,\ 
\rangle_{P_1}$ turns $\KK^0$ into a pre--\HSP. We assume that the
completion \KK is separable. By continuity, the involution ``$*$'' extends
to a conjugation on \KK, $P_1$ and $P_2$ to orthogonal projections,
$C$ to a self--adjoint unitary, and $\kappa$ to a nondegenerate
hermitian form. These extensions will be denoted by the same symbols.
Setting
$$\KK_n\DEF P_n(\KK),\qquad n=1,2,$$
we get a direct sum decomposition
$\KK=\KK_1\oplus\KK_2$ which is orthogonal \WRT both $\kappa$ and
$\langle\ ,\ \rangle_{P_1}$. The following notations will frequently
be used for $A\in\BB(\KK)$:
\begin{align*}
  A_{mn} & \DEF P_mAP_n,\quad m,n=1,2,\\
  A\+    & \DEF CA^*C\\
  A^\tau & \DEF \4{A^*}.
\end{align*}
The components $A_{mn}$ of $A$ are regarded as operators from $\KK_n$
to $\KK_m$, and $A$ will sometimes be written as a matrix
$\bigl(\begin{smallmatrix}A_{11}&A_{12}\\
  A_{21}&A_{22}\end{smallmatrix}\bigr)$. $A\+$ is the adjoint of $A$
relative to $\kappa$, while $A^*$ is the \HSP adjoint. Thus one has
relations like
$$\4{P_2}=P_1=P_1\+=P_1^*,\quad\4{C}=-C,\quad{A_{12}}\+={A\+}_{21}=
-{A_{12}}^*,\quad\4{A_{11}}=\4{A}_{22},\quad\text{etc.}$$
The {\em\FS}
$\omega_{P_1}$ is the unique state\footnote{A \emph{state} $\omega$
  over \CKG is a linear functional with $\omega(\1)=1$ and
  $\omega(a^*a)\geq0,\ a\in\CKG$.} which is annihilated by all
$f\in\RAN P_2$:
$$\omega_{P_1}(f^*f)=0\quad\text{if }P_1f=0.$$
(In the conventional
setting mentioned above, $\omega_{P_1}$ is the \FS corresponding to
the complex structure $iC$ on $\RE\KK$.) Let $\mathcal{F}_s(\KK_1)$ be
the symmetric Fock space over $\KK_1$ and let \DD be the dense
subspace of algebraic tensors. A GNS representation $\pi_{P_1}$ for
$\omega_{P_1}$ is provided by
$$\pi_{P_1}(f)=a^*(P_1f)+a(P_1f^*),\quad f\in\KK$$
where $a^*(g)$ and
$a(g)$, $g\in\KK_1$, are the usual (Bosonic) creation and annihilation
operators on \DD. The cyclic vector inducing the state $\omega_{P_1}$
is $\Omega_{P_1}$, the Fock vacuum. The operators $\pi_{P_1}(a)$,
$a\in\CKG$, have invariant domain \DD, are closable, and one has
$\pi_{P_1}(a^*)\subset\pi_{P_1}(a)^*$. In particular, if $f\in\RE\KK$,
then $\pi_{P_1}(f)$ is essentially self--adjoint on \DD, and the
unitary \emph{Weyl operator} $w(f)$ is defined as the exponential of
the closure of $i\pi_{P_1}(f)$. Its vacuum expectation value is
$$\omega_{P_1}(w(f))\DEF\langle{\Omega_{P_1},
  w(f)\Omega_{P_1}}\rangle=e^{-\frac{1}{4}\NORM{f}^2_{P_1}},$$
and the \emph{Weyl relations} hold
$$w(f)w(g)=e^{-\2\kappa(f,g)}w(f+g),\quad f,g\in\RE\KK.$$
The Weyl
operators generate a simple \CAL $\WK$ which acts irreducibly on
$\mathcal{F}_s(\KK_1)$. If \HH is a subspace of \KK with $\HH=\HH^*$,
then the \CAL generated by all $w(f)$ with $f\in\RE\HH$ will be
denoted by $\WW(\HH)$. If $\HH\SC$ is the orthogonal complement of
\HH \WRT $\kappa$, then \emph{duality} holds \cite{A63,AS71,A82}:
\begin{equation}
  \label{ccr:DUAL}
  \WW(\HH)'=\WW(\HH\SC)''.
\end{equation}
\begin{Lem}
  \label{lem:AFF}\hspace*{\fill}\\
  For $f\in\KK$, let $\HH_f$ be the subspace spanned by $f$ and $f^*$.
  Then the closure of $\pi_{P_1}(f)$ is affiliated with
  $\WW(\HH_f)''$.
\end{Lem}
\begin{proof}
  Let $T$ be the closure of $\pi_{P_1}(f)$, with domain $D(T)$. We
  have to show that, for any $A\in\WW(\HH_f)'$
  $$A(D(T))\subset D(T),\qquad AT=TA\text{ on }D(T).$$
  Now by virtue
  of the CCR \eqref{ccr:CCR}, $\NORM{T\phi}^2=\NORM{T^*\phi}^2
  +\kappa(f,f)\NORM{\phi}^2$ for $\phi\in\DD$. Hence, for a given
  Cauchy sequence $\phi_n\in\DD$, $T\phi_n$ converges \IFF $T^*\phi_n$
  does. This implies that
  $$D(T)=D(T^*).$$
  Let $f^\pm\in\RE\HH_f$ be defined as
  $f^+\DEF\2(f+f^*),\ f^-\DEF\frac{i}{2}(f-f^*)$, and let $T^\pm$ be
  the (self--adjoint) closure of $\pi_{P_1}(f^\pm)$. We claim that
  $$D(T)=D(T^+)\cap D(T^-),\qquad T=T^+-iT^-\text{ on }D(T).$$
  For if
  $\phi\in D(T)$, then there exists a sequence $\phi_n\in\DD$
  converging to $\phi$ such that $\pi_{P_1}(f)\phi_n$ and
  $\pi_{P_1}(f^*)\phi_n$ converge. Thus $\phi$ belongs to the domain
  of the closure of $\pi_{P_1}(f^\pm)$. Conversely, if $\phi\in
  D(T^+)\cap D(T^-)$, then there exists a sequence $\phi_n\in\DD$
  converging to $\phi$ such that both $\pi_{P_1}(f+f^*)\phi_n$ and
  $\pi_{P_1}(f-f^*)\phi_n$ converge (this follows from the detailed
  description of the domains of such operators given in \cite{R78}).
  Therefore $\pi_{P_1}(f)\phi_n$ is also convergent, i.e.\ $\phi$ is
  contained in $D(T)$, and $T\phi=(T^+-iT^-)\phi$.
  
  Now if $A\in\WW(\HH_f)'$, then $A$ commutes with the one--parameter
  unitary groups $w(tf^\pm)=\exp(itT^\pm)$. As a consequence, $A$
  leaves $D(T^\pm)$ invariant and commutes with $T^\pm$ on $D(T^\pm)$.
  It follows that $A(D(T))\subset D(T)$ and $AT=TA$ on $D(T)$ as was
  to be shown.
\end{proof}
\subsection{Implementability of quasi--free endomorphisms}
\label{sec:IMP}
\emph{Quasi--free endomorphisms} are the unital *--endomorphisms of
\CKG which map \KK, viewed as a subspace of \CKG, into itself. They
are completely determined by their restrictions to \KK which are
called \emph{Bogoliubov operators.} Hence $V\in\BB(\KK)$ is a
Bogoliubov operator\footnote{We may disregard unbounded Bogoliubov
  operators $V$ (defined on $\KK^0$) since the topologies induced by
  the corresponding states $\omega_{P_1}\0\rho_V$ on $\KK^0$ differ
  from the one induced by $\omega_{P_1}$. Hence these states cannot be
  quasi--equivalent to $\omega_{P_1}$ (cf.~\cite{A71,A82}), and
  $\rho_V$ cannot be implemented.} \IFF it commutes with complex
conjugation and preserves the hermitian form $\kappa$. Bogoliubov
operators form a unital semigroup which we denote by
$$\SP{}{}\DEF\SET{V\in\BB(\KK)}{\4{V}=V,\ V\+V=\1}.$$
Each
$V\in\SP{}{}$ extends to a unique quasi--free endomorphism of \CKG and
to a unique *--endomorphism of $\WK$. By abuse of notation, both
endomorphisms are denoted by $\rho_V$, so that $\rho_V(f)=Vf,\ 
f\in\KK$, and $\rho_V(w(g))=w(Vg),\ g\in\RE\KK$.

The condition $V\+V=\1$ entails that $V$ is injective and $V^*$
surjective; hence $\RAN V$ is closed, and $V$ is a semi--Fredholm
operator \cite{K}. We claim that the Fredholm index $-\IND V=\dim\ker
V\+$ cannot be odd, in contrast to the CAR case (cf.~\eqref{INK}). For
let $f\in\ker V\+$ such that $0=\kappa(f,g)\equiv\langle{f,Cg}
\rangle_{P_1}\ \forall g\in\ker V\+$.  Then $f\in(C\ker
V\+)^\perp=(\ker V^*)^\perp=\RAN V$, but $\RAN V\cap\ker V\+=\{0\}$
due to $V\+V=\1$, so $f$ has to vanish. This shows that the
restriction of $\kappa$ to $\ker V\+$ stays nondegenerate. It follows
that $\dim\ker V\+$ cannot be odd (there is no nondegenerate
symplectic form on an odd dimensional space).

On the other hand, each even number (and $\infty$) occurs as $\dim\ker
V\+$ for some $V$. Hence we have an epimorphism of semigroups
$$\SP{}{}\to\NN\cup\{\infty\},\quad V\mapsto-\3\IND V=\3\dim\ker V\+$$
(remember that $0\in\NN$). Let
$$\SP{}{n}\DEF\{V\in\SP{}{}\ |\ \IND V=-2n\},\quad
n\in\NN\cup\{\infty\}.$$
$\SP{}{0}$ is the group of quasi--free
automorphisms (isomorphic to the symplectic group of $\RE\KK$). It
acts on $\SP{}{}$ by left multiplication. Analogous to the CAR case,
the orbits under this action are the subsets $\SP{}{n}$, and the
stabilizer of $V\in\SP{}{n}$ is isomorphic to the symplectic group
Sp($n$).

We are interested in endomorphisms $\rho_V$ which can be implemented
by \HSPs of isometries on $\mathcal{F}_s(\KK_1)$. This means that
there exist isometries $\Psi_j$ on $\mathcal{F}_s(\KK_1)$ which
fulfill the \CA relations \eqref{CUNTZ1} and implement
$\rho_V$ according to \eqref{CUNTZ2}:
$$\rho_V(w(f))=\sum_j\Psi_jw(f)\Psi_j^*,\quad f\in\RE\KK.$$
As
explained in Section~\ref{sec:REP}, such isometries exist \IFF
$\rho_V$, viewed as a representation of $\WK$ on
$\mathcal{F}_s(\KK_1)$, is quasi--equivalent to the defining (Fock)
representation.

To study $\rho_V$ as a representation, for fixed $V\in\SP{}{}$, let us
decompose it into cyclic subrepresentations. Let $f_1,f_2,\dotsc$ be
an \ONB in $\KK_1\cap\ker V\+$ and let
$\alpha=(\alpha_1,\dots,\alpha_l)$ be a multi--index with
$\alpha_j\leq\alpha_{j+1}$. Such $\alpha$ has the form
\begin{equation}
  \label{ccr:ALPHA}
  \alpha=(\underbrace{\alpha'_1,\dots,\alpha'_1}_{l_1},
  \underbrace{\alpha'_2,\dots,\alpha'_2}_{l_2},\dots,
  \underbrace{\alpha'_r,\dots,\alpha'_r}_{l_r})
\end{equation}
with $\alpha'_1<\alpha'_2<\dotsb<\alpha'_r$ and $l_1+\dots+l_r=l$. Let
\begin{equation}
  \label{ccr:DEF}
  \begin{split}
    c_\alpha &\DEF(l_1!\dotsm l_r!)^{-\2}\\
    \phi_\alpha &\DEF c_\alpha a^*(f_{\alpha_1})\dotsm
      a^*(f_{\alpha_l})\Omega_{P_1},\\ 
    \mathcal{F}_\alpha &\DEF \4{\WW(\RAN V)\phi_\alpha},\\ 
    \pi_\alpha &\DEF \rho_V|_{\mathcal{F}_\alpha}.
  \end{split}
\end{equation}
The Bosonic analogue of Proposition~\ref{prop:DECCAR} is
\begin{Prop}
  \label{ccr:prop:CYC}\hspace*{\fill}\\
  One has $\rho_V=\oplus_\alpha\pi_\alpha$, where the sum extends
  over all multi--indices $\alpha$ as above, including $\alpha=0$
  ($\phi_0\DEF\Omega_{P_1}$). Each
  $(\pi_\alpha,\mathcal{F}_\alpha,\phi_\alpha)$ is a GNS
  representation for $\omega_{P_1}\0\rho_V$ \textup{(}regarded as a
  state over $\WK$\textup{)}.
\end{Prop}
\begin{proof}
  By definition, the $\phi_\alpha$ constitute an \ONB for
  $\mathcal{F}_s(\KK_1\cap\ker V\+)$, and
  $(\pi_\alpha,\mathcal{F}_\alpha,\phi_\alpha)$ is a cyclic
  representation of $\WK$.  Since the closures of $a^*(f_j)$ and
  $a(f_j)$ are affiliated with $\WW(\ker V\+)''=\WW(\RAN V)'$ (see
  Lemma~\ref{lem:AFF} and \eqref{ccr:DUAL}), one obtains for
  $f\in\RE\KK$
  \begin{equation*}
    \begin{split}
      \langle{\phi_\alpha,\pi_\alpha(w(f))\phi_\alpha}\rangle 
        &=c_\alpha^2\langle{a^*(f_{\alpha_1})\dotsm
        a^*(f_{\alpha_l})\Omega_{P_1}, w(Vf)a^*(f_{\alpha_1})\dotsm
        a^*(f_{\alpha_l})\Omega_{P_1}}\rangle\\ 
      &=c_\alpha^2\langle\Omega_{P_1},
        w(Vf)\underbrace{a(f_{\alpha_l})\dotsm
        a(f_{\alpha_1})a^*(f_{\alpha_1})\dotsm a^*(f_{\alpha_l})
        \Omega_{P_1}}_{c_\alpha^{-2}\Omega_{P_1}}\rangle\\
      &=\langle{\Omega_{P_1},w(Vf)\Omega_{P_1}}\rangle.
    \end{split}
  \end{equation*}
  This proves that $(\pi_\alpha,\mathcal{F}_\alpha,\phi_\alpha)$ is a
  GNS representation for $\omega_{P_1}\0\rho_V$. Similarly, one finds
  that $\langle{\phi_\alpha,w(Vf)\phi_{\alpha'}}\rangle=0$ for
  $\alpha\neq\alpha'$, so the $\mathcal{F}_\alpha$ are mutually
  orthogonal.
  
  It remains to show that $\oplus_\alpha\mathcal{F}_\alpha=
  \mathcal{F}_s(\KK_1)$. We claim that $\mathcal{F}_0$ equals
  $\mathcal{F}_s(\4{\RAN P_1V})$, the symmetric Fock space over the
  closure of $\RAN P_1V$. The inclusion $\mathcal{F}_0\subset
  \mathcal{F}_s(\4{\RAN P_1V})$ holds because vectors of the form
  $w(Vf)\Omega_{P_1}=\exp i\bigl(a^*(P_1Vf)+ a(P_1Vf)\bigr)
  \Omega_{P_1}\in\mathcal{F}_s(\4{\RAN P_1V})$ are total in
  $\mathcal{F}_0$. The converse inclusion may be proved inductively.
  Assume that $a^*(g_1)\dotsm a^*(g_m)\Omega_{P_1}$ is contained in
  $\mathcal{F}_0$ for all $m\leq n,\ g_1,\dots,g_m\in\RAN P_1V$.
  Then, for $f\in V(\RE\KK)$ and $g_1,\dots,g_n\in\RAN P_1V$,
  $\tfrac{1}{i}\tfrac{w(tf)-\1}{t}a^*(g_1)\dotsm a^*(g_n)\Omega_{P_1}$
  has a limit $a^*(P_1f)a^*(g_1)\dotsm a^*(g_n)\Omega_{P_1}+
  a(P_1f)a^*(g_1)\dotsm a^*(g_n)\Omega_{P_1}$ in $\mathcal{F}_0$ as
  $t\searrow0$. By assumption, the second term lies in
  $\mathcal{F}_0$, and so does the first. Since each $g\in\RAN P_1V$
  is a linear combination of such $P_1f$, it follows that
  $a^*(g_1)\dotsm a^*(g_{n+1})\Omega_{P_1}$ is contained in
  $\mathcal{F}_0$ for arbitrary $g_j\in\RAN P_1V$, and, by induction,
  for arbitrary $n\in\NN$. But such vectors span a dense subspace in
  $\mathcal{F}_s(\4{\RAN P_1V})$, so
  $\mathcal{F}_0=\mathcal{F}_s(\4{\RAN P_1V})$ as claimed.
  
  Finally, $\KK_1\cap\ker V\+$ equals $\ker V^*P_1$, where $V^*P_1$ is
  regarded as an operator from $\KK_1$ to \KK. Thus we have
  $\KK_1=\4{\RAN P_1V}\oplus(\KK_1\cap\ker V\+)$ and
  $\mathcal{F}_s(\KK_1)\cong\mathcal{F}_0\otimes\mathcal{F}_s
  (\KK_1\cap\ker V\+)$. Under this isomorphism, $\mathcal{F}_\alpha$ is
  identified with $\mathcal{F}_0\otimes(\CC\phi_\alpha)$. Since the
  $\phi_\alpha$ form an \ONB for $\mathcal{F}_s(\KK_1\cap\ker V\+)$,
  the desired result
  $\oplus_\alpha\mathcal{F}_\alpha=\mathcal{F}_s(\KK_1)$ follows.
\end{proof}
As a consequence, the representation $\rho_V$ is quasi--equivalent to
the GNS representation associated with the quasi--free state
$\omega_{P_1}\0\rho_V$. So $\rho_V$ is implementable \IFF
$\omega_{P_1}\0\rho_V$ and $\omega_{P_1}$ are quasi--equivalent. Now
the two--point function of $\omega_{P_1}\0\rho_V$ (as a state over
\CKG) is given by
$$\omega_{P_1}\0\rho_V(f^*g)=\kappa(f,Sg)=
\langle{f,\tilde{S}g}\rangle_{P_1},\quad f,g\in\KK,$$ 
with
$$S\DEF V\+P_1V,\qquad\tilde{S}\DEF V^*P_1V.$$
The latter operators
contain valuable information about $\omega_{P_1}\0\rho_V$. For
example, it can be shown (cf.~\cite{MV}) that $\omega_{P_1}\0\rho_V$
is a \emph{pure} state over $\WK$ \IFF $S$ is a \BP, that is, \IFF $S$
is idempotent (the remaining conditions in \eqref{ccr:BP} are
automatically fulfilled). This is further equivalent to
$[P_1,VV\+]=0$, by the following argument:
$$\begin{array}{rcll}
  S^2=S &\Leftrightarrow& 0=S\4{S} & (\text{since }\4{S}=\1-S)\\
  &\Leftrightarrow& 0=V^*P_1VCV^*P_2V & \\
  &\Leftrightarrow& 0=P_1VCV^*P_2 & 
    (\text{since }\RAN V^*P_2V=\RAN V^*P_2\\
  &&& \text{ and }\ker V^*P_1V=\ker P_1V)\\
  &\Leftrightarrow& 0=P_1VV\+P_2 & \\
  &\Leftrightarrow& 0=[P_1,VV\+].&
\end{array}$$
On the other hand, the criterion for quasi--equivalence of quasi--free
states, in the form given by Araki and Yamagami \cite{A82} (see also
\cite{AS71,A71,vD}), yields that $\omega_{P_1}\0\rho_V$ is
quasi--equivalent to $\omega_{P_1}$ \IFF $P_1-\tilde{S}^{\2}$ is a
Hilbert--Schmidt operator on \KK. Using Theorem~\ref{th:IMPCAR}, this
condition can be simplified:
\begin{Th}
  \label{ccr:th:IMP}\hspace*{\fill}\\
  Let a Bogoliubov operator $V\in\SP{}{}$ be given. Then there exists
  a \HSP of isometries $H(\rho_V)$ which implements the endomorphism
  $\rho_V$ in the Fock representation determined by the \BP $P_1$ \IFF
  $[P_1,V]$ \textup{(}or, equivalently, $V_{12}$\textup{)} is a
  Hilbert--Schmidt operator. The dimension of $H(\rho_V)$ is $1$ if
  $\,\IND V=0$, otherwise $\infty$.
\end{Th}
\begin{proof}
  First note that $[P_1,V]=V_{12}-V_{21}=V_{12}-\4{V_{12}}$ is
  Hilbert--Schmidt (HS) \IFF $V_{12}$ is.

  By the preceding discussion, $\rho_V$ is implementable \IFF
  $P_1-\tilde{S}^{\2}$ is HS\@. In this case,
  $P_2(P_1-\tilde{S}^{\2})^2P_2=P_2\tilde{S}P_2= {V_{12}}^*V_{12}$ is
  of trace class, hence $V_{12}$ is HS\@.

  Conversely, assume $V_{12}$ to be HS\@. Let $V=V'\ABS{V}$ be the
  polar decomposition of $V$. Then $\ABS{V}=\4{\ABS{V}}$ is a bounded
  bijection with a bounded inverse, and $\ABS{V}-\1=
  (\ABS{V}^2-\1)(\ABS{V}+\1)^{-1}=(V^*-V\+)V(\ABS{V}+\1)^{-1}
  =2({V_{12}}^*+{V_{21}}^*)V(\ABS{V}+\1)^{-1}$ is HS\@. Thus, by a
  corollary \cite{A82} of an inequality of Araki and
  Yamagami \cite{AY}, $(\ABS{V}A\ABS{V})^\2-A^\2$ is HS for any
  positive $A\in\BB(\KK)$. Applying this to $A={V'}^*P_1V'$, we get
  that
  \begin{equation}
    \label{ccr:HS}
    \tilde{S}^\2-({V'}^*P_1V')^\2\text{ is HS\@.}
  \end{equation}
  Now $V'$ is an isometry with $\4{V'}=V'$, i.e.\ a CAR Bogoliubov
  operator (see \eqref{IK}). Since $[P_1,V]$ and $[P_1,\ABS{V}^{-1}]=
  \ABS{V}^{-1}\bigl[\ABS{V},P_1\bigr]\ABS{V}^{-1}=
  \ABS{V}^{-1}\bigl[\ABS{V}-\1,P_1\bigr]\ABS{V}^{-1}$ are HS, the same
  holds true for $[P_1,V']=[P_1,V\ABS{V}^{-1}]$. So $V'$ fulfills the
  implementability condition for CAR Bogoliubov operators derived in
  Theorem~\ref{th:IMPCAR}, and, as shown there, this forces
  $P_1-({V'}^*P_1V')^\2$ to be HS\@. This, together with
  \eqref{ccr:HS}, implies that $P_1-\tilde{S}^{\2}$ is HS as claimed.
  
  It remains to prove the statement about $\dim H(\rho_V)$. Let
  $\tilde\varrho_V$ be the normal extension of $\rho_V$ to
  $\BB(\mathcal{F}_s(\KK_1))$. Then
  $\BB(H(\rho_V))\cong\tilde\varrho_V(\BB(\mathcal{F}_s(\KK_1)))'
  =\rho_V(\WK)' =\WW(\RAN V)'=\WW(\ker V\+)''$. The latter (and hence
  $H(\rho_V)$) is one--dimensional if $\ker V\+=\{0\}$ and
  infinite dimensional if $\ker V\+\not=\{0\}$.
\end{proof}
\begin{Rem}
  Shale's original result \cite{S} asserts that a quasi--free
  automorphism $\rho_V$, $V\in\SP{}{0}$, is implementable \IFF
  $\ABS{V}-\1$ is HS\@. This condition is equivalent to $[P_1,V]$
  being HS, not only for $V\in\SP{}{0}$, but for all $V\in\SP{}{}$
  with $-\IND V<\infty$. However, the two conditions are \emph{not}
  equivalent for $V\in\SP{}{\infty}$, as the following example shows.
  Let $\KK_1=\HH\oplus\HH'$ be a decomposition into
  infinite dimensional subspaces. Choose an operator $V_{12}$ from
  $\KK_2$ to \HH with $\TR\ABS{V_{12}}^4<\infty$, but
  $\TR\ABS{V_{12}}^2=\infty$. Let $V_{21}\DEF\4{V_{12}}$ and
  $\ABS{V_{11}}\DEF(P_1+ \ABS{V_{21}}^2)^\2$. Choose an isometry
  $v_{11}$ from $\KK_1$ to $\HH'$ and set $V_{11}\DEF
  v_{11}\ABS{V_{11}},\ V_{22}\DEF\4{V_{11}}$. These components define
  a Bogoliubov operator $V\in\SP{}{\infty}$ (cf.\ 
  \eqref{ccr:REL1}--\eqref{ccr:REL4} below) which violates the
  condition of Theorem~\ref{ccr:th:IMP}. But it fulfills Shale's
  condition since $\ABS{V}^2-\1=2(\ABS{V_{12}}^2+ \ABS{V_{21}}^2)$ is
  HS and since $\ABS{V}-\1=(\ABS{V}^2-\1) (\ABS{V}+\1)^{-1}$.
\end{Rem}
Let $V\in\SP{}{}$ with $V_{12}$ compact. Due to stability under
compact perturbations \cite{K}, $V_{11}$ and $V_{22}=\4{V_{11}}$ are
semi--Fredholm with
\begin{equation}
  \label{ccr:IND}
  \IND V_{11}=\IND V_{22}=\2\IND V.
\end{equation}
We will occasionally use the relation $V\+V=\1$ componentwise:
\begin{subequations}\label{ccr:REL}
  \begin{align}
    {V_{11}}^*V_{11}-{V_{21}}^*V_{21} &= P_1, \label{ccr:REL1} \\
    {V_{22}}^*V_{22}-{V_{12}}^*V_{12} &= P_2, \label{ccr:REL2} \\
    {V_{11}}^*V_{12}-{V_{21}}^*V_{22} &= 0, \label{ccr:REL3} \\
    {V_{22}}^*V_{21}-{V_{12}}^*V_{11} &= 0. \label{ccr:REL4}
  \end{align}
\end{subequations}
Since $V_{11}$ is injective by \eqref{ccr:REL1} and has closed range, we
may define a bounded operator ${V_{11}}^{-1}$ as the inverse of ${V_{11}}$
on $\RAN{V_{11}}$ and as zero on $\ker {V_{11}}^*$ (the same applies
to $V_{22}$). These operators will be needed later. Note that $\dim\ker
{V_{11}}^*=-\2\IND V$.
\subsection{The semigroup of implementable endomorphisms}
\label{sec:SG}
According to Theorem~\ref{ccr:th:IMP}, the semigroup of implementable
quasi--free endomorphisms is isomorphic to the following semigroup of
Bogoliubov operators:
$$\SP{P_1}{}\DEF\SET{V\in\SP{}{}}{V_{12}\text{ is
    Hilbert--Schmidt}}.$$
$\SP{P_1}{}$ is a topological semigroup \WRT
the metric $\delta_{P_1}(V,V')\DEF\NORM{V-V'}+\HSNORM{V_{12}-
  V'_{12}}$, where $\HSNORM{\ }$ denotes Hilbert--Schmidt norm. It
contains the closed sub--semigroup of diagonal Bogoliubov operators
$$\SP{\text{diag}}{}=\SET{V\in\SP{}{}}{[P_1,V]=0}.$$
One has (cf.~\eqref{IDIAG})
$$\SP{\text{diag}}{}=\IDIAG(\KK)=\SP{}{}\cap\II(\KK).$$
The Fredholm index yields a decomposition
$$\SP{P_1}{}=\bigcup_{n\in\NN\cup\{\infty\}}\SP{P_1}{n},\qquad
\SP{P_1}{n}\DEF\SP{P_1}{}\cap\SP{}{n}.$$
The group $\SP{P_1}{0}$ is
usually called the \emph{restricted symplectic group} \cite{S,Se}. It
has a natural normal subgroup
$$\SP{\text{HS}}{}\DEF\SET{V\in\SP{}{}}{V-\1\text{ is
    Hilbert--Schmidt}} \subset\SP{P_1}{0}.$$
As in the CAR case, we
will eventually show that each $V\in\SP{P_1}{}$ can be written as a
product $V=UW$ with $U\in\SP{\text{HS}}{}$ and
$W\in\SP{\text{diag}}{}$. Assume that such $U$ and $W$ exist. Then
$P_V\DEF UP_1U\+$ is a \BP extending the ``partial \BP'' $VP_1V\+$
such that
\begin{equation}
  \label{ccr:PV}
  P_1-P_V\text{ is Hilbert--Schmidt},\qquad V\+P_VV=P_1,
\end{equation}
so the corresponding \FS $\omega_{P_V}$ is unitarily equivalent to
$\omega_{P_1}$ and fulfills $\omega_{P_V}\0\rho_V=\omega_{P_1}$. In
order to construct such \BPs, let us investigate the set $\PP_{P_1}$
(not to be confused with $\PP_{P_1}$ from Section~\ref{sec:CAR}) of
\BPs of $(\KK,\kappa)$ which differ from $P_1$ only by a
Hilbert--Schmidt operator:
$$\PP_{P_1}\DEF\SET{P}{P\text{ is a \BP, } P_1-P 
  \text{ is Hilbert--Schmidt}}.$$
$\PP_{P_1}$ is isomorphic to the set of all
\FSs which are unitarily equivalent to $\omega_{P_1}$.  Further let
$\EE_{P_1}$ be the infinite dimensional analogue of the open unit
disk \cite{Si,Se}, consisting of all symmetric Hilbert--Schmidt
operators $Z$ from $\KK_1$ to $\KK_2$ with norm less than $1$
\begin{equation}
  \label{ccr:SYM}
  \EE_{P_1}\DEF\SET{Z\in\BB(\KK_1,\KK_2)}{Z=Z^\tau,\ 
    \NORM{Z}<1,\ Z\text{ is Hilbert--Schmidt}}
\end{equation}
(the condition $\NORM{Z}<1$ is equivalent to $P_1+Z\+Z$ being
positive definite on $\KK_1$). Then the following is more or less
well--known (cf.~\cite{Se}).
\begin{Prop}
  \label{ccr:prop:BP}\hspace*{\fill}\\
  The map $P\mapsto P_{21}{P_{11}}^{-1}$ defines a bijection from
  $\PP_{P_1}$ onto $\EE_{P_1}$, with inverse 
  \begin{equation}
    \label{ccr:PZ}
    Z\mapsto P_Z\DEF(P_1+Z)(P_1+Z\+Z)^{-1}(P_1+Z\+).
  \end{equation}
  The restricted symplectic group $\SP{P_1}{0}$ acts transitively on
  either set, in a way compatible with the above bijection, through
  the formulas
  \begin{align}
    P&\mapsto UPU\+\label{ccr:OP1}\\
    Z&\mapsto(U_{21}+U_{22}Z)(U_{11}+U_{12}Z)^{-1}. 
      \tag{\ref{ccr:OP1}$'$}\label{ccr:OP2}
  \end{align}
  The restrictions of these actions to the subgroup $\SP{\text{\rm
      HS}}{}$ remain transitive, as follows from the fact that, for
  $Z\in\EE_{P_1}$,
  \begin{equation}
    \label{ccr:UZ}
    U_Z\DEF(P_1+Z)(P_1+Z\+Z)^{-\2}+(P_2-Z\+)(P_2+ZZ\+)^{-\2}    
  \end{equation}
  lies in $\SP{\text{\rm HS}}{}$ and fulfills $U_ZP_1U_Z\+=P_Z$
  \textup{(}equivalently, under the action \eqref{ccr:OP2}, $U_Z$
  takes $0\in\EE_{P_1}$ to $Z$\textup{)}.
\end{Prop}
\begin{proof}
  Having made \KK into a \HSP, the conditions \eqref{ccr:BP} on $P$ to
  be a \BP may be rewritten as
  \begin{equation}
    \label{ccr:BP0}
    P=P\+=\1-\4{P}=P^2,\qquad CP\text{ is positive definite on }\RAN P;  
  \end{equation}
  or, in components:
  \begin{subequations}
    \begin{alignat}{2}
      P_{11}&={P_{11}}^*&&=P_1-\4{P_{22}},\label{ccr:BP1}\\
      P_{22}&={P_{22}}^*&&=P_2-\4{P_{11}},\label{ccr:BP2}\\
      P_{21}&=\4{{P_{21}}^*}&&=-{P_{12}}^*,\label{ccr:BP3}
    \end{alignat}
    \begin{align}
      {P_{11}}^2-P_{11}&={P_{21}}^*P_{21},\label{ccr:BP4}\\
      {P_{22}}^2-P_{22}&={P_{12}}^*P_{12},\label{ccr:BP5}\\
      (P_1-P_{11})P_{12}&=P_{12}P_{22},\label{ccr:BP6}\\
      (P_2-P_{22})P_{21}&=P_{21}P_{11},\label{ccr:BP7}
    \end{align}
    \begin{equation}
      \label{ccr:BP8}
      \begin{pmatrix}
        P_{11} & P_{12} \\ -P_{21} & -P_{22}
      \end{pmatrix}
      \text{ is positive definite on\,}\RAN P.
    \end{equation}
  \end{subequations}
  Moreover, $P_1-P$ is Hilbert--Schmidt \IFF $P_2P$ is.
  
  Now let $P\in\PP_{P_1}$. Then $P_{22}\leq0$ by
  \eqref{ccr:BP8}, hence, by \eqref{ccr:BP1},
  $$P_{11}=P_1-\4{P_{22}}\geq P_1,$$
  so that $P_{11}$ has a bounded
  inverse. Thus $Z\DEF P_{21}{P_{11}}^{-1}$ is a well--defined
  Hilbert--Schmidt operator.  By \eqref{ccr:BP1}--\eqref{ccr:BP3} and
  \eqref{ccr:BP7},
  \begin{equation*}
    \begin{split}
      Z-\4{Z^*}&=P_{21}{P_{11}}^{-1}-\4{{P_{11}}^{-1}{P_{21}}^*}\\
      &=\4{{P_{11}}^{-1}}\bigl((P_2-P_{22})P_{21}
      -P_{21}P_{11}\bigr){P_{11}}^{-1}\\
      &=0,
    \end{split} 
  \end{equation*}
  so $Z$ is symmetric in the sense of \eqref{ccr:SYM}. Furthermore,
  by~\eqref{ccr:BP4},
  \begin{equation}
    \label{ccr:Z}
    \begin{split}
      P_1-Z^*Z&=P_1-{P_{11}}^{-1}{P_{21}}^*P_{21}{P_{11}}^{-1}\\
      &=P_1-{P_{11}}^{-1}({P_{11}}^2-P_{11}){P_{11}}^{-1}\\
      &={P_{11}}^{-1}
    \end{split}
  \end{equation}
  is positive definite on $\KK_1$, which proves $Z\in\EE_{P_1}$.

  Next let $Z\in\EE_{P_1}$ and let $P_Z$ be given by \eqref{ccr:PZ}. We
  associate with $Z$ an operator 
  \begin{equation}
    \label{ccr:Y}
    Y\DEF(P_1+Z\+Z)^{-1}=(P_1-Z^*Z)^{-1}
  \end{equation}
  which is bounded by assumption. Then $P_Z=P_Z\+=P_Z^2$ since
  $(P_1+Z\+)(P_1+Z)=Y^{-1}$. To prove that $P_Z+\4{P_Z}=\1$ holds,
  note that $ZY^{-1}=\4{Y}^{\,-1}Z$ and therefore $\4{Y}Z=ZY$,
  $YZ\+=Z\+\4{Y}$. It follows that
  \begin{equation*}
    \begin{split}
      P_Z+\4{P_Z}&=(P_1+Z)Y(P_1+Z\+)+(P_2-Z\+)\4{Y}(P_2-Z)\\
      &=Y+ZY+YZ\++ZZ\+\4{Y}+\4{Y}-YZ\+-ZY+Z\+ZY\\
      &=Y^{-1}Y+\4{Y}^{\,-1}\4{Y}\\
      &=P_1+P_2\\
      &=\1. 
    \end{split} 
  \end{equation*}
  Since $P_2P_Z$ is clearly HS and since 
  \begin{equation}
    \label{ccr:CPZ}
    CP_Z=(P_1-Z)Y(P_1-Z^*)    
  \end{equation}
  is positive definite on $\RAN P_Z=\RAN(P_1+Z)$, we get that
  $P_Z\in\PP_{P_1}$ as desired.
  
  To show that these two maps are mutually inverse, let first
  $Z\in\EE_{P_1}$. Then $(P_Z)_{21}{(P_Z)_{11}}^{-1}=ZYY^{-1}=Z$.
  Conversely, let $P\in\PP_{P_1}$ be given and set $Z\DEF
  P_{21}{P_{11}}^{-1}$. Then $ZP_{11}=P_{21}$ and
  $P_{11}Z\+={P_{21}}\+=P_{12}$. By \eqref{ccr:Z} and \eqref{ccr:Y},
  $Y=P_{11}$, hence $P_{11}Z\+=Z\+\4{P_{11}}$. Thus we get
  \begin{equation*}
    \begin{split}
      P-P_Z&=P-(P_1+Z)P_{11}(P_1+Z\+)\\
      &=P-P_{11}-ZP_{11}-P_{11}Z\+-ZP_{11}Z\+\\
      &=P-P_{11}-P_{21}-P_{12}-ZZ\+\4{P_{11}}\\
      &=P_{22}-ZZ\+\4{P_{11}}\\
      &=P_2-(P_2+ZZ\+)\4{P_{11}}\text{ (by \eqref{ccr:BP2})}\\
      &=0.
    \end{split} 
  \end{equation*}
  
  It remains to prove the statements about the group actions. It is
  fairly obvious that $\SP{P_1}{0}$ acts on $\PP_{P_1}$
  via \eqref{ccr:OP1}. The proof that $U_Z$ is a Bogoliubov operator
  which takes $P_1$ to $P_Z$ is also straightforward. To show that
  $U_Z\in\SP{\text{HS}}{}$, let $Y$ be given by \eqref{ccr:Y}. Then
  $$Y^\2-P_1=Y^\2(P_1-Y^{-1})(P_1+Y^{-\2})^{-1}
  =Y^\2Z^*Z(P_1+Y^{-\2})^{-1}$$
  is of trace class. Therefore
  $(U_Z-\1)P_1=(P_1+Z)Y^\2-P_1=Y^\2-P_1+ZY^\2$ is HS, which implies
  $U_Z\in\SP{\text{HS}}{}$.
  
  Finally we have to show that the action \eqref{ccr:OP1} on
  $\PP_{P_1}$ carries over to the action \eqref{ccr:OP2} on
  $\EE_{P_1}$. Thus, for given $Z\in\EE_{P_1}$ and $U\in\SP{P_1}{0}$,
  we have to compute the operator $Z'=P'_{21}{P'_{11}}^{-1}$ which
  corresponds to $P'=UP_ZU\+$. By definition,
  \begin{equation}
    \label{ccr:P'}
    \begin{aligned}
      P'_{21}&=(U_{21}+U_{22}Z)Y(U_{11}+U_{12}Z)^*,\\
      P'_{11}&=(U_{11}+U_{12}Z)Y(U_{11}+U_{12}Z)^*.
    \end{aligned}
  \end{equation}
  Suppose that $(U_{11}+U_{12}Z)f=0$ for some $f\in\KK_1$. Then
  $\NORM{f}_{P_1}=\NORM{{U_{11}}^{-1}U_{12}Zf}_{P_1}$. Since
  $\NORM{{U_{11}}^{-1}U_{12}}^2 =\NORM{{U_{12}}^*{U_{11}}^{-1*}
    {U_{11}}^{-1}U_{12}}
  =\NORM{{U_{12}}^*(P_1+U_{12}{U_{12}}^*)^{-1}U_{12}}
  =\NORM{U_{12}}^2/(1+\NORM{U_{12}}^2)<1$ and $\NORM{Z}<1$, it follows
  that $f=0$. Hence $U_{11}+U_{12}Z$ is injective, and, as a Fredholm
  operator with vanishing index \eqref{ccr:IND}, it has a bounded
  inverse. So we get from \eqref{ccr:P'} that
  $Z'=P'_{21}{P'_{11}}^{-1}=(U_{21}+U_{22}Z)(U_{11}+U_{12}Z)^{-1}$ as
  claimed.
\end{proof}
\begin{Rem}
  It is known that the unique (up to a phase) cyclic vector
  in $\mathcal{F}_s(\KK_1)$ inducing the \FS $\omega_{P_Z}$ is
  proportional to $\exp(\2Z\+a^*a^*)\Omega_{P_1}$.
\end{Rem}
The following construction will enable us to assign, in an unambiguous
way, to each Bogoliubov operator $V\in\SP{P_1}{}$ a \BP $P_V$ such
that \eqref{ccr:PV} holds.
\begin{Lem}
  \label{ccr:lem:H}\hspace*{\fill}\\
  Let $\HH\subset\KK$ be a closed *--invariant subspace such that
  $\kappa|_{\HH\times \HH}$ is nondegenerate and such that $[P_1,E]$
  is Hilbert--Schmidt, where $E$ is the orthogonal projection onto
  \HH.  Let $A\DEF ECE$ be the self--adjoint operator, invertible on
  \HH, such that $\kappa(f,g)=\langle{f,Ag}\rangle_{P_1},\ f,g\in\HH$,
  and let $A_\pm$ be the unique positive operators such that
  $A=A_+-A_-$ and $A_+A_-=0$. Further let $A^{-1}$ be defined as the
  inverse of $A$ on \HH and as zero on $\HH^\bot$, and similarly for
  $A_\pm^{-1}$. Then $A^{-1}C$ is the $\kappa$--orthogonal projection
  onto \HH, $p_+\DEF A_+^{-1}C$ is a \BP of
  $(\HH,\kappa|_{\HH\times\HH})$, and $P_2p_+$ is Hilbert--Schmidt.
  Moreover, $p_+=P_1E$ \IFF $[P_1,E]=0$.
\end{Lem}
\begin{proof}
  Let $E'\DEF \1-E$. Since $ECE'$ and $E'CE$ are compact by
  assumption, $C-ECE'-E'CE=A+E'CE'$ is a Fredholm operator on \KK with
  vanishing index. Hence $A$ is Fredholm on \HH with $\IND A=0$. $A$
  is injective since $\kappa$ is nondegenerate on \HH. It is therefore
  a bounded bijection on \HH with a bounded inverse (the same holds
  true for $A_\pm$ as operators on $\RAN A_\pm$). Thus $Q\DEF A^{-1}C$
  is well--defined. It fulfills $Q^2=A^{-1}(ECE)A^{-1}C=Q$ and
  $Q\+=C(CA^{-1})C=Q$. So $Q$ is a projection, self--adjoint \WRT
  $\kappa$. Since its range equals $\RAN A^{-1}=\HH$, it is the
  $\kappa$--orthogonal projection onto \HH.
  
  By a similar argument, $p_+$ is also a $\kappa$--orthogonal
  projection.  It is straightforward to see that $p_+=P_1E$ \IFF
  $[P_1,E]=0$. To show that $p_+$ is actually a \BP of \HH
  (cf.~\eqref{ccr:BP0}), note that $\4{A_+}=A_-$ because of $\4{A}=-A$
  (and uniqueness of $A_\pm$). This implies
  $p_++\4{p_+}=A_+^{-1}C-A_-^{-1}C=A^{-1}C=\1_\HH$. Positive
  definiteness of $Cp_+$ on $\RAN p_+$ follows from
  $\langle{f,Cp_+f}\rangle_{P_1}=\NORM{A_+^{-1/2}Cf}_{P_1}^2$.
  
  To prove that $P_2p_+$ is HS, let $D\DEF EP_1E-A_+$. Since
  $EP_1E-EP_2E=A=A_+-A_-$, we have $D=\4{D}$. We claim that $D$ is of
  trace class. Since $ECE'$ is HS,
  \begin{equation*}
    \begin{split}
      ECE'CE &= EC(\1-E)CE\\
        &= E-(ECE)^2\\
        &= E-A^2\\
        &= (E+\ABS{A})(E-\ABS{A})
    \end{split}
  \end{equation*}
  is of trace class. Since $E+\ABS{A}$ has a bounded inverse (as an
  operator on \HH) and since $\ABS{A}=A_++A_-$, it follows that
  $E-\ABS{A}=EP_1E+EP_2E-A_+-A_-=D+\4{D}=2D$ is of trace class as
  claimed. As a consequence, $A_+P_2=(EP_1E-D)P_2$ is HS ($P_1EP_2$
  is HS by assumption). By boundedness of $A_+^{-1}$,
  $p_+P_2=-A_+^{-2}(A_+P_2)$ and $P_2p_+=(p_+P_2)\+$ are also HS. This
  completes the proof.
\end{proof}
Now let $V\in\SP{P_1}{}$. We already showed in Section~\ref{sec:IMP}
that the restriction of $\kappa$ to $\ker V\+$ is nondegenerate. We
also showed in the proof of Theorem~\ref{ccr:th:IMP} that $[P_1,V']$
is Hilbert--Schmidt where $V'$ is the isometry arising from polar
decomposition of $V$. Hence $[P_1,E]$ is Hilbert--Schmidt where
\begin{equation}
  \label{ccr:E}
  E\DEF C(\1-V'{V'}^*)C
\end{equation}
is the orthogonal projection onto $\ker V\+$. Thus
Lemma~\ref{ccr:lem:H} applies to $\HH=\ker V\+$.
\begin{Def}
  \label{ccr:def:UW}\hspace*{\fill}\\
  For $V\in\SP{P_1}{}$, let $p_V\DEF p_+$ be the \BP of $(\ker V\+,
  \kappa|_{\ker V\+\times\ker V\+})$ given by Lemma~\ref{ccr:lem:H},
  and set
  \begin{alignat}{2}
    P_V&\DEF VP_1V\++p_V&&\in\PP_{P_1},\label{ccr:PV2}\\
    Z_V&\DEF(P_V)_{21}{(P_V)_{11}}^{-1}&&\in\EE_{P_1}
    \label{ccr:ZV2}
  \end{alignat}
  (cf.\ Proposition~\ref{ccr:prop:BP}). Further let
  $U_V\in\SP{\text{\rm HS}}{}$ be the Bogoliubov operator associated
  with $Z_V$ according to \eqref{ccr:UZ}, and define $W_V\DEF
  U_V\+V\in\SP{\text{\rm diag}}{}$.
\end{Def}\noindent
$P_V$ clearly is a \BP which satisfies \eqref{ccr:PV}. Actually, any
\BP $P$ fulfilling $V\+PV=P_1$ or, equivalently, $PV=VP_1$, is of the
form $P=VP_1V\++q$ where $q$ is some \BP of $(\ker V\+,\kappa|_{\ker
  V\+\times\ker V\+})$. What had to be proved above is that $q$ can be
chosen such that $P_2q$ is Hilbert--Schmidt, which is not obvious in
the case $\dim\ker V\+=\infty$. In fact, any such extension of
$VP_1V\+$ would suffice for what follows.

The condition $V\+P_VV=P_1$ translates into the condition 
\begin{equation}
  \label{ccr:ZV}
  Z_VV_{11}=V_{21}
\end{equation}
for $Z_V$. Again, each $Z\in\EE_{P_1}$ fulfilling \eqref{ccr:ZV} would
do, but we prefer to have a definite choice. It follows from symmetry
\eqref{ccr:SYM} that any $Z$ which solves \eqref{ccr:ZV} must have the
form
\begin{equation}
  \label{ccr:ZV1}
  Z=V_{21}{V_{11}}^{-1}+{V_{22}}^{-1*}{V_{12}}^*p_{\ker{V_{11}}^*}+Z'
\end{equation}
where $p_\HH$ denotes the orthogonal projection onto some closed
subspace $\HH\subset\KK$, ${V_{11}}^{-1}$ and ${V_{22}}^{-1}$ have
been defined below \eqref{ccr:REL}, and $Z'$ is a symmetric
Hilbert--Schmidt operator from $\ker{V_{11}}^*$ to $\ker{V_{22}}^*$.
The freedom in the choice of $Z'$ corresponds to the freedom in the
choice of $q$. Note that $Z$ can be written, \WRT the decompositions
$\KK_1=\RAN V_{11}\oplus\ker{V_{11}}^*,\ \KK_2=\RAN
V_{22}\oplus\ker{V_{22}}^*$, as
\begin{equation}
  Z=\begin{pmatrix}
    p_{\RAN{V_{22}}}V_{21}{V_{11}}^{-1} & 
    {V_{22}}^{-1*}{V_{12}}^*p_{\ker{V_{11}}^*} \\
    p_{\ker{V_{22}}^*}V_{21}{V_{11}}^{-1} & Z'
  \end{pmatrix}.
\end{equation}
The Hilbert--Schmidt norm of $Z$ is minimized by choosing $Z'=0$, but
there are examples in which this choice violates the condition
$\NORM{Z}<1$, i.e.\ it does not always define an element of
$\EE_{P_1}$. This is in contrast to the CAR case where the choice
analogous to $Z'=0$ appeared to be natural (cf.~\eqref{TV}). 

The operators $U_V$ and $W_V$ constitute the product decomposition of
$V$ that was announced earlier, generalizing a construction given by
Maa\ss\ \cite{M} to the infinite dimensional case. $W_V$ is diagonal
because $P_1W_V=P_1U_V\+V=U_V\+P_VV=U_V\+VP_1=W_VP_1$. Explicitly,
one computes that
$$W_V=\begin{pmatrix}
  (P_1+Z_V\+Z_V)^\2V_{11} & 0 \\ 0 & (P_2+Z_VZ_V\+)^\2V_{22}
\end{pmatrix}$$
\WRT the decomposition $\KK=\KK_1\oplus\KK_2$. Let us summarize the
properties of these operators.
\begin{Prop}
  \label{ccr:prop:UW}\hspace*{\fill}\\
  Definition \ref{ccr:def:UW} establishes a product decomposition of
  $V\in\SP{P_1}{}$,
  $$V=U_VW_V,$$
  where $U_V\in\SP{\text{\rm HS}}{}$ and
  $W_V\in\SP{\text{\rm diag}}{}$ have the properties
  \begin{align*}
    \IND U_V&=0, & Z_{U_V}&=Z_V, & P_{U_V}&=P_V,\\ 
    \IND W_V&=\IND V, & Z_{W_V}&=0, & P_{W_V}&=P_1.
  \end{align*}
  In particular, if $V\in\SP{P_1}{0}$, then 
  $$U_V=\begin{pmatrix}
    \ABS{{V_{11}}^*} & V_{12}{v_{22}}^* \\ 
    V_{21}{v_{11}}^* & \ABS{{V_{22}}^*}
  \end{pmatrix},\qquad 
  W_V=\begin{pmatrix}
    v_{11} & 0 \\ 0 & v_{22}
  \end{pmatrix}$$ 
  where $v_{11}\DEF V_{11}\ABS{V_{11}}^{-1}$ and $v_{22}=\4{v_{11}}$
  are the unitary parts of $V_{11}$ and $V_{22}$; whereas if
  $V\in\SP{\text{\rm diag}}{}$, then $U_V=\1$ and $W_V=V$.
\end{Prop}
The well--known fact that the restricted symplectic group
$\SP{P_1}{0}$ is connected \cite{Se,Ca84} entails for $\SP{P_1}{}$
\begin{Cor}
  \label{ccr:cor:SP}\hspace*{\fill}\\
  $\SP{P_1}{}=\SP{\text{\rm HS}}{}\cdot\SP{\text{\rm diag}}{}$.
  The orbits of the action of $\SP{P_1}{0}$ on $\SP{P_1}{}$ are the
  subsets $\SP{P_1}{n},\ n\in\NN\cup\{\infty\}$. They coincide with
  the connected components of $\SP{P_1}{}$.
\end{Cor}
\subsection{Normal form of implementers}
\label{sec:CON}
The first step in the construction of implementers consists in a
generalization of the definition of ``bilinear Hamiltonians'' \cite{A71}
from the finite rank case to the case of bounded operators. If $H$ is
a finite rank operator on \KK such that $H=H^\tau=-H^*$, then $e^{HC}$
belongs to $\SP{\text{HS}}{}$. Expanding $H=\sum f_j\langle
g_j,.\rangle_{P_1}$, 
one obtains a skew--adjoint element $b_0(H)\DEF\sum f_jg_j^*$ in \CKG
which is a linear function of $H$, independent of the choice of
$f_j,g_j\in\KK$.  Then $\pi_{P_1}\bigl(b_0(H)\bigr)$ is essentially
skew--adjoint on \DD, and, if $b(H)$ denotes its closure,
$\exp\bigl(\2b(H)\bigr)$ is a unitary which implements the
automorphism induced by $e^{HC}$ \cite{A71,A82}.

Using Wick ordering, the definition of bilinear Hamiltonians can be
extended to arbitrary bounded symmetric\footnote{The bilinear
  Hamiltonian corresponding to an antisymmetric operator ($H=-H^\tau$)
  vanishes.} operators $H$:
\begin{equation}
  \label{ccr:SYMH}
  H_{11}={H_{22}}^\tau,\qquad H_{12}={H_{12}}^\tau,\qquad
  H_{21}={H_{21}}^\tau.
\end{equation}
Without loss of generality, we henceforth assume that
$\KK_1=L^2(\RR^d)$. Then let $\SS\subset\mathcal{F}_s(\KK_1)$ be the
dense subspace consisting of finite particle vectors $\phi$ with
$n$--particle wave functions $\phi^{(n)}$ in the Schwartz space
$\SS(\RR^{dn})$. The unsmeared annihilation operator $a(p)$ with
(invariant) domain $\SS$ is defined as usual
$$(a(p)\phi)^{(n)}(p_1,\dots,p_n)\DEF\sqrt{n+1}\,\phi^{(n+1)}
(p,p_1,\dots,p_n).$$
Let $a^*(p)$ be its \QF adjoint on
$\SS\times\SS$. Then Wick ordered monomials $a^*(q_1)\dotsm
a^*(q_m)a(p_1)\dotsm a(p_n)$ make sense as \QFs on $\SS\times\SS$
\cite{GJ,RS2}, and, for $\phi,\phi'\in\SS$,
$$\langle\phi,a^*(q_1)\dotsm a^*(q_m)a(p_1)\dotsm a(p_n)\phi'\rangle
\DEF\langle a(q_1)\dotsm a(q_m)\phi,a(p_1)\dotsm a(p_n)\phi'\rangle$$
is a Schwartz function to which tempered distributions can be applied.
In particular, the distributions $H_{jk}(p,q),\ j,k=1,2$, given by
\begin{align*}
  \langle{f,H_{11}g}\rangle_{P_1}&=\int\4{f(p)}H_{11}(p,q)g(q)\,
    dp\,dq,\\
  \langle{f,H_{12}g^*}\rangle_{P_1}&=\int\4{f(p)}H_{12}(p,q)\4{g(q)}\,
    dp\,dq,\\
  \langle{f^*,H_{21}g}\rangle_{P_1}&=\int f(p)H_{21}(p,q)g(q)\,
    dp\,dq,\\
  \langle{f^*,H_{22}g^*}\rangle_{P_1}&=\int f(p)H_{22}(p,q)\4{g(q)}\,
    dp\,dq
\end{align*}
for $f,g\in\SS(\RR^d)\subset\KK_1$, give rise to the following \QFs
on $\SS\times\SS$: 
\begin{align*}
  H_{11}a^*a &\DEF \int a(p)^*H_{11}(p,q)a(q)\,dp\,dq \\
  H_{12}a^*a^* &\DEF \int a(p)^*H_{12}(p,q)a(q)^*\,dp\,dq \\
  H_{21}aa &\DEF \int a(p)H_{21}(p,q)a(q)\,dp\,dq\\
  \WO{H_{22}aa^*} &\DEF\int a(q)^*H_{22}(p,q)a(p)\,dp\,dq=H_{11}a^*a.
\end{align*}
Wick ordering of $H_{22}aa^*$ is necessary to make this expression
well--defined. The last equality follows from symmetry of $H$:
$$H_{11}(p,q)=H_{22}(q,p),\quad H_{12}(p,q)=H_{12}(q,p),\quad
H_{21}(p,q)=H_{21}(q,p).$$ 
We next define $\WO{b(H)}$ and its Wick ordered powers as \QFs on
$\SS\times\SS$:
\begin{align*}
  \WO{b(H)} &\DEF H_{12}a^*a^*+2H_{11}a^*a+H_{21}aa,\\
  \WO{b(H)^l} &\DEF l!\sum_{\substack{l_1,l_2,l_3=0 \\ l_1+l_2+l_3=l}}^l
    \frac{2^{l_2}}{l_1!l_2!l_3!}H_{l_1,l_2,l_3},\qquad l\in\NN,\\
  \text{with}\qquad H_{l_1,l_2,l_3}\DEF\int &H_{12}(p_1,q_1)\dotsm 
    H_{12}(p_{l_1},q_{l_1})
    H_{11}(p_1',q_1')\dotsm H_{11}(p_{l_2}',q_{l_2}')\\
  \cdot\, &H_{21}(p_1'',q_1'')\dotsm H_{21}(p_{l_3}'',q_{l_3}'')
    a^*(p_1)\dotsm a^*(p_{l_1})a^*(q_1)\dotsm a^*(q_{l_1})\\
  \cdot\, &a^*(p_1')\dotsm a^*(p_{l_2}') a(q_1')\dotsm a(q_{l_2}')
    a(p_1'')\dotsm a(p_{l_3}'')a(q_1'')\dotsm a(q_{l_3}'')\\
  \cdot\, &dp_1\,dq_1\dots dp_{l_1}\,dq_{l_1}\,dp_1'\,dq_1'\dots
    dp_{l_2}'\,dq_{l_2}'\,dp_1''\,dq_1''\dots dp_{l_3}''\,dq_{l_3}''.
\end{align*}
The Wick ordered exponential of $\2b(H)$ is also well-defined on
$\SS\times\SS$, since only a finite number of terms contributes when
applied to vectors from $\SS$:
$$\EH\DEF\sum_{l=0}^\infty\frac{1}{l!2^l}\WO{b(H)^l}.$$
The important point is
that these \QFs are actually the forms of uniquely determined linear
operators, defined on the dense subspace \DD and mapping \DD into the
domain of (the closure of) any creation or annihilation operator,
provided that \cite{R78}
\begin{equation}
  \label{ccr:OPH}
  \NORM{H_{12}}<1,\qquad H_{12}\text{ is Hilbert--Schmidt.}
\end{equation}
These operators will be denoted by the same symbols as the \QFs. The
analogue of Lemma~\ref{lem:REL} is
\begin{Lem}
  \label{ccr:lem:CR}\hspace*{\fill}\\
  Let $H\in\BB(\KK)$ satisfy \eqref{ccr:SYMH} and \eqref{ccr:OPH}.
  Then the following commutation relations hold on \DD, for
  $f\in\KK_1$:
  \begin{align*}
    [H_{l_1,l_2,l_3},a(f)^*] &=l_2a(H_{11}f)^*H_{l_1,l_2-1,l_3}
      +2l_3H_{l_1,l_2,l_3-1}a\bigl((H_{21}f)^*\bigr),\\
    [a(f),H_{l_1,l_2,l_3}] &=2l_1a(H_{12}f^*)^*H_{l_1-1,l_2,l_3}+
      l_2H_{l_1,l_2-1,l_3}a({H_{11}}^*f),
  \end{align*}
  implying that
  \begin{align*}
    [\EH,a(f)^*] &=a(H_{11}f)^*\EH+\EH a\bigl((H_{21}f)^*\bigr),\\
    [a(f),\EH] &=a(H_{12}f^*)^*\EH+\EH a({H_{11}}^*f).
  \end{align*}
\end{Lem}
\begin{proof}
  Compute as in \cite{R78,CB1}.
\end{proof}
For given $V\in\SP{P_1}{}$, we are now looking for bounded symmetric
operators $H$ which satisfy \eqref{ccr:OPH} and the following
intertwiner relation on \DD
\begin{equation}
  \label{ccr:INT}
  \EH\pi_{P_1}(f)=\pi_{P_1}(Vf)\EH,\quad f\in\KK
\end{equation}
(taking the closure of $\pi_{P_1}(Vf)$ is tacitly assumed here). This
problem turns out to be equivalent to the determination of the
operators $Z$ done in \eqref{ccr:ZV}, \eqref{ccr:ZV1}. 
\begin{Lem}
  \label{ccr:lem:HZ}\hspace*{\fill}\\
  Each $Z\in\EE_{P_1}$ fulfilling \eqref{ccr:ZV} gives rise to a
  unique solution $H$ of the above problem through the formula
  $$H=\begin{pmatrix}
    V_{11}-P_1+Z\+V_{21} & Z\+ \\
    ({V_{22}}^*+{V_{12}}^*Z\+)V_{21} & {V_{22}}^*-P_2+{V_{12}}^*Z\+
  \end{pmatrix},$$ 
  and each solution arises in this way.
\end{Lem}
\begin{proof}
  Let us abbreviate $\eta_H\DEF\EH$. Choosing $f\in\KK_2$ resp.\ 
  $f\in\KK_1$ and inserting the definition of $\pi_{P_1}$, one finds
  that \eqref{ccr:INT} is equivalent to
  $$\eta_H a(g)=\bigl(a(V_{11}g)+a^*(V_{12}g^*)\bigr)\eta_H,\quad
  \eta_H a^*(g)=\bigl(a^*(V_{11}g)+a(V_{12}g^*)\bigr)\eta_H$$
  for
  $g\in\KK_1$. Using the commutation relations from
  Lemma~\ref{ccr:lem:CR}, these equations may be brought into Wick
  ordered form:
  \begin{align*}
    0 &= a^*\bigl((V_{12}+H_{12}V_{22})g^*\bigr)\eta_H + 
      \eta_H a\Bigl(\bigl((P_1+{H_{11}}^*)V_{11}-P_1\bigr)g\Bigr),\\
    0 &= a^*\bigl((P_1+H_{11}-V_{11}-H_{12}V_{21})g\bigr)\eta_H +
      \eta_H a\Bigl(\bigl(\4{H_{21}}-(P_1+{H_{11}}^*)V_{12}\bigr)g^*\Bigr).
  \end{align*}
  As in the CAR case (see the proof of Lemma~\ref{lem:HT}), these
  equations hold for all $g\in\KK_1$ \IFF
  \begin{subequations}
    \begin{align}
      0&= V_{12}+H_{12}V_{22},\label{ccr:H1}\\
      0&= P_1+H_{11}-V_{11}-H_{12}V_{21},\label{ccr:H2}\\
      0&= H_{21}-(P_2+H_{22})V_{21},\label{ccr:H3}\\
      0&= P_2-(P_2+H_{22})V_{22}\label{ccr:H4}
    \end{align}
  \end{subequations}
  (we applied complex conjugation and used ${H_{11}}^\tau=H_{22}$).
  
  Now assume that $H$ solves the above problem. It is then obvious
  from \eqref{ccr:SYMH}, \eqref{ccr:OPH} and \eqref{ccr:H1} that
  $Z\DEF{H_{12}}\+$ belongs to $\EE_{P_1}$ and fulfills
  \eqref{ccr:ZV}.
  
  Conversely, let $Z\in\EE_{P_1}$ satisfy \eqref{ccr:ZV}. If there
  exists a solution $H$ with $H_{12}=Z\+$, then $H_{11}$ is fixed by
  \eqref{ccr:H2}, $H_{22}$ must equal ${H_{11}}^\tau$, and $H_{21}$ is
  determined by \eqref{ccr:H3}. Thus there can be at most one solution
  corresponding to $Z$, and it is necessarily of the form stated in
  the proposition.

  It remains to prove that the so--defined $H$ has all desired
  properties, i.e.\ that $H_{21}$ is  symmetric and that \eqref{ccr:H4}
  holds, the rest being clear by construction. The first claim follows
  from \eqref{ccr:REL4}:
  $$H_{21}-\4{{H_{21}}^*}=({V_{22}}^*+{V_{12}}^*Z\+)V_{21}-
  {V_{12}}^*(V_{11}+Z\+V_{21})=0,$$
  and the second from \eqref{ccr:ZV}
  and \eqref{ccr:REL2}:
  $$(P_2+H_{22})V_{22}=({V_{22}}^*-{V_{12}}^*\4{Z})V_{22}=
  {V_{22}}^*V_{22}-{V_{12}}^*V_{12}=P_2.$$
\end{proof}
Inserting the formula~\eqref{ccr:ZV1} for $Z$, one obtains
\begin{align*}
  H_{11} &= {V_{11}}^{-1*}-P_1-p_{{\ker V_{11}}^*}
    V_{12}{V_{22}}^{-1}V_{21} + {Z'}\+V_{21},\\
  H_{12} &= -V_{12}{V_{22}}^{-1}-{V_{11}}^{-1*}{V_{21}}^* 
    p_{{\ker V_{22}}^*} + {Z'}\+,\\
  H_{21} &= ({V_{22}}^{-1}-{V_{12}}^*{V_{11}}^{-1*}{V_{21}}^*
    p_{{\ker V_{22}}^*})V_{21} + {V_{12}}^*{Z'}\+V_{21},\\
  H_{22} &= {V_{22}}^{-1}-P_2-{V_{12}}^*{V_{11}}^{-1*}{V_{21}}^*
    p_{{\ker V_{22}}^*} + {V_{12}}^*{Z'}\+.
\end{align*}
$H$ corresponds to Ruijsenaars' ``associate'' $\Lambda$ \cite{R78}. If
one compares the above formula for $H$ with Ruijsenaars' formula for
$\Lambda$ in the case of automorphisms ($\ker {V_{jj}}^*=\{0\},\ 
j=1,2;\ Z'=0$), one finds that the off--diagonal components carry
opposite signs. This is due to the fact that Ruijsenaars constructs
implementers for the transformation induced by $CVC$ rather than $V$,
cf.\ (3.27) and (3.29) in \cite{R78}.

Note that
\begin{equation}
  \label{ccr:EOM}
  \EH\Omega_{P_1}=\exp(\3H_{12}a^*a^*)\Omega_{P_1}.
\end{equation}
By Ruijsenaars' computation \cite{R78} (see also \cite{Se}), 
the norm of such vectors is
$$\BigNORM{\EH\Omega_{P_1}}=
\left({\det}_{\KK_1}(P_1+H_{12}{H_{12}}\+)\right)^{-1/4}.$$
\begin{Def}
  \label{ccr:def:PSIA}\hspace*{\fill}\\
  Let $V\in\SP{P_1}{}$, and let $P_V,\ Z_V$ and $H_V$ be the operators
  associated with $V$ according to Definition~\ref{ccr:def:UW} and
  Lemma~\ref{ccr:lem:HZ}. Choose a $\kappa$--\ONB $g_1,g_2,\dotsc$ in
  \begin{equation}
    \label{ccr:KV}
    \mathfrak{k}_V\DEF P_V(\ker V\+),
  \end{equation}
  i.e.\ a basis such that $\kappa(g_j,g_k)=\delta_{jk}$ (this is
  possible because the restriction of $\kappa$ to $\mathfrak{k}_V$ is
  positive definite). Note that $\dim\mathfrak{k}_V=-\2\IND V$. Let
  $\psi_j$ be the isometry obtained by polar decomposition of the
  closure of $\pi_{P_1}(g_j)$. Then define operators $\PV{\alpha}$ on
  \DD, for any multi--index $\alpha=(\alpha_1,\dots,\alpha_l)$ with
  $\alpha_j\leq\alpha_{j+1}$ (or $\alpha=0$) as in \eqref{ccr:ALPHA},
  as
  \begin{equation}
    \label{ccr:PSIA}
    \PV{\alpha}\DEF\left({\det}_{\KK_1}(P_1+Z_V\+Z_V)
    \right)^{\frac{1}{4}}\psi_{\alpha_1}\dotsm\psi_{\alpha_l}
    \WO{\exp\left(\3b(H_V)\right)}.
  \end{equation}
\end{Def}
\begin{Th}
  \label{ccr:th:PSIA}\hspace*{\fill}\\
  The $\PV{\alpha}$ extend continuously to isometries
  \textup{(}denoted by the same symbols\textup{)} on the symmetric
  Fock space $\mathcal{F}_s(\KK_1)$ such that
  \begin{equation}
    \label{ccr:CUNTZ}
    \PV{\alpha}^*\PV{\beta}=\delta_{\alpha\beta}\1,\qquad
    \sum_\alpha\PV{\alpha}\PV{\alpha}^*=\1,
  \end{equation}
  and such that, for any element $w$ of the Weyl algebra $\WK$,
  \begin{equation}
    \label{ccr:IMP}
    \rho_V(w)=\sum_\alpha\PV{\alpha}w\PV{\alpha}^*.    
  \end{equation}
  The infinite sums converge in the strong topology.
\end{Th}
\begin{proof}
  By \eqref{ccr:CCR} we have
  $\pi_{P_1}(g_j)^*\pi_{P_1}(g_j)=\1+\pi_{P_1}(g_j)\pi_{P_1}(g_j)^*$
  on \DD, so the closure of $\pi_{P_1}(g_j)$ is injective, and
  $\psi_j$ is isometric. It is also easy to see, using \eqref{ccr:INT},
  the CCR and $\NORM{\PV{\alpha}\Omega_{P_1}}=1$, that for
  $f_1,\dots,f_m,h_1,\dots,h_n\in\KK$ 
  \begin{multline*}
    \langle\PV{\alpha}\pi_{P_1}(f_1\dotsm f_m)\Omega_{P_1},
    \PV{\alpha}\pi_{P_1}(h_1\dotsm h_n)\Omega_{P_1}\rangle\\
    =\langle\pi_{P_1}(f_1\dotsm f_m)\Omega_{P_1}, \pi_{P_1}(h_1\dotsm
    h_n)\Omega_{P_1}\rangle.
  \end{multline*}
  Hence $\PV{\alpha}$ is isometric on \DD and has a continuous
  extension to an isometry on $\mathcal{F}_s(\KK_1)$.
  
  Let $\HH_j\DEF{\rm{span}}(g_j,g_j^*)$, so that
  $\psi_j\in\WW(\HH_j)''$ by virtue of Lemma~\ref{lem:AFF}. Since
  $\HH_j\subset\ker V\+$, the duality relation \eqref{ccr:DUAL}
  implies that $\WW(\HH_j)\subset\WW(\RAN V)'$. Now let $f\in\RE\KK$
  and $\phi\in\DD$.  Since $\phi$ is an entire analytic vector for
  $\pi_{P_1}(f)$ \cite{AS71}, since \DD is invariant under
  $\pi_{P_1}(f)$, and since $\4{\pi_{P_1}(Vf)}$ is affiliated with
  $\WW(\RAN V)$ by Lemma~\ref{lem:AFF} (the bar denotes closure), it
  follows from \eqref{ccr:INT} that
  \begin{eqnarray*}
    \PV{\alpha}w(f)\phi &=&
    \sum_{n=0}^\infty\frac{i^n}{n!}\PV{\alpha}(\pi_{P_1}(f))^n\phi\\ 
    &=&\sum_{n=0}^\infty\frac{i^n}{n!}\psi_{\alpha_1}\dotsm
      \psi_{\alpha_l}\bigl(\4{\pi_{P_1}(Vf)}\bigr)^n\PV{0}\phi\\ 
    &=&\sum_{n=0}^\infty\frac{i^n}{n!}\bigl(\4{\pi_{P_1}(Vf)}\bigr)^n
      \PV{\alpha}\phi\\  
    &=&w(Vf)\PV{\alpha}\phi.
  \end{eqnarray*}
  By continuity, this entails
  \begin{equation}
    \label{ccr:INTA}
    \PV{\alpha}w=\rho_V(w)\PV{\alpha},\qquad w\in\WK.    
  \end{equation}
  We next claim that
  \begin{equation}
    \label{ccr:PSI0}
    \psi_j^*\PV{0}=0
  \end{equation}
  or, equivalently, that $\pi_{P_1}(g_j)^*\PV{0}=0$. To see this,
  apply Lemma~\ref{ccr:lem:CR} and write $\pi_{P_1}(g_j)^*\PV{0}$ in
  Wick ordered form:
  $$\pi_{P_1}(g_j)^*\PV{0}=a\bigl((P_1+H_{12})g_j^*\bigr)^*\PV{0}+
  \PV{0}a\bigl((P_1+{H_{11}}^*)g_j\bigr)$$
  on \DD, with $H\DEF H_V$.
  Then \eqref{ccr:PSI0} holds \IFF
  \begin{equation}
    (P_1+H_{12})g_j^*=0,\qquad(P_1+{H_{11}}^*)g_j=0.
    \tag{\ref{ccr:PSI0}$'$}\label{ccr:PSI1}
  \end{equation}
  Now $g_j\in\RAN P_V$ is equivalent to $g_j^*\in\ker P_V=\ker
  CP_V=\ker(P_1+H_{12})$ (we used \eqref{ccr:CPZ}). This proves the
  first equation in \eqref{ccr:PSI1}. It also shows that
  ${H_{12}}^*g_j=-P_2g_j$. Hence by Lemma~\ref{ccr:lem:HZ},
  $$(P_1+{H_{11}}^*)g_j=({V_{11}}^*+{V_{21}}^*{H_{12}}^*)g_j
  =({V_{11}}^*-{V_{21}}^*)g_j=P_1V\+g_j=0$$
  which proves the second
  equation in \eqref{ccr:PSI1} and therefore \eqref{ccr:PSI0}.
  
  The orthogonality relation $\PV{\alpha}^*\PV{\beta}=0$
  ($\alpha\not=\beta$) now follows from \eqref{ccr:PSI0} and from
  $\WW(\HH_j)\subset\WW(\HH_k)'$ ($j\not=k$) which in turn is a
  consequence of $\kappa(\HH_j,\HH_k)=0$ and \eqref{ccr:DUAL}.
  
  The proof of the completeness relation $\sum\PV{\alpha}
  \PV{\alpha}^*=\1$ is facilitated by invoking the product
  decomposition $V=U_VW_V$ from Proposition~\ref{ccr:prop:UW}. Set
  $f_j\DEF U_V\+g_j$ to obtain a $\kappa$--\ONB $f_1,f_2,\dotsc$ in
  $\mathfrak{k}_{W_V}=P_1(\ker W_V\+)$. Let $\psi_j'$ be the isometric
  part of $a(f_j)^*$. An application of Definition~\ref{ccr:def:PSIA}
  to $W_V$ yields implementers
  $\PW{\alpha}=\psi'_{\alpha_1}\dotsm\psi'_{\alpha_l}\PW{0}$ for
  $W_V$. $Z_{W_V}=0$ entails that
  \begin{equation}
    \label{ccr:PWO}
    \PW{\alpha}\Omega_{P_1}=\psi'_{\alpha_1}\dotsm
    \psi'_{\alpha_l}\Omega_{P_1}.
  \end{equation}
  One computes, using the CCR, that
  \begin{equation}
    \label{ccr:P'O}
    \psi'_{\alpha_1}\dotsm\psi'_{\alpha_l}\Omega_{P_1}=\phi'_\alpha,
  \end{equation}
  where the $\phi'_\alpha$ are the cyclic vectors associated with the
  pure state $\omega_{P_1}\0\rho_{W_V}=\omega_{P_1}$ as in
  Proposition~\ref{ccr:prop:CYC}. Let $\mathcal{F}'_\alpha$ be the
  closure of $\WW(\RAN W_V)\phi'_\alpha$. Since the
  $\mathcal{F}'_\alpha$ are irreducible subspaces for $\WW(\RAN W_V)$
  by Proposition~\ref{ccr:prop:CYC}, they must coincide with the
  irreducible subspaces $\RAN\PW{\alpha}$.
  $\oplus\mathcal{F}'_\alpha=\mathcal{F}_s(\KK_1)$ then implies
  completeness of the $\PW{\alpha}$.
 
  The proof will be completed by showing that
  \begin{equation}
    \label{ccr:COMP}
    \PV{\alpha}=\PUV\PW{\alpha}
  \end{equation}
  holds where \PUV is the unitary implementer for $U_V$
  given by Definition~\ref{ccr:def:PSIA}. It suffices to show that
  \eqref{ccr:COMP} holds on $\Omega_{P_1}$ since any bounded operator
  fulfilling \eqref{ccr:INTA} is already determined by its value on
  $\Omega_{P_1}$.  Because of $Z_{U_V}=Z_V$ we have
  \begin{equation}
    \label{ccr:PSIU}
    \PV{0}\Omega_{P_1}=\PUV\Omega_{P_1},    
  \end{equation}
  so it remains to show that
  $\psi_{\alpha_1}\dotsm\psi_{\alpha_l}\PUV\Omega_{P_1}=
  \PUV\psi'_{\alpha_1}\dotsm\psi'_{\alpha_l}\Omega_{P_1}$. We
  claim that
  \begin{equation}
    \label{ccr:PJPU}
    \psi_j\PUV=\PUV\psi'_j.
  \end{equation}
  For let $T$ (resp.\ $T'$) be the closure of $\pi_{P_1}(g_j)$ (resp.\ 
  $\pi_{P_1}(f_j)$). It suffices to prove that
  $\PUV\bigl(D(T')\bigr)=D(T)$ and that
  \begin{equation}
    \label{ccr:PUT}
    T\PUV=\PUV T'.
  \end{equation}
  \eqref{ccr:PUT} clearly holds on \DD. Now let $\phi\in D(T')$, and
  choose $\phi_n\in\DD$ with $\phi_n\to\phi$ and $T'\phi_n\to T'\phi$.
  Then $\PUV\phi_n\in\DD$ converges to $\PUV\phi$, and
  $T\PUV\phi_n=\PUV T'\phi_n$ converges to $\PUV T'\phi$. It follows
  that $\PUV\phi\in D(T)$ and $T\PUV\phi=\PUV T'\phi$, i.e.\ that
  $T\PUV\supset\PUV T'$. In the same way one obtains that
  $T'\PUV^*\supset\PUV^* T$, so that \eqref{ccr:PUT} and
  \eqref{ccr:PJPU} hold. (An alternative proof of \eqref{ccr:PJPU}
  goes as follows. Let $T^\pm$ (resp.\ ${T'}^\pm$) be the
  self--adjoint operators corresponding to $T$ (resp.\ $T'$) as in the
  proof of Lemma~\ref{lem:AFF}. Then one has $D(T)=D(T^+)\cap D(T^-)$
  and $T=T^+-iT^-$, and similar for $T'$.  There holds
  $\PUV\exp(it{T'}^\pm)\PUV^*=\exp(itT^\pm),\ t\in\RR$.  Therefore
  \PUV maps $D({T'}^\pm)$ onto $D(T^\pm)$, and
  $\PUV{T'}^\pm\PUV^*=T^\pm$. Consequently,
  $\PUV\bigl(D(T')\bigr)=D(T)$ and $\PUV T'\PUV^*=T$.  This implies
  that $\PUV\psi'_j\PUV^*=\psi_j$ as claimed.)
  
  The proof is complete since \eqref{ccr:CUNTZ} and \eqref{ccr:INTA}
  together imply \eqref{ccr:IMP}.
\end{proof}
\begin{Cor}
  \label{ccr:cor:HRV}\hspace*{\fill}\\
  There is a unitary isomorphism from $H(\rho_V)$, the \HSP generated
  by the $\PV{\alpha}$, onto the symmetric Fock space
  $\mathcal{F}_s(\mathfrak{k}_V)$ over $\mathfrak{k}_V$, which maps
  $\PV{\alpha}$ to $c_\alpha a^*(g_{\alpha_1})\dotsm a^*(g_{\alpha_l})
  \Omega$, where the normalization factor $c_\alpha$ is defined in
  \eqref{ccr:DEF}, and $a^*(g_j)$ and $\Omega$ are now creation
  operators and the Fock vacuum in $\mathcal{F}_s(\mathfrak{k}_V)$.
\end{Cor}
We shall see in Section~\ref{sec:CCRCH} that, for gauge invariant $V$,
the isomorphism described above is not only an isomorphism of graded
\HSPs but also of modules of the \GG.


%% file: sectors.tex
\specialsection[Sectors Reached by Quasi--free
  Endomorphisms]{Superselection Sectors Reached by Gauge Invariant
  Quasi--free Endomorphisms}
\label{sec:SECTORS}
In the present section we will apply our results from
Sections~\ref{sec:CAR} and \ref{sec:CCR} to the theory of \SSSs. We
are especially interested in the possible ``charge quantum numbers''
that can be realized by quasi--free endomorphisms. We will consider
situations where the theory of \DR applies, i.e.\ where the observable
algebra \AA consists of the invariant elements of a \FA (given in its
vacuum representation) under a group $G$ of gauge automorphisms of the
first kind. As mentioned in the introduction, the charge quantum
numbers of a localized endomorphism \rho then are labels for the
unitary representation of $G$ which is realized on the \HSP $H(\rho)$.

The CAR and \CCRs will play the r\^ole of the \FA, so that quasi--free
endomorphisms are from the outset endomorphisms of the \FA rather than
the observable algebra. The following simple observation shows that
one has to restrict attention to \emph{gauge invariant} endomorphisms,
i.e.\ to endomorphisms which commute with all gauge transformations.
\begin{Prop}
  \hspace*{\fill}\\
  Let \rho be an endomorphism of the \FA which is implemented by a
  \HSP $H(\rho)$ of isometries. Then $H(\rho)$ is invariant under $G$
  \IFF \rho is gauge invariant.
\end{Prop}
\begin{proof}
  Assume first that $H(\rho)$ is invariant under $G$. Let $R$ be the
  representation of $G$ on $H(\rho)$:
  $$R(\gamma)\DEF\gamma|_{H(\rho)},\quad\gamma\in G.$$
  $R$ is clearly unitary because one has for any $\gamma\in G$
  $$\langle R(\gamma)\Psi,R(\gamma)\Psi'\rangle\1=\gamma(\Psi^*\Psi')
  =\gamma(\langle\Psi,\Psi'\rangle\1)=\langle\Psi,\Psi'\rangle\1,
  \quad\Psi,\Psi'\in H(\rho).$$
  Writing $R(\gamma)$ as a matrix \WRT
  the \ONB $(\Psi_j)$, one gets $\gamma(\Psi_j)=\sum_k R(\gamma)_{kj}
  \Psi_k$ and $\sum_j R(\gamma)_{kj}\4{R(\gamma)_{lj}}=\delta_{kl}$,
  so that
  \begin{equation*}
    \begin{split}
      \gamma(\rho(F)) &=\sum_j\gamma(\Psi_j)\gamma(F)\gamma(\Psi_j)^*\\
      &=\sum_{k,l}\Big(\sum_j R(\gamma)_{kj}\4{R(\gamma)_{lj}}\Big)
        \Psi_k\gamma(F)\Psi_l^*\\
      &=\sum_k\Psi_k\gamma(F)\Psi_k^*\\
      &=\rho(\gamma(F))
    \end{split}
  \end{equation*}
  for any field $F$. Therefore \rho is gauge invariant.
  
  Conversely, assume that \rho is gauge invariant. Since the \FA is
  irreducibly represented, the \HSP $H(\rho)$ consists of all field
  operators $\Psi$ which intertwine between the vacuum representation
  and the representation induced by \rho:
  \begin{equation}
    \label{INTER1}
    \Psi\in H(\rho)\iff\Psi F=\rho(F)\Psi
    \text{ for all local fields }F
  \end{equation}
  (cf.\ \eqref{INTER}). Now let $\Psi\in H(\rho)$ and $\gamma\in G$.
  Then one has for any $F$
  $$\gamma(\Psi)F=\gamma\big(\Psi\gamma^{-1}(F)\big)
  =\gamma\big(\rho(\gamma^{-1}(F))\Psi\big)=\rho(F)\gamma(\Psi),$$
  so that $\gamma(\Psi)\in H(\rho)$ by \eqref{INTER1}.
\end{proof}
We are thus led to consider the following setting. We assume that a
distinguished \BP $P_1$ of \KK (CAR) resp.\ of $(\KK,\kappa)$ (CCR) is
given. The global \GG $G$ will consist of diagonal Bogoliubov
operators (resp.\ of the corresponding automorphisms)
\begin{alignat*}{2}
  G&\subset\II^0_{P_1}(\KK)\cap\IDIAG(\KK)\cong U(\KK_1)
    &\quad&\text{(CAR)}\\
  G&\subset\SP{P_1}{0}\cap\SP{\text{\rm diag}}{}\cong U(\KK_1)
    &\quad&\text{(CCR)}.
\end{alignat*}
Gauge transformations leave the \FS $\omega_{P_1}$ invariant. The
usual second quantization of $U\in G$ (or, more precisely, of
$U_{11}$) will be denoted by $\Gamma(U)$. This is the same as
$\Psi(U)$ defined by \eqref{PAV} resp.\ \eqref{ccr:PSIA}; the map
$U\mapsto\Gamma(U)$ is strongly continuous. The gauge automorphism
corresponding to $U$ on Fock space will be denoted by $\gamma_U$:
$$\gamma_U(F)\DEF \Gamma(U)F\Gamma(U)^*,$$
where $F$ can now be any
bounded operator on Fock space. The charge structure of all
implementable gauge invariant quasi--free endomorphisms $\rho_V$ will
be unraveled in the next sections, i.e.\ of all $\rho_V$ with $V$
contained in either of the following semigroups
\begin{alignat*}{2}
  \II_{P_1}(\KK)^G&\DEF\SET{V\in\II_{P_1}(\KK)}{[V,U]=0
    \text{ for all }U\in G}&\quad&\text{(CAR)}\\
  \SP{P_1}{}^G&\DEF\SET{V\in\SP{P_1}{}}{[V,U]=0\text{ for all }U\in G}
    &\quad&\text{(CCR)}.
\end{alignat*}
Bogoliubov operators commuting with $G$ will be called \emph{gauge
  invariant.} The notations $\II_{P_1}^n(\KK)^G$ resp.\ 
$\SP{P_1}{n}^G$ will be used to denote the subsets of $V$ with $\IND
V=-n$. At this stage of generality, it is not necessary to assume that
$G$ is a (strongly) compact topological group. However, if $G$ is
``too large'', then it can happen that $\II_{P_1}(\KK)^G$ and
$\SP{P_1}{}^G$ become trivial and consist only of the operators of the
form $e^{i\lambda}P_1+e^{-i\lambda}P_2,\ \lambda\in\RR$. On the other
hand, if $G$ is compact, then $\II_{P_1}(\KK)^G$ and $\SP{P_1}{}^G$
can be described more explicitly as follows. The representation of $G$
on \KK can be brought into the form
\begin{equation}
  \label{GDEC}
  \KK=\bigoplus_\xi(\KK_\xi\otimes\mathfrak{h}_\xi)
\end{equation}
where the sum extends over all equivalence classes $\xi$ of
irreducible representations of $G$ realized on \KK, $\mathfrak{h}_\xi$
is a finite dimensional subspace carrying a representation
$\mathcal{U}_\xi$ of class $\xi$, and $\KK_\xi$ is a \HSP with
dimension equal to the multiplicity of $\xi$. $G$ acts on
$\KK_\xi\otimes\mathfrak{h}_\xi$ like $\1_{\KK_\xi}\otimes
\mathcal{U}_\xi$. Gauge invariant Bogoliubov operators then have the
form $V=\oplus_\xi(V_\xi\otimes\1_{\mathfrak{h}_\xi})$, and there
exist non--surjective gauge invariant Bogoliubov operators \IFF at
least one $\KK_\xi$ is infinite dimensional. Also note that one should
require on physical grounds to have $-\1\in G$, but we do not need
this assumption at this point.

The analysis of the representations of $G$ on the implementing \HSPs
$H(\rho_V)$ is facilitated by the following lemma.
\begin{Lem}
  \label{lem:REPHRV}\hspace*{\fill}\\
  Let $V$ be an element of $\II_{P_1}(\KK)^G$ or $\SP{P_1}{}^G$. Then
  the representation of $G$ on $H(\rho_V)$ \textup{(}obtained by
  restricting $\gamma_U$ to $H(\rho_V)$\textup{)} is canonically
  unitarily equivalent to the representation on
  $H(\rho_V)\Omega_{P_1}$ \textup{(}obtained by restricting
  $\Gamma(U)$ to $H(\rho_V)\Omega_{P_1}$\textup{)}, via the map
  $\Psi\mapsto\Psi\Omega_{P_1}$.
\end{Lem}
\begin{proof}
  Obvious from gauge invariance of the vacuum.
\end{proof}
\begin{Ex}[{\sc The free Dirac field with gauge symmetry}]
  \label{ex:DIRAC}\hspace*{\fill}\\
  The setting specified above is abstracted from field theoretic
  examples of the following kind. Let $\HH=L^2(\RR^{2n-1},\CC^{2^n})$
  be the single particle space of the free time--zero Dirac field in
  $2n$ spacetime dimensions. Let $H=-i\vec\alpha\vec\nabla+\beta m$ be
  the free Dirac Hamiltonian, with spectral projections $p_\pm$
  corresponding to the positive resp.\ negative part of the spectrum
  of $H$. Tensored with $\1_N$, these operators act on the space
  $\HH'\DEF\HH\otimes\CC^N$. The gauge group $U(N)$ also acts
  naturally on $\HH'$. In the selfdual CAR formalism, one sets
  \begin{equation}
    \label{H'}
    \KK\DEF\HH'\oplus{\HH'}^*,
  \end{equation}
  where ${\HH'}^*$ is the \HSP conjugate to $\HH'$. There is a natural
  conjugation $f\mapsto f^*$ on \KK which is inherited from the
  antiunitary identity map $\HH'\to{\HH'}^*$. The \BP $P_1$
  corresponding to the vacuum representation of the field is then
  given by
  $$P_1\DEF p_+'\oplus\4{p_-'}$$
  with $p_\pm'=p_\pm\otimes\1_N$. Gauge transformations act like
  $U=(\1_\HH\otimes u)\oplus(\1_{\HH^*}\otimes\4{u}),\ u\in U(N)$, on
  \KK. They commute with $P_1$. With respect to the decomposition
  $\KK=\KK_1\oplus\KK_2$ induced by $P_1$, $U$ has the form
  \begin{equation}
    \label{GFORM}
    U=\begin{pmatrix}
      p_+\otimes u+\4{p_-}\otimes\4{u} & 0\\
      0 & \4{p_+}\otimes\4{u}+p_-\otimes u
    \end{pmatrix}.
  \end{equation}
  The field operators $\varphi_t$ at time $t$ are given by 
  $$\varphi_t(f)\DEF\pi_{P_1}(e^{itH'}f)=
  a(p_+'e^{itH'}f)^*+a(\4{p_-'}e^{-it\4{H'}}f^*)$$
  where $H'\DEF
  H\otimes\1_N,\ f\in\HH'$. They are solutions of the
  Dirac--Schr\"odinger equation
  $$-i\frac{d}{dt}\varphi_t(f)=\varphi_t(H'f),\qquad f\in D(H')$$
  (the
  minus sign is due to the fact that our field operators
  $\varphi_t(f)$ are complex linear in $f$). If $O$ is a double cone
  with base $B\subset\RR^{2n-1}$ at time $t$, then the local \FA
  associated with $O$ is generated by all $\varphi_t(f)$ with $\SUPP
  f\subset B$. The local observable algebra $\AA(O)$ is the fixed
  point subalgebra of this local \FA under the gauge action. (The
  whole net of local algebras is generated from these special ones by
  applying Lorentz transformations.) Bogoliubov operators in
  $\II_{P_1}(\KK)^{U(N)}$ which act like the identity on functions
  $f$ with $\SUPP f\cap B=\emptyset$ induce endomorphisms of \AA which
  are localized in the double cone $O$, cf.\ Prop.~\ref{prop:LOC}
  below. (The charge carried by such localized endomorphisms can be
  read off from our formulas in the following section.) All gauge
  invariant Bogoliubov operators $V$ have the form
  \begin{equation}
    \label{VFORM}
    V=(v\otimes\1_N)\oplus(\4{v}\otimes\4{\1_N})
  \end{equation}
  \WRT the decomposition \eqref{H'} where $v$ is some isometry of \HH.
  This holds because \eqref{H'}, with $\HH'=\HH\otimes\CC^N$, is the
  decomposition of \KK analogous to \eqref{GDEC}. Thus
  $\II_{P_1}(\KK)^{U(N)}$ is isomorphic to the semigroup of all
  isometries of \HH. This fact remains true if the group $U(N)$ is
  replaced by $SU(N)$, except for the case of $SU(2)$. In the latter
  case (and only in that case), the defining representation of the
  group is equivalent to its complex conjugate representation, since
  one has $JuJ^*=\4{u},\ u\in SU(2)$, with $J=\MAT{0}{-1}{1}{0}$.
  Therefore \eqref{H'} is \emph{not} the decomposition of \KK as in
  \eqref{GDEC}, and there exist Bogoliubov operators in
  $\II_{P_1}(\KK)^{SU(2)}$ which do not have the form \eqref{VFORM}.
  
  Similar constructions work for the free charged Klein--Gordon field.
\end{Ex}
\subsection{The charge of gauge invariant endomorphisms of the \CAR}
\label{sec:CARCH}
In this section we will compute the behaviour of the implementers
\PV{\alpha} defined by \eqref{PAV} under gauge transformations for
arbitrary gauge invariant $V$. Recall from \eqref{PSIAOM} that the
value of \PV{\alpha} on the Fock vacuum is given by
\begin{equation}
  \label{PAVOM}
  \PV{\alpha}\Omega_{P_1}=D_V\psi_{P_1}(g_{\alpha_1}\dotsm 
  g_{\alpha_l}e_1\dotsm e_{L_V})\exp(\3\4{T_V}a^*a^*)\Omega_{P_1}
\end{equation}
where $D_V$ is a numerical constant, $\{g_1,\dots,g_{M_V}\}$ is an
\ONB in $\mathfrak{k}_V=P_V(\ker V^*)$, $\{e_1,\dots,e_{L_V}\}$ is an
\ONB in $\mathfrak{h}_V=V_{12}(\ker V_{22})$, and $T_V$ is the
antisymmetric Hilbert--Schmidt operator defined in \eqref{TV}. We have
to calculate the transformed vectors
$\Gamma(U)\PV{\alpha}\Omega_{P_1},\ U\in G$.
\begin{Lem}
  \label{lem:GP}\hspace*{\fill}\\
  $\Gamma(U)$ implements the gauge automorphism $\rho_U$ in the
  twisted Fock representation $\psi_{P_1}$, i.e.\ one has
  $$\Gamma(U)\psi_{P_1}(a)=\psi_{P_1}(\rho_U(a))\Gamma(U),\qquad a\in\CK.$$
\end{Lem}
\begin{proof}
  It suffices to consider the case that $a$ belongs to $\KK$.  Recall
  from \eqref{PSI1} that $\psi_{P_1}(f)=i\pi_{P_1}(f)\PME$ for
  $f\in\KK$. Since $U$ is diagonal, the implementer $\Gamma(U)$ is
  even:
  $$[\Gamma(U),\PME]=0.$$
  This implies $\Gamma(U)\psi_{P_1}(f)\Gamma(U)^*=\psi_{P_1}(Uf)$.
\end{proof}
It turns out that the transformation properties of the exponential
term in \eqref{PAVOM} are also easily obtained.
\begin{Lem}
  \label{lem:EXP}\hspace*{\fill}\\
  If $V\in\II_{P_1}(\KK)^G$, then $\exp(\3\4{T_V}a^*a^*)\Omega_{P_1}$
  is invariant under all gauge transformations $\Gamma(U),\ U\in G$.
\end{Lem}
\begin{proof}
  If $T\in\mathfrak{H}_{P_1}$ (see \eqref{HP1}) has finite rank, then
  one readily verifies that
  $$\Gamma(U)(\3\4{T}a^*a^*)\Gamma(U)^*=\3(U\4{T}U^*)a^*a^*.$$
  Approximating $T_V$ by finite rank operators from
  $\mathfrak{H}_{P_1}$ relative to Hilbert--Schmidt norm
  (cf.~\cite{CR}), one convinces oneself that
  $$\Gamma(U)(\3\4{T_V}a^*a^*)^n\Omega_{P_1}=
  \big(\3(U\4{T_V}U^*)a^*a^*\big)^n\Omega_{P_1},\qquad n\in\NN,$$
  because $\Gamma(U)\Omega_{P_1}=\Omega_{P_1}$. It follows that
  \begin{equation*}
    \begin{split}
      \Gamma(U)\exp(\3\4{T_V}a^*a^*)\Omega_{P_1}&=\sum_{n=0}^\infty
        \frac{1}{n!}\Gamma(U)(\3\4{T_V}a^*a^*)^n\Omega_{P_1}\\
      &=\exp\big(\3(U\4{T_V}U^*)a^*a^*\big)\Omega_{P_1}.
    \end{split}
  \end{equation*}
  Since $U$ commutes with $P_1,\ P_2$ and $V$, it also commutes with
  all components of $V$ and $V^*$, including the operators
  ${V_{11}}^{-1},\ p_{\ker{V_{11}}^*}$ etc. $T_V$ is by definition
  \eqref{TV} a bounded function of these operators, so that
  \begin{equation}
    \label{GT}
    [U,T_V]=0,\qquad U\in G.
  \end{equation}
  Hence we get
  $\Gamma(U)\exp(\3\4{T_V}a^*a^*)\Omega_{P_1}=
  \exp(\3\4{T_V}a^*a^*)\Omega_{P_1}$ as claimed.
\end{proof}
Thus we arrive at the following formula
\begin{equation}
  \label{GPSIOM}
  \Gamma(U)\PV{\alpha}\Omega_{P_1}=D_V\psi_{P_1}(Ug_{\alpha_1}\dotsm 
  Ug_{\alpha_l})\psi_{P_1}(Ue_1\dotsm Ue_{L_V})\exp(\3\4{T_V}a^*a^*)
  \Omega_{P_1}
\end{equation}
which enables us to derive the ``charge'' carried by $\rho_V$. This is
best described by using the Fock space structure on $H(\rho_V)$
established in Corollary~\ref{cor:HRFOCK}.
\begin{Th}
  \label{th:CHARGE}\hspace*{\fill}\\
  Let $P_1$ be a \BP of \KK, let $G$ be a group consisting of diagonal
  Bogoliubov operators, and let $V\in\II_{P_1}(\KK)^G$. Further let
  $\mathfrak{h}_V\subset\KK_1$ be the $L_V$--dimensional subspace
  defined in \eqref{HV}, $L_V<\infty$, and let
  $\mathfrak{k}_V\subset\KK$ be the $M_V$--dimensional subspace
  defined in \eqref{KV}, $M_V=-\2\IND V$. Then both $\mathfrak{h}_V$
  and $\mathfrak{k}_V$ are invariant under $G$. Let
  $\Lambda_{\mathfrak{k}_V}$ be the unitary representation of $G$ on
  the antisymmetric Fock space $\mathcal{F}_a(\mathfrak{k}_V)$ over
  $\mathfrak{k}_V$ that is obtained by taking antisymmetric tensor
  powers of the representation on $\mathfrak{k}_V$. Then the unitary
  representation $\mathcal{U}_V$ of $G$ on the \HSP of isometries
  $H(\rho_V)$ which implements $\rho_V$ in the Fock representation
  $\pi_{P_1}$ is unitarily equivalent to $\Lambda_{\mathfrak{k}_V}$,
  tensored with the one--dimensional representation
  ${\det}_{\mathfrak{h}_V}(U)\DEF\det(U|_{\mathfrak{h}_V})$:
  \begin{equation}
    \label{REPG}
    \mathcal{U}_V\simeq{\det}_{\mathfrak{h}_V}\otimes
    \Lambda_{\mathfrak{k}_V}.
  \end{equation}
\end{Th}
\begin{proof}
  The finite dimensional subspace $\mathfrak{h}_V=V_{12}(\ker V_{22})$
  is invariant under $G$ because the elements of $G$ commute with the
  components of $V$. Since the $e_r$ form an \ONB in $\mathfrak{h}_V$,
  it follows from the CAR that 
  $$U(e_1)\dotsm U(e_{L_V})=\det\big(U|_{\mathfrak{h}_V}\big)\cdot
  e_1\dotsm e_{L_V},\qquad U\in G.$$
  
  Similarly, the elements of $G$ commute with the \BP $P_V$ (cf.\ 
  \eqref{GT} and \eqref{PV1}) and leave $\ker V^*$ invariant, so that
  $\mathfrak{k}_V=P_V(\ker V^*)$ is also left invariant. It then
  follows from the CAR that $g_{\alpha_1}\dotsm g_{\alpha_l}$
  transforms like the $l$--fold antisymmetric tensor product of
  $g_{\alpha_1},\dots,g_{\alpha_l}$ under $G$.
  
  Thus we see from \eqref{GPSIOM} that the representation of $G$ on
  $H(\rho_V)\Omega_{P_1}$ is unitarily equivalent to
  ${\det}_{\mathfrak{h}_V}\otimes\Lambda_{\mathfrak{k}_V}$. By
  Lemma~\ref{lem:REPHRV}, the same holds true for the representation
  on $H(\rho_V)$.
\end{proof}
Theorem~\ref{th:CHARGE} shows that genuine (i.e.\ non--surjective)
quasi--free endomorphisms $\rho_V$ are always \emph{reducible} in the
sense that the representation $\mathcal{U}_V$ of $G$ on $H(\rho_V)$
or, equivalently, the representation of the gauge invariant
``observable'' algebra $\CK^G$ induced by $\rho_V$ on the subspace of
$\Gamma(G)$--invariant vectors in $\mathcal{F}_a(\KK_1)$, is
reducible. In fact, each ``$n$--particle'' subspace of $H(\rho_V)$
(i.e.\ the closed linear span of all $\PV{\alpha}$ with $\alpha$ of
length $n$) is invariant under $G$.  Let $\mathcal{U}_V^{(n)}$ be the
restriction of $\mathcal{U}_V$ to this subspace. \label{page:RED}
Closest to irreducibility is the case that at least
$\mathcal{U}_V^{(1)}$ is irreducible. In typical situations, the
remaining representations $\mathcal{U}_V^{(n)}$ will then also be
irreducible. This happens for instance if $G\cong U(N)$ or $G\cong
SU(N)$, and $\mathfrak{k}_V$ carries the defining representation of
$G$. In the $U(N)$ case, the $\mathcal{U}_V^{(n)}$ are not only
irreducible, but also mutually inequivalent. In the $SU(N)$ case, the
representations $\mathcal{U}_V^{(0)},\dots,\mathcal{U}_V^{(N-1)}$ are
mutually inequivalent, but $\mathcal{U}_V^{(N)}$ is equivalent to
$\mathcal{U}_V^{(0)}$. In general, it can nevertheless happen that
$\mathcal{U}_V^{(1)}$ is irreducible but some $\mathcal{U}_V^{(n)}$
are not, as is the case if $G\cong SO(N)$ ($N>2$ even) and
$\mathfrak{k}_V$ carries the defining representation of $G$ (cf.\ 
\cite{Wey,Boe}). If $\mathcal{U}_V^{(1)}$ is reducible, then one has
an additional Clebsch--Gordan type splitting according to the unitary
isomorphism $\mathcal{F}_a(\mathfrak{k}_1\oplus\mathfrak{k}_2)
\cong\mathcal{F}_a(\mathfrak{k}_1)\otimes
\mathcal{F}_a(\mathfrak{k}_2)$ where $\mathfrak{k}_1$ and
$\mathfrak{k}_2$ are $G$--invariant subspaces of $\mathfrak{k}_V$.
\medskip

In the special case that $\rho_V$ is an automorphism ($\IND V=0$),
there survives only the factor ${\det}_{\mathfrak{h}_V}$ in
Theorem~\ref{th:CHARGE}. This is consistent with Matsui's result
\eqref{MATSUI1}, \eqref{MATSUI2} on the equivalence of \FSs over
$\CK^G$ ($G$ compact). For let $V\in\II_{P_1}^0(\KK)^G$ and set $P\DEF
VP_1V^*$.  Then the GNS representation $\pi_P$ for the \FS $\omega_P$
can be realized on $\mathcal{F}_a(\KK_1)$ as
$\pi_P=\pi_{P_1}\0{\rho_V}^{-1}$, with cyclic vector $\Omega_{P_1}$.
The unitary implementer \PV{} intertwines the representations
$\pi_P|_{\CK^G}$ and $\pi_{P_1}|_{\CK^G}$. \PV{} restricts to a
unitary isomorphism between the closed cyclic subspaces
$(\pi_P|_{\CK^G}\Omega_{P_1})^{\bar{\;}}$ and
$(\pi_{P_1}|_{\CK^G}\Omega_{P_1})^{\bar{\;}}$ \IFF
$[\PV{},\Gamma(G)]=0$, i.e.\ \IFF ${\det}_{\mathfrak{h}_V}(G)=1$. Now
one has $\mathfrak{h}_V=V_{12}(\ker V_{22})=\ker {V_{11}}^*=\ker
P_{11}= \KK_1\cap\4{P}(\KK)$, where we used \eqref{DIMKER} and
\eqref{P11}. Therefore ${\det}_{\mathfrak{h}_V}(G)=1$ is equivalent to
Matsui's condition \eqref{MATSUI2}. However, Matsui's result applies
more generally because a \BP $P\in\PP_{P_1}$ (see \eqref{PP1DEF})
which commutes with $G$ is not necessarily of the form $P=VP_1V^*$
with $V\in\II_{P_1}^0(\KK)^G$. In general, there exist
$V\in\II_{P_1}^0(\KK)$ such that $[VP_1V^*,G]=0$ but $[V,G]\neq0$. For
such $V$, one does not have $\Gamma(U)\PV{}\Gamma(U)^*\sim\PV{}$ but
only $\gamma(U)\PV{}\Omega_{P_1}\sim\PV{}\Omega_{P_1}$. If, for a
given \BP $P$, the representation $\mathcal{U}$ of $G$ on $\ker
P_{11}$ is equivalent to its complex conjugate representation
$\mathcal{U}^*$, then there exists a gauge invariant $V$ with
$VP_1V^*=P$. If $\mathcal{U}$ and $\mathcal{U}^*$ are disjoint but
both contained in $\KK_1$ with infinite multiplicity then $V$ with
$VP_1V^*=P$ can also be chosen to be gauge invariant. (These
statements follow from part (a) of Proposition~\ref{prop:HVKV}
together with Proposition~\ref{prop:PP1}.)
\medskip

\label{page:U1}
Another case of interest is obtained by specializing further to the
case of automorphisms and $G\cong\TT$. Assume that \TT acts on \KK via
$U_\lambda=e^{i\lambda}P+e^{-i\lambda}\4{P},\ \lambda\in\RR$, where
$P$ is a \BP of \KK which commutes with the given \BP $P_1$. Then all
gauge invariant Bogoliubov operators commute with $P$, and the
semigroup of (not necessarily implementable) gauge invariant
Bogoliubov operators is isomorphic to the semigroup of all isometries
of $P(\KK)$. Let 
\begin{equation}
  \label{p+-}
  p_+\DEF PP_1,\qquad p_-\DEF PP_2.
\end{equation}
(In the situation of the free Dirac field outlined in
Example~\ref{ex:DIRAC}, $p_\pm$ are just the spectral projections of
the Dirac Hamiltonian.)  If $V$ is gauge invariant, then the
implementability condition ($V_{12}$ Hilbert--Schmidt) is equivalent
to the condition that the components $V_{+-}$ and $V_{-+}$ be
Hilbert--Schmidt, because one has
\begin{equation}
  \label{V+-}
  \begin{aligned}
    V_{11}&=V_{++}+\4{V_{--}}, &\qquad V_{12}&=V_{+-}+\4{V_{-+}},\\
    V_{21}&=V_{-+}+\4{V_{+-}}, &\qquad V_{22}&=V_{--}+\4{V_{++}},
  \end{aligned}
\end{equation}
with $V_{\epsilon\epsilon'}\DEF p_\epsilon Vp_{\epsilon'},\ 
\epsilon,\epsilon'=\pm$.  These relations also entail that, for
$V\in\II_{P_1}(\KK)^\TT$,
\begin{equation}
  \label{HV1}
  \mathfrak{h}_V=V_{12}(\ker V_{22})=V_{+-}(\ker V_{--})\oplus
  \4{V_{-+}}(\ker \4{V_{++}}). 
\end{equation}
The subgroup $\II_{P_1}^0(\KK)^\TT$ of gauge invariant implementable
Bogoliubov operators with index zero can be identified with the {\em
  restricted unitary group} of $P(\KK)$, i.e.\ with the group of all
unitaries on $P(\KK)$ whose $(+-)$ and $(-+)$ components are
Hilbert--Schmidt, through the restriction map $V\mapsto PVP$. To
describe the charge corresponding to $V\in\II_{P_1}^0(\KK)^\TT$, we
have to compute the determinant of $U_\lambda|_{\mathfrak{h}_V},\
\lambda\in\RR$. By $\IND V=0$, $V_{+-}$ maps $\ker V_{--}$
isometrically onto $\ker{V_{++}}^*$, and $V_{-+}$ maps $\ker V_{++}$
isometrically onto $\ker{V_{--}}^*$ (cf.\ \eqref{KER},
\eqref{DIMKER}). Hence we get from \eqref{HV1}
\begin{align*}
  {\det}_{\mathfrak{h}_V}(U_\lambda)&=\exp\big(i\lambda(\dim\ker
    V_{--}-\dim\ker V_{++})\big)\\ 
  &=\exp(i\lambda\IND V_{--})\\
  &=\exp(-i\lambda\IND V_{++}).
\end{align*}
The charge corresponding to an element $V$ of the restricted unitary
group is therefore equal to $-\IND V_{++}$. \label{page:U1CHARGE} The
implementer \PV{} maps the charge--$q$ sector in
$\mathcal{F}_a(\KK_1)$, $q\in\ZZ$, onto the sector with charge $q-\IND
V_{++}$. Of course, $\IND V_{++}$ can only be nonzero if $p_+$ and
$p_-$ both have infinite rank. The fact that the charge created by a
gauge invariant quasi--free automorphism $\rho_V$ is in this way
related to the Fredholm index of $V_{++}$ was implicit in the
literature of the 1970s on the external field problem (see e.g.\ 
\cite{Lab,Lab75,F77,KS,R77a,R78,Sei78}) but has apparently first been
pointed out explicitly by Carey, Hurst and O'Brien \cite{CHOB}. These
authors showed that the connected components of the restricted unitary
group $\II_{P_1}^0(\KK)^\TT$ are precisely labelled by $\IND V_{++}$.
Computations of $\IND V_{++}$ for certain classes of unitary operators
$V$ can be found in more recent publications
\cite{CR,Mat87,R89,R89a,Mat90,BH}.
\medskip

Let us return to the general situation. So far we have analyzed the
representation of $G$ on $H(\rho_V)$ in terms of the given
representations on $\mathfrak{h}_V$ and $\mathfrak{k}_V$. But we can
also characterize the representations of $G$ which can possibly occur
on $\mathfrak{h}_V$ and $\mathfrak{k}_V$. Note that, if
$\KK_1'\subset\KK_1$ is a $G$--invariant subspace carrying a
representation of class $\xi$, then the complex conjugate space
$(\KK_1')^*\subset\KK_2$ carries a representation of the complex
conjugate class $\xi^*$.
\begin{Prop}
  \label{prop:HVKV}\hspace*{\fill}\\\indent
  a) Let $\mathfrak{h}\subset\KK_1$ be a finite dimensional
  $G$--invariant subspace carrying a representation of class $\xi$. If
  $\xi$ is self--conjugate \textup{(}i.e.\ $\xi=\xi^*$\textup{)}, then
  there exists $V\in\II_{P_1}^0(\KK)^G$ with
  $\mathfrak{h}_V=\mathfrak{h}$ such that $V-\1$ has finite rank. If
  $\xi$ and $\xi^*$ are disjoint, then there exists
  $V\in\II_{P_1}^0(\KK)^G$ with $\mathfrak{h}_V=\mathfrak{h}$ \IFF
  both $\xi$ and $\xi^*$ are contained in $\KK_1$ with infinite
  multiplicity.
  
  b) Let $\mathfrak{k}\subset\KK_1$ be a closed $G$--invariant
  subspace carrying a representation of class $\xi$. If $\xi$ is
  contained in $\KK_1$ with infinite multiplicity, then there exists a
  diagonal Bogoliubov operator $V\in\II_{P_1}(\KK)^G$ with
  $\mathfrak{k}_V=\mathfrak{k}$. If $\mathfrak{k}$ is irreducible and
  if $\xi$ is contained in $\KK_1$ with finite multiplicity, then
  there exists $V\in\II_{P_1}(\KK)^G$ with
  $\mathfrak{k}_V=\mathfrak{k}$ \IFF $\mathfrak{k}$ is finite
  dimensional and $\xi^*$ is contained in $\KK_1$ with infinite
  multiplicity. A class $\xi'$ which is only contained in $\KK_2$ but
  not in $\KK_1$ cannot be realized on any $\mathfrak{k}_V$,
  $V\in\II_{P_1}(\KK)^G$.
\end{Prop}
\begin{proof}
  The proof will be constructive.

  a) 1. Assume first that $\xi$ is self--conjugate. Then
  there exists a partial isometry $u$ from $\KK_2$ to $\KK_1$ with
  initial space $\mathfrak{h}^*$ and final space $\mathfrak{h}$ which
  commutes with $G$. Let $p_{\mathfrak{h}^\bot}$ be the orthogonal
  projection onto $\mathfrak{h}^\bot\subset\KK_1$. Then
  $$V\DEF\begin{pmatrix}
    p_{\mathfrak{h}^\bot} & u \\
    \4{u} & \4{p_{\mathfrak{h}^\bot}}
  \end{pmatrix}$$ 
  is a unitary Bogoliubov operator with $[V,G]=0$,
  $\mathfrak{h}_V=V_{12}(\ker V_{22})=u(\mathfrak{h}^*)=\mathfrak{h}$,
  and $V-\1$ has finite rank.

  2. Next assume that $\xi$ and $\xi^*$ are disjoint. If $\xi$ and
  $\xi^*$ are contained in $\KK_1$ with infinite multiplicity, one has
  a decomposition
  $$\KK_1\cong\ell^2(\mathfrak{h})\oplus\ell^2(\mathfrak{h}')\oplus\HH$$
  where $\mathfrak{h}$ is the given subspace carrying a representation
  of class $\xi$, $\mathfrak{h}'$ is a subspace carrying a
  representation of the complex conjugate class $\xi^*$,
  $\ell^2(\mathfrak{h})$ is the \HSP of square summable
  sequences of elements of $\mathfrak{h}$, and \HH is the
  orthogonal complement of
  $\ell^2(\mathfrak{h})\oplus\ell^2(\mathfrak{h}')$ in $\KK_1$.
  Then let $s$ resp.\ $s'$ be the ($G$--invariant) unilateral shift on
  $\ell^2(\mathfrak{h})$ resp.\ $\ell^2(\mathfrak{h}')$, given by
  $(f_1,f_2,\dotsc)\mapsto(0,f_1,f_2,\dotsc)$. Define a Fredholm
  operator $V_{11}$ on $\KK_1$ with index zero by
  $$V_{11}\DEF s\oplus{s'}^*\oplus p_\HH$$
  so that $\ker
  V_{11}=\mathfrak{h}',\ \ker{V_{11}}^*=\mathfrak{h}$. (Here
  $\mathfrak{h}$ is identified with $(\mathfrak{h},0,0,\dotsc)$, and
  similar for $\mathfrak{h}'$.)  Furthermore let $V_{12}$ be a
  $G$--invariant partial isometry from $\KK_2$ to $\KK_1$ with initial
  space ${\mathfrak{h}'}^*$ and final space $\mathfrak{h}$ (such
  $V_{12}$ exists because ${\mathfrak{h}'}^*$ and $\mathfrak{h}$ both
  belong to the class $\xi$). Then $V_{11}$ and $V_{12}$ define a
  Bogoliubov operator $V\in\II_{P_1}^0(\KK)^G$ with
  $\mathfrak{h}_V=\mathfrak{h}$. (One can show that $V$ fulfills in
  addition $VP_1V^*=P_{(0,\mathfrak{h})}$, in the notation of
  Proposition~\ref{prop:PP1}.) 
  
  Conversely, if there exists $V\in\II_{P_1}^0(\KK)^G$ with
  $\mathfrak{h}_V=\mathfrak{h}$, then $\xi^*$ must be contained in
  $\KK_1$ since $\ker{V_{11}}$ belongs to this class (recall that
  $V_{12}$ restricts to a unitary intertwiner between
  $\ker{V_{22}}=(\ker{V_{11}})^*$ and $\mathfrak{h}_V$). Since $\xi$
  and $\xi^*$ are disjoint, one has a decomposition
  \begin{equation}
    \label{K1DEC}
    \KK_1=\ker{V_{11}}^*\oplus\ker V_{11}\oplus\HH
  \end{equation}
  with $\HH\DEF\RAN V_{11}\cap\RAN{V_{11}}^*$. Viewing $V_{11}$ as a
  bounded bijection from $\RAN{V_{11}}^*$ onto $\RAN V_{11}$, one gets
  a unitary intertwiner from $\RAN{V_{11}}^*$ onto $\RAN V_{11}$ by
  polar decomposition of $V_{11}$. Therefore the representations of
  $G$ on $\RAN{V_{11}}^*$ and $\RAN V_{11}$ must be unitarily
  equivalent. One has by \eqref{K1DEC}
  $$\RAN{V_{11}}^*=\ker{V_{11}}^*\oplus\HH,\qquad \RAN
  V_{11}=\ker V_{11}\oplus\HH.$$
  Since $\ker{V_{11}}^*$ belongs to the class $\xi$ and $\ker V_{11}$
  to the (disjoint) class $\xi^*$, the representations on
  $\RAN{V_{11}}^*$ and $\RAN V_{11}$ can only be equivalent if $\xi$
  and $\xi^*$ are both contained in \HH with infinite multiplicity.

  b) 1. Assume first that $\xi$ appears in $\KK_1$ infinitely
  often. Then $\KK_1$ has the form
  $$\KK_1\cong\ell^2(\mathfrak{k})\oplus\HH$$
  where $\mathfrak{k}$ belongs to the class $\xi$. If $s$ is the shift
  on $\ell^2(\mathfrak{k})$ as above, then
  $$V\DEF
  \begin{pmatrix}
    s\oplus p_\HH & 0 \\ 0 & \4{s\oplus p_\HH}
  \end{pmatrix}$$
  is an element of $\II_{P_1}(\KK)^G$ with
  $\mathfrak{k}_V=\ker{V_{11}}^*=\mathfrak{k}$.
  
  2. Next assume that $\mathfrak{k}\subset\KK_1$ is irreducible of
  class $\xi$ and that $\xi$ is contained in $\KK_1$ with multiplicity
  $M<\infty$. Let first $\mathfrak{k}$ be finite dimensional and
  $\xi^*$ be contained in $\KK_1$ with multiplicity $\infty$. Then
  $\xi$ and $\xi^*$ are necessarily disjoint. Without loss of
  generality, we may restrict attention to the closed *--invariant
  subspace of \KK which comprises all representations of the classes
  $\xi$ and $\xi^*$, because a Bogoliubov operator having the desired
  properties on this subspace can be canonically extended to \KK by
  letting it act like the identity on the complement of that subspace
  (and because all other relevant operators leave this space
  invariant). Thus we will assume that \KK is of the form
  $$\KK=\ell^2(\mathfrak{k})\oplus\ell^2(\mathfrak{k})^*$$
  with $\ell^2(\mathfrak{k})\cap\KK_1=
  \SET{(f_1,\dots,f_M,0,0,\dotsc)}{f_j\in\mathfrak{k}}$. If $P$
  denotes the \BP onto $\ell^2(\mathfrak{k})$, then we are in a
  situation similar to the one discussed on p.~\pageref{page:U1}. $P$
  commutes with $P_1$, and if we define $p_\pm$ as in \eqref{p+-}:
  $p_+\DEF PP_1,\ p_-\DEF PP_2$, then $p_+$ has finite rank since
  $M<\infty$ by assumption. Any gauge invariant Bogoliubov operator
  $V$ has to commute with $P$ and is therefore of the form
  \eqref{V+-}. Such $V$ is automatically implementable because $p_+$
  has finite rank. Now let $s$ again be the unilateral shift on
  $\ell^2(\mathfrak{k})$, and set
  $$V\DEF s\oplus\4{s}.$$
  Then $V$ lies in $\II_{P_1}(\KK)^G$. We
  claim that $\mathfrak{k}_V=\mathfrak{k}$. Note that, by
  \eqref{PVQV}, $\mathfrak{k}_V\DEF P_V(\ker V^*)=P_{T_V}(\ker V^*)$
  where the \BP $P_{T_V}$ is defined in \eqref{PT} and \eqref{TV}. By
  definition of $V$, one finds that $V_{21}$ is a partial isometry so
  that $T_V=0$ by the remark on p.~\pageref{page:TV}. As a
  consequence, $P_{T_V}=P_1$ and $\mathfrak{k}_V=P_{T_V}(\ker V^*)=
  P_1(\mathfrak{k}\oplus\mathfrak{k}^*)=\mathfrak{k}$ as claimed.
  
  Conversely, assume that $V\in\II_{P_1}(\KK)^G$ exists with
  $\mathfrak{k}_V=\mathfrak{k}$. Since all $V^n(\mathfrak{k}_V),\ 
  n\in\NN$, are mutually orthogonal and belong to the class $\xi$,
  $\xi$ must be contained in \KK with infinite multiplicity. This
  means that $\xi^*$ must appear in $\KK_1$ with infinite
  multiplicity, because the multiplicity of $\xi$ in $\KK_1$ is finite
  by assumption. Then consider again the closed *--invariant subspace
  $\KK'$ of \KK which contains all the representations of class $\xi$
  and $\xi^*$. Since $\xi$ and $\xi^*$ are disjoint, we can write
  $$\KK'=\ell^2(\mathfrak{k})\oplus\ell^2(\mathfrak{k})^*
  \cong\big(\ell^2(\CC)\otimes\mathfrak{k}\big)\oplus
  \big(\ell^2(\CC)\otimes\mathfrak{k}\big)^*.$$
  Any operator $A$ on \KK
  which commutes with $G$ leaves $\KK'$ invariant, and its restriction
  to $\KK'$ has the form $(A'\otimes\1_\mathfrak{k})\oplus
  (\4{A''\otimes\1_\mathfrak{k}})$. In particular, the projections
  $p_\pm$ introduced above can be written as
  $p_\pm=p_\pm'\otimes\1_\mathfrak{k}$, and the component $V_{12}$ of
  a gauge invariant Bogoliubov operator restricts to an operator of
  the form $(V_{+-}'\otimes\1_\mathfrak{k})\oplus
  (\4{V_{-+}'\otimes\1_\mathfrak{k}})$ on $\KK'$ (cf.\
  \eqref{V+-}). Since $V$ fulfills the implementability condition,
  $V_{12}$ is Hilbert--Schmidt. If $\mathfrak{k}$ were infinite
  dimensional, this would entail that $V_{+-}'$ and $V_{-+}'$ had to
  vanish. But then the restrictions of $P_1$ and $V$ to $\KK'$ would
  commute, so that the powers of $V$ acting on $\mathfrak{k}_V$ would
  produce infinitely many mutually orthogonal copies of $\mathfrak{k}$
  in $\KK_1$. This contradicts the assumption, and therefore
  $\mathfrak{k}$ has to be finite dimensional.
  
  3. Finally, $\mathfrak{k}_V$ cannot be a subspace of $\KK_2$ because
  one has $\mathfrak{k}_V=P_{T_V}(\ker V^*)$ (see above), but
  $\KK_2\cap\RAN P_T=\{0\}$ for any operator $T\in\mathfrak{H}_{P_1}$
  (cf.\ the remark on p.~\pageref{page:KERP}).
\end{proof}
One could also ask which combinations of representations can occur
simultaneously on $\mathfrak{h}_V$ and $\mathfrak{k}_V$, for a single
$V$. Let us only remark here that, if a fixed $V\in\II_{P_1}(\KK)^G$
is multiplied with a unitary $V'\in\II_{P_1}^0(\KK)^G$, then the
corresponding representation changes according to
$\mathcal{U}_{V'V}\simeq\det_{\mathfrak{h}_{V'}}\cdot\mathcal{U}_V$ so
that the determinant factor in \eqref{REPG} can be modified almost
arbitrarily without affecting the factor $\Lambda_{\mathfrak{k}_V}$.

Also note that, in the generic field theoretic situation described in
Example~\ref{ex:DIRAC}, $\KK_1$ has the form
$\KK_1=\KK_+\oplus{\KK_-}^*$ where $\KK_+$ and $\KK_-$ both carry the
defining representation of $G$ with multiplicity $\infty$. Hence the
defining representation and its complex conjugate are both realized on
some $\mathfrak{h}_V$ and $\mathfrak{k}_V$.
\medskip

These remarks complete our general discussion of the charge structure
of quasi--free endomorphisms of the \CAR. We would like to devote the
rest of this section to a comparison between the semigroup of gauge
invariant quasi--free endomorphisms discussed above and the semigroup
of localized endomorphisms appearing in the theory of \SSSs. First of
all, \SSSs are by definition \emph{irreducible} so that the question
arises how to obtain the irreducible {\em``subobjects''} of gauge
invariant quasi--free endomorphisms $\rho_V$. Suppose that
$\{\Psi_1,\dots,\Psi_m\}$ is an orthonormal set in $H(\rho_V)$ which
transforms irreducibly under $G$. According to the general theory,
there should exist a gauge invariant isometry $\Phi$ on Fock space
with $\RAN\Phi=\oplus_{j=1}^m\RAN\Psi_j$ (cf.\ \eqref{BORCHERS}).  The
corresponding irreducible endomorphism \rho (which is not quasi--free)
would then be given by
$\rho(a)\DEF\sum_{j=1}^m\Phi^*\Psi_ja\Psi_j^*\Phi$. However, the
construction of such gauge invariant isometries $\Phi$ is at present
unclear.

Similarly, the \emph{direct sum} of quasi--free endomorphisms (which
should be defined as in \eqref{DIRSUM}) can in general not be
quasi--free. This is evident from the index formula \eqref{DV} and
from the additivity of the statistics dimension on direct sums.  The
statistics dimension and the Bosonized statistics operator from
\eqref{BSO} and \eqref{BSO1} are the only invariants related to the
\emph{statistics} of a \SSS which can be unambiguously ascribed to a
quasi--free endomorphism in this general setting.

Thus $\II_{P_1}(\KK)^G$ is only closed under composition and not under
taking direct sums or subobjects. Furthermore, the existence of {\em
  conjugates} is only guaranteed if one makes additional assumptions
on the action of $G$ on \KK. Specifically, one needs charge
conjugation already on the level of first quantization:
\begin{Prop}
  \label{prop:CONJ}\hspace*{\fill}\\
  Assume that there exists a further \BP $P$ of \KK which commutes
  with $P_1$ and with $G$, and let $p_\pm$ be defined as in
  \eqref{p+-}:
  $$p_+\DEF PP_1,\qquad p_-\DEF PP_2.$$
  Assume further that there
  exists a unitary operator $\mathcal{C}_{+-}$ from $\KK_-\DEF
  p_-(\KK)$ onto $\KK_+\DEF p_+(\KK)$ which commutes with $G$, and let
  $\mathcal{C}$ be the unique Bogoliubov operator which commutes with
  $P$ and which is given on $P(\KK)=\KK_+\oplus\KK_-$ by the matrix
  $$P\mathcal{C}P\DEF
  \begin{pmatrix}
    0 & \mathcal{C}_{+-} \\ {\mathcal{C}_{+-}}^* & 0
  \end{pmatrix}.$$
  Then $\mathcal{C}$ is gauge invariant, unitary and self--adjoint,
  and the map
  \begin{equation}
    \label{CHCONJ}
    V\mapsto V^c\DEF \mathcal{C}V\mathcal{C}
  \end{equation}
  is an involutive automorphism of $\II_{P_1}(\KK)^G$ which preserves
  the statistics dimension. In addition, one has
  $$\mathfrak{h}_{V^c}=\mathcal{C}({\mathfrak{h}_V}^*),\qquad
    \mathfrak{k}_{V^c}=\mathcal{C}({\mathfrak{k}_V}^*)$$
  so that the
  representation $\mathcal{U}_{V^c}$ associated with $V^c$ according
  to \eqref{REPG} is unitarily equivalent to the complex conjugate of
  $\mathcal{U}_V$ \emph{(``charge conjugation'')}.
\end{Prop}
\begin{proof}
  $\mathcal{C}=P\mathcal{C}P+\4{P\mathcal{C}P}$ is clearly a gauge
  invariant self--adjoint element of $\II^0(\KK)$. Its components \WRT
  $P_1,P_2$ are
  $$\mathcal{C}_{11}=0,\qquad \mathcal{C}_{12}=
  \begin{pmatrix}
    0 & \mathcal{C}_{+-} \\ {\mathcal{C}_{+-}}^\tau & 0
  \end{pmatrix}.$$
  Thus one has ${V^c}_{12}=\mathcal{C}_{12}V_{21}\mathcal{C}_{12}$ so
  that $V^c$ belongs to $\II_{P_1}(\KK)^G$ if $V$ does. $(V^c)^c=V$
  follows from $\mathcal{C}^2=\1$. The map \eqref{CHCONJ} is
  consequently an involutive automorphism of $\II_{P_1}(\KK)^G$ (it is
  unital and multiplicative). It preserves the statistics dimension
  because $\IND \mathcal{C}=0$.
  
  Now let $V\in\II_{P_1}(\KK)^G$. One has
  ${V^c}_{22}=\mathcal{C}_{21}V_{11}\mathcal{C}_{12}$ so that
  \begin{equation*}
    \begin{split}
      \mathfrak{h}_{V^c} &\DEF{V^c}_{12}(\ker{V^c}_{22})\\
      &=\mathcal{C}_{12}V_{21}\mathcal{C}_{12}(\ker\mathcal{C}_{21}V_{11}
        \mathcal{C}_{12})\\
      &=\mathcal{C}_{12}V_{21}(\ker V_{11})\\
      &=\mathcal{C}_{12}({\mathfrak{h}_V}^*)\\
      &=\mathcal{C}({\mathfrak{h}_V}^*).
    \end{split}
  \end{equation*}
  Similarly, one finds that the antisymmetric Hilbert--Schmidt
  operators $T_V,\ T_{V^c}$ associated with $V,\ V^c$ through
  \eqref{TV} are related to each other by
  $$T_{V^c}=\mathcal{C}_{21}\4{T_V}\mathcal{C}_{21}$$
  since
  \begin{equation*}
    \begin{split}
      T_{V^c} &\DEF {V^c}_{21}{{V^c}_{11}}^{-1}-{{V^c}_{22}}^{-1*}
        {{V^c}_{12}}^*p_{{\ker {V^c}_{11}}^*}\\
      &=\mathcal{C}_{21}V_{12}\mathcal{C}_{21}\cdot\mathcal{C}_{12}
        {V_{22}}^{-1}\mathcal{C}_{21}-\mathcal{C}_{21}{V_{11}}^{-1*}
        \mathcal{C}_{12}\cdot \mathcal{C}_{21}{V_{21}}^*
        \mathcal{C}_{21}\cdot\mathcal{C}_{12}p_{{\ker V_{22}}^*}
        \mathcal{C}_{21}\\
      &=\mathcal{C}_{21}(V_{12}{V_{22}}^{-1}-{V_{11}}^{-1*}{V_{21}}^*
        p_{{\ker V_{22}}^*})\mathcal{C}_{21}\\
      &=\mathcal{C}_{21}\4{T_V}\mathcal{C}_{21}.
    \end{split}
  \end{equation*}
  One obtains for the corresponding \BPs (cf.~\eqref{PT})
  \begin{equation*}
    \begin{split}
      P_{T_{V^c}} &\DEF(P_1+T_{V^c})(P_1+{T_{V^c}}^*T_{V^c})^{-1}
        (P_1+{T_{V^c}}^*)\\
      &=(P_1+\mathcal{C}_{21}\4{T_V}\mathcal{C}_{21})
        (P_1+\mathcal{C}_{12}T_V{T_V}^*\mathcal{C}_{21})^{-1}
        (P_1+\mathcal{C}_{12}{T_V}^\tau \mathcal{C}_{12})\\
      &=\mathcal{C}(P_2+\4{T_V})\mathcal{C}_{21}\cdot
        \mathcal{C}_{12}(P_2+T_V{T_V}^*)^{-1}\mathcal{C}_{21}\cdot
        \mathcal{C}_{12}(P_2+{T_V}^\tau)\mathcal{C}\\
      &=\mathcal{C}\4{P_{T_V}}\mathcal{C}.
    \end{split}
  \end{equation*}
  This entails further that 
  $$\mathfrak{k}_{V^c}\DEF P_{T_{V^c}}(\ker V^{c*})
    =\mathcal{C}\4{P_{T_V}}\mathcal{C}(\mathcal{C}\ker V^*)
    =\mathcal{C}({\mathfrak{k}_V}^*)$$
  and finally
  $$\mathcal{U}_{V^c}\simeq{\det}_{\mathfrak{h}_{V^c}}\otimes
    \Lambda_{\mathfrak{k}_{V^c}}\simeq({\det}_{\mathfrak{h}_V})^*
    \otimes(\Lambda_{\mathfrak{k}_V})^*\simeq{\mathcal{U}_V}^*.$$
\end{proof}
Note that the assumptions of the proposition amount to a decomposition
of the single--particle space $\KK_1$ into the direct sum of two
antiunitarily equivalent $G$--modules $\KK_+$ and ${\KK_-}^*$. These
assumptions are of course satisfied
in Example~\ref{ex:DIRAC}. It is also fairly obvious that, if these
assumptions are violated, there will in general exist operators $V$ in
$\II_{P_1}(\KK)^G$ which do \emph{not} admit conjugates in
$\II_{P_1}(\KK)^G$ (i.e.\ there is no $V'\in\II_{P_1}(\KK)^G$ with
$\mathcal{U}_{V'}\simeq {\mathcal{U}_V}^*$;
cf.~Proposition~\ref{prop:HVKV}). But also note that, if
$V\in\II_{P_1}(\KK)^G$ has finite index and if $G$ acts on both
subspaces $\mathfrak{h}_V$ and $\mathfrak{k}_V$ with determinant~$1$,
then $\rho_V$ is automatically \emph{self--conjugate} in the sense
that $\mathcal{U}_V\simeq{\mathcal{U}_V}^*$. Recall from
Theorem~\ref{th:CHARGE} that ${\mathcal{U}_V}$ is then equivalent to
the representation $\Lambda_{\mathfrak{k}_V}$ on the antisymmetric
Fock space $\mathcal{F}_a(\mathfrak{k}_V)$ over $\mathfrak{k}_V$. Let
$\Lambda_{\mathfrak{k}_V}^{(n)}$ be the restriction of
$\Lambda_{\mathfrak{k}_V}$ to the $n$--particle subspace of
$\mathcal{F}_a(\mathfrak{k}_V)$ so that
$$\Lambda_{\mathfrak{k}_V}=
\bigoplus_{n=0}^{M_V}\Lambda_{\mathfrak{k}_V}^{(n)}.$$
The character
$\chi^{(n)}(U),\ U\in G$, of $\Lambda_{\mathfrak{k}_V}^{(n)}$ is equal
to the $n$th elementary symmetric function
$\sum_{j_1<\dots<j_n}\lambda_{j_1} \cdots\lambda_{j_n}$ of the
eigenvalues $\lambda_1,\dots, \lambda_{M_V}$ of $U|_{\mathfrak{k}_V}$
(cf.\ \cite{Mur}). It follows that
$$\chi^{(n)}=\chi^{(M_V)}\4{\chi^{(M_V-n)}}$$
and that
$$\Lambda_{\mathfrak{k}_V}^{(n)}\simeq{\det}_{\mathfrak{k}_V}\otimes
\Lambda_{\mathfrak{k}_V}^{(M_V-n)*},\qquad n=0,\dots,M_V$$
(this holds for any unitary representation on $\mathfrak{k}_V$, the
condition on the determinant is not needed here). Assuming now that
$\det_{\mathfrak{h}_V}=\det_{\mathfrak{k}_V}=1$, we get that
$$\mathcal{U}_V^{(n)}\simeq\mathcal{U}_V^{(M_V-n)*}$$
so that $\mathcal{U}_V$ is self--conjugate as claimed. (This remark
applies in particular to the Dirac field with \GG $SU(N)$, see
Example~\ref{ex:DIRAC}.)
\medskip

The main reason for our inability to mimic the generic superselection
structure more closely is of course the complete lack of {\em
  locality} in our preceding considerations. If one could find, in a
specific model, a localized implementable quasi--free endomorphism,
then it is clear that our methods would apply and could be used to
construct the corresponding charged local fields and to determine
their charge quantum numbers. It is however \emph{not} clear from the
outset that localization and implementability are compatible with each
other. Known results in this direction only deal with the case of
automorphisms. Building on the work of Carey and Ruijsenaars \cite{CR}
and others, we constructed in \cite{CB0} a family of (implementable
and transportable) localized automorphisms, carrying arbitrary
\TT--charges, for the free Dirac field in two spacetime dimensions
with arbitrary mass. The Bogoliubov operators $V$ belonging to these
automorphisms are given by two \TT--valued functions which are equal
to $1$ at spacelike infinity, and the charge $-\IND V_{++}$ (cf.\ 
p.~\pageref{page:U1CHARGE}) created by $\rho_V$ is equal to the
difference of the winding numbers of these functions. However, in
contrast to the two--dimensional case, there are no known examples of
implementable charge--carrying (\WRT $G=\TT$) Bogoliubov automorphisms
in four spacetime dimensions. (Non--zero charge seems to be compatible
with $V_{12}$ compact, but not with $V_{12}$ Hilbert--Schmidt, cf.\ 
\cite{Mat87,R89,Mat90,BH}.)

As a slight generalization of a result in \cite{CB0}, we can
characterize localized gauge invariant endomorphisms of the free Dirac
field as follows. Recall from Example~\ref{ex:DIRAC} that all gauge
invariant Bogoliubov operators for the $N$--component Dirac field with
$U(N)$ gauge symmetry have the form \eqref{VFORM}
\begin{equation}
  \label{VFORM1}
  V=(v\otimes\1_N)\oplus(\4{v}\otimes\4{\1_N})
\end{equation}
where $\KK=\HH'\oplus{\HH'}^*,\ \HH'=\HH\otimes\CC^N,\ 
\HH=L^2(\RR^{2n-1},\CC^{2^n})$, and $v$ is an isometry of \HH. We
restrict to the time zero situation.
\begin{Prop}
  \label{prop:LOC}\hspace*{\fill}\\
  Let $O$ be a double cone with base $B\subset\RR^{2n-1}$ at time
  zero.  Let $V\in\II_{P_1}(\KK)^{U(N)}$ and let $v$ be the isometry
  of \HH associated with $V$ via \eqref{VFORM1}. Then $\rho_V$ is
  localized in $O$ in the sense of \eqref{LOCEND}\footnote{More
    precisely, the normal extension of $\rho_V$ in the representation
    $\pi_{P_1}$ is localized in $O$.} \IFF there exists, for each
  connected component $\Delta$ of $\RR^{2n-1}\setminus B$, a phase
  factor $\tau_\Delta\in\TT$ such that
  \begin{equation}
    \label{TAUDELTA}
    v(f)=\tau_\Delta f\qquad\text{for all }f\in\HH\text{ with }
    \SUPP f\subset\Delta.
  \end{equation}
\end{Prop}
\begin{proof}
  Assume that $\rho_V$ is localized in $O$. Let $b_1,\dots,b_N$ be the
  standard basis in $\CC^N$, let $\Delta$ be a component of the
  complement of $B$, and let $f,g\in\HH$ with $\SUPP
  f,g\subset\Delta$. Then
  $$a(f,g)\DEF \sum_{j=1}^N(f\otimes b_j)(g\otimes b_j)^*\in
  \CK^{U(N)}$$
  is gauge invariant, and $\pi_{P_1}(a(f,g))$ is a local
  observable in $\AA(O')$. Since $\rho_V$ is localized in $O$, one has
  $a(f,g)=\rho_V(a(f,g))=\sum_j(vf\otimes b_j)(vg\otimes b_j)^*$.
  Since the $b_j$ are linearly independent, it follows that
  \begin{equation}
    \label{FG}
    (f\otimes b_j)(g\otimes b_j)^*=(vf\otimes b_j)(vg\otimes b_j)^*,
    \qquad j=1,\dots,N.
  \end{equation}
  Now let $P'$ be the \BP onto $\HH'\subset\KK$, and let $\omega_{P'}$
  be the corresponding Fock state. One has 
  $$\omega_{P'}\Big((f\otimes b_j)^*(f\otimes b_j)(f\otimes b_j)^*
  (f\otimes b_j)\Big)=\NORM{f}_\HH^4,$$
  and, since $(vf\otimes b_j)^*$
  belongs to the Gelfand (or annihilator) ideal of $\omega_{P'}$,
  $$\omega_{P'}\Big((f\otimes b_j)^*(vf\otimes b_j)(vf\otimes b_j)^*
  (f\otimes b_j)\Big)=\langle vf,f\rangle\omega_{P'}\big((f\otimes
  b_j)^*(vf\otimes b_j)\big)=\ABS{\langle vf,f\rangle}^2.$$
  Therefore
  one gets from \eqref{FG}, in the special case $f=g$, that
  $\NORM{f}^2_\HH=\ABS{\langle vf,f\rangle}$. It follows that there
  exists $\tau_f\in\TT$ such that $v(f)=\tau_f f$. By the same
  argument, $v(g)=\tau_g g$ for some $\tau_g\in\TT$. Then \eqref{FG}
  yields that $\tau_f=\tau_g$. Therefore these phase factors depend
  only on $\Delta$ and not on the functions.
  
  Conversely, assume that \eqref{TAUDELTA} holds. Let $\tilde O\subset
  O'$ be a symmetric double cone with base $\tilde B$ at time zero,
  and let $\Delta$ be the connected component of $\RR^{2n-1}\setminus
  B$ which contains $\tilde B$. Then $\rho_V$ acts on the local \FA
  belonging to $\tilde O$ like the gauge transformation induced by
  $\tau_\Delta\in\TT\subset U(N)$, and it consequently acts like the
  identity on $\AA(\tilde O)$. Since the algebra $\AA(O')$ is
  generated by such local algebras $\AA(\tilde O)$, it follows that
  $\rho_V$ is localized in $O$.
\end{proof}
Of course, $\RR^{2n-1}\setminus B$ is connected if $n>1$, but has two
connected components if $n=1$. Recall that this is the basic reason
for the generic violation of Haag duality and for the possible
occurrence of braid group statistics and soliton sectors in
two--dimensional Minkowski space.
\medskip

We would like to close this section with a demonstration that at least
the free massless Dirac field in two spacetime dimensions possesses
non--surjective implementable localized quasi--free endomorphisms. It
suffices to consider one chiral component of the field. Thus consider
the \HSP $\HH=L^2(\RR)$ with Dirac Hamiltonian $-i\tfrac{d}{dx}$. It
is convenient to transform to the \HSP $L^2(\TT)$ via the Cayley
transform $\varkappa$
$$\varkappa: \RR\cup\{\infty\}\to\TT,\qquad x\mapsto -e^{2i\arctan
  x}=\frac{x-i}{x+i}$$
(cf.\ \cite{CR,R89a}). $\varkappa$ induces a
unitary transformation $\tilde\varkappa$
$$\tilde\varkappa:L^2(\TT)\to L^2(\RR),\qquad(\tilde\varkappa f)(x)=
\pi^{-\2}\frac{f(\varkappa(x))}{x+i}.$$
The important point is that
the spectral projections $p_\pm$ of $-i\tfrac{d}{dx}$ are transformed
into the Hardy space projections (cf.\ \cite{D}): Set
$q_\pm\DEF\tilde\varkappa^{-1}p_\pm\tilde\varkappa$, then
\begin{equation}
  \label{q+-}
  q_+=\sum_{n\geq0}e_n\langle e_n,.\,\rangle,\quad 
  q_-=\sum_{n<0}e_n\langle e_n,.\,\rangle,\quad e_n(z)\DEF z^n\
  (z\in\TT,n\in\ZZ). 
\end{equation}
Our task is to construct an isometry $v$ of $L^2(\TT)$ with
$q_+vq_-,q_-vq_+$ Hilbert--Schmidt (implementability), with $\IND
v=-1$ (close to irreducibility, cf.\ p.~\pageref{page:RED}), and such
that $vf=f$ for all $f\in L^2(\TT)$ with $\SUPP f\subset\TT\setminus
I$ where $I\subset\TT$ is a fixed interval (localization). As
localization region we shall choose the interval
$$I\DEF\SET{e^{i\lambda}}{\tfrac{\pi}{2}\leq\lambda\leq\tfrac{3\pi}{2}}$$
which corresponds, by the inverse Cayley transform, to the interval
$\varkappa^{-1}(I)=[-1,1]$ in \RR. We need the following \ONB
$(f_m)_{m\in\ZZ}$ in $L^2(I)\subset L^2(\TT)$
$$f_m(z)\DEF\sqrt{2}(-1)^mz^{2m}\chi_I(z),\quad z\in\TT$$
where $\chi_I$ is the characteristic function of $I$. The isometry $v$
is now defined by
\begin{equation}
  \label{v}
  v\DEF\1+\sum_{m\geq0}(f_{m+1}-f_m)\langle f_m,.\,\rangle.
\end{equation}
Note that $v$ acts like the identity on functions with support in
$\TT\setminus I$, that $vf_m=f_m$ if $m<0$, and that $v$ acts like the
unilateral shift on the remaining $f_m$: $vf_m=f_{m+1}$ if $m\geq0$.

The author gratefully acknowledges some private lessons in estimating
infinite series given to him by P.~Grinevich which were indispensable
for the proof of the following lemma.
\begin{Lem}
  \label{lem:HS}\hspace*{\fill}\\
  $q_+vq_-$ and $q_-vq_+$ are Hilbert--Schmidt.
\end{Lem}
\begin{proof}
  \allowdisplaybreaks[4]
  The following inner products are easily computed
  $$\langle e_l,f_m\rangle=\frac{(-1)^m}{\pi\sqrt{2}}
  \int_{\frac{\pi}{2}}^\frac{3\pi}{2}e^{i(2m-l)\lambda}d\lambda=
  \begin{cases}
    \frac{(-1)^m}{\sqrt{2}}\delta_{2m,l}, &l\text{ even}\\
    \frac{\sqrt{2}}{\pi}\frac{(-1)^{\frac{l-1}{2}}}{2m-l}, &l\text{ odd.}
  \end{cases}$$
  1. This yields for the Hilbert--Schmidt norm of $q_+vq_-$, using
  \eqref{q+-} and \eqref{v}
  \begin{align*}
    \HSNORM{q_+vq_-}^2 &=\sum_{n<0}\bigNORM{q_+ve_n}^2\\
    &=\sum_{n<0}\biggNORM{\sum_{l,m\geq0}e_l\langle
      e_l,f_{m+1}-f_m\rangle\langle f_m,e_n\rangle}^2\\
    &=\sum_{n<0}\sum_{l\geq0}\biggABS{\sum_{m\geq0}\langle
      e_l,f_{m+1}-f_m\rangle\langle f_m,e_n\rangle}^2\\
    &=\frac{4}{\pi^4}\sum_{\substack{n<0,\\n\,\text{odd}}}
      \sum_{\substack{l>0,\\l\,\text{odd}}}\left(\sum_{m\geq0}
        \left(\frac{1}{2m+2-l}-\frac{1}{2m-l}\right)
        \frac{1}{2m-n}\right)^2\\
    &\quad+\frac{1}{\pi^2}\sum_{\substack{n<0,\\n\,\text{odd}}}
      \sum_{\substack{l\geq0,\\l\,\text{even}}}\left(\sum_{m\geq0}
        \Bigl((-1)^{m+1}\delta_{2m+2,l}-(-1)^m\delta_{2m,l}\Bigr)
        \frac{1}{2m-n}\right)^2\\
    &=\frac{16}{\pi^4}\sum_{\substack{n>0,\\n\,\text{odd}}}
      \sum_{\substack{l>0,\\l\,\text{odd}}}\left(\sum_{m\geq0}
        \frac{1}{(2m+2-l)(2m-l)(2m+n)}\right)^2\\
    &\quad+\frac{1}{\pi^2}\sum_{\substack{n>0,\\n\,\text{odd}}}
      \sum_{\substack{l\geq0,\\l\,\text{even}}}\left(\frac{1}{l-2+n}-
        \frac{1}{l+n}\right)^2\\
    &=\frac{16}{\pi^4}\sum_{n,l\geq0}\left(\sum_{m\geq0}
      \frac{1}{(2m-2l+1)(2m-2l-1)(2m+2n+1)}\right)^2\\
    &\quad+\frac{4}{\pi^2}\sum_{n,l\geq0}\left(
      \frac{1}{(2l+2n-1)(2l+2n+1)}\right)^2\\
    &=\frac{16}{\pi^4}\sum_{n,l\geq0}\left(\sum_{m\geq0}\frac{1}{
        \bigl(4(m-l)^2-1\bigr)(2m+2n+1)}\right)^2\\
    &\quad+\frac{4}{\pi^2}\sum_{n,l\geq0}\frac{1}{\bigl(4(l+n)^2
      -1\bigr)^2}\\
    &=\frac{16}{\pi^4}\sum_{n,l\geq0}\left(\sum_{m\geq-l}
      \frac{1}{(4m^2-1)(2m+a_{ln})}\right)^2
      +\frac{4}{\pi^2}\sum_{n\geq0}\frac{1}{(4n^2-1)^2}\\
    &\quad+\frac{4}{\pi^2}\sum_{n\geq0}\sum_{l\geq1}
      \frac{1}{\bigl(4(l+n)^2-1\bigr)^2},
      \qquad\text{with }a_{ln}\DEF 2l+2n+1\\
    &<\frac{16}{\pi^4}\sum_{n,l\geq0}\left(\sum_{m\geq-l}
      \frac{1}{4m^2-1}\Biggl(\frac{1}{a_{ln}}-\left(\frac{1}{a_{ln}}-
        \frac{1}{2m+a_{ln}}\right)\Biggr)\right)^2\\
    &\quad+\frac{4}{\pi^2}\Biggl(1+\underbrace{\sum_{n\geq1}
      \frac{1}{4n^2-1}}_{=\2}+\underbrace{\sum_{n\geq0}
      \frac{1}{4(n+1)^2-1}}_{=\2}\underbrace{\sum_{l\geq1}
      \frac{1}{4l^2-1}}_{=\2}\Biggr)\\
    &=\frac{7}{\pi^2}+\frac{16}{\pi^4}\sum_{n,l\geq0}\Biggl(
      \frac{1}{a_{ln}}\biggl(\underbrace{\sum_{m\in\ZZ}
          \frac{1}{4m^2-1}}_{=0}-\sum_{m<-l}\frac{1}{4m^2-1}\biggr)\\
    &\phantom{=\frac{7}{\pi^2}+\frac{16}{\pi^4}\sum_{n,l\geq0}\Biggl(}
      -\frac{1}{a_{ln}}\sum_{m\geq-l}\frac{2m}{(4m^2-1)(2m+a_{ln})}
      \Biggr)^2\\
    &=\frac{7}{\pi^2}+\frac{16}{\pi^4}\sum_{n,l\geq0}\frac{1}{{a_{ln}}^2}
      \Biggl(\sum_{m>l}\frac{1}{4m^2-1}+\sum_{m>l}\frac{2m}{(4m^2-1)
          (2m+a_{ln})}\\
    &\phantom{=\frac{7}{\pi^2}+\frac{16}{\pi^4}\sum_{n,l\geq0}
      \frac{1}{{a_{ln}}^2}\Biggl(}+\sum_{m=-l}^l\frac{2m}{(4m^2-1)
        (2m+a_{ln})}\Biggr)^2\\
    &=\frac{7}{\pi^2}+\frac{16}{\pi^4}\sum_{n,l\geq0}\frac{1}{{a_{ln}}^2}
      \Biggl(\sum_{m>l}\frac{1}{4m^2-1}\biggl(\underbrace{1+
        \frac{2m}{2m+a_{ln}}}_{<2}\biggr)\\
    &\phantom{=\frac{7}{\pi^2}+\frac{16}{\pi^4}\sum_{n,l\geq0}
      \frac{1}{{a_{ln}}^2}\Biggl(}+\sum_{m=1}^l
      \underbrace{\frac{2m}{4m^2-1}\biggl(\frac{1}{2m+a_{ln}}+
        \frac{1}{2m-a_{ln}}\biggr)}_{=\underbrace{\frac{8m^2}{4m^2-1}}_{
          \leq\frac{8}{3}}\underbrace{\frac{1}{4m^2-{a_{ln}}^2}}_{<0}}
      \Biggr)^2\\
    &<\frac{7}{\pi^2}+\frac{64}{\pi^4}\sum_{n,l\geq0}\frac{1}{{a_{ln}}^2}
      \Biggl(\sum_{m>l}\frac{1}{4m^2-1}+\frac{4}{3}\sum_{m=1}^l
      \frac{1}{{a_{ln}}^2-4m^2}\Biggr)^2\\
    &<\frac{7}{\pi^2}+\frac{64}{\pi^4}\underbrace{\sum_{n\geq0}
      \frac{1}{(2n+1)^2}}_{=\frac{\pi^2}{8}}
      \Biggl(\underbrace{\sum_{m>0}\frac{1}{4m^2-1}}_{=\2}\Biggr)^2\\
    &\quad+\frac{64}{\pi^4}\sum_{n\geq0}\sum_{l\geq1}\frac{1}{{a_{ln}}^2}
      \biggl(\int_l^\infty\frac{dm}{(2m-1)^2}+\frac{4}{3}\int_1^l
      \frac{dm}{(2m-a_{ln})^2}\biggr)^2\\
    &=\frac{9}{\pi^2}+\frac{64}{\pi^4}\sum_{n\geq0}
      \sum_{l\geq1}\frac{1}{{a_{ln}}^2}\biggl(\frac{1}{2(2l-1)}+
      \frac{4}{3}\frac{l-1}{(a_{ln}-2l)(a_{ln}-2)}\biggr)^2\\
    &=\frac{9}{\pi^2}+\frac{16}{\pi^4}\sum_{n,l\geq0}
      \biggl(\underbrace{\frac{1}{(a_{ln}+2)(2l+1)}}_{<\frac{1}{(2n+1)
          (2l+1)}}+\frac{8}{3}\underbrace{\frac{l}{(a_{ln}+2)
          (2n+1)a_{ln}}}_{<\frac{1}{(2n+1)(2l+1)}}\biggr)^2\\
    &<\frac{9}{\pi^2}+\frac{16}{\pi^4}\biggl(
      \frac{11}{3}\biggr)^2\Biggl(\underbrace{\sum_{n\geq0}
        \frac{1}{(2n+1)^2}}_{=\frac{\pi^2}{8}}\Biggr)^2\\
    &=\frac{9}{\pi^2}+\biggl(\frac{11}{6}\biggr)^2<\infty.
  \end{align*}
  This proves that $q_+vq_-$ is Hilbert--Schmidt. The Hilbert--Schmidt
  norm of $q_-vq_+$ can be estimated as follows:
  \begin{align*}
    \HSNORM{q_-vq_+}^2 &=\sum_{n\geq0}\bigNORM{q_-ve_n}^2\\
    &=\sum_{n\geq0}\sum_{l<0}\biggABS{\sum_{m\geq0}\langle
      e_l,f_{m+1}-f_m\rangle\langle f_m,e_n\rangle}^2\\
    &=\frac{1}{\pi^2}\sum_{\substack{n\geq0,\\n\,\text{even}}}
      \sum_{\substack{l<0,\\l\,\text{odd}}}\left(\sum_{m\geq0}
        \left(\frac{1}{2m+2-l}-\frac{1}{2m-l}\right)
        (-1)^m\delta_{2m,n}\right)^2\\
    &\quad+\frac{4}{\pi^4}\sum_{\substack{n>0,\\n\,\text{odd}}}
      \sum_{\substack{l<0,\\l\,\text{odd}}}\left(\sum_{m\geq0}
        \left(\frac{1}{2m+2-l}-\frac{1}{2m-l}\right)
        \frac{1}{2m-n}\right)^2\\
    &=\frac{1}{\pi^2}\sum_{n,l\geq0}\biggl(\frac{1}{a_{ln}+2}-
      \frac{1}{a_{ln}}\biggr)^2\qquad
      \text{(with $a_{ln}$ as above)}\\
    &\quad+\frac{4}{\pi^4}\sum_{n,l\geq0}\Biggl(\sum_{m\geq0}\biggl(
      \frac{1}{2m+2l+3}-\frac{1}{2m+2l+1}\biggr)
      \frac{1}{2m-2n-1}\Biggr)^2\\
    &=\frac{4}{\pi^2}\sum_{n,l\geq0}\biggl(\frac{1}{a_{ln}(a_{ln}+2)}
      \biggr)^2\\
    &\quad+\frac{16}{\pi^4}\sum_{n,l\geq0}\Biggl(\sum_{m\geq0}
      \frac{1}{(2m+2l+1)(2m+2l+3)(2m-2n-1)}\Biggr)^2\\
    &<\frac{4}{\pi^2}\sum_{n,l\geq0}\frac{1}{(2n+1)^2}
      \frac{1}{(2l+1)^2}\\
    &\quad+\frac{16}{\pi^4}\sum_{n\geq0}\sum_{l\geq1}\Biggl(
      \sum_{m\geq0}\frac{1}{\bigl(4(l+m)^2-1\bigr)(2m-2n-1)}
      \Biggr)^2\\
    &<\frac{\pi^2}{16}+\frac{16}{\pi^4}\sum_{n\geq0}\sum_{l\geq1}
      \Biggl(\biggABS{\sum_{m=-n}^{n+1}\frac{1}{\bigl(4(l+m+n)^2
          -1\bigr)(2m-1)}}\\
    &\phantom{<\frac{\pi^2}{16}+\frac{16}{\pi^4}\sum_{n\geq0}
      \sum_{l\geq1}\Biggl(}
      +\sum_{m\geq2n+2}\frac{1}{\bigl(4(l+m)^2-1\bigr)(2n+3)}\Biggr)^2\\
    &=\frac{\pi^2}{16}+\frac{16}{\pi^4}\sum_{n\geq0}\sum_{l\geq1}
      \Biggl(\biggl\lvert\sum_{m=1}^{n+1}\frac{1}{\bigl(4(l+m+n)^2
        -1\bigr)(2m-1)}\\
    &\phantom{=\frac{\pi^2}{16}+\frac{16}{\pi^4}\sum_{n\geq0}
      \sum_{l\geq1}\Biggl(}
      +\sum_{m=1}^{n+1}\frac{1}{\bigl(4(l+1-m+n)^2-1\bigr)(1-2m)}
      \biggr\rvert\\
    &\phantom{=\frac{\pi^2}{16}+\frac{16}{\pi^4}\sum_{n\geq0}
      \sum_{l\geq1}\Biggl(}
      +\frac{1}{2n+3}\sum_{m\geq2n+l+2}\frac{1}{4m^2-1}\Biggr)^2\\
    &<\frac{\pi^2}{16}+\frac{16}{\pi^4}\sum_{n\geq0}\sum_{l\geq1}
      \Biggl(\sum_{m=1}^{n+1}\biggl\lvert\frac{1}{4(l+m+n)^2-1}\\
    &\phantom{<\frac{\pi^2}{16}+\frac{16}{\pi^4}\sum_{n\geq0}
      \sum_{l\geq1}\Biggl(\sum_{m=1}^{n+1}\biggl(}
      -\frac{1}{4(l+1-m+n)^2-1}\biggr\rvert\cdot\frac{1}{2m-1}\\
    &\phantom{<\frac{\pi^2}{16}+\frac{16}{\pi^4}\sum_{n\geq0}
      \sum_{l\geq1}\Biggl(}
      +\frac{1}{2n+3}\int_{2n+l+1}^\infty\frac{dm}{(2m-1)^2}\Biggr)^2\\
    &=\frac{16}{\pi^4}\sum_{n\geq0}\sum_{l\geq1}
      \Biggl(\sum_{m=1}^{n+1}\frac{4a_{ln}}{\underbrace{\bigl(
          4(l+m+n)^2-1\bigr)}_{=\underbrace{\scriptstyle
          (2(l+m+n)-1)}_{>a_{ln}}(2(l+m+n)+1)}
      \bigl(4(l+1-m+n)^2-1 \bigr)}\\
    &\phantom{=\frac{16}{\pi^4}\sum_{n\geq0}\sum_{l\geq1}\Biggl(}
      +\frac{1}{2(2n+3)(4n+2l+1)}\Biggr)^2+\frac{\pi^2}{16}\\
    &<\frac{\pi^2}{16}+\frac{16}{\pi^4}\sum_{n\geq0}\sum_{l\geq1}
      \Biggl(\sum_{m=0}^n\frac{4}{\bigl(2(l+m+n)+3\bigr)
        \bigl(4(l-m+n)^2-1\bigr)}\\
    &\phantom{<\frac{\pi^2}{16}+\frac{16}{\pi^4}\sum_{n\geq0}
      \sum_{l\geq1}\Biggl(}
      +\frac{1}{2(2n+3)(4n+2l+1)}\Biggr)^2\\
    &<\frac{\pi^2}{16}+\frac{16}{\pi^4}\sum_{n\geq0}\sum_{l\geq1}
      \Biggl(\frac{4}{2l+2n+3}\sum_{m=l}^{l+n}\frac{1}{4m^2-1}\\
    &\phantom{<\frac{\pi^2}{16}+\frac{16}{\pi^4}\sum_{n\geq0}
      \sum_{l\geq1}\Biggl(}
      +\frac{1}{2(2n+3)(4n+2l+1)}\Biggr)^2\\
    &<\frac{\pi^2}{16}+\frac{16}{\pi^4}\sum_{n\geq0}\sum_{l\geq1}
      \Biggl(\frac{4}{2l+2n+3}\biggl(\frac{1}{4l^2-1}+\underbrace{
        \int_l^{l+n}\frac{dm}{4m^2-1}}_{
        \begin{subarray}{l}
          <\int_l^{l+n}\frac{dm}{(2m-1)^2}\\
          =\frac{n}{(2l+2n-1)(2l-1)}\\
          <\frac{1}{2(2l-1)}
        \end{subarray}}
        \biggr)\\
    &\phantom{<\frac{\pi^2}{16}+\frac{16}{\pi^4}\sum_{n\geq0}
      \sum_{l\geq1}\Biggl(}
      +\frac{1}{2(2n+3)(4n+2l+1)}\Biggr)^2\\
    &<\frac{\pi^2}{16}+\frac{16}{\pi^4}\sum_{n\geq0}\sum_{l\geq1}
      \biggl(\frac{4}{2l+2n+3}\cdot\underbrace{
        \frac{2l+3}{2(2l+1)}}_{<1}\cdot
      \frac{1}{2l-1}\\
    &\phantom{<\frac{\pi^2}{16}+\frac{16}{\pi^4}\sum_{n\geq0}
      \sum_{l\geq1}\Biggl(}
      +\frac{1}{2(2n+3)(4n+2l+1)}\biggr)^2\\
    &<\frac{\pi^2}{16}+\frac{16}{\pi^4}\sum_{n,l\geq1}\biggl(
      \frac{4}{(2l+2n+1)(2l-1)}+\frac{1}{2(2n+1)(4n+2l-3)}\biggr)^2\\
    &<\frac{\pi^2}{16}+\frac{16}{\pi^4}\sum_{n,l\geq1}\biggl(
      \frac{2}{(n+l)l}+\frac{2}{n(n+l)}\biggr)^2\\
    &=\frac{\pi^2}{16}+\frac{64}{\pi^4}\sum_{n,l\geq1}\biggl(
      \frac{n+l}{(n+l)nl}\biggr)^2\\
    &=\frac{\pi^2}{16}+\frac{64}{\pi^4}\Biggl(\sum_{n\geq1}
      \frac{1}{n^2}\Biggr)^2\\
    &=\frac{\pi^2}{16}+\frac{16}{9}.
  \end{align*}
\end{proof}
Thus the Bogoliubov operator $V$ induced by $v$ by \eqref{VFORM1}
yields a $U(N)$--gauge invariant implementable localized endomorphism
$\rho_V$ of the chiral Dirac field. Since $M_V=\dim\mathfrak{k}_V=N$
by construction, it is clear that $\mathfrak{k}_V$ carries either the
defining representation of $U(N)$ or its complex conjugate.
By the discussion on p.~\pageref{page:RED}, the irreducible
constituents of $\rho_V$ correspond to the irreducible, mutually
inequivalent representations $\mathcal{U}_V^{(n)},\ n=0,\dots,N$. It
would be interesting to find a manageable description of the
Cayley--transformed operator $\tilde\varkappa v\tilde\varkappa^{-1}$
on $L^2(\RR)$, which would perhaps give an idea how to obtain
implementable localized endomorphisms in the two--dimensional massive
case.
\subsection{The charge of gauge invariant endomorphisms of the \CCR}
\label{sec:CCRCH}
Our discussion of the charge structure of gauge invariant
endomorphisms of the \CCR will be short compared to the CAR case.
Recall from the beginning of Section~\ref{sec:SECTORS} that we
consider the following situation: We have an infinite dimensional
vector space \KK together with a nondegenerate hermitian sesquilinear
form $\kappa$ and a fixed \BP $P_1$. \KK is assumed to be complete
\WRT the inner product induced by $P_1$. The gauge group $G$ acts by
diagonal Bogoliubov operators on \KK and can be identified with a
subgroup of the unitary group of $\KK_1\DEF P_1(\KK)$.  Our interest
is in the representations of $G$ on the \HSPs $H(\rho_V)$ which
implement gauge invariant quasi--free endomorphisms $\rho_V,\ 
V\in\SP{P_1}{}^G$.

By Lemma~\ref{lem:REPHRV}, it suffices to look at the values of the
implementers on the Fock vacuum. We have by \eqref{ccr:EOM},
\eqref{ccr:PSIA}
\begin{equation}
  \label{ccr:PAVOM}
  \PV{\alpha}\Omega_{P_1}=D_V\psi_{\alpha_1}\dotsm\psi_{\alpha_l}
  \exp\Bigl(\3Z_V\+a^*a^*\Bigr)\Omega_{P_1}
\end{equation}
where $D_V$ is a numerical constant, $Z_V$ is the symmetric
Hilbert--Schmidt operator defined in \eqref{ccr:ZV2}, $g_1,g_2,\dotsc$
is a $\kappa$--\ONB in $\mathfrak{k}_V\DEF P_V(\ker V\+)$, and
$\psi_j$ is the isometry obtained by polar decomposition of the
closure of $\pi_{P_1}(g_j)$.  In contrast to the Fermionic case (cf.\ 
Lemma~\ref{lem:GP}), there is no simple transformation law for the
$\psi_j$ under $G$. This difficulty can however be circumvented by
rewriting $\PV{\alpha}\Omega_{P_1}$ with the help of
\eqref{ccr:PWO}--\eqref{ccr:PSIU} and \eqref{ccr:DEF} as follows
\begin{equation}
  \label{ccr:PAVOM1}
  \begin{split}
    \PV{\alpha}\Omega_{P_1} &=\PUV\PW{\alpha}\Omega_{P_1}\\
    &=c_\alpha\PUV\pi_{P_1}\bigl(U_V\+(g_{\alpha_1})\dotsm U_V\+
      (g_{\alpha_l})\bigr)\Omega_{P_1}\\
    &=c_\alpha\4{\pi_{P_1}(g_{\alpha_1})}\dotsm
      \4{\pi_{P_1}(g_{\alpha_l})}\PV{0}\Omega_{P_1}\\
    &=c_\alpha D_V\4{\pi_{P_1}(g_{\alpha_1})}\dotsm
      \4{\pi_{P_1}(g_{\alpha_l})}\exp\Bigl(\3Z_V\+a^*a^*\Bigr)
      \Omega_{P_1}
  \end{split}
\end{equation}
(the bar denotes closure).
The behaviour of the $\4{\pi_{P_1}(g_j)}$ under gauge transformations
is obvious. It remains to consider the exponential term. It is a
salient feature of Definition~\ref{ccr:def:UW} that this term is gauge
invariant:
\begin{Lem}
  \label{ccr:lem:GEXP}\hspace*{\fill}\\
  Let $V\in\SP{P_1}{}^G$ be given, and let $Z_V\in\EE_{P_1}$ be
  defined by \eqref{ccr:ZV2}. Then $\exp(\3Z_V\+a^*a^*)\Omega_{P_1}$
  is invariant under all gauge transformations $\Gamma(U)$, $U\in G$.
\end{Lem}
\begin{proof}
  Let $E$ be the orthogonal projection onto $\HH\DEF\ker V\+=C\ker
  V^*,\ C=P_1-P_2$, as in \eqref{ccr:E}. Then $E$ commutes with $G$
  since $V$ and $P_1$ do so. Let $A\DEF ECE$ be the self--adjoint
  operator introduced in Lemma~\ref{ccr:lem:H}. Then $A$ and all its
  spectral projections commute with $G$. It follows that the positive
  part $A_+$ of $A$ and the operator ${A_+}^{-1}$ defined in
  Lemma~\ref{ccr:lem:H} also commute with $G$. Therefore $P_V\DEF
  VP_1V\++{A_+}^{-1}C$ and $Z_V\DEF(P_V)_{21}{(P_V)_{11}}^{-1}$
  commute with $G$ as well.  
  
  Arguing as in the proof of Lemma~\ref{lem:EXP}, one finds for $U\in
  G$
  $$\Gamma(U)\Bigl(\3Z_V\+a^*a^*\Bigr)^n\Omega_{P_1}
  =\Bigl(\3(UZ_V\+U\+)a^*a^*\Bigr)^n\Omega_{P_1}$$ 
  and finally
  $$\Gamma(U)\exp\Bigl(\3Z_V\+a^*a^*\Bigr)\Omega_{P_1}=
  \exp\Bigl(\3(UZ_V\+U\+)a^*a^*\Bigr)\Omega_{P_1}=
  \exp\Bigl(\3Z_V\+a^*a^*\Bigr)\Omega_{P_1}.$$
\end{proof}
Invoking the isomorphism $H(\rho_V)\cong\mathcal{F}_s
(\mathfrak{k}_V)$ from Corollary~\ref{ccr:cor:HRV}, we thus deduce
\begin{Th}
  \label{ccr:th:CHARGE}\hspace*{\fill}\\
  Let $P_1$ be a \BP of $(\KK,\kappa)$, let $G$ be a group consisting
  of diagonal Bogoliubov operators, and let $V\in\SP{P_1}{}^G$. Let
  $\mathfrak{k}_V=P_V(\ker V\+)$ be the subspace of \KK defined in
  Definition~\ref{ccr:def:PSIA}, with $\dim\mathfrak{k}_V=-\2\IND V$.
  Then $\mathfrak{k}_V$ is invariant under $G$, and the unitary
  representation $\mathcal{U}_V$ of $G$ on the \HSP of isometries
  $H(\rho_V)$ which implements $\rho_V$ in the Fock representation
  $\pi_{P_1}$ is unitarily equivalent to the representation on
  $\mathcal{F}_s(\mathfrak{k}_V)$ that is obtained by taking symmetric
  tensor powers of the representation on $\mathfrak{k}_V$.
\end{Th}
\begin{proof}
  $\mathfrak{k}_V$ is invariant under $G$ because $\ker V\+$ is
  invariant and because $G$ commutes with $P_V$ (see the proof of
  Lemma~\ref{ccr:lem:GEXP}). The assertion hence follows from
  \eqref{ccr:PAVOM1}, Lemma~\ref{ccr:lem:GEXP} and
  Lemma~\ref{lem:REPHRV}.
\end{proof}
Theorem~\ref{ccr:th:CHARGE} shows that genuine quasi--free
endomorphisms of the \CCR are even ``more reducible'' than
endomorphisms of the \CAR in that they are always \emph{infinite}
direct sums of irreducibles, a fact which explains the generic
occurrence of infinite statistics in the CCR case (cf.\ 
Theorem~\ref{ccr:th:IMP}). Again, each closed subspace of $H(\rho_V)$
spanned by all $\PV{\alpha}$ with length of $\alpha$ fixed is
invariant under $G$.  The special case of quasi--free automorphisms is
of little interest in the CCR case because they are all neutral
($\mathcal{U}_V$ is the trivial representation of $G$ if $\IND V=0$).

There is also less freedom in the choice of the representation of $G$
on $\mathfrak{k}_V$ (cf.\ Proposition~\ref{prop:HVKV}):
\begin{Prop}
  \hspace*{\fill}\\\indent
  a) If $\mathfrak{k}\subset\KK_1$ is a closed $G$--invariant subspace
  carrying a representation of class $\xi$ and if $\xi$ is contained
  in $\KK_1$ with infinite multiplicity, then there exists a diagonal
  Bogoliubov operator $V\in\SP{P_1}{}^G$ with
  $\mathfrak{k}_V=\mathfrak{k}$.
  
  b) A class $\xi$ of irreducible representations of $G$ is realized
  on some $\mathfrak{k}_V$, $V\in\SP{P_1}{}^G$, \IFF $\xi$ is contained
  in $\KK_1$ with infinite multiplicity.
\end{Prop}
\begin{proof}
  The proof of a) is the same as the proof of
  Proposition~\ref{prop:HVKV} b) 1.\ because
  $\SP{\text{diag}}{}=\IDIAG(\KK)$.
  
  To prove b), assume that $\xi$ is contained in $\KK_1$ with finite
  (or zero) multiplicity. (Recall that $G$ commutes with $P_1$ so that
  any irreducible class $\xi$ contained in \KK is contained in $\KK_1$
  or $\KK_2$.) Assume further that $V\in\SP{P_1}{}^G$ exists such that
  $\mathfrak{k}_V$ belongs to the class $\xi$. Then let $q_n$ be the
  $\kappa$--orthogonal projection onto $V^n(\mathfrak{k}_V)$
  $$q_n\DEF V^nP_V(\1-VV\+)V^{\dagger n}.$$
  Since $q_mq_n=0\ (m\neq
  n)$, the $V^n(\mathfrak{k}_V)$ are mutually $\kappa$--orthogonal,
  and their direct sum $\KK_V\DEF\oplus_{n=0}^\infty
  V^n(\mathfrak{k}_V)\subset\KK$ carries the representation $\xi$ with
  infinite multiplicity. Therefore $\xi$ must be contained in $\KK_2$
  with multiplicity $\infty$. However, the restriction of $\kappa$ to
  $\KK_V$ is positive definite, whereas the restriction to $\KK_2$ is
  negative definite. This is a contradiction and shows that no $V$
  with $\mathfrak{k}_V$ of class $\xi$ can exist. On the other hand,
  if $\xi$ is contained in $\KK_1$ with infinite multiplicity, then we
  are back in the situation of a).
\end{proof}

The remarks made in Section~\ref{sec:CARCH} after
Proposition~\ref{prop:HVKV} on the relation to the generic
superselection structure apply here as well. In particular, one needs
additional assumptions to ensure the existence of conjugates (cf.\ 
Proposition~\ref{prop:CONJ}):
\begin{Prop}
  \hspace*{\fill}\\
  Assume that there exists an orthogonal projection\footnote{Such $P$
    is a CAR \BP of \KK, but not a \BP of $(\KK,\kappa)$ since
    $\kappa$ is not positive definite on $\RAN P$.} $P$ of \KK with
  $P+\4{P}=\1$ which commutes with $P_1$ and $G$. Let $p_\pm$ be
  defined by
  $$p_+\DEF PP_1,\qquad p_-\DEF PP_2,$$
  and assume that there exists a
  unitary operator $\mathcal{C}_{+-}$ from $\KK_-\DEF p_-(\KK)$ onto
  $\KK_+\DEF p_+(\KK)$ which commutes with $G$. Let $\mathcal{C}$ be
  the unique operator on \KK which commutes with $P$, which fulfills
  $\4{\mathcal{C}}=\mathcal{C}$, and which is given on
  $P(\KK)=\KK_+\oplus\KK_-$ by the matrix
  $$P\mathcal{C}P\DEF
  \begin{pmatrix}
    0 & \mathcal{C}_{+-} \\ {\mathcal{C}_{+-}}^* & 0
  \end{pmatrix}.$$
  Then $\mathcal{C}$ is gauge invariant, unitary and self--adjoint, it
  fulfills $\{\mathcal{C},C\}=0$ and $\mathcal{C}\+=-\mathcal{C}$, and
  the map
  \begin{equation}
    \label{ccr:CHCONJ}
    V\mapsto V^c\DEF \mathcal{C}V\mathcal{C}
  \end{equation}
  is an involutive automorphism of $\SP{P_1}{}^G$. One has 
  \begin{equation}
    \label{ccr:KVC}
    \mathfrak{k}_{V^c}=\mathcal{C}({\mathfrak{k}_V}^*)
  \end{equation}
  so that the
  representation $\mathcal{U}_{V^c}$ of $G$ on $H(\rho_{V^c})$ is
  unitarily equivalent to the complex conjugate of $\mathcal{U}_V$.
\end{Prop}
\begin{proof}
  The asserted properties of $\mathcal{C}$ are readily verified. It
  only remains to prove \eqref{ccr:KVC}. Let us first show that
  \begin{equation}
    \label{ccr:PVC}
    P_{V^c}=\mathcal{C}\4{P_V}\mathcal{C}
  \end{equation}
  where $P_V,\ P_{V^c}$ are the \BPs associated with $V,\ V^c$
  according to Definition~\ref{ccr:def:UW}:
  \begin{equation}
    \label{ccr:PVPVC}
    P_V=VP_1V\++p_V,\qquad P_{V^c}=V^cP_1V^{c\dagger}+p_{V^c}.
  \end{equation}
  Let $E$ resp.\ $E^c$ be the orthogonal projections onto $\ker V\+$
  resp.\ $\ker V^{c\dagger}$, and let $A\DEF ECE,\ A^c\DEF E^cCE^c$ be
  the corresponding operators as in Lemma~\ref{ccr:lem:H}, with
  positive/negative parts $A_\pm,\ {A^c}_\pm$. Since $\ker
  V^{c\dagger}=\mathcal{C}(\ker V\+)$, one has
  $E^c=\mathcal{C}E\mathcal{C}$ and
  $$A^c=\mathcal{C}E(\mathcal{C}C\mathcal{C})E\mathcal{C}=
    -\mathcal{C}ECE\mathcal{C}=-\mathcal{C}A\mathcal{C}$$
  so that
  ${A^c}_+=\mathcal{C}A_-\mathcal{C}$. Since $A_-=\4{A_+}$ (see the
  proof of Lemma~\ref{ccr:lem:H}), since $\{\mathcal{C},C\}=0$ and
  $\4{C}=-C$, this implies
  $$p_{V^c}\DEF{{A^c}_+}^{-1}C=\mathcal{C}\4{A_+}^{\,-1}\mathcal{C}C=
    \mathcal{C}\4{{A_+}^{-1}C}\mathcal{C}=
    \mathcal{C}\4{p_V}\mathcal{C}.$$
  Thus we get from \eqref{ccr:PVPVC}, using
  $\mathcal{C}P_1\mathcal{C}=P_2$ and $\mathcal{C}\+=-\mathcal{C}$
  $$P_{V^c}=V^cP_1V^{c\dagger}+p_{V^c}=\mathcal{C}VP_2V\+\mathcal{C}+
    \mathcal{C}\4{p_V}\mathcal{C}=\mathcal{C}\4{P_V}\mathcal{C}$$
  which proves \eqref{ccr:PVC}. It follows that \eqref{ccr:KVC} holds:
  $$\mathfrak{k}_{V^c}=P_{V^c}(\ker V^{c\dagger})=
    \mathcal{C}\4{P_V}\mathcal{C}(\mathcal{C}\ker V\+)=\mathcal{C}
    ({\mathfrak{k}_V}^*).$$
  Therefore the representation of $G$ on
  $\mathfrak{k}_{V^c}$ is unitarily equivalent to the complex
  conjugate of the representation on $\mathfrak{k}_V$, so that
  $\mathcal{U}_{V^c}\simeq{\mathcal{U}_V}^*$ by
  Theorem~\ref{ccr:th:CHARGE}.
\end{proof}


%% file: not.tex
\chapter*{What We Have Not Learnt Yet}
\addtocontents{toc}{\vspace{-1ex}}
Let us comment on some open questions and possible generalizations.
As we have seen, implementable quasi--free endomorphisms of the CAR
and \CCRs are always accompanied by \emph{Fock spaces of isometries.}
This fact implies, among other things, that these endomorphisms
restrict to reducible endomorphisms of the observable algebra. Gauge
invariant endomorphisms behave in some respect like ``master
endomorphisms'', because they carry, together with some ``basic''
representation of the \GG, all higher (anti)symmetric tensor powers
thereof. In typical cases, e.g.\ for the classical compact Lie groups,
any irreducible representation of the group is equivalent to a
subrepresentation of some tensor power of the defining representation.
In such cases there will exist quasi--free endomorphisms which contain
each \SSS as a subrepresentation.

It is an interesting problem how to obtain irreducible subobjects of
quasi--free endomorphisms. According to the general theory, one needs
gauge invariant isometries fulfilling relation \eqref{BORCHERS} on
Fock space for this purpose. Such isometries would also permit to
define direct sums of quasi--free endomorphisms, so that one would get
the full Doplicher--Roberts category generated by quasi--free
endomorphisms.

Another important question is whether one can exhibit classes of
localized Bogoliubov operators with finite nonzero index, say on the
single--particle space of the time--zero Dirac field, and such that
the implementability condition holds. Recall that our construction of
such Bogoliubov operators made essential use of the existence of local
Fourier bases on the circle. It would be desirable to have a
basis--independent characterization of such Bogoliubov operators in
Minkowski space, but this is an unsolved functional analytic problem.
Of particular interest would be the massive case in two dimensions,
where one might hope to find localized quasi--free endomorphisms with
non--Abelian (``plektonic'') braid group statistics. But preliminary
calculations based on our explicit formulas indicate that the
implementers corresponding to irreducible subobjects of quasi--free
endomorphisms can only have ``anyonic'' \CRs (i.e.\ Abelian braid
group statistics). A related question is what kind of quantum
symmetries beyond compact groups can be realized by quasi--free
endomorphisms. 

Let us mention that our results might be relevant for the discussion
of the superselection structure of the massless Thirring model. The
sectors of this model have been described by Streater \cite{Str74}.
Later, Carey, Ruijsenaars and Wright constructed approximate cutoff
fields which reproduce Klaiber's $n$--point functions of the Thirring
model \cite{CRW}. The Thirring fields can be obtained as strong limits
of the approximate fields. Whereas the latter are always (multiples
of) unitary operators, it could well be that the Thirring fields are
only (multiples of) isometries. They would then fall into the scope of
our work. The \HSP chosen by Streater would then be too ``big'', in
that the Thirring fields would not reach all sectors from the vacuum.
\medskip

On the mathematical side, our investigations can be extended in
several interesting directions. One could enlarge the class of
representations under consideration, by studying the implementation of
quasi--free endomorphisms in arbitrary quasi--free states, and not
just in \FSs. Such an analysis should be possible because the main
technical tools, the 
quasi--equivalence criteria, apply to arbitrary quasi--free states.

One could also enlarge the class of transformations and consider e.g.\ 
the implementation of completely positive quasi--free maps. Some work
has been done in this context, notably by Evans
\cite{DEE77,DEE79,DEE80a}. The ``Bogoliubov operators'' corresponding
to completely positive quasi--free maps are in general not isometric,
but only contractive. As an analogue of implementation, Evans
constructed ``dissipators'' for completely positive quasi--free maps,
i.e.\ nonzero contractions on Fock space which are ``dominated'' by
the completely positive map. But, as shown by Arveson in his
generalization of Powers' index theory of $E_0$--semigroups to
semigroups of completely positive maps \cite{Arv96a,Arv96b}, such maps
also admit an implementation by ``metric operator spaces'', and this
notion is very close to the implementation of endomorphisms by \HSPs
of isometries.

Let us finally mention that it would be interesting to study the
representations of the \CAs that are generated by the implementers of
quasi--free endomorphisms, and that one could examine whether the
Fredholm index of CCR Bogoliubov operators, which is a finer invariant
than the Jones--Kosaki index, has a counterpart on the algebraic
level.


%% file: master.bbl
\newcommand{\etalchar}[1]{$^{#1}$}
\providecommand{\bysame}{\leavevmode\hbox to3em{\hrulefill}\thinspace}
\begin{thebibliography}{BDM{\etalchar{+}}96}

\bibitem[ACE84]{ACE}
H.~Araki, A.~L. Carey, and D.~E. Evans, \emph{On $\mathcal{O}_{n+1}$}, J.
  Operator Theory \textbf{12} (1984), 247.

\bibitem[Adl96]{Ad}
C.~Adler, \emph{Braid group statistics in two--dimensional quantum field
  theory}, Rev. Math. Phys. \textbf{8} (1996), 907.

\bibitem[AE83]{AE}
H.~Araki and D.~E. Evans, \emph{On a {\CAL} approach to phase transition in the
  two--dimensional {I}sing model}, Commun. Math. Phys. \textbf{91} (1983), 489.

\bibitem[Ara63]{A63}
H.~Araki, \emph{A lattice of von {N}eumann algebras associated with the quantum
  theory of a free {B}ose field}, J. Math. Phys. \textbf{4} (1963), 1343.

\bibitem[Ara64]{A64}
H.~Araki, \emph{Type of {\VN} algebras associated with the free field}, Progr.
  Theoret. Phys. \textbf{32} (1964), 956.

\bibitem[Ara68]{A68}
H.~Araki, \emph{On the diagonalization of a bilinear {H}amiltonian by a
  {B}ogoliubov transformation}, Publ. Res. Inst. Math. Sci. \textbf{4} (1968),
  387.

\bibitem[Ara87]{A87}
H.~Araki, \emph{Bogoliubov automorphisms and {F}ock representations of
  canonical anticommutation relations}, Operator Algebras and Mathematical
  Physics, vol.~62, Am. Math. Soc., 1987, p.~23.

\bibitem[Ara71]{A70}
H.~Araki, \emph{On quasifree states of {CAR} and {B}ogoliubov automorphisms},
  Publ. Res. Inst. Math. Sci. \textbf{6} (1970/71), 385.

\bibitem[Ara72]{A71}
H.~Araki, \emph{On quasifree states of the canonical commutation relations
  ({II})}, Publ. Res. Inst. Math. Sci. \textbf{7} (1971/72), 121.

\bibitem[Arv94]{A94}
W.~Arveson, \emph{${E_0}$--semigroups in quantum field theory}, 1994, Lecture
  presented at the von {N}eumann symposium on quantization and nonlinear wave
  equations held at {MIT}.

\bibitem[Arv96a]{Arv96a}
W.~Arveson, \emph{The index of a quantum dynamical semigroup}, preprint, Univ.
  California, 1996.

\bibitem[Arv96b]{Arv96b}
W.~Arveson, \emph{On the index and dilations of completely positive
  semigroups}, preprint, Univ. California, 1996.

\bibitem[AS72]{AS71}
H.~Araki and M.~Shiraishi, \emph{On quasifree states of the canonical
  commutation relations ({I})}, Publ. Res. Inst. Math. Sci. \textbf{7}
  (1971/72), 105.

\bibitem[AY81]{AY}
H.~Araki and S.~Yamagami, \emph{An inequality for {H}ilbert--{S}chmidt norm},
  Commun. Math. Phys. \textbf{81} (1981), 89.

\bibitem[AY82]{A82}
H.~Araki and S.~Yamagami, \emph{On quasi--equivalence of quasifree states of
  the canonical commutation relations}, Publ. Res. Inst. Math. Sci. \textbf{18}
  (1982), 283.

\bibitem[BCL97]{BCL}
P.~Bertozzini, R.~Conti, and R.~Longo, \emph{Covariant sectors with infinite
  dimension and positivity of the energy}, preprint, Univ. Roma ``Tor
  Vergata'', 1997.

\bibitem[BDF87]{BAF}
D.~Buchholz, C.~D'Antoni, and K.~Fredenhagen, \emph{The universal structure of
  local algebras}, Commun. Math. Phys. \textbf{111} (1987), 123.

\bibitem[BDLR92]{BDLR}
D.~Buchholz, S.~Doplicher, R.~Longo, and J.~E. Roberts, \emph{A new look at
  {G}oldstone's {T}heorem}, Rev. Math. Phys. \textbf{{\rm special issue}}
  (1992), 49.

\bibitem[BDM{\etalchar{+}}96]{BDMRS}
D.~Buchholz, S.~Doplicher, G.~Morchio, J.~E. Roberts, and F.~Strocchi, \emph{A
  model for charges of electromagnetic type}, Operator Algebras and Quantum
  Field Theory (S.~Doplicher, R.~Longo, J.~E. Roberts, and L.~Zsido, eds.),
  Accademia Nazionale dei Lincei, Roma, International Press, 1996, p.~647.

\bibitem[Ber66]{Ber}
F.~A. Berezin, \emph{The method of second quantization}, Academic Press, New
  York, 1966.

\bibitem[BF82]{BF}
D.~Buchholz and K.~Fredenhagen, \emph{Locality and the structure of particle
  states}, Commun. Math. Phys. \textbf{84} (1982), 1.

\bibitem[BFK96]{BFK}
R.~Brunetti, K.~Fredenhagen, and M.~K{\"o}hler, \emph{The microlocal spectrum
  condition and the {W}ick's polynomials of free fields}, Commun. Math. Phys.
  \textbf{180} (1996), 312.

\bibitem[BGL93]{BGL}
R.~Brunetti, D.~Guido, and R.~Longo, \emph{Modular structure and duality in
  conformal quantum field theory}, Commun. Math. Phys. \textbf{156} (1993),
  201.

\bibitem[BH92]{BH}
U.~Bunke and T.~Hirschmann, \emph{The index of the scattering operator on the
  positive spectral subspace}, Commun. Math. Phys. \textbf{148} (1992), 487.

\bibitem[BHJ26]{BHJ}
M.~Born, W.~Heisenberg, and P.~Jordan, \emph{Zur {Q}uantenmechanik {II}}, Z.
  Phys. \textbf{35} (1926), 557.

\bibitem[Bin93]{CB0}
C.~Binnenhei, \emph{Bogoliubov--{T}ransformationen und lokalisierte
  {M}orphismen f{\"u}r freie {D}iracfelder}, diploma thesis BONN-IR-93-33,
  Phys. {I}nst. {U}niv. {B}onn, 1993.

\bibitem[Bin95]{CB1}
C.~Binnenhei, \emph{Implementation of endomorphisms of the {CAR} algebra}, Rev.
  Math. Phys. \textbf{7} (1995), 833.

\bibitem[Bin97]{CB2}
C.~Binnenhei, \emph{On the even {CAR} algebra}, Lett. Math. Phys. \textbf{40}
  (1997), 91.

\bibitem[Bin98]{CB3}
C.~Binnenhei, \emph{{${\mathcal O}_\infty$} realized on {B}ose {F}ock space},
  Commun. Math. Phys. \textbf{195} (1998), 353.

\bibitem[BJ25]{BJ}
M.~Born and P.~Jordan, \emph{On quantum mechanics}, Z. Phys. \textbf{34}
  (1925), 858.

\bibitem[BJ96]{BJ96}
O.~Bratteli and P.~E.~T. J{\o}rgensen, \emph{A connection between
  multiresolution wavelet theory of scale ${N}$ and representations of the
  {C}untz algebra $\mathcal{O}_{N}$}, 1996, preprint, funct-an/9612006.

\bibitem[BJKW97]{BJKW}
O.~Bratteli, P.~E.~T. J{\o}rgensen, A.~Kishimoto, and R.~F. Werner, \emph{Pure
  states on $\mathcal{O}_{d}$}, 1997, preprint, funct-an/9711004.

\bibitem[BJP96]{BJP}
O.~Bratteli, P.~E.~T. J{\o}rgensen, and G.~L. Price, \emph{Endomorphisms of
  ${\BB(\HH)}$}, Quantization, {N}onlinear {P}artial {D}ifferential
  {E}quations, and {O}perator {A}lgebra (W.~Arveson, T.~Branson, and I.~Segal,
  eds.), Proc. Sympos. Pure Math., vol.~59, Amer. Math. Soc., 1996, p.~93.

\bibitem[Bla58]{Blat}
R.~J. Blattner, \emph{Automorphic group representations}, Pacific J. Math.
  \textbf{8} (1958), 665.

\bibitem[BLOT90]{BLOT}
N.~N. Bogolubov, A.~A. Logunov, A.~I. Oksak, and I.~T. Todorov, \emph{General
  principles of quantum field theory}, Kluwer Acad. Publ., Dordrecht, Boston,
  London, 1990.

\bibitem[BMT88]{BMT}
D.~Buchholz, G.~Mack, and I.~Todorov, \emph{The current algebra on the circle
  as a germ of local field theories}, Nuclear Phys. B \textbf{5B} (1988), 20.

\bibitem[BNS98]{BNS}
G.~B{\"o}hm, F.~Nill, and K.~Szlach{\'a}nyi, \emph{Weak {H}opf algebras: {I}.
  {I}ntegral theory and {$C^*$}--structure}, preprint, Sfb 288, 1998.

\bibitem[B{\"o}c94]{JMB1}
J.~M. B{\"o}ckenhauer, \emph{Localized endomorphisms of the chiral {I}sing
  model}, Commun. Math. Phys. \textbf{177} (1994), 265.

\bibitem[B{\"o}c96]{JMB2}
J.~M. B{\"o}ckenhauer, \emph{An algebraic formulation of level one
  {W}ess--{Z}umino--{W}itten models}, Rev. Math. Phys. \textbf{8} (1996), 925.

\bibitem[Boe70]{Boe}
H.~Boerner, \emph{Representations of groups}, second ed., North--Holland
  Publishing Company, Amsterdam, London, 1970.

\bibitem[Bor65]{B65}
H.-J. Borchers, \emph{Local rings and the connection of spin with statistics},
  Commun. Math. Phys. \textbf{1} (1965), 281.

\bibitem[Bor67a]{B67}
H.-J. Borchers, \emph{On the theory of local observables}, Carg\`ese Lectures
  in Theoretical Physics 1965 (F.~Lur\c{c}at, ed.), Gordon and Breach, New
  York, London, Paris, 1967, p.~3.

\bibitem[Bor67b]{B67a}
H.-J. Borchers, \emph{A remark on a theorem of {B}. {M}isra}, Commun. Math.
  Phys. \textbf{4} (1967), 315.

\bibitem[Bor95]{B95}
H.-J. Borchers, \emph{On the use of modular groups in quantum field theory},
  Ann. Inst. H. Poincar{\'e} \textbf{63} (1995), 331.

\bibitem[BP83a]{BPb}
B.~M. Baker and R.~T. Powers, \emph{Product states and {$C^*$}--dynamical
  systems of product type}, J. Funct. Anal. \textbf{50} (1983), 229.

\bibitem[BP83b]{BPa}
B.~M. Baker and R.~T. Powers, \emph{Product states of the gauge invariant and
  rotationally invariant {\CARs}}, J. Operator Theory \textbf{10} (1983), 365.

\bibitem[BR81]{BR}
O.~Bratteli and D.~W. Robinson, \emph{Operator algebras and quantum statistical
  mechanics}, vol. 1,2, Springer--Verlag, Berlin, Heidelberg, New York,
  1979--1981.

\bibitem[BSM90]{BSM}
D.~Buchholz and H.~Schulz-Mirbach, \emph{Haag duality in conformal quantum
  field theory}, Rev. Math. Phys. \textbf{2} (1990), 105.

\bibitem[Buc82]{Buch82}
D.~Buchholz, \emph{The physical state space of quantum electrodynamics},
  Commun. Math. Phys. \textbf{85} (1982), 49.

\bibitem[BW75]{BW}
J.~J. Bisognano and E.~H. Wichmann, \emph{On the duality condition for a
  {H}ermitian scalar field}, J. Math. Phys. \textbf{16} (1975), 985.

\bibitem[Car84]{Ca84}
A.~L. Carey, \emph{Some infinite dimensional groups and bundles}, Publ. Res.
  Inst. Math. Sci. \textbf{20} (1984), 1103.

\bibitem[CHO82]{CHOB}
A.~L. Carey, C.~A. Hurst, and D.~M. O'Brien, \emph{Automorphisms of the
  canonical anticommutation relations and index theory}, J. Funct. Anal.
  \textbf{48} (1982), 360.

\bibitem[CO83]{COB}
A.~L. Carey and D.~M. O'Brien, \emph{Automorphisms of the infinite dimensional
  {C}lifford algebra and the {A}tiyah--{S}inger mod 2 index}, Topology
  \textbf{22} (1983), 437.

\bibitem[Cob67]{Co}
L.~A. Coburn, \emph{The {\CAL} generated by an isometry}, Bull. Amer. Math.
  Soc. (N.S.) \textbf{73} (1967), 722.

\bibitem[Con94]{Con}
A.~Connes, \emph{Noncommutative geometry}, Academic Press, San Diego, New York,
  Boston, 1994.

\bibitem[CR87]{CR}
A.~L. Carey and S.~N.~M. Ruijsenaars, \emph{On {F}ermion gauge groups, current
  algebras and {K}ac--{M}oody algebras}, Acta Appl. Math. \textbf{10} (1987),
  1.

\bibitem[CRW85]{CRW}
A.~L. Carey, S.~N.~M. Ruijsenaars, and J.~D. Wright, \emph{The massless
  {T}hirring model: {P}ositivity of {K}laiber's $n$--point functions}, Commun.
  Math. Phys. \textbf{99} (1985), 347.

\bibitem[Cun77]{C77}
J.~Cuntz, \emph{Simple {\CALs} generated by isometries}, Commun. Math. Phys.
  \textbf{57} (1977), 173.

\bibitem[Cun80]{C80}
J.~Cuntz, \emph{Automorphisms of certain simple {\CALs}}, Quantum Fields ---
  Algebras, Processes. Proceedings of the Symposium "Bielefeld Encounters in
  Physics and Mathematics II" 1978 (Wien, New York) (L.~Streit, ed.),
  Springer--Verlag, 1980, p.~187.

\bibitem[Cun81]{C81}
J.~Cuntz, \emph{K--theory for certain {\CALs}}, Ann. of Math. (2) \textbf{113}
  (1981), 181.

\bibitem[Cun91]{C91}
J.~Cuntz, \emph{Regular actions of {H}opf algebras on the {\CAL} generated by a
  {H}ilbert space}, preprint, Math. {I}nst. {U}niv. {H}eidelberg, 1991.

\bibitem[Del68]{DA}
G.~F. Dell'Antonio, \emph{Structure of the algebras of some free systems},
  Commun. Math. Phys. \textbf{9} (1968), 81.

\bibitem[DHR69a]{DHR1}
S.~Doplicher, R.~Haag, and J.~E. Roberts, \emph{Fields, observables and gauge
  transformations {I}}, Commun. Math. Phys. \textbf{13} (1969), 1.

\bibitem[DHR69b]{DHR2}
S.~Doplicher, R.~Haag, and J.~E. Roberts, \emph{Fields, observables and gauge
  transformations {II}}, Commun. Math. Phys. \textbf{15} (1969), 173.

\bibitem[DHR71]{DHR3}
S.~Doplicher, R.~Haag, and J.~E. Roberts, \emph{Local observables and particle
  statistics {I}}, Commun. Math. Phys. \textbf{23} (1971), 199.

\bibitem[DHR74]{DHR4}
S.~Doplicher, R.~Haag, and J.~E. Roberts, \emph{Local observables and particle
  statistics {II}}, Commun. Math. Phys. \textbf{35} (1974), 49.

\bibitem[Dir26]{Dir26}
P.~A.~M. Dirac, \emph{The fundamental equations of quantum mechanics}, Proc.
  Roy. Soc. London Ser. A \textbf{109} (1926), 642.

\bibitem[Dir27]{Dir27}
P.~A.~M. Dirac, \emph{The quantum theory of the emission and absorption of
  radiation}, Proc. Roy. Soc. London Ser. A \textbf{114} (1927), 243.

\bibitem[Dou72]{D}
R.~D. Douglas, \emph{Banach algebra techniques in operator theory}, Academic
  Press, New York, San Francisco, London, 1972.

\bibitem[DP68]{DP}
S.~Doplicher and R.~T. Powers, \emph{On the simplicity of the even {\CAR} and
  free field models}, Commun. Math. Phys. \textbf{7} (1968), 77.

\bibitem[DR72]{DR72}
S.~Doplicher and J.~E. Roberts, \emph{Fields, statistics and non--abelian gauge
  groups}, Commun. Math. Phys. \textbf{28} (1972), 331.

\bibitem[DR87]{DR87}
S.~Doplicher and J.~E. Roberts, \emph{Duals of compact {L}ie groups realized in
  the {C}untz algebras and their actions on {C*}--algebras}, J. Funct. Anal.
  \textbf{74} (1987), 96.

\bibitem[DR88]{DR88}
S.~Doplicher and J.~E. Roberts, \emph{Compact group actions on {\CALs}}, J.
  Operator Theory \textbf{19} (1988), 283.

\bibitem[DR89a]{DR89a}
S.~Doplicher and J.~E. Roberts, \emph{Endomorphisms of {\CALs}, cross products
  and duality for compact groups}, Ann. of Math. (2) \textbf{130} (1989), 75.

\bibitem[DR89b]{DR89}
S.~Doplicher and J.~E. Roberts, \emph{A new duality theory for compact groups},
  Invent. Math. \textbf{98} (1989), 157.

\bibitem[DR90]{DR90}
S.~Doplicher and J.~E. Roberts, \emph{Why there is a field algebra with a
  compact gauge group describing the superselection structure in particle
  physics}, Commun. Math. Phys. \textbf{131} (1990), 51, see further references
  herein.

\bibitem[EK98]{EK}
D.~E. Evans and Y.~Kawahigashi, \emph{Quantum symmetries on operator algebras},
  Oxford Mathematical Monographs, Oxford Univ. Press, 1998.

\bibitem[Ell95]{Ell}
G.~A. Elliott, \emph{The classification problem for amenable {\CALs}}, Proc.
  {ICM} {Z}{\"u}rich 1994, Birkh{\"a}user {V}erlag, {B}asel, 1995, p.~922.

\bibitem[Eva77]{DEE77}
D.~E. Evans, \emph{Some semigroups of completely positive maps on the {\CCR}},
  J. Funct. Anal. \textbf{26} (1977), 369.

\bibitem[Eva79]{DEE79}
D.~E. Evans, \emph{Completely positive quasi--free maps on the {\CAR}}, Commun.
  Math. Phys. \textbf{70} (1979), 53.

\bibitem[Eva80a]{DEE80a}
D.~E. Evans, \emph{Dissipators for symmetric quasi--free dynamical semigroups
  on the {\CAR}}, J. Funct. Anal. \textbf{37} (1980), 318.

\bibitem[Eva80b]{DEE80}
D.~E. Evans, \emph{On $\mathcal{O}_n$}, Publ. Res. Inst. Math. Sci. \textbf{16}
  (1980), 915.

\bibitem[FG93]{FG}
J.~Fr{\"o}hlich and F.~Gabbiani, \emph{Operator algebras and conformal quantum
  field theory}, Commun. Math. Phys. \textbf{155} (1993), 569.

\bibitem[FGM90]{FGM}
J.~Fr{\"o}hlich, F.~Gabbiani, and P.-A. Marchetti, \emph{Braid statistics in
  three--dimensional local quantum theory}, The Algebraic Theory of
  Superselection Sectors. {I}ntroduction and Recent Results (I.S.I. Guccia,
  Palermo 1989) (D.~Kastler, ed.), World Scientific, Singapore, 1990, p.~259.

\bibitem[FGR96]{FGR}
K.~Fredenhagen, M.~R. Gaberdiel, and S.~M. R{\"u}ger, \emph{Scattering states
  of {P}lektons (particles with braid group statistics) in $2+1$ dimensional
  quantum field theory}, Commun. Math. Phys. \textbf{175} (1996), 319.

\bibitem[FGV92]{FGV}
J.~Fuchs, A.~Ganchev, and P.~Vecserny\'es, \emph{Level 1 {WZW} superselection
  sectors}, Commun. Math. Phys. \textbf{146} (1992), 553.

\bibitem[FK93]{FK}
J.~Fr{\"o}hlich and T.~Kerler, \emph{Quantum groups, quantum categories and
  quantum field theory}, Lecture Notes in Mathematics, no. 1542,
  Springer--Verlag, 1993.

\bibitem[FL97]{FL}
N.~J. Fowler and M.~Laca, \emph{Endomorphisms of ${\BB(\HH)}$, extensions of
  pure states, and a class of representations of $\mathcal{O}_n$}, preprint,
  Univ. of Newcastle, Australia, 1997, funct-an/9709004.

\bibitem[FM83]{FM}
K.~Fredenhagen and M.~Marcu, \emph{Charged states in ${\ZZ_2}$ gauge theories},
  Commun. Math. Phys. \textbf{92} (1983), 81.

\bibitem[Foi83]{Fo}
J.~J. Foit, \emph{Abstract twisted duality for quantum free {F}ermi fields},
  Publ. Res. Inst. Math. Sci. \textbf{19} (1983), 729.

\bibitem[Fre73]{F73}
K.~Fredenhagen, \emph{Lokale {S}truktur und {P}olarzerlegung des geladenen
  {S}kalarfeldes}, diploma thesis, Univ. Hamburg, 1973.

\bibitem[Fre77]{F77}
K.~Fredenhagen, \emph{Implementation of automorphisms and derivations of the
  {CAR}--algebra}, Commun. Math. Phys. \textbf{52} (1977), 255.

\bibitem[Fre85]{F85}
K.~Fredenhagen, \emph{On the modular structure of local algebras of
  observables}, Commun. Math. Phys. \textbf{97} (1985), 79.

\bibitem[Fre90]{Fre90}
K.~Fredenhagen, \emph{Generalizations of the theory of {\SSSs}}, The Algebraic
  Theory of Superselection Sectors. {I}ntroduction and Recent Results (I.S.I.
  Guccia, Palermo 1989) (D.~Kastler, ed.), World Scientific, Singapore, 1990,
  p.~379.

\bibitem[Fre92]{F92}
K.~Fredenhagen, \emph{Product of states}, Quantum Groups and Related Topics
  (R.~Gielerak, J.~Lukierski, and Z.~Popwicz, eds.), Max Born Symposium in
  Theoretical Physics, Wroclaw 1991, Kluwer Acad. Publ., 1992, p.~199.

\bibitem[Fre93]{Fre93}
K.~Fredenhagen, \emph{Superselection sectors in low dimensional quantum field
  theory}, J. Geom. Phys. \textbf{11} (1993), 337.

\bibitem[Fre94]{F94}
K.~Fredenhagen, \emph{Superselection sectors with infinite statistical
  dimension}, preprint DESY--94--071, Hamburg Univ., 1994.

\bibitem[Fre95]{F}
K.~Fredenhagen, \emph{Superselection sectors}, lectures given at Hamburg Univ.,
  available at http://www.desy.de/algebra/notes.html, 1995.

\bibitem[Fri53]{Fri}
K.~O. Friedrichs, \emph{Mathematical aspects of the quantum theory of fields},
  Interscience Publ., New York, 1953.

\bibitem[Fr{\"o}76]{Froe}
J.~Fr{\"o}hlich, \emph{New superselection sectors (`soliton states') in
  two--dimensional {B}ose quantum field models}, Commun. Math. Phys.
  \textbf{47} (1976), 269.

\bibitem[FRS89]{FRS1}
K.~Fredenhagen, K.-H. Rehren, and B.~Schroer, \emph{Superselection sectors with
  braid group statistics and exchange algebras {I}. {G}eneral theory}, Commun.
  Math. Phys. \textbf{125} (1989), 201.

\bibitem[FRS92]{FRS2}
K.~Fredenhagen, K.-H. Rehren, and B.~Schroer, \emph{Superselection sectors with
  braid group statistics and exchange algebras {II}. {G}eometric aspects and
  conformal covariance}, Rev. Math. Phys. \textbf{{\rm special issue}} (1992),
  113.

\bibitem[GJ71]{GJ}
J.~Glimm and A.~Jaffe, \emph{Quantum field theory models}, Statistical
  Mechanics and Quantum Field Theory, Les Houches 1970 (C.~DeWitt and R.~Stora,
  eds.), Gordon and Breach, New York, 1971.

\bibitem[GL92]{GL}
D.~Guido and R.~Longo, \emph{Relativistic invariance and charge conjugation in
  quantum field theory}, Commun. Math. Phys. \textbf{148} (1992), 521.

\bibitem[GLW97]{GLW}
D.~Guido, R.~Longo, and H.-W. Wiesbrock, \emph{Extensions of conformal nets and
  superselection structure}, preprint 255, Sfb 288, 1997.

\bibitem[Haa58]{H58}
R.~Haag, \emph{Quantum field theories with composite particles and asymptotic
  conditions}, Phys. Rev. \textbf{112} (1958), 669.

\bibitem[Haa59]{H57}
R.~Haag, \emph{Discussion des ``axiomes'' et des propri\'et\'es asymptotiques
  d'une th\'eorie des champs locales avec particules compos\'ees}, Les
  Probl\`emes Math\'ematiques de la Th\'eorie Quantique des Champs, Colloque
  Internationaux du CNRS LXXV (Lille 1957), CNRS Paris, 1959, p.~151.

\bibitem[Haa63]{H63}
R.~Haag, \emph{Bemerkungen zum {N}ahwirkungsprinzip in der {Q}uantenphysik},
  Ann. Physik (7) \textbf{11} (1963), 29.

\bibitem[Haa96]{H}
R.~Haag, \emph{Local quantum physics. {F}ields, particles, algebras}, second
  ed., Springer--Verlag, Berlin, Heidelberg, New York, 1996.

\bibitem[HK64]{HK}
R.~Haag and D.~Kastler, \emph{An algebraic approach to quantum field theory},
  J. Math. Phys. \textbf{5} (1964), 848.

\bibitem[HS62]{HS}
R.~Haag and B.~Schroer, \emph{Postulates of quantum field theory}, J. Math.
  Phys. \textbf{3} (1962), 248.

\bibitem[Izu93]{Iz}
M.~Izumi, \emph{Subalgebras of infinite {\CALs} with finite {W}atatani indices.
  {I}: {C}untz algebras}, Commun. Math. Phys. \textbf{155} (1993), 157.

\bibitem[Jon83]{J}
V.~F.~R. Jones, \emph{Index for subfactors}, Invent. Math. \textbf{72} (1983),
  1.

\bibitem[JW28]{JW}
P.~Jordan and E.~P. Wigner, \emph{\"{U}ber das {P}aulische
  \"{A}quivalenzverbot}, Z. Phys. \textbf{47} (1928), 631.

\bibitem[Kat66]{K}
T.~Kato, \emph{Perturbation theory for linear operators}, Springer--Verlag,
  Berlin, Heidelberg, New York, 1966.

\bibitem[Kir94a]{Ki}
E.~Kirchberg, Lecture at {ICM} satellite meeting, Geneva, 1994.

\bibitem[Kir94b]{Ki94}
E.~Kirchberg, \emph{The classification of purely infinite {\CALs} using
  {K}asparov's theory}, preprint, Humboldt {U}niv. {B}erlin, 1994.

\bibitem[Kir95]{Ki95}
E.~Kirchberg, \emph{On subalgebras of the {CAR}--algebra}, J. Funct. Anal.
  \textbf{129} (1995), 35.

\bibitem[KMR90]{KMR}
D.~Kastler, M.~Mebkhout, and K.-H. Rehren, \emph{Introduction to the algebraic
  theory of {\SSSs} (space--time dimension = 2 --- strictly localized
  morphisms)}, The Algebraic Theory of Superselection Sectors. {I}ntroduction
  and Recent Results (I.S.I. Guccia, Palermo 1989) (D.~Kastler, ed.), World
  Scientific, Singapore, 1990, p.~113.

\bibitem[Kos86]{Ko}
H.~Kosaki, \emph{Extension of {J}ones' theory on index to arbitrary factors},
  J. Funct. Anal. \textbf{66} (1986), 123.

\bibitem[KP97]{KP}
E.~Kirchberg and N.~C. Phillips, \emph{Embedding of exact {\CALs} and
  continuous fields in the {C}untz algebra {$\mathcal{O}_2$}}, preprint, 1997,
  funct-an/9712002.

\bibitem[KS77]{KS}
M.~Klaus and G.~Scharf, \emph{The regular external field problem in quantum
  electrodynamics}, Helv. Phys. Acta \textbf{50} (1977), 779.

\bibitem[Kun97]{Walter}
W.~Kunhardt, \emph{On infravacua and superselection theory}, preprint, Univ.
  G{\"o}ttingen, 1997.

\bibitem[Lab74]{Lab}
G.~Labont{\'e}, \emph{On the nature of `strong' {B}ogoliubov transformations
  for {F}ermions}, Commun. Math. Phys. \textbf{36} (1974), 59.

\bibitem[Lab75]{Lab75}
G.~Labont{\'e}, \emph{On the particle interpretation for the quantized fields
  in external field problems}, Canad. J. Phys. \textbf{53} (1975), 1533.

\bibitem[Lac93a]{La93a}
M.~Laca, \emph{Endomorphisms of {$\BB(\HH)$} and {C}untz algebras}, J. Operator
  Theory \textbf{30} (1993), 85.

\bibitem[Lac93b]{La93b}
M.~Laca, \emph{Gauge invariant states of {$\mathcal{O}_\infty$}}, J. Operator
  Theory \textbf{30} (1993), 381.

\bibitem[Lon89]{L}
R.~Longo, \emph{Index of subfactors and statistics of quantum fields. {I}},
  Commun. Math. Phys. \textbf{126} (1989), 217.

\bibitem[Lon90]{L2}
R.~Longo, \emph{Index of subfactors and statistics of quantum fields {II}.
  {C}orrespondences, {B}raid group statistics}, Commun. Math. Phys.
  \textbf{130} (1990), 285.

\bibitem[Maa71]{M}
H.~Maa{\ss}, \emph{Siegel's modular forms and {D}irichlet series},
  Springer--Verlag, Berlin, Heidelberg, New York, 1971.

\bibitem[Man70]{Man}
J.~Manuceau, \emph{{\'E}tude alg\'ebrique des \'etats quasi--libres}, Carg\`ese
  Lectures in Physics (D.~Kastler, ed.), vol.~4, Gordon and Breach, New York,
  London, Paris, 1970.

\bibitem[Mat87a]{Mat87}
T.~Matsui, \emph{The index of scattering operators of {D}irac equations},
  Commun. Math. Phys. \textbf{110} (1987), 553.

\bibitem[Mat87b]{M87a}
T.~Matsui, \emph{On quasi--equivalence of quasifree states of gauge invariant
  {CAR} algebras}, J. Operator Theory \textbf{17} (1987), 281.

\bibitem[Mat90]{Mat90}
T.~Matsui, \emph{The index of scattering operators of {D}irac equations, {II}},
  J. Funct. Anal. \textbf{94} (1990), 93.

\bibitem[MRT69]{MRT}
J.~Manuceau, F.~Rocca, and D.~Testard, \emph{On the product form of quasi--free
  states}, Commun. Math. Phys. \textbf{12} (1969), 43.

\bibitem[MS90]{MS}
G.~Mack and V.~Schomerus, \emph{Conformal field algebras with quantum symmetry
  from the theory of superselection sectors}, Commun. Math. Phys. \textbf{134}
  (1990), 139.

\bibitem[MS92]{MS92}
G.~Mack and V.~Schomerus, \emph{Quasi {H}opf quantum symmetry in quantum
  theory}, Nuclear Phys. B \textbf{370} (1992), 185.

\bibitem[MS95]{MSch}
J.~Mund and R.~Schrader, \emph{{\HSPs} for nonrelativistic and relativistic
  "free" plektons (particles with braid group statistics)}, Advances in
  Dynamical Systems and Quantum Physics (Capri 1993) (S.~Albeverio, R.~Figari,
  E.~Orlandi, and A.~Teta, eds.), World Scientific, Singapore, 1995, p.~235.

\bibitem[M{\"u}g97]{MM2}
M.~M{\"u}ger, \emph{Superselection structure of massive quantum field theories
  in $1+1$ dimensions}, preprint 97--081, DESY, Hamburg, 1997.

\bibitem[M{\"u}g98]{MM1}
M.~M{\"u}ger, \emph{Quantum double actions on operator algebras and orbifold
  quantum field theories}, Commun. Math. Phys. \textbf{191} (1998), 137.

\bibitem[Mur62]{Mur}
F.~D. Murnaghan, \emph{The unitary and rotation groups}, Spartan Books,
  Washington {D.~C.}, 1962.

\bibitem[MV68]{MV}
J.~Manuceau and A.~Verbeure, \emph{Quasi--free states of the {C.C.R.}--algebra
  and {B}ogoliubov transformations}, Commun. Math. Phys. \textbf{9} (1968),
  293.

\bibitem[MY95]{MY}
T.~Murakami and S.~Yamagami, \emph{On types of quasifree representations of
  {C}lifford algebras}, Publ. Res. Inst. Math. Sci. \textbf{31} (1995), 33.

\bibitem[NSW98]{NSW}
F.~Nill, K.~Szlach{\'a}nyi, and H.-W. Wiesbrock, \emph{Weak {H}opf algebras and
  reducible {J}ones inclusions of depth $2$. {I}: {F}rom crossed products to
  {J}ones towers}, preprint, Sfb 288, 1998.

\bibitem[Pet90]{Pe}
D.~Petz, \emph{An invitation to the algebra of canonical commutation
  relations}, Leuven Notes in Mathematical and Theoretical Physics, Series A:
  Mathematical Physics, 2, Leuven Univ. Press, 1990.

\bibitem[Pow67]{P67}
R.~T. Powers, \emph{Representations of the canonical anticommutation
  relations}, Ph.D. thesis, Princeton Univ., 1967.

\bibitem[Pow87]{P87}
R.~T. Powers, \emph{A non spatial continuous semigroup of *--endomorphisms of
  {$\BB(\HH)$}}, Publ. Res. Inst. Math. Sci. \textbf{23} (1987), 1053.

\bibitem[Pow88]{P88}
R.~T. Powers, \emph{An index theory for semigroups of *--endomorphisms of
  ${\BB(\HH)}$ and type {${\rm II}_1$} factors}, Canad. J. Math. \textbf{{XL}}
  (1988), no.~1, 86.

\bibitem[PR94]{PR}
R.~J. Plymen and P.~L. Robinson, \emph{Spinors in {H}ilbert space}, Cambridge
  tracts in mathematics, no. 114, Cambridge University Press, 1994.

\bibitem[PS70]{PS}
R.~T. Powers and E.~St{\o}rmer, \emph{Free states of the canonical
  anticommutation relations}, Commun. Math. Phys. \textbf{16} (1970), 1.

\bibitem[PS86]{PrSe}
A.~Pressley and G.~Segal, \emph{Loop groups}, Oxford mathematical monographs,
  Oxford Univ. Press, 1986.

\bibitem[Rec93]{Rec93}
A.~Recknagel, \emph{Fusion rules from algebraic {K}--theory}, Int. J. Mod.
  Phys. A \textbf{8} (1993), 1345.

\bibitem[Rec96]{Rec96}
A.~Recknagel, \emph{From path representations to global morphisms for a class
  of minimal models}, preprint ETH--TH/96--44, ETH Z\"urich, 1996.

\bibitem[Reh90]{Re90}
K.-H. Rehren, \emph{Braid group statistics and their superselection rules}, The
  Algebraic Theory of Superselection Sectors. {I}ntroduction and Recent Results
  (I.S.I. Guccia, Palermo 1989) (D.~Kastler, ed.), World Scientific, Singapore,
  1990, p.~333.

\bibitem[Reh96]{Re96}
K.-H. Rehren, \emph{Weak {$C^*$} {H}opf symmetry}, preprint 96--231, DESY,
  1996.

\bibitem[Reh97]{Re97}
K.-H. Rehren, \emph{Spin--statistics and {CPT} for solitons}, preprint 501,
  ESI, 1997.

\bibitem[Rid68]{Ri}
G.~Rideau, \emph{On some representations of the anticommutation relations},
  Commun. Math. Phys. \textbf{9} (1968), 229.

\bibitem[Rob70]{R70}
J.~E. Roberts, \emph{The structure of sectors reached by a field algebra},
  Carg\`ese Lectures in Physics 1969 (D.~Kastler, ed.), vol.~4, Gordon and
  Breach, New York, London, Paris, 1970, p.~61.

\bibitem[Rob76a]{R76}
J.~E. Roberts, \emph{Cross products of von {N}eumann algebras by group duals},
  Sympos. Math. \textbf{20} (1976), 335.

\bibitem[Rob76b]{R74}
J.~E. Roberts, \emph{Spontaneously broken gauge symmetries and superselection
  rules}, Proceedings Camerino 1974 (G.~Gallavotti, ed.), 1976.

\bibitem[Rob80]{R80}
J.~E. Roberts, \emph{Net cohomology and its applications to field theory},
  Quantum Fields --- Algebras, Processes. Proceedings of the Symposium
  "Bielefeld Encounters in Physics and Mathematics II" 1978 (Wien, New York)
  (L.~Streit, ed.), Springer--Verlag, 1980, p.~239.

\bibitem[Rob90]{R90}
J.~E. Roberts, \emph{Lectures on algebraic quantum field theory}, The Algebraic
  Theory of Superselection Sectors. {I}ntroduction and Recent Results (I.S.I.
  Guccia, Palermo 1989) (D.~Kastler, ed.), World Scientific, Singapore, 1990,
  p.~1.

\bibitem[Rob93]{Ro}
P.~L. Robinson, \emph{The even {$C^*$} {C}lifford algebra}, Proc. Amer. Math.
  Soc. \textbf{118} (1993), 713.

\bibitem[R{\o}r93]{R93}
M.~R{\o}rdam, \emph{Classification of inductive limits of {C}untz algebras}, J.
  Reine Angew. Math. \textbf{440} (1993), 175.

\bibitem[RS75]{RS2}
M.~Reed and B.~Simon, \emph{Methods of modern mathematical physics}, vol. II:
  Fourier analysis, self-adjointness, Academic Press, New York, San Francisco,
  London, 1975.

\bibitem[RS89]{RSch}
K.-H. Rehren and B.~Schroer, \emph{Einstein causality and {A}rtin braids},
  Nuclear Phys. B \textbf{312} (1989), 715.

\bibitem[Rue62]{Ru}
D.~Ruelle, \emph{On the asymptotic condition in quantum field theory}, Helv.
  Phys. Acta \textbf{35} (1962), 147.

\bibitem[Rui77]{R77a}
S.~N.~M. Ruijsenaars, \emph{On {B}ogoliubov transformations for systems of
  relativistic charged particles}, J. Math. Phys. \textbf{18} (1977), 517.

\bibitem[Rui78]{R78}
S.~N.~M. Ruijsenaars, \emph{On {B}ogoliubov transformations. {II}. {T}he
  general case}, Ann. Physics \textbf{116} (1978), 105.

\bibitem[Rui89a]{R89}
S.~N.~M. Ruijsenaars, \emph{Index formulas for generalized {W}iener--{H}opf
  operators and {B}oson-{F}ermion correspondence in ${2N}$ dimensions}, Commun.
  Math. Phys. \textbf{124} (1989), 553.

\bibitem[Rui89b]{R89a}
S.~N.~M. Ruijsenaars, \emph{Index theorems and anomalies: a common playground
  for mathematicians and physicists}, preprint, Stichting Mathematisch Centrum,
  Amsterdam, 1989.

\bibitem[Sch95]{VSch}
V.~Schomerus, \emph{Construction of {\FAs} with quantum symmetry from local
  observables}, Commun. Math. Phys. \textbf{169} (1995), 193.

\bibitem[Sch96a]{Schl}
D.~Schlingemann, \emph{On the existence of kink (soliton) states}, Rev. Math.
  Phys. \textbf{8} (1996), 1187.

\bibitem[Sch96b]{BSch96}
B.~Schroer, \emph{Motivations and physical aims of algebraic {QFT}}, Ann.
  Physics \textbf{255} (1996), 270.

\bibitem[Sch97a]{BSch97b}
B.~Schroer, \emph{Modular localization and the bootstrap--formfactor program},
  Nuclear Phys. B \textbf{499} (1997), 547.

\bibitem[Sch97b]{BSch97a}
B.~Schroer, \emph{Wigner representation theory of the {P}oincar{\'e} group,
  localization, statistics and the ${S}$--matrix}, Nuclear Phys. B \textbf{499}
  (1997), 519.

\bibitem[Sch98a]{BSch98}
B.~Schroer, \emph{Coincidences between {M}(atrix) theory and algebraic {QFT}?},
  preprint 305, Sfb 288, 1998.

\bibitem[Sch98b]{BSch}
B.~Schroer, \emph{A {C}ourse on: ``{M}odular localization and nonperturbative
  local quantum physics}, Lecture notes, available at
  http://www.lqp.uni-goettingen.de/lqp/books/, 1998.

\bibitem[Seg81]{Se}
G.~Segal, \emph{Unitary representations of some infinite dimensional groups},
  Commun. Math. Phys. \textbf{80} (1981), 301.

\bibitem[Sei78]{Sei78}
R.~Seiler, \emph{Particles with spin ${S}\leq1$ in an external field},
  Proceedings of the 1977 Erice Summer School on Invariant Wave Equations
  (G.~Velo and A.~S. Wightman, eds.), Lecture Notes in Physics, no.~73,
  Springer--Verlag, Berlin, Heidelberg, New York, 1978, p.~165.

\bibitem[Sha62]{S}
D.~Shale, \emph{Linear symmetries of free boson fields}, Trans. Amer. Math.
  Soc. \textbf{103} (1962), 149.

\bibitem[Sie64]{Si}
C.~L. Siegel, \emph{Symplectic geometry}, Academic Press, New York, London,
  1964.

\bibitem[SS64]{SS64}
D.~Shale and W.~F. Stinespring, \emph{States of the {C}lifford algebra}, Ann.
  of Math. (2) \textbf{80} (1964), 365.

\bibitem[SS65]{SS65}
D.~Shale and W.~F. Stinespring, \emph{Spinor representations of infinite
  orthogonal groups}, J. Math. Mech. \textbf{14} (1965), 315.

\bibitem[SSS70]{SSS}
B.~Schroer, R.~Seiler, and J.~A. Swieca, \emph{Problems of stability for
  quantum fields in external time--dependent potentials}, Phys. Rev. D (3)
  \textbf{2} (1970), 2927.

\bibitem[St{\o}70]{St}
E.~St{\o}rmer, \emph{The even {CAR} algebra}, Commun. Math. Phys. \textbf{16}
  (1970), 136.

\bibitem[Str74]{Str74}
R.~F. Streater, \emph{Charges and currents in the {T}hirring model}, Physical
  Reality and Mathematical Description (C.~P. Enz and J.~Mehra, eds.), Reidel,
  Dordrecht, 1974, p.~375.

\bibitem[SW68]{PCT}
R.~F. Streater and A.~S. Wightman, \emph{{PCT}, spin \& statistics, and all
  that}, second ed., W. A. Benjamin, New York, Amsterdam, 1968.

\bibitem[SW70]{StrW}
R.~F. Streater and I.~F. Wilde, \emph{Fermion states of a {B}oson field},
  Nuclear Phys. B \textbf{24} (1970), 561.

\bibitem[Szl94]{Sz}
K.~Szlach{\'a}nyi, \emph{Chiral decomposition as a source of quantum symmetry
  in the {I}sing model}, Rev. Math. Phys. \textbf{6} (1994), 649.

\bibitem[vD71]{vD}
A.~van Daele, \emph{Quasi--equivalence of quasi--free states on the {W}eyl
  algebra}, Commun. Math. Phys. \textbf{21} (1971), 171.

\bibitem[vN39]{vN}
J.~von Neumann, \emph{On infinite direct products}, Compositio Math. \textbf{6}
  (1939), 1.

\bibitem[Was]{Was}
A.~J. Wassermann, \emph{Operator algebras and conformal field theory {III}},
  preprint, Univ. Cambridge, available at
  http://www.pmms.cam.ac.uk/Staff/A.J.Wassermann.html.

\bibitem[Wat90]{Wat}
Y.~Watatani, \emph{Index for ${C^*}$--subalgebras}, Memoirs of the {A}m.
  {M}ath. {S}oc. {V}ol. 83 {N}o. 424, 1990.

\bibitem[Wey28]{Wey28}
H.~Weyl, \emph{Quantenmechanik und {G}ruppentheorie}, Z. Phys. \textbf{46}
  (1928), 1.

\bibitem[Wey46]{Wey}
H.~Weyl, \emph{The classical groups}, second ed., Princeton University Press,
  1946.

\bibitem[Wie96]{HWW96}
H.-W. Wiesbrock, \emph{Modular inclusions and intersections of algebras in
  {QFT}}, Operator Algebras and Quantum Field Theory (S.~Doplicher, R.~Longo,
  J.~E. Roberts, and L.~Zsido, eds.), Accademia Nazionale dei Lincei, Roma,
  International Press, 1996, p.~609.

\bibitem[Wig95]{W95}
A.~S. Wightman, \emph{Superselection rules; old and new}, Nuovo Cimento B (11)
  \textbf{110} (1995), 751.

\bibitem[Wol75]{W75}
J.~C. Wolfe, \emph{Free states and automorphisms of the {C}lifford algebra},
  Commun. Math. Phys. \textbf{45} (1975), 53.

\bibitem[WWW52]{WWW}
G.~C. Wick, A.~S. Wightman, and E.~P. Wigner, \emph{The intrinsic parity of
  elementary particles}, Phys. Rev. \textbf{88} (1952), 101.

\end{thebibliography}
